\ifdefined\DoubleSided
  \documentclass[twoside,openright,a4paper]{nusthesis}
\else
  \documentclass[a4paper]{nusthesis}
\fi

\usepackage[utf8]{inputenc}
\usepackage[T1]{fontenc}

\usepackage{bookmark}
\usepackage[
  backend=biber,
  style=ieee,
  citestyle=numeric,
  giveninits=true,
  sorting=nyt,
  maxbibnames=99,
  maxcitenames=99,
  dashed=false,
  doi=false]{biblatex}
\DeclareFieldFormat*{title}{``#1''\newunitpunct} % move comma outside the quotation mark
\usepackage{microtype} % Better typography
\usepackage{tabulary} % Automatic table sizing
\usepackage{hhline}

\usepackage{hyperref}

\usepackage{listings}

\usepackage{mathtools}

\usepackage{interval}
\intervalconfig{soft open fences}

\usepackage[linesnumbered,ruled,vlined]{algorithm2e}
\SetAlgorithmName{Algorithm}{Algorithm}{Algorithm}
\IncMargin{0.5em}
\SetCommentSty{textnormal}
\SetNlSty{}{}{:}
\SetAlgoNlRelativeSize{0}
\SetKwInput{KwGlobal}{Global}
\SetKwInput{KwPrecondition}{Precondition}
\SetKwProg{Proc}{Procedure}{:}{}
\SetKwProg{Func}{Function}{:}{}
\SetKw{And}{and}
\SetKw{Or}{or}
\SetKw{To}{to}
\SetKw{DownTo}{downto}
\SetKw{Break}{break}
\SetKw{Continue}{continue}
\SetKw{SuchThat}{\textit{s.t.}}
\SetKw{WithRespectTo}{\textit{wrt}}
\SetKw{Iff}{\textit{iff.}}
\SetKw{MaxOf}{\textit{max of}}
\SetKw{MinOf}{\textit{min of}}
\SetKwBlock{Match}{match}{}{}

\usepackage{array}
\newcolumntype{L}[1]{>{\raggedright\let\newline\\\arraybackslash\hspace{0pt}}m{#1}}
\newcolumntype{C}[1]{>{\centering\let\newline\\\arraybackslash\hspace{0pt}}m{#1}}
\newcolumntype{R}[1]{>{\raggedleft\let\newline\\\arraybackslash\hspace{0pt}}m{#1}}

\usepackage{color}
\usepackage[dvipsnames]{xcolor}
\usepackage{soul}
\soulregister\cite7
\soulregister\ref7
\soulregister\pageref7
\soulregister\autoref7
\soulregister\eqref7

\usepackage{lipsum} % for generating dummy text
\usepackage{amsmath, amsfonts, amssymb, amsthm, amscd}
\usepackage{graphicx}
\usepackage{psfrag}
\usepackage{perpage}
\usepackage{url}
\usepackage{mathrsfs}
\usepackage{scalerel}
\usepackage{enumerate}
\usepackage{comment}
\usepackage[normalem]{ulem}
\usepackage{tikz}\usepackage{caption,subcaption}
\usetikzlibrary{arrows,decorations.pathmorphing,backgrounds,positioning,fit,petri}
\usepackage{tikz-cd} % Note taking
\usepackage{bm,bbm}
\usepackage{mathtools}
\usepackage{subdepth}

\usepackage{nomencl} % List of figures
\makenomenclature

\usepackage{dsfont} % nice indicator function symbol (blackboard 1)

\usepackage{microtype}
\allowdisplaybreaks[0] % [0]disable align page breaks

\numberwithin{equation}{section}

\newtheorem{theorem}{Theorem}[section]
\newtheorem*{theorem*}{Theorem}
\newtheorem{lemma}[theorem]{Lemma}
\newtheorem*{lemma*}{Lemma}
\newtheorem{proposition}[theorem]{Proposition}
\newtheorem*{proposition*}{Proposition}
\newtheorem{corollary}[theorem]{Corollary}
\newtheorem*{corollary*}{Corollary}

\newtheorem{conjecture}[theorem]{Conjecture}
\newtheorem*{conjecture*}{Conjecture}

\theoremstyle{definition}
\newtheorem{definition}[theorem]{Definition}
\newtheorem{question}[theorem]{Question}
\newtheorem{remark}[theorem]{Remark}

\newtheorem{assumption}[theorem]{Assumption}
\newtheorem{example}[theorem]{Example}
\newcommand{\E}{\mathbb{E}}
\newcommand{\N}{\mathbb{N}}
\renewcommand{\P}{\mathbb{P}}

\newcommand{\R}{\mathbb{R}}
\newcommand{\Z}{\mathbb{Z}}
\newcommand{\T}{\mathbb{T}}

\def\bs{\boldsymbol}

\newcommand\bP{\ensuremath{\bs{\mathrm{P}}}}
\newcommand\bE{\ensuremath{\bs{\mathrm{E}}}}

\newcommand{\cC}{{\ensuremath{\mathcal C}} }

\newcommand{\cT}{{\ensuremath{\mathcal T}} }

\newcommand{\ClCp}[2]{ \cT_{\cC_{#1}(#2)}} % \overline{\mathcal{C}_{#1}(#2)}}

\newcommand{\bigO}[1]{O\Big(#1\Big)}
\newcommand{\Ker}{{\mathrm{Ker}}}
\newcommand{\ErdosRenyi}{Erd\H{o}s--R\'enyi }
\newcommand{\LE}{\mathrm{LE}}

\newcommand{\e}{{\mathrm{e}}}
\newcommand{\weight}{\mathrm{w}}
\newcommand{\diam}{\mathrm{diam}}
\renewcommand{\theta}{\vartheta}
\renewcommand{\rho}{\varrho}

\newcommand{\effR}[3]{{\mathrm{R}}_{\mathrm{eff}}^{#1}(#2 \leftrightarrow #3)} % Effective resistance
\newcommand{\tmix}{\ensuremath{t_{\text{\normalfont mix}}}}

\newcommand{\limn}{\lim\limits_{n \rightarrow \infty}}

\makeatletter
\renewcommand\xleftrightarrow[2][]{%
  \ext@arrow 9999{\longleftrightarrowfill@}{#1}{#2}}
\newcommand\longleftrightarrowfill@{%
  \arrowfill@\leftarrow\relbar\rightarrow}
\makeatother

\begin{document}

\title{Observables of random spanning trees in random environment}

\author{Luca Makowiec}
\prevdegrees{%
BSc.\ Technische Universit\"at Berlin, MSc.\ Utrecht University}
\degree{Doctor of Philosophy}
\field{Mathematics}
\degreeyear{2025}
\supervisor{Professor Rongfeng Sun}

\maketitle

%\declaredate{\today}

\begin{frontmatter}

  %\dedicate{To my family and partner.}
  
  \begin{acknowledgments}

I am sincerely grateful to my supervisor, Rongfeng Sun, for their guidance and patience during the course of this work. His advice has been very helpful and significantly enhanced my mathematical skills, and without his support, much of this would not have been possible. I also thank Michele Salvi for the many interesting and fruitful discussions that helped shape this thesis.

I would like to thank my family and partner for their support and encouragement throughout this journey. Their unwavering belief in me has provided invaluable motivation and reassurance.

My gratitude extends to the National University of Singapore for providing a supportive academic environment and for offering resources that enabled me to attend overseas conferences and workshops. These opportunities have greatly contributed to my development as a researcher.

Finally, I wish to express my appreciation to the researchers and fellow students in mathematics, who have always been friendly and with whom it has been a joy to discuss ideas.

\end{acknowledgments}

  \tableofcontents 
  \begin{abstract}

In this thesis, we study a new disordered system called random spanning tree in random environment (RSTRE) across different families of graphs with varying disorder distributions. %This model is defined by assigning random disorder variables to the edges, and by considering a Gibbs measure with disorder strength (inverse temperature) $\beta \geq 0$ and Hamiltonian given by the sum of the disorder variables. 
We examine several observables as functions of the disorder strength (inverse temperature) $\beta \geq 0$, and compare their values to the extreme cases $\beta = 0$ and $\beta \to \infty$, which correspond to the uniform spanning tree (UST) and the minimum spanning tree (MST), respectively.

After introducing the relevant background information about USTs, random walks, and electric networks, we show the existence of a limiting RSTRE measure on the Euclidean lattice in $\mathbb{Z}^d$ by considering the RSTRE on finite graphs in an exhaustion of $\mathbb{Z}^d$. We then analyze the number of connected components and study the overlap density of two trees sampled under the same environment.

We then turn our focus to local observables of the RSTRE on the complete graph on $n$ vertices, with the disorder variables being distributed uniformly on $[0,1]$. We show that when $\beta = o(n/\log n)$ the length (Hamiltonian) is $(1 + o(1)) \beta/n$ and the edge overlap is $(1+o(1)) \beta$, while for $\beta$ much larger than $n \log^2 n$, the edge overlap is $(1-o(1))n$ and the length is approximately that of the MST. Furthermore, there is a sharp transition of the local limit when $\beta = n^{\gamma}$ and $\gamma$ crosses the critical threshold $\gamma_c = 1$: when $\gamma < \gamma_c$ the RSTRE locally converges to the same limit as the uniform spanning tree, whereas for $\gamma > \gamma_c$ the local limit of the RSTRE is the same as that of the MST.

Finally, we consider the global observable of the diameter, which is of particular interest as it is required to establish a non-trivial scaling limit of the random spanning trees. When the underlying graphs are either bounded degree expanders or boxes in $\Z^d$ ($d \geq 5)$ on $n$ vertices, then, regardless of the edge disorder distribution, the diameter of the RSTRE is of the order $\sqrt{n}$ (up to logarithmic terms), the same order as the diameter of the UST. For the complete graph with the disorder variables being distributed uniformly on $[0,1]$, we show that for $\beta = O(n/\log n)$ the diameter of the RSTRE is of order $\sqrt{n}$, while for $\beta \geq n^{4/3} \log n$ the order of the diameter is $n^{1/3}$, the same order as the diameter of the MST.

\end{abstract}

  % \listoffigures
  % \listoftables
   \renewcommand{\nomname}{List of symbols}

\nomenclature[0001]{$G$}{Undirected graph $G = (V,E)$}
\nomenclature[0002]{$V(G)$}{Vertices of the graph $G$}
\nomenclature[0003]{$E(G)$}{Edges of the graph $G$}
\nomenclature[001]{$\T(G)$}{Set of spanning trees of $G$}
\nomenclature[003]{$\mu$}{Distribution of the edge disorder}
\nomenclature[002]{$\omega$}{Edge disorder (random environment)}
\nomenclature[004]{$\P$}{Joint probability distribution of the environment}
\nomenclature[005]{$\E$}{Expectation of $\P$}
\nomenclature[006]{$\beta$}{Disorder strength}
\nomenclature[007]{$\weight$}{Edge weights (conductances)}
\nomenclature[0071]{$(G,\weight)$}{Weighted graph}
\nomenclature[008]{$\bP^\weight$}{UST probability measure given weights on the edges}
\nomenclature[009]{$\bE^\weight$}{UST expectation given weights on the edges}
\nomenclature[0081]{$\bP^\omega_{G, \beta}$}{RSTRE probability measure given the environment}
\nomenclature[0091]{$\bE^\omega_{G, \beta}$}{RSTRE expectation given the environment}
\nomenclature[010]{$\widehat{\P}$}{Joint law of the environment and the RSTRE}
\nomenclature[011]{$\widehat{\E}$}{Joint expectation of the environment and the RSTRE}
\nomenclature[012]{$\cT$}{Random spanning tree sampled from either $\bP^\weight$ or $\bP^\omega_{G, \beta}$}
\nomenclature[013]{$H(T, \omega)$}{Hamiltonian $\sum_{e \in T} \omega_e$ of a spanning tree $T$}

\nomenclature[11]{$q(u,v)$}{One-step transition kernel of the lazy random walk on $(G, \weight)$}
\nomenclature[12]{$Q_v( \cdot )$}{Probability measure of the lazy random walk on $(G, \weight)$ started at $v$}
\nomenclature[13]{$\pi(\cdot)$}{Stationary distribution of the lazy random walk on $(G, \weight)$} 
\nomenclature[14]{$\tmix$}{Mixing time of the lazy random walk on $(G, \weight)$}
\nomenclature[15]{$\Phi_G$}{Bottleneck ratio of the lazy random walk on $(G, \weight)$}
\nomenclature[16]{$\Phi_G(r)$}{Bottleneck profile of the lazy random walk on $(G, \weight)$}

\nomenclature[21]{$G_{n,p}$}{\ErdosRenyi random graph}
\nomenclature[220]{$\cC_\ell(p)$}{The $\ell$-th largest component of $G_{n,p}$}
\nomenclature[221]{$\cC_\ell(G, p)$}{The $\ell$-th largest component of percolation with parameter $p$ on $G$}
\nomenclature[23]{$\mathcal{O}(\beta)$}{Edge overlap of two independent RSTREs under the same environment}
\nomenclature[014]{$\effR{G,\weight}{A}{B}$}{Effective resistance between $A$ and $B$ in $(G,\weight)$}

% \nomenclature{$a$}{}

\addcontentsline{toc}{chapter}{List of symbols}
\printnomenclature

\end{frontmatter}

\chapter{Introduction}

Spanning trees play a central role in both combinatorial optimization and probability theory, emerging in the study of graph structures and stochastic processes. In this thesis, we study a new model called \textit{random spanning trees in random environment} (RSTRE), which interpolates between the \textit{uniform spanning tree} (UST) and \textit{minimum spanning tree} (MST) measures through a strength parameter $\beta$, which governs the influence of the random edge disorder.

\section{Uniform spanning trees}

The study of \textit{uniform spanning trees} (USTs) on finite graphs dates back to Borchardt \cite{Bor60} and Cayley \cite{Cay89}, who analyzed the number of such trees on the complete graph. Since then, a substantial amount of research has been done for the UST, both in the context of combinatorial enumeration arguments \cite{Moo70, MM78, Gri80, FO82}, and using techniques from probability theory. Central to the probabilistic analysis of USTs is the elegant connection to random walks discovered by Broder \cite{Bro89} and Aldous \cite{Ald90}, and later extended by Wilson \cite{Wil96} to loop-erased random walks (LERW), leading to the much celebrated Wilson's algorithm. This tool has led to significant developments such as the introduction of Schramm–Loewner evolutions \cite{Sch00}, which describe the scaling limits of LERWs in 2D  (corresponding to $\kappa = 2$) and the paths separating the UST from its dual (corresponding to $\kappa = 8$) on the 2D lattice \cite{LSW04}.

By taking an infinite volume limit, a natural extension of the UST measure to the lattice $\Z^d$ was introduced and analyzed by Pemantle \cite{Pem91}. In this work, the existence of a \textit{uniform spanning forest} (USF) measure on $\Z^d$ was established, and it was shown that the measure concentrates on connected graphs if and only if $d \leq 4$. %That is, $d_c = 4$ is the critical dimension above which the UST behaves in a \textit{mean-field} manner. 
Later the authors of \cite{BLPS01} further extended these results to general graphs using Wilson's algorithm rooted at infinity, and analyzed topological properties such as the number of ends of the connected components. More recently, Hutchcroft \cite{Hut18} used the interlacement process by Sznitman \cite{Szn10} to construct the USF, which gives rise to yet another tool useful for the study of both the UST and USF. 

Returning to finite graphs, Aldous \cite{Ald91a, Ald91b, Ald93} proved the existence of a scaling limit for the UST on the complete graph $K_n$, when distances are rescaled by a factor of $n^{-1/2}$. The fascinating fractal limiting object, known as Aldous' \textit{Brownian continuum random tree} (Brownian CRT), emerges as the limit of various combinatorial objects, such as Galton-Watson trees with finite variance conditioned on their size \cite{Ald93} and Pólya trees \cite{MM11}, and may be encoded as a Brownian excursion on $[0,1]$. This has given rise to a large amount of research into general continuum random trees, see e.g. the survey by Le Gall \cite{LeG06} and the references therein. The convergence of finite-dimensional distributions has been extended to the torus in $\Z^d$ (under proper rescaling) in dimension 4 \cite{Sch09} and above \cite{PR05}, and is believed to hold for ``high dimensional'' graphs with certain regularity structures. For instance, in \cite{ANS21} and \cite{ANS22}, it is shown that the UST on either transitive high-dimensional graphs or dense graphs converges to Aldous' Brownian CRT with respect to the Gromov-Hausdorff-Prokhorov distance. Additionally, in \cite{MNS21} it is proven that the diameter and typical distances of the UST in high-dimensional graphs, with some regularity structure but not requiring transitivity, are of order $\sqrt{n}$.

%%%%%%%%%%%%%%%%%%%%%%%%%% NEW SECTION %%%%%%%%%%%%%%%%%%%%%%%%%%

\section{Minimum spanning trees} \label{S:MSTintro}

Another important spanning tree model, stemming from combinatorial optimization, is the \textit{minimum spanning tree} (MST). Introduced by Borůvka \cite{Bor26}, the MST has many practical applications and is also connected to invasion percolation and the ground states of a spin glass model \cite{NS96}. We define the MST in the following way: given distinct weights $\weight(\cdot)$ on the edges of a finite connected graph $G = (V,E)$, the MST is the unique tree $M = M(\weight)$ that minimizes
\begin{equation*}
    \weight(T) = \sum_{e \in T} \weight(e)
\end{equation*}
over all $T$ in the set of spanning trees on $G$. 

Finding the MST given the weights is an algorithmically easy task that is straightforward to describe and is computationally efficient. The two most well-known methods are Kruskal's algorithm \cite{Kru56} and Prim's algorithm \cite{Pri57}, which we very briefly recall now. Note that both these methods only depend on the ordering of the weights and can run in $O(n \log n)$ many steps.

\medskip
    \textbf{Kruskal's algorithm:} Start with an empty forest and iteratively add the unique lowest weighted edge in $E$ that does not create a cycle.

\smallskip \textbf{Prim's algorithm:} Start at any vertex and iteratively add the lowest outgoing edge of the current tree that does not create a cycle. 

\medskip

A natural probabilistic model is to choose the weights independently from some common continuous distribution and study the \textit{random} MST. A frequent choice for the weight distribution is the uniform distribution on $[0,1]$, which allows for an easy coupling with \ErdosRenyi random graphs. For a historical overview of the research of the MST on finite graphs, see the introduction of \cite{ABGM17}, whereas for the problem of the \textit{Euclidean} MST, we refer to the book of \cite{Pen03}. Similarly to the UST and USF, the MST has also been extended to infinite graphs, see e.g.\ \cite{LPS06} or \cite{NTW17} for some recent results. While much of the research on finite graphs focuses on the total weight of the MST and its dependence on the weight distribution, as initially investigated by Frieze \cite{Fri85}, the geometric properties of the MST have been studied much less intensively. Recently, in \cite{ABGM17} it is shown that, under appropriate rescaling, the MST on the complete graph converges to a random real tree, and this limiting object is different than Aldous' Brownian CRT. It is believed that this continuum object should be a universal limit of the MST for any high-dimensional graphs, with this being confirmed for random $3$-regular graphs in \cite{AS21}.

%%%%%%%%%%%%%%%%%%%%%%%%%% NEW SECTION %%%%%%%%%%%%%%%%%%%%%%%%%%

\section{Definition of RSTRE} \label{S:def_RSTRE}

Let $G = (V,E)$ be a finite connected undirected graph with positive weights $\weight : E \rightarrow \R^+$ assigned to the edges of $G$. A spanning tree of $G$ is a subgraph that is both connected and cycle-free, and the (finite) set of all spanning trees of $G$ is denoted by $\T(G)$. The weighted uniform spanning tree measure on $\T(G)$ is the probability measure that picks a tree proportional to the product of the weights of the edges on the tree, in other words, for $T \in \T(G)$
\begin{equation} \label{eq:def_w_UST}
    \bP^\weight(\cT = T) = \frac{1}{Z} \prod_{e \in T} \weight(e),
\end{equation}
where $Z$ is the normalization constant (also called partition function) given by
\begin{equation*}
    Z = \sum_{T \in \T(G)} \prod_{e \in T} \weight(e).
\end{equation*}
When the weights are equal for each edge, we recover the (unweighted) UST.

In this thesis, the weights on the edges will be chosen randomly, creating an underlying random environment for the spanning trees. More precisely, let $\omega := (\omega_e)_{e \in E}$ be a collection of i.i.d.\ random variables distributed according to some common law $\mu$ supported on $\R$. Let $\E$ and $\P$ denote the (joint over all edges) expectation and probability measure associated with $\omega$, where we occasionally enlarge the probability space to capture couplings. Given any environment $(\omega_e)_{e \in E}$ and disorder strength (inverse temperature)  $\beta \geq 0$, we define the Gibbs measure
\begin{equation}\label{eq:PomegaT}
    \bP_{G, \beta}^{\omega}(\mathcal{T} = T) := \frac{1}{Z_{G, \beta}^{\omega}} \prod_{e \in T} e^{-\beta \omega_e} = \frac{1}{Z_{G, \beta}^{\omega}} \, e^{-\beta H(T, \omega)}\,,
\end{equation}
where the Hamiltonian is defined as $H(T, \omega):= \sum_{e\in T}\omega_e$ and $Z_{G, \beta}^{\omega}$ is the normalization constant.  Notice that this measure corresponds to the weighted UST measure, as in \eqref{eq:def_w_UST}, by letting $\weight(e) = \exp(-\beta \omega_e)$. We call $\bP_{G, \beta}^{\omega}$ the law of a \textit{random spanning tree in random environment} (RSTRE) with expectation $\bE^\omega_{G,\beta}$, and when there is no ambiguity, we shall drop sub- and superscripts. The random variable $\cT = \cT^\omega_{G, \beta}$ will be used to denote a random spanning tree distributed as an RSTRE.

Furthermore, we denote by $\widehat \P$ (with expectation $\widehat \E$) the joint law of $\omega$ and $\cT$, where $\omega$ has marginal law $\P$, and conditional on $\omega$, the  random spanning tree $\cT$ has law $\bP^\omega_{G, \beta}$. The marginal law 
\begin{equation*}
    \widehat \P(\cT \in \cdot) := \E \big[ \bP^\omega_{G, \beta}(\cT\in \cdot) \big] 
\end{equation*}
is known as the \textit{averaged law}, while $\bP^\omega_{G, \beta}$ is known as the \textit{quenched law}. We are interested in the behavior of observables of $\cT$ depending on the choice of $\mu$ and $\beta$, and in many parts of this thesis, we fix $\mu$ to be the uniform distribution on $[0,1]$. It will often be easier to study the RSTRE under the averaged law $\widehat{\P}$ instead of under the quenched law $\bP^\omega_{G, \beta}$.

We remark that in the extreme cases $\beta=0$ and $\beta \uparrow \infty$ (for a fixed environment and size of the graph), the RSTRE measure is simply given by the UST and MST measure, respectively. The latter follows by the observation that the MST is the unique tree that minimizes the Hamiltonian $H(T, \omega)$. Consequently, as $\beta$ varies, the RSTRE measure interpolates between the UST and MST measures. We will therefore need to identify observables that can differentiate between these measures for different choices of $\beta$. 

\begin{figure}[htp]
  \centering
  %\hspace*{-0.3cm}%
  \subcaptionbox{RSTRE with $\weight(e) = 1/U_e$ for $U_e$ uniform on $[0,1]$. \label{fig:UST}}{\includegraphics[width=0.52\textwidth]{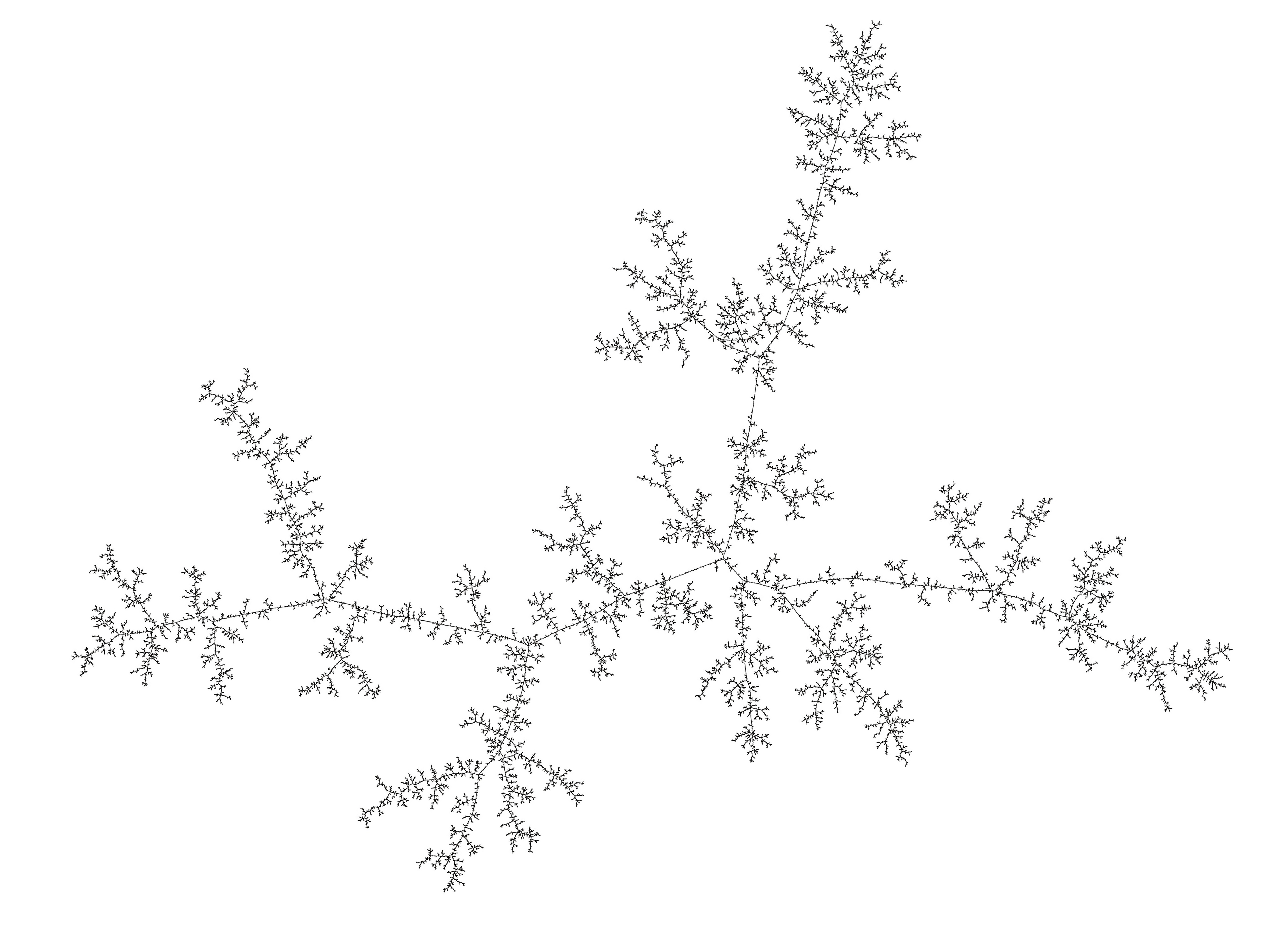}}\hspace*{-0.15cm}%
  \subcaptionbox{Random MST\label{fig:MST}}{\includegraphics[width=0.52\textwidth]{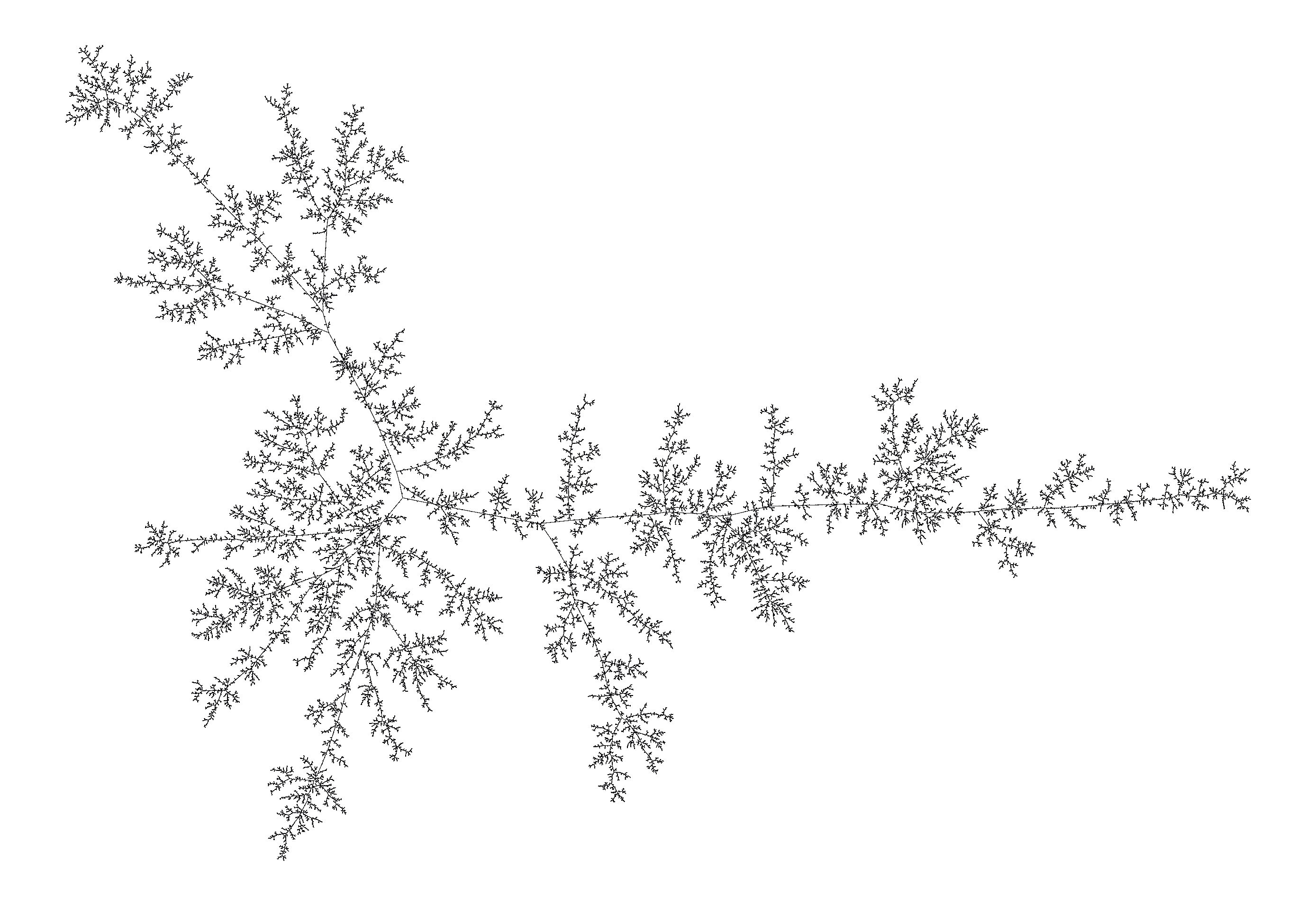}}
  \caption{A simulation of a RSTRE and random MST on a random $3$-regular graph with $100,000$ vertices. In Chapter \ref{ch:diam}, we will show that for bounded degree graphs the RSTRE looks more comparable to an unweighted UST than to a random MST.}
  \label{fig:RSTREvsMST}
\end{figure}

The RSTRE was first studied in \cite{MSS23} (under a slightly different parametrization) on expander graphs and boxes in $\Z^d$, whereas the paper \cite{MSS24} explicitly studied the dependence of the RSTRE on $\beta$ on the complete graph. Furthermore, around the same time, Kúsz \cite{K24} independently studied the RSTRE model (under a different name) on the complete graph and additionally analyzed the behavior of the algorithms used to generate the spanning trees. In particular, using similar arguments, both \cite{K24} and \cite{MSS24} manage to independently bound the diameter of the RSTRE in the low and high disorder regime of the complete graph (see Theorem \ref{T:main}, \cite[Theorems 1.6 and 1.7]{K24} and \cite[Theorem 1.1]{MSS24}). We refer to Figure \ref{fig:RSTREvsMST} for a simulation of the RSTRE and random MST.

\subsection*{Similar models} % CHECK WITH CHAT GPT FROM HERE>>>>

One motivation for the RSTRE comes from the {\em directed polymer
model in random environment} (DPRE). In the DPRE, random walk paths in $\Z^d$ are given a weight according to a random time-space environment $\big( \omega(n, x) \big)_{n \in \N, x \in \Z^d}$ similar to the underlying environment $(\omega_e)_{e \in E}$ for the RSTRE. The corresponding Gibbs measure at inverse temperature $\beta$ is then defined as
\begin{equation*}
    \bP_{N, \beta}^{\omega}(S) := \frac{1}{Z_{N, \beta}^{\omega}} e^{ \beta \sum_{n=1}^N \omega(n, S_n)} \bP(S),
\end{equation*}
where $S$ is a simple symmetric random walk (SRW) on $\Z^d$ with law $\bP$, and $Z_{N, \beta}^\omega$ is the partition function. We refer to the monograph of \cite{Com17} for an in-depth introduction. When $\beta = 0$, the DPRE is given by the simple random walk measure, while for $\beta \uparrow \infty$, the trajectory concentrates on the unique path minimizing Hamiltonian $H^\omega_N(S) = -\sum_{n=1}^N \omega(n, S_n)$. The interpolation of the RSTRE between the UST and MST mimics the dependence on $\beta$ in this model. In Chapter \ref{ch:Zd}, we approach the RSTRE from a statistical physics viewpoint in a similar fashion as for the DPRE. % In dimension $d=1$, the latter model is known as the {\em last passage percolation model} and belongs to the celebrated KPZ universality class. 

\smallskip

Yet another similar model is first passage percolation with (random) edge weights $\weight(e):=\exp(\beta \omega_e)$, where $\omega_e$ are i.i.d.\ according to some common distribution. The distance between any two vertices $u$ and $v$ is given by 
\begin{equation*}
\min_{\Gamma :u \rightarrow v} \sum_{e \in \Gamma} \weight(e),
\end{equation*}
where the minimum is over (acyclic) paths $\Gamma$ connecting $u$ and $v$. In the extreme cases $\beta = 0$ and $\beta \uparrow \infty$, the paths achieving the above minimum are given by the shortest path in the (unweighted) graph and by the unique path in the MST connecting $u$ and $v$. We refer to \cite{BvdH12, EGHN13,EGH20} for further information.

\smallskip

Lastly, we mention the model of random walks in random environment (RWRE) and refer to \cite{Zei04} for a survey of some results on $\Z^d$ for $d \geq 1$. In this model, we consider (random) weights on the edges and define a random walk with transition probabilities proportional to the edge weights. Topics studied in this model include the speed of the random walk and its convergence to Brownian motion. We note that, due to the connections of loop-erased random walks to (weighted) USTs, one may view paths in the UST as LERWs in a random environment. To the best of our knowledge, there is no existing literature on LERWs in random environments.

%%%%%%%%%%%%%%%%%%%%%%%%%% NEW SECTION %%%%%%%%%%%%%%%%%%%%%%%%%%

\section{Observables and main results}

To distinguish the RSTRE for different choices of random environments and disorder strength $\beta$, we study properties of the resulting trees in the hope that they differ from each other. In this thesis, the observables we focus on will be the (expected) edge overlap, the local limit, and the diameter of the tree. 

\subsection*{Edge overlap and length}

Consider a graph with $n$ vertices, with edge weights $\weight(e) = \exp(-\beta \omega_e)$ induced by the random environment. Now sample two independent trees $\cT$ and $\cT'$ according to the weighted UST measure in \eqref{eq:PomegaT}. How many edges do the trees $\cT$ and $\cT'$ share, either typically or in expectation? We call this quantity the edge overlap, and we shall study the expectation (over the trees)
\begin{equation}
    \mathcal{O}(\beta) := \bE^\omega_{G, \beta} \otimes \bE^\omega_{G, \beta}\big[ | \cT \cap \cT'| \big],
\end{equation}
where $\cT$ and $\cT'$ are sampled independently under the same law $\bP^\omega_{G, \beta}$. For two MSTs under the same environment, the overlap is trivially equal to $n-1$, whereas for the unweighted UST on the complete graph, simple calculations show that the expected overlap is $2 (n-1)/n$. The following is one of our results for the RSTRE with disorder strength $\beta$, and is discussed in Chapter \ref{ch:local}.

\begin{theorem*}[cf.\ Theorem \ref{T:Overlap} and {\cite[Theorem 1.1]{Mak24}}]
    Let $K_n$ be the complete graph and assume that $(\omega_e)_{e \in E}$ are independently uniformly distributed on $[0,1]$. Then for $\beta = \beta(n)$ the disorder strength, we have with high probability that
    \begin{equation*}
        \mathcal{O}(\beta) = 
        \begin{dcases}
            (1 + o(1)) \beta \frac{e^{2\beta} - 1}{(e^\beta - 1)^2} & \text{ if } 0 < \beta \ll n/\log n, \\
            (1 - o(1)) n & \text{ if } \beta \gg n (\log n)^2.
        \end{dcases}
    \end{equation*}
\end{theorem*}
\noindent We remark that if $\beta$ is much larger than $1$ but much smaller than $n/\log n$, then the expected overlap is typically of order $\beta$. Once $\beta$ is much larger than $n$ times a log factor, most of the edges of the RSTRE agree with the MST.

Closely related to the edge overlap is the length 
\begin{equation*}
    L(\cT) = H(\cT, \omega) = \sum_{e \in \cT} \omega_e,
\end{equation*}
as first studied by Kúsz in \cite{K24} for the RSTRE. Using our techniques developed for the edge overlap, we will be able to prove the following.

\begin{theorem*}[cf. Theorem \ref{T:length} and {\cite[Theorem 1.4]{Mak24}}]
Let $K_n$ be the complete graph and assume that $(\omega_e)_{e \in E}$ are independently uniformly distributed on $[0,1]$. Then for $\beta = \beta(n)$, we have with high probability (and in $\P$-expectation) that
\begin{equation*}
        \bE^\omega_{n,\beta}\big[ L(\cT) \big] = 
        \begin{dcases}
             (1 + o(1)) \frac{n}{\beta} \cdot \frac{1 - \beta e^{-\beta} - e^{-\beta}}{1- e^{-\beta}} & \text{ if } 0 < \beta \ll n/\log n, \\
            (1 + o(1)) \zeta(3) & \text{ if } \beta \gg n (\log n)^5,
        \end{dcases}
    \end{equation*}
    where $\zeta(3) = \sum_{k=1}^\infty k^{-3} = 1.202\dots$ is Apéry's constant. 
\end{theorem*}
\noindent We remark that the value $\zeta(3)$ is exactly the limiting weight of the MST, see \cite{Fri85, Jan95}. 

\subsection*{Local limit}

Assume that we have sampled a tree $\cT$ and now fix a root $o \in \cT$ uniformly at random from all the vertices. The local limit studies what neighborhoods of finite radius $r \in \N$ look like when viewed from the root $o$. More precisely, we study the probabilities
\begin{equation*}
    \widehat{\P}( B_\cT(o, r) \simeq t) \qquad r \in \N,
\end{equation*}
for finite trees $t$ of height $r$. In Chapter \ref{ch:local}, we will prove the following theorem that says there is a sharp transition for the local limit when $\beta = n^\gamma$ and $\gamma$ crosses the critical threshold $\gamma_c = 1$.
\begin{theorem*}[cf.\ Theorem \ref{T:local_limit} and {\cite[Theorem 1.3]{Mak24}}]
    Let $K_n$ be the complete graph and assume that $(\omega_e)_{e \in E}$ are independently uniformly distributed on $[0,1]$. Then for $\beta = \beta(n)$ with $\beta \ll n/\log n$, the local limit of the RSTRE is the same as the local limit of the unweighted UST. If on the other hand $\beta = \beta(n) \gg n (\log n)^\lambda$, for $\lambda = \lambda(n) \rightarrow \infty$, then the local limit of the RSTRE agrees with the local limit of the MST.
\end{theorem*}

\subsection*{Diameter}

Any finite tree can be regarded as a metric space on the vertices by letting the distance between two vertices be the length of the (unique) path in the tree connecting the vertices. The diameter measures the maximum distance between any two vertices and is used to rescale the metric space to obtain a scaling limit. It turns out that on the complete graph the diameters of the UST and MST scale as $n^{1/2}$ and $n^{1/3}$, respectively. Large parts of this thesis are aimed at finding the correct exponent for the diameter of the RSTRE in terms of $\beta$. Namely, we conjecture the following smooth interpolation between the UST and MST on the complete graph.
\begin{conjecture*}[cf. Conjecture \ref{C:Intermediate} and {\cite[Conjecture 1.3]{MSS24}}]
    Let $K_n$ be the complete graph and denote by $\cT^\omega_{n, \beta}$ the RSTRE with $(\omega_e)_{e \in E}$ being independently uniformly distributed on $[0,1]$. Then with high probability as $n \to\infty$, we have that
	\begin{equation*}
	{\diam}(\cT^\omega_{n, \beta})=
	\begin{dcases}
	n^{\frac{1}{2}+o(1)} & \beta = O(n), \\
	n^{\frac{1}{2} - \frac{\gamma}{2}+o(1)} & \beta  = n^{1 + \gamma+o(1)}, \quad 0 < \gamma < \frac{1}{3}, \\
	n^{\frac{1}{3}+o(1)} & \beta =  \Omega(n^{\frac{4}{3}}).
	\end{dcases}
	\end{equation*}
\end{conjecture*}
The cases $\beta \ll n/\log n$ and $\beta \geq n^{4/3} \log n$ are covered in Theorem \ref{T:main} (cf. \cite[Theorem 1.1]{MSS24}) of Chapter \ref{ch:diam}, while Chapter \ref{ch:conj} discusses the conjectured intermediate regime. We emphasize here that we expect to see a smooth transition in terms of the diameter compared to the sharp transition as seen in the local limit. That is, local observables do not give us enough information to discern global properties.

%%%%%%%%%%%%%%%%%%%%%%%%%% NEW SECTION %%%%%%%%%%%%%%%%%%%%%%%%%%

\section{Outline of the thesis}

In Chapter \ref{ch:prelims}, we will give a brief introduction to electric networks and their relation to random walks and USTs. Furthermore, we will describe a few useful properties of the UST measure in Section \ref{S:UST_properties}, discuss mixing time and bottleneck ratio of a random walk in Section \ref{S:Mixing}, and give a coupling of the RSTRE and MST processes with an \ErdosRenyi random graph in Section \ref{S:ERintro}. Chapter \ref{ch:Zd} will discuss the RSTRE on the Euclidean lattice $\Z^d$ viewed from a statistical physics standpoint.

\smallskip

In Chapter \ref{ch:local}, we briefly review the notion of local convergence, and then we give sharp asymptotics for the edge overlap and length of the RSTRE and sketch the proof of the local convergence for $\beta$ in the low ($\beta = n^\gamma$, $\gamma <1$) and high disorder ($\beta = n^\gamma$, $\gamma > 1$) regimes. This will mainly be a summary of the preprint:
\begin{quote}
    \cite{Mak24} \citeauthor{Mak24} \citetitle{Mak24} \citeyear{Mak24}, available online at \href{https://arxiv.org/abs/2410.16836}{https://arxiv.org/abs/2410.16836}.
\end{quote}

\smallskip

\sloppy %
Next, in Chapter \ref{ch:diam}, we study the diameter of the RSTRE first for bounded-degree graphs, then for the low ($\beta = n^\gamma$, $\gamma <1$) and high disorder ($\beta = n^\gamma$, $\gamma > 4/3$) regimes on the complete graph, where the high disorder regime differs from that of the local limit. We only sketch out our proof arguments and we refer to the following two papers for the full versions of the arguments: %
\begin{quote}
    \cite{MSS23} \citeauthor{MSS23} \citetitle{MSS23} \citeyear{MSS23}, available online at \href{https://arxiv.org/abs/2311.01808}{https://arxiv.org/abs/2311.01808}.
\end{quote}
\begin{quote}
    \cite{MSS24} \citeauthor{MSS24} \citetitle{MSS24} \citeyear{MSS24}, available online at \href{https://arxiv.org/abs/2410.16830}{https://arxiv.org/abs/2410.16830}. 
\end{quote} %
We remark that the results of Chapter \ref{ch:diam} include many slight improvements over those in \cite{MSS23} and \cite{MSS24}. This chapter represents the main contribution of this thesis and aims to provide initial steps that might hopefully contribute to establishing a scaling limit of the RSTRE in different regimes.

\fussy %

\smallskip

Lastly, in Chapter \ref{ch:conj}, we give an argument that reduces the study of the spanning tree on the complete graph to that of a spanning tree on a slightly supercritical giant component. The goal of this chapter is to give heuristics for Conjecture \ref{C:Intermediate}, and as a byproduct, we will prove that the diameter of the UST on the giant component, say with $m$ many vertices, of a slightly supercritical graph is of order $\sqrt{m}$, up to some polylogarithmic errors. The contents of Chapter \ref{ch:Zd} and Chapter \ref{ch:conj} have not appeared in the literature.

\subsection{Asymptotic notation}

As most of our results only hold along a sequence of graphs $G_n$, we briefly recall some standard notation. For two functions $f$ and $g$, we say that $f = O(g)$ if there exists some constants $C > 0$ and $n_0$ such that for $n \geq n_0$ it holds that $f(n) \leq C g(n)$. Conversely, we write $f = \Omega(g)$ if $f(n) \geq C g(n)$ for some $C > 0$ and all $n \geq n_0$, for some $n_0 \in \N$. Alternatively, we also write $f \lesssim g$ (or $f \gtrsim g$) if we have $f = O(g)$ (or $f = \Omega(g)$). If both $f = O(g)$ and $f = \Omega(g)$ hold, then we denote this as either $f = \Theta(g)$ or $f \asymp g$, meaning the order of the two functions is the same. Furthermore, we will write $f \ll g$ (and sometimes $f = o(g)$) if $\lim\limits_{n \rightarrow \infty} f(n)/g(n) = 0$ and use an analogous notation for $f \gg g$. Lastly, we say that a sequence of events $A_n$ holds with high $\P_n$-probability if, under a sequence of laws $\P_n$, we have $\limn \P_n(A_n) = 1$.
\chapter{UST, electric networks and MST} \label{ch:prelims}

In this chapter, we gather some preliminary material. In Section \ref{S:electric}, we recall some results about electric networks and random walks; in Section \ref{S:Wilson}, we describe Wilson's algorithm to construct the UST; in Section \ref{S:UST_properties}, we state some properties about the UST measure, especially in relation to electric networks; Section \ref{S:Mixing} covers the bottleneck ratio and its application to the mixing time of the random walk. Lastly, in Section \ref{S:ERintro}, we collect several results about \ErdosRenyi random graphs and define a coupling of the RSTRE and MST to it.

%%%%%%%%%%%%%%%%%%%%%%%%%% NEW SECTION %%%%%%%%%%%%%%%%%%%%%%%%%%

\section{Random walks and electric networks } \label{S:electric}

In this section, we introduce the definition of electric networks and we will mostly follow the exposition of Chapter 2 in \cite{LP16}. Let $(G,\weight)$ be an electric network consisting of a finite connected graph $G=(V, E)$ and non-negative weights (also called conductances) $\weight = (\weight(e))_{e \in E}$ indexed by the undirected edges $e = (x,y) = (y,x) \in E$. Occasionally, we will write $V(G)$ and $E(G)$ for the vertex and edge sets of $G$, respectively.  In our setting, we shall choose weights given by 
\begin{equation*}
\weight(e) = \weight_\beta(e) = \exp(-\beta \omega_e),
\end{equation*}
where $(\omega_e)_{e \in E}$ are i.i.d.\ random variables with common distribution $\mu$. 
On $(G,\weight)$ we define a lazy random walk with law $Q$, by letting the transition probabilities be
\begin{equation*}
    \forall u,v \in V, \quad q(u,v) =
    \begin{dcases}
        \frac{1}{2} & \text{if } u=v, \\
       \frac{1}{2} \frac{\weight(u,v)}{\sum_{x \sim u} \weight(u,x)} =: \frac{1}{2} \frac{\weight(u,v)}{\weight(u)} & \text{if } u\neq v,
    \end{dcases}
\end{equation*}
with corresponding $t$ step probabilities $q_t(\cdot, \cdot)$. Here $x \sim u$ denotes that there is an edge $(x,u)$ in $G$ between $x$ and $u$, and for a vertex $u \in V$ we write $\weight(u) := \sum_{x \sim u} \weight(u,x)$. When the random walk starts at a vertex $x \in V$, we denote the law by $Q_x = Q^\omega_x$. We remark that taking the lazy version of the random walk avoids periodicity issues and will not change the law of UST on $(G, \weight)$ as generated by the random walk (see Section \ref{S:Wilson}). The stationary distribution $\pi$ satisfies
\begin{equation*}
    \pi(v) = \frac{\weight(v)}{2 \sum_{e \in E} \weight(e)} \qquad \text{and} \qquad \pi(S) := \sum_{v \in S} \pi(v) = \frac{\sum_{v \in S}\weight(v)}{2 \sum_{e \in E} \weight(e)} .
\end{equation*}
We let $\pi_{\mathrm{min}} := \min_{v \in V} \pi(v)$ and $\pi_{\mathrm{max}} := \max_{v \in V} \pi(v)$ be the minimum and maximum values of the stationary distribution.

\begin{definition} \label{D:effR}
    For disjoint sets $A, B \subset V$, we define the effective resistance between $A$ and $B$ by 
\begin{equation} \label{eq:def_effR}
    \effR{G, \weight}{A}{B} := \frac{1}{\sum_{a \in A} \weight(a) Q_{a}(\tau_B < \tau_A^+)}, 
\end{equation}
where $Q_a$ is the law of the lazy random walk started at $a$, $\tau_B$ is the first hitting time of $B$, and $\tau_A^+$ is the first return time to $A$. When $A=\{a\}$ and $B = \{b\}$ are singletons, we will just write
$\effR{G, \weight}{a}{b}$ instead. When the underlying graphs and weights are clear we will often drop either one or both of the superscripts.
\end{definition}

\subsection{Flows}

There is an equivalent definition of the effective resistance using the notion of flows. For each edge $e \in E$ fix an arbitrary orientation and let $e^-$ and $e^+$ be the starting and end points of $e$. We denote by $E^\rightarrow$ the set of oriented edges where we include each edge with both orientations $(e^-, e^+)$ and $(e^+, e^-)$. A flow $\theta : E^\rightarrow \rightarrow \R$ from $A \subset V$ to $B \subset V$, with $A \cap B = \emptyset$, is a function on oriented edges that is
\begin{enumerate}
    \item anti-symmetric:
    \begin{equation*}
        \theta(x,y) = - \theta(y,x) \quad \forall x,y \in V,
    \end{equation*}
    \item has zero net outgoing flow for vertices not in $A \cup B$:
    \begin{equation*}
        \sum_{y \sim x} \theta(x,y) = 0, \quad \forall x \not\in A \cup B,
    \end{equation*}
    \item and has non-negative (resp.\ non-positive) flow out of A (resp.\ B):
        \begin{equation*}
        \sum_{y \sim a} \theta(a,y) \geq 0 \quad \forall a \in A, \quad \qquad \sum_{y \sim b} \theta(b,y) \leq 0 \quad \forall b \in B. 
    \end{equation*}
\end{enumerate}
The strength $|\theta|$ of a flow is defined as
\begin{equation*}
    |\theta| := \sum_{a \in A} \sum_{y \sim a} \theta(a,y) \geq 0,
\end{equation*}
and we say $\theta$ is a unit flow if $|\theta| = 1$.

Given disjoint vertex sets $A, B \subset V$, let $v(\cdot)$ be the unique harmonic function (called \textit{voltage}) on the vertices with boundary values $v|_A \equiv 1$ and $v|_B \equiv 0$. That is, for non-boundary vertices $x \not\in A \cup B$
\begin{equation*}
    v(x) = \sum_{y \sim x} \frac{\weight(x,y)}{\weight(x)} v(y).
\end{equation*}
Given $v(\cdot)$, let $i : E^\rightarrow \rightarrow \R$ be the associated \textit{current} satisfying 
\begin{equation*}
    i(x,y) = \weight(x,y) \big( v(x) - v(y) \big) \qquad \forall (x,y) \in E.
\end{equation*}
It is an easy exercise to check that the current is indeed a flow. As the function $Q_{v}(\tau_A < \tau_B)$ is the unique harmonic function with the required boundary conditions, one may verify that for the associated current it holds that
\begin{equation} \label{eq:def_effR_current}
    \effR{}{A}{B} = \frac{1}{|i|}.
\end{equation}
Lastly, we define the energy $\Vert\theta\Vert^2_r$ of a flow as 
\begin{equation}
    \Vert\theta\Vert^2_r : = \frac{1}{2} \sum_{e \in E} \frac{1}{\weight(e)} \theta(e)^2.
\end{equation}

\subsection{Properties of the effective resistance}

One main advantage of using the flow interpretation of the effective resistance in \eqref{eq:def_effR_current} over the random walk definition in \eqref{eq:def_effR}, is that it is often easy to provide tight upper bounds for the effective resistance by providing a suitable flow. The following theorem formalizes this, and we refer to \cite[Section 2.4]{LP16} for a proof.
% - Exercise 2.13 in \cite{LP16}
\begin{theorem}[Thompson's principle] \label{T:thompson}
    Let $A, B \subset V$ be two disjoint sets of vertices. Then
    \begin{equation*}
        \effR{}{A}{B} = \min \left\{ \Vert \theta \Vert^2_r \, : \, \theta \text{ is a unit flow from } A \text{ to } B \right\}.
    \end{equation*}
\end{theorem}

Next, we give two more tools for bounding the effective resistance.

\begin{theorem}[Rayleigh’s monotonicity principle] \label{T:Rayleigh}
Let $\weight, \weight'$ be weights on $G$ with $\weight(e) \leq \weight'(e)$ for all edges $e \in E$. If $A,B \subset V$ are two disjoint sets, then
\begin{equation*}
    \effR{\weight'}{A}{B} \leq \effR{\weight}{A}{B},
\end{equation*}
i.e.\ the effective resistance is a non-increasing function of the weights.
\end{theorem}

\begin{lemma}[Nash-Williams inequality] \label{L:nash_williams}
    Let $\Pi_1, \ldots \Pi_k \subset E$ be disjoint cutsets of $A$ and $B$, i.e.\ any path between $A$ and $B$ uses at least one edge in $\Pi_i$ for every $i \in \{1, \ldots, k\}$. Then
    \begin{equation}
        \effR{}{A}{B} \geq \sum_{i=1}^k \Big( \sum_{e \in \Pi_i} \weight(e) \Big)^{-1}.
    \end{equation}
\end{lemma}
% \begin{proof}
%     This is simply the Nash-Williams inequality with only one cutset. See e.g.\ Section 2.5 of \cite{LP16}.
% \end{proof}

Another method of estimating the effective resistance is to consider a (smaller) simpler graph with the same effective resistances via network reductions. The two main transforms are the series law and parallel law, but we note the existence of another (family) of transforms called $Y$-$\Delta$ transformations.

\begin{lemma}[Series law] \label{L:SeriesLaw}
   Suppose $v \in V$ is a vertex of degree $2$ with incident edges $e_1 = (u_1, v)$ and $e_2 = (v, u_2)$ that have resistances $r_1 := 1/\weight(e_1)$ and $r_2 := 1/\weight(e_2)$, respectively. Consider the  graph $G'$
   on $V \setminus \{v\}$  where $e_1$ and $e_2$ are replaced by a single edge $e$ between $u_1$ and $u_2$, with resistance $r(e) := r_1 + r_2$. Then, for any disjoint $A,B\subseteq V\setminus\{v\}$, the effective resistance between $A$ and $B$ in $G$ is equal to that in $G'$.
\end{lemma}

\begin{lemma}[Parallel law]
    Suppose $u,v \in V$ are such that there exist two parallel edges $e_1$ and $e_2$ between $u$ and $v$ in $G$. Consider a new graph $G'$ by replacing $e_1$ and $e_2$ by a single edge $e$ with weight $\weight(e) = \weight(e_1) + \weight(e_2)$. Then the effective resistance in $G$ between $A$ and $B$ in $G$ is equal to that in $G'$.
\end{lemma}

%%%%%%%%%%%%%%%%%%%%%%%%%% NEW SECTION %%%%%%%%%%%%%%%%%%%%%%%%%%

\section{UST and electric networks} \label{S:UST_properties}

Recall that given weights $\weight(e)$ on the edges, the weighted UST measure $\bP^\weight_G = \bP$ on the set of spanning trees $\T(G)$ is defined by
\begin{equation} \label{eq:defUST}
    \bP(\mathcal{T} = T) = \frac{1}{Z^{\weight}_G} \prod_{e \in T} \weight(e), \qquad Z^{\weight}_G:= \sum_{T\in \T(G)} \prod_{e \in T} \weight(e).
\end{equation}
In this thesis, the environment $\omega = (\omega_e)_{e \in E}$ consists of i.i.d.\ random variables with common distribution $\mu$, and the edge weights are given by $\weight(e) = \exp(-\beta \omega_e)$ for some $\beta \geq 0$. Consequently, the spanning tree measure can be parameterized as
\begin{equation} \label{eq:defRSTRE}
    \bP(\mathcal{T} = T) = \bP^\omega_{G, \beta}(\mathcal{T} = T) = \frac{1}{Z^{\omega}_{G, \beta}} \prod_{e \in T} \exp( - \beta \omega_e) = \frac{1}{Z^{\omega}_{G, \beta}} \exp \big( - \beta \sum_{e\in T} \omega_e \big),
\end{equation}
with corresponding normalization constant $Z^{\omega}_{G, \beta}$. As we shall see in Section \ref{S:Wilson}, there is a strong relationship between random walks, electric networks, and USTs. We now recall some basic properties of the weighted UST measure and refer to Chapter 2 and Chapter 4 of \cite{LP16} for more details.

\begin{theorem}[Kirchhoff's Formula] \label{T:Kirchhoff}
Given any finite electric network $(G, \weight)$ and corresponding weighted UST measure $\bP^\weight_G$, we have
\begin{equation} \label{eq:p_in_tree_r_eff}
    \bP^\weight_G \big( (u,v) \in \mathcal{T} \big) = \weight(u,v) \effR{G, \weight}{u}{v}.
\end{equation}
\end{theorem}
\begin{definition}[Transfer-impedance matrix]
    Recall that we fixed an (arbitrary) orientation of the edges with $e^-$ and $e^+$ denoting the starting and ending points of an edge $e$. Let $v_e(\cdot)$ be the voltage on the vertices when a unit current flow $\theta_e$ is sent from $e^-$ to $e^+$, i.e.\ $v_e(\cdot)$ is the (unique) harmonic function on $V \setminus \{e^-, e^+\}$ with $v_e(e^-) = \effR{}{e^-}{e^+}$ and $v_e(e^+) = 0$. The transfer-impedance matrix $Y$ is a $|E| \times |E|$ matrix with entries given by
    \begin{equation*}
        Y(e,f) := w(f) \big[ v_e(f^-) - v_e(f^+) \big] = \theta_e(f^-, f^+).
    \end{equation*}
\end{definition}

Note that we have
\begin{equation*}
    Y(e,f) \weight(e) = Y(f,e) \weight(f),
\end{equation*}
implying that $Y(e,f)$ is symmetric in the unweighted case. The following theorem is due to \cite{BP93}.

\begin{theorem}[Transfer-impedance theorem] \label{T:transfer-impedance}
    For distinct edges $e_1, \ldots, e_k$
    \begin{equation}
        \bP(e_1, \ldots, e_k \in \cT) = \det \big( Y(e_i, e_j)_{1 \leq i,j\leq k} \big),
    \end{equation}
    where $Y(e_i, e_j)_{1 \leq i,j\leq k}$ are entries of the transfer-impedance matrix $Y$.
\end{theorem}

% Section 4.2 in lp - E -

A direct consequence of the transfer-impedance theorem is that for distinct edges $e$ and $f$
% \begin{equation}
%     \bP(e \in \cT) \bP(f \in \cT) - \bP(e,f \in \cT) = \frac{w(e)}{w(f)}Y(e,f)^2. 
% \end{equation}
\begin{equation} \label{eq:Yef_as_UST}
    Y(e,f) = \pm \sqrt{\frac{w(f)}{w(e)}} \sqrt{\bP(e \in \cT) \bP(f \in \cT) - \bP(e,f \in \cT)} \, .
\end{equation}
This means that the UST measure is determined, up to signs, by its $1$ and $2$ dimensional marginal distribution on the edges. This will be used to prove the convergence of the local limit in Chapter \ref{ch:local}.

\subsection{Spatial Markov property}

We now introduce another useful property of the UST measure, known as the spatial Markov property, and refer to \cite[Sec.\ 2.2.1]{HN19} for further details. This will be one of the main tools we use in Chapter \ref{ch:diam}, and in particular in Section \ref{S:toy_model}. The deletion of an edge $e \in E$ produces a new graph $G - \{e\}$, formed by removing $e$ from the edge set. Similarly, the contraction of $e$ results in the graph $G / \{e\}$, which is obtained by first removing $e$ and then contracting its endpoints into a single vertex. The deletion and contraction operations commute, meaning that for disjoint sets of edges $A, B \subset E$, we may write $(G - B) / A$ for the graph that is obtained by first removing the edges in $B$ and then contracting those in $A$.

Notice that $(G - B) / A$ may contain loops and multiple edges between vertices. However, the weighted UST measure $\bP^{\weight}_{(G-B)/A}$ can still be defined similarly to \eqref{eq:defUST}, by retaining the original edge weights for edges that are not removed, and by considering all spanning trees in $\mathbb{T}{(G - B) / A}$. In this case, multiple edges between two vertices are treated as distinct edges in different spanning trees. 

\begin{lemma}[Spatial Markov Property] \label{L:USTmarkov}
	Let $(G, \weight)$ be a finite connected weighted graph and let $A, B\subset E$ be two disjoint sets of edges such that $\bP^\weight(A\subset \cT, B\cap \cT=\emptyset)>0$. Then for any set of edges $F\subset E$,
	\begin{equation*}
	    \bP^\weight_G (\cT = F \mid A\subset \cT, B\cap \cT =\emptyset) = \bP^{\weight}_{(G-B)/A}(\cT \cup A = F).
	\end{equation*}
\end{lemma}
\noindent
Notably, conditioned on $\mathcal{T}$ containing all edges in $A$ and none in $B$, the law of $\bP^\weight_G$ on the remaining edges in $E \setminus (A \cup B)$ is given by the weighted UST measure on $(G - B) / A$.

\begin{remark} \label{R:correlation}
    The spatial Markov property together with Rayleigh's monotonicity principle also implies a negative correlation of edge occupancies, as for two edges $e \neq f$ we have
    \begin{equation*}
        \bP_G(e, f \in \cT) = \bP_{G / f}(e \in \cT) \bP_G(f \in \cT) \leq \bP_{G}(e \in \cT) \bP_G(f \in \cT).
    \end{equation*}
    The same is however not always true for the MST measure and the averaged law $\widehat{\P}$. For instance, in Section 5 of \cite{LPS03}, a (multi) graph with 4 vertices is constructed where there are two edges $e,f$ such that
    \begin{equation*}
        \P(e,f \in \textrm{MST}) > \P(e \in \textrm{MST}) \P(f \in \textrm{MST}).
    \end{equation*}
    As the number of vertices is fixed, the averaged RSTRE law converges to the MST law for $\beta \rightarrow \infty$ whenever $\mu$ is atomless. The positive correlation thus also holds for the averaged law if $\beta$ is large enough and, say, $\mu$ is the uniform distribution on $[0,1]$.  
\end{remark}

\begin{question}
    Assume that $e,f \in E(G)$ are such that the edge occupancies of $e$ and $f$ are positively correlated in the random MST. For what values of $\beta$ and choices of $\mu$ does the same hold for the averaged law of the RSTRE?
\end{question}

\section{Sampling USTs} \label{S:Wilson}

As mentioned in the introduction, there are several methods for generating a sample of a (weighted) uniform spanning tree. For instance, the Aldous-Broder algorithm (see, e.g., \cite[Corollary 4.9]{LP16}) starts a random walk from some vertex and constructs the tree by including edges that are used to visit a vertex for the first time. This process continues until every vertex is visited. A related approach, called the Aldous-Broder interlacement process \cite{Hut18}, considers random walk trajectories arriving in time as a Poisson point process. Although the main part of this thesis will not use either of these constructions, we note that they are an important ingredient in the proof of Theorem \ref{T:Nach} in \cite{MNS21}, which is one of our main tools for establishing bounds on the diameter of the RSTRE.

\subsection*{Wilson's Algorithm}

A computationally more efficient algorithm is Wilson's algorithm, which uses \textit{loop erased random walks} (LERW) to generate paths in the (weighted) UST. The loop erasure of a trajectory $\Vec{X} = (X_0, \ldots, X_m)$ is the path in which cycles of $\Vec{X}$ are removed in chronological order. More precisely, start at $X_0$ and let $0 \leq i \leq m$ be the largest index such that $X_i = X_0$. Then, set $Y_0 = X_i = X_0$. Next, consider the largest index $j$, with $i+1 \leq j \leq m$ (if it exists), such that $X_j = X_{i+1}$ and set $Y_1 = X_j$ so that $(Y_0, Y_1) = (X_i, X_{i+1}) = (X_0, X_j)$. Repeating this procedure until $Y_k = X_m$ gives the loop erasure $Y =: \LE (\Vec{X})$. %,we also refer to Appendix \ref{ch:LERWapp} for more information. 
If $X$ is a random walk, then we call the law of $Y = \LE (X)$ the law of a loop erased random walk.

We now briefly explain Wilson's algorithm. Fix an arbitrary ordering\footnote{The ordering of the vertices may be dynamically adjusted between different steps of the algorithm without affecting the law of the UST.} $\{v_1, \ldots, v_n\}$ of the vertices and let $v_1$ be the root. The initial tree consists only of the root and we set $\cT(1) = (\{v_1\}, \emptyset)$. In the next step, run a random walk $X^{2}$ started at $v_2$ until $X^{2}$ hits $\cT(1)$. As $G$ is finite and connected, this happens after a finite number of random walk steps. Update the tree by setting $\cT(2)$ as the union of $\cT(1)$ and the trajectory (including both vertices and edges) of the loop erased trajectory $\LE(X^{2})$. Repeat this procedure for each vertex $v_i$, where for the $i$-th step we run a random walk $X^{i}$ started at $v_i$ until $X$ hits $\cT(i-1)$ and update the tree by letting $\cT(i)$ be the union of $\cT(i-1)$ and the (possibly empty if $v_i \in \cT(i-1)$) path $\LE(X^{i})$. After all steps, the algorithm outputs a spanning tree $\cT = \cT(n)$. The following theorem is due to Wilson and first appeared in \cite{Wil96}.

\begin{theorem}
Let $(G, \weight)$ be a finite connected weighted graph. Wilson's Algorithm samples a spanning tree $\mathcal{T}$ according to the weighted UST measure $\bP^\weight(\mathcal{T})$.
\end{theorem}

\subsection*{Wilson's algorithm rooted at infinity}

Suppose now the graph $G = (V,E)$ is infinite with locally finite weights $(\weight(e))_{e \in E}$, i.e.\ $\weight(v) < \infty$ for all $v \in V$. If $G$ is recurrent, then the length of a random walk started at some vertex $v_1$ and ended at another vertex $v_2$ will almost surely be finite. Therefore, we may define Wilson's algorithm on this graph similarly to the finite graph case. However, if $G$ is transient, then the random walk may never hit $v_2$. Nonetheless, in \cite{BLPS01}, Wilson's algorithm was adapted to infinite transient graphs by ``treating infinity as the root''. We defer the exact definitions of the limiting uniform spanning forest measure to Chapter \ref{ch:Zd}.

Consider again some arbitrary ordering $V = \{v_1, v_2, \ldots \}$ of the vertices. Set $\cT(0) = (\emptyset, \emptyset)$ and think of $\cT(0)$ containing only the ``point at infinity''. Inductively, for each $i \geq 1$ run a random walk $X^i$ started at $v_i$ and stopped if $X^i$ hits $\cT(i-1)$. As the graph is transient, there is a positive probability that $X^i$ never hits $\cT(i-1)$, however, as each vertex is visited only a finite number of times, the loop erasure $\LE(X^i)$ is well defined. We then set $\cT(i)$ as the union of $\cT(i-1)$ with the possibly infinite length path $\LE(X^i)$, and define $\cT = \cup_{i \in \N} \cT(i)$. We call this procedure \textit{Wilson's algorithm rooted at infinity} and refer to \cite[Proposition 10.1]{LP16} or \cite{BLPS01} for the following theorem.

\begin{theorem} \label{T:wilson_infty}
    If $(G, \weight)$ is a transient weighted graph, then the forest produced by Wilson's algorithm rooted at infinity is distributed as a (wired) uniform spanning forest on $(G, \weight)$.
\end{theorem}

%%%%%%%%%%%%%%%%%%%%%%%%%% NEW SECTION %%%%%%%%%%%%%%%%%%%%%%%%%%

%%%%%%%%%%%%%%%%%%%%%%%%%% NEW SECTION %%%%%%%%%%%%%%%%%%%%%%%%%%

\section{Bottleneck ratio and mixing time} \label{S:Mixing}

We now recall a few definitions related to random walks that will play a crucial role in our proofs. Given weights $\weight = (\weight(e))_{e \in E}$ on the edges, define the {\em bottleneck ratio} of $S \subset V$ by
\begin{equation*}
    \Phi_{(G, \weight)}(S)=  \frac{\sum_{x\in S, y\in S^c} \pi(x) q(x,y)}{\pi(S)},
\end{equation*}
where $q(\cdot, \cdot)$ and $\pi(\cdot)$ are respectively the transition probabilities and stationary distribution of the lazy random walk on $(G, \weight)$ (or more generally of a Markov chain) introduced in Section \ref{S:electric}. In the literature, this quantity is sometimes also called the \textit{conductance} of $S$, and we remark that $\Phi_{(G, \weight)}(S)$ can also be written in the form of
\begin{equation}\label{eq:PhiGS}
    \Phi_{(G, \weight)}(S) := \frac{1}{2} \frac{\sum_{e \in E(S, S^c)} \weight(e)}{\sum_{e \in E(S,V) } \weight(e)},
\end{equation}
where $E(A,B)$ denotes the set of edges between the vertex sets $A$ and $B$, and the factor of $1/2$ arises from using the lazy version of the random walk. %\lu?ca{TO%DO: double-check constants coming from 1/2 factor... In toy model check the factor 1/2...}

\begin{definition}[Bottleneck ratio] \label{D:bottleneck}
    The bottleneck ratio of $(G,\weight)$ is defined as
    \begin{equation*}
          \Phi_{(G, \weight)} := \min_{0 < \pi(S) \leq \frac{1}{2}}\Phi_{(G, \weight)}(S).
    \end{equation*}
\end{definition}
\noindent In the unweighted case $\weight \equiv 1$, this is closely related to the notion of edge expansion.

\begin{definition}[Edge expansion and expander graphs] \label{D:expander}
    The edge expansion $h_G$ of a graph is defined as
    \begin{equation*}
        h_G := \min_{1 \leq |S| \leq n/2} h_G(S) = \min_{1 \leq |S| \leq n/2} \frac{|E(S,S^c)|}{|E(S,V)|},
    \end{equation*}
    where $E(A,B)$ denotes all edges between vertex sets $A$ and $B$.
    Given $b > 0$, we say a graph $G$ on $n$ vertices is a $b$-{\em expander} if
\begin{equation}
	h_G =  \min_{1 \leq |S| \leq n/2} \frac{|E(S,S^c)|}{|E(S,V)|} \geq b, \label{eq:expander}
\end{equation}
or equivalently if
\begin{equation*}
    |E(S,S^c)| \geq b \min \big\{ |E(S,V)|, |E(S^c, V)| \big\} \quad \forall S \subseteq V.
\end{equation*}
We further say that $(G_n)_{n \geq 1}$ is a sequence of $b$-expander graphs if each $G_n$ is a $b$-expander on $n$ vertices.
\end{definition}

We note that for bounded-degree graphs, as covered in \cite{MSS23}, occasionally the expansion $|E(S,S^c)|/|S|$ is used instead. In this work, we will only use the definition provided in \eqref{eq:expander}. Moreover, in the unweighted case $\weight \equiv 1$, the ratio in \eqref{eq:expander} is equal to that in \eqref{eq:PhiGS} up to a constant multiple of $1/2$ that is necessary due to random walk being lazy. The following is a version of the bottleneck ratio on a more refined scale.
% We remark that there are several slightly different definitions of expanders in the literature, for the current work it is easier to deal with the version defined in \eqref{eq:expander}.

\begin{definition}[Bottleneck profile] \label{D:bottleneck_profile}
\sloppy % 
Let $(G, \weight)$ be an electric network with weights $\weight = (\weight(e))_{e \in E}$ and consider the corresponding lazy random walk with stationary measure $\pi$. For $\pi_{\mathrm{min}} \leq r \leq 1/2$ define the quantity
\begin{equation*}
    \Phi_{(G, \weight)}(r) := \min\limits_{0 < \pi(S) \leq r} \Phi_{(G, \weight)}(S),
\end{equation*}
and for $r \geq 1/2$, let $\Phi_{(G, \weight)}(r) := \Phi_{(G, \weight)}(1/2) = \Phi_{(G, \weight)}$.
\fussy %
\end{definition}

\noindent There is a strong connection between the mixing time and the bottleneck profile of a graph which goes back to the results of \cite{SJ89}.

\begin{definition}[Mixing time]
    Define the (uniform) mixing time as
    \begin{equation*}
    \tmix(G, \weight) = \tmix := \min \left\{ t \geq 0 : \max\limits_{u,v \in V} \left| \frac{q_t(u,v)}{\pi(v)} - 1\right| \leq \frac{1}{2}\right\}.
    \end{equation*}
\end{definition}

We shall use the following theorem, which gives bounds both on the mixing time and on the transition probabilities, see also \cite[Chapter 6]{LP16}. 
\begin{theorem}[Theorem 1 in \cite{MP05}] \label{T:HeatCheeger}
    For $u,v \in V(G)$, let $q_t(u,v)$ be the $t$-step transition kernel of the lazy random walk on $(G, \weight)$. 
    For any $t>1$ and $\chi \in (0,1)$ such that 
        \begin{equation*}
        t \geq 1 + \int_{4 \min\{\pi(u), \pi(v)\}}^{\frac{4}{\chi}} \frac{4}{r \Phi_{(G, \weight)}^2(r)} dr,
    \end{equation*}
    we have 
        \begin{equation*}
        \Big| \frac{q_t(u,v)}{\pi(v)} -1 \Big| \leq \chi.
    \end{equation*}
    Furthermore,
    \begin{equation*}
        \tmix \leq 1 + \int_{4 \pi_{\mathrm{min}}}^{8} \frac{4}{r \Phi_{(G, \weight)}^2(r)} dr.
    \end{equation*}
\end{theorem}

\noindent The main application of this will be in for the proof of Theorem \ref{T:toy_model} and Theorem \ref{T:GeneralLow} in Chapter \ref{ch:diam}.

% \lu?ca{[JUST FOR US: This is actually a bit of an overkill. We do not necessairly need the isoperimetric profile, only the bottleneck ratio, i.e.\ we only use $\Phi(r)$ when $r=1/2$. There is an easier (than the paper above) argument using the bottleneck ratio that
% 	\begin{equation*}
% 		|\frac{p_t(x,y)}{\pi(y)} - 1 | \leq \frac{\exp(  -t \frac{1}{4?} \Phi(r)^2 ) }{\sqrt{\pi(x) \pi(y)}},
% 	\end{equation*}
% 	see e.g. equation 12.13 of \cite{LP16}. On the other hand, we already wrote everything down for the profile.... We would have to change a few lines in the proof. DISCUSS]}

%%%%%%%%%%%%%%%%%%%%%%%%%% NEW SECTION %%%%%%%%%%%%%%%%%%%%%%%%%%

\section{Random graphs} \label{S:ERintro} %\label{S:coupleRG}

For $p \in [0,1]$, let $G_{n,p}$ denote an \ErdosRenyi random graph on $n$ vertices obtained by independently keeping each edge of the complete graph $K_n$ with probability $p$ and removing it otherwise. If the edge $e$ is kept, we say it is $p$-open; otherwise, we say it is $p$-closed. We denote $\{u \xleftrightarrow{p} v\}$ as the event that two vertices $u$ and $v$ are connected by a path along which each edge is $p$-open.  This results in a set of disjoint, connected $p$-clusters, $\mathcal{C}_1(p), \mathcal{C}_2(p), \dots$, each viewed as a subgraph of $K_n$. We order these clusters by decreasing size $|\mathcal{C}_1(p)| \geq |\mathcal{C}_2(p)| \geq \ldots$, where $|\mathcal{C}_i(p)|$ represents the number of vertices in the $i$-th cluster. The largest cluster, $\mathcal{C}_1(p)$, is referred to as the \textit{giant component}, and bounds on its size play a key role in our proofs.

\subsection{Component sizes}

The seminal paper by Erd\H{o}s and R\'enyi \cite{ER60} describes the phase transition of $G_{n,p}$ around the critical point $p = 1/n$. Since then, the behavior near the critical point has become much better understood; see e.g.\ \cite{Luc90, NP07} for detailed results, or, more generally, either of the books \cite{FK16, vdH17} for broader results on \ErdosRenyi random graphs. We will often consider a parameter of the form $p = (1 + \epsilon)/n$, where $\epsilon:= \epsilon(n)$ is a sequence with absolute value tending to zero. Specifically, if
\begin{equation*}
   p = \frac{1 + \epsilon}{n} \, \text{ then we say } G_{n,p} \text{ is } \,
    \begin{dcases}
        \text{ critical} & \text{if } |\epsilon| = O(n^{-1/3}), \\
        \text{ slightly supercritical} & \text{if } \epsilon \gg  n^{-1/3}, \\
        \text{ slightly subcritical} & \text{if } -\epsilon \gg  n^{-1/3}.
    \end{dcases}
\end{equation*}

The component sizes for different choices of $\epsilon = o(1)$ have been well studied. We state below the qualitative versions for the sizes in the different regimes and remark that in our proofs we will implicitly use more quantitative bounds (see e.g. Section 2 of \cite{MSS24}). 

\begin{theorem}[Slightly subcritical]
    Let $p= \frac{1 - \epsilon}{n}$ with $\epsilon^3 n \rightarrow \infty$. Then for any $\ell \in \N$ and $\delta > 0$
    \begin{equation} \label{eq:s-sub-critical}
        \P \Bigg( \Big| \frac{|\cC_\ell(p)|}{2 \epsilon \log( \epsilon^3 n) } -1 \Big| \geq \delta \Bigg) \rightarrow 0.
    \end{equation}
\end{theorem}

\begin{theorem}[Slightly supercritical]
Let $p= \frac{1 + \epsilon}{n}$ with $\epsilon^3 n \rightarrow \infty$. Then for any $\delta > 0$
\begin{equation*}
        \P \Bigg( \Big| \frac{|\cC_1(p)|}{2 \epsilon n } -1 \Big| \geq \delta \Bigg) \rightarrow 0,
    \end{equation*}
    and $\cC_\ell$, for any $\ell \geq 2$, satisfies \eqref{eq:s-sub-critical}.
\end{theorem}

For the proof of the above theorems, we refer to either \cite{Luc90} or \cite{NP07}. For completeness' sake, we also list the component sizes in the critical, supercritical, and subcritical regimes, and refer to Chapters 4 and 5 of \cite{vdH17} for the proofs.

\begin{theorem}[Critical] \label{T:C1critical}
Let $p= \frac{1 + \epsilon}{n}$ with $\epsilon = \lambda n^{-1/3}$, $\lambda \in \R$. Then for all $\delta < 1$
\begin{equation*} 
        \P \Big( \delta n^{1/3} \leq |\cC_1(p)| \leq \delta^{-1} n^{1/3}  \Big) \geq 1 - C(\lambda) \delta,
\end{equation*}
where $C(\lambda) > 0$ is some constant.
\end{theorem}

\begin{theorem}[Supercritical]
    Let $p= \frac{\lambda}{n}$ with $\lambda > 1$. Then for any $\delta > 0$
    \begin{equation*}
        \P \Bigg( \Big| \frac{|\cC_1(p)|}{(1-\eta_\lambda) n } -1 \Big| \geq \delta \Bigg) \rightarrow 0,
    \end{equation*}
    where $\eta_\lambda > 0 $ satisfies $\eta_\lambda = \exp(\lambda(\eta_\lambda - 1))$.
\end{theorem}

\begin{theorem}[Subcritical]
    Let $p= \frac{\lambda}{n}$ with $\lambda < 1$. Then there exists a constant $C(\lambda) > 0$ such that
    \begin{equation*}
        \P \Big( C(\lambda)^{-1} \log n \leq |\cC_1(p)| \leq C(\lambda) \log n \Big) \geq 1 - o(1).
    \end{equation*}
\end{theorem}

In Chapter \ref{ch:local}, we will make use of the following theorem, which follows as a corollary of the (slightly) sub- and supercritical random graphs bounds as above.
\begin{theorem}[cf.\ {\cite[Theorem 2.5]{Mak24}}] \label{T:C2_size}
    There exists constants $L, \eta > 0$ such that for $\epsilon \geq n^{-1/4}$ with $\epsilon =o(1)$ and
    \begin{equation*}
        p \not\in \Big[ \frac{1 - \epsilon}{n}, \frac{1+\epsilon}{n}\Big],
    \end{equation*}
    we have that
    \begin{equation}
        \P \Big( \big| \mathcal{C}_2(p) \big| \geq L \log(n) \epsilon^{-2} \Big) = O(n^{-\eta}).
    \end{equation}
\end{theorem}

Finally, we recall the connectivity threshold for $G_{n,p}$. When $pn = \log n + c$, the probability that the random graph is connected approaches $1$ as $c$ increases. We will use the following theorem and refer to Theorem 9 in Chapter 7 of \cite{Bol98} for its proof.
\begin{theorem}
\label{T:Gnp_connect}
    Let $pn = \log n + \lambda(n)$ for $\lambda(n) \rightarrow \infty$. Then for $n$ large enough
    \begin{equation*}
         \P( G_{n,p} \text{ is not connected} \,) \leq 4 e^{-\lambda(n)} .
    \end{equation*}
\end{theorem}

\subsection{Percolation on general graphs} \label{SS:perc}

We may repeat the random graph process (also called percolation), as defined at the start of Section \ref{S:ERintro}, for arbitrary locally finite graphs $G = (V,E)$: for $p \in [0,1]$ each edge in $E$ is independently kept with probability $p$ and removed otherwise. We denote the resulting graph by $G_p$, and view $(G_p)_{p \in [0,1]}$ as a random graph process in terms of $p$. For finite graphs, this again gives rise to connected $p$-clusters $\cC_\ell(G, p)$, which we order in decreasing size. We will require the following lemma on the component sizes of subcritical percolation on finite degree graphs, and refer to e.g.\ Chapter 3 and Chapter 4.2 of \cite{vdH17} for an introduction to branching processes and their relation to random graphs.

\begin{lemma} \label{L:perc_subcrit}
    Let $G = (V,E)$, $|V| = n$, be a graph with maximum degree $\Delta \in \N$, and let $p = \frac{1-s}{\Delta - 1} $ for some $0 < s \leq 1$. Then there exists a constant $C = C(s) > 0$ such that for $k\geq 1$
    \begin{equation} \label{eq:perc_subcrit_eq}
		\P \big( |\mathcal C_1(G, p)| \geq C k \log n \big) \leq \frac{1}{n^k}  .
	\end{equation}
\end{lemma}
\begin{proof}
    Denote by $\cC(v,p)$ the $p$-open connected component of $v \in V$. The size of the following branching process stochastically bounds the size of $\cC(v,p)$: let $X_v - 1$ be distributed as a $\mathrm{Bin}(\Delta-1, p)$ random variable, attach $X_v$ many children to $v$, and from each child start an independent a Galton-Watson branching process with offspring distribution $\mathrm{Bin}(\Delta - 1, p)$. In other words, we create a Galton-Watson tree with offspring distribution $\mathrm{Bin}(\Delta - 1, p)$ with one additional branch attached to the root (the root may have $\Delta$ many children, while each child has at most $\Delta - 1$ further children). Denote by $T$ the size of such a process without the additional branch at the root, and further let $(X_i)_{i \geq 1}$ be i.i.d.\ $\mathrm{Bin}(\Delta - 1, p)$ random variables. If we let $S_t = X_1 + \ldots X_t$, then by a random walk representation of the branching process (see e.g.\ Chapter 3.3 or Chapter 4.2 of \cite{vdH17}) we have the stochastic domination
    \begin{equation*}
        \P( T > t) \leq \P( S_t > t-1).
    \end{equation*}
    As each $X_i$ is binomially distributed with the same parameter, we see that $S_t$ is distributed as a $\mathrm{Bin}(t(\Delta - 1), p)$ random variable. Using a standard large deviation inequality on binomial random variables (cf.\ \cite[Corollary 2.20]{vdH17}) gives that
    \begin{equation*}
        \P( T > t) \leq \P(S_t \geq t) \leq \exp \big(- t I(1-s) \big),
    \end{equation*}
    where for $a < 1$
    \begin{equation*}
        I(a) = a - 1 -  \log (a ) > 0.
    \end{equation*}
    
    By independence, the size of the additional branch attached to the root has the same distribution of $T$, so that
    \begin{equation*}
        \P \big( |\cC(v,p)| > 2t \big) \leq 2 \exp \big(- t I(1-s) \big) .
    \end{equation*}
    Using a union bound over all vertices then gives that 
    \begin{equation*}
       \P \Big( |\mathcal C_1(G, p)| \geq \frac{k}{I(1-s)} \log n \Big) \leq 2 n^{1 - k/2}, 
    \end{equation*}
    from which one can easily deduce \eqref{eq:perc_subcrit_eq}.
\end{proof}

\subsection{Coupling to random graphs} \label{SS:CoupleRG}

It will be useful to couple both the RSTRE and MST to a random graph process $(G_p)_{p \in [0,1]}$ on a finite graph $G$ as defined above. First, we couple the MST (denoted by $M$) to the random environment $\omega$ by letting $M$ be the tree that minimizes
\begin{equation*}
    H(T, \omega) = \sum_{e \in T} \omega_e 
\end{equation*}
over all $T \in \T(G)$. One can observe that the MST is the unique tree that maximizes the probability in \eqref{eq:PomegaT} (and \eqref{eq:defRSTRE}). 

Then, we couple the environment $\omega$ to $(G_{p})_{p \in [0,1]}$ by letting 
\begin{equation*}
    e \in G_{p} \iff \omega_e \leq p,
\end{equation*}
if $\mu$ is the uniform distribution on $[0,1]$, or more generally
\begin{equation*}
    e \in G_{p} \iff \omega_e \leq F^{-1}_\mu(p),
\end{equation*}
where $F_\mu$ is the cumulative distribution function (CDF) of $\mu$. The probability of $\omega_e$ being at most $F^{-1}(p)$ is for continuous distributions, of course, equal to $p$. 

Denote by $E(p)$ the $p$-open edges of $G_{p}$. One important aspect of this coupling, which can for instance be observed by using Kruskal's algorithm, is that the connected components of $M \cap E(p)$ are precisely the connected components of $G_{p}$. In other words, $M \cap E(p)$ is equal to $G_{p}$ with all cycles removed in some specific way. In particular, if $G = K_n$ is the complete graph and $G_p = G_{n,p}$ is almost cycle-free, say when $p = 1/n$, then $M \cap E(p)$ almost agrees with $G_{n,p}$. Another consequence is the following lemma. 

\begin{lemma} \label{L:max_deg_MST}
    Let $M$ be the random MST of the complete graph on $n$ vertices, then
    \begin{equation*}
        \P \big( \exists v \in V : deg_M(v) > 60 \log n \big) = O(n^{-4}).
    \end{equation*}
\end{lemma}
\begin{proof}
    This follows by Theorem \ref{T:Gnp_connect} and standard bounds on the degrees in $G_{n,p}$, see e.g. \cite[Lemma 2.7]{Mak24}.
\end{proof}

\chapter[RSTRE on the Euclidean lattice]{RSTRE on the Euclidean lattice} \label{ch:Zd}

The main goal of this chapter is to study the RSTRE on the $d$-dimensional Euclidean lattice $\mathbb{L}^d = (\mathbb{Z}^d,E^d)$. In particular, we are interested in the process on $\Lambda_n:=[-n, n]^d \cap \Z^d$ as $ n \rightarrow \infty$. First, in Section \ref{S:boundary}, we introduce new notation to address different boundary conditions on $\Lambda_n$, and then we define the limiting \textit{uniform spanning forest} (USF) and \textit{minimal spanning forest} (MSF) measures. Large parts of this section are similar in spirit to \cite{Pem91}. In Section \ref{S:RSTREZd}, we prove the existence of a limiting RSTRE measure on $\Z^d$ and study some properties of the resulting measure. In Section \ref{S:free_energy}, we take a statistical physics point of view and study the normalized free energy $\log Z / |\Lambda_n|$ of our model. Finally, in Section \ref{S:Zd_overlap}, we consider the edge overlap of two spanning trees in the same random environment and relate this quantity to the free energy.

\section{Boundary conditions and weak convergence} \label{S:boundary}
Consider the $d$-dimensional Euclidean lattice $\mathbb{L}^d = (\mathbb{Z}^d,E^d)$, $d \geq 2$, where $E^d = \{ (u,v)  \ : \ u,v \in \Z^d, \Vert u - v \Vert_1 = 1 \}$. %As for $d=1$ the graph $\mathbb{L}^d$ does not contain any cycles, the RSTRE is simply is equal to $\mathbb{L}^1$, we will implicitly assume that $d \geq 2$. 
By abuse of notation, for $A \subset \mathbb{R}^d$ denote by $A \cap E^d$ the set of edges with both endpoints in $A \cap \Z^d$. For the rest of the chapter, we will use the following assumption, but remark that the atomless requirement is only needed when discussing the minimum spanning tree (or forest).

\begin{assumption} \label{as:atomless}
    We consider environment distributions $\mu$ that are supported on $(-\infty, \infty)$ and have no atoms, i.e.\ for any $a \in \R$ we assume that $\mu(a) = 0$. We let the environment $\omega = (\omega_e)_{e \in E^d}$ be an i.i.d.\ collection of random variables distributed according to $\mu$, and define weights $\weight(e) = \exp(-\beta \omega_e)$ for $\beta \geq 0$.
\end{assumption}

Let $A \subset \Z^d$ be a finite set and denote by
\begin{equation*}
    \partial A = \big\{u \in A : \exists v \in \Z^d \setminus A \text{ with } (u,v) \in E^d \big\}
\end{equation*}
the (inner) vertex boundary of $A$. We consider graphs $G^\#_A$ with the following two boundary conditions:
\begin{itemize}
    \item \textbf{Free} $G^F_A$: Remove all edges which do not have both endpoints in $A$. In other words, $G_A$ has vertex set $A$ and edge set
    \begin{equation*}
        E( G^F_A) = \big\{ (u,v) \in E^d : u,v \in A \big\}.
    \end{equation*}
    
    \item \textbf{Wired} $G^W_A$: Contract all vertices outside of $A$ into a single vertex. To ensure that conductances stay finite, we remove any self-loops that formed this way. That is, by denoting $\dagger$ for the contracted vertex, we have
    \begin{align*}
        V(G^W_A) &= A \cup \{ \dagger \}, \\
        E(G^W_A) &= \big\{ (u,v) \in E^d : u,v \in A\} \cup \{ (u, \dagger) : u \in \partial A \big\}.
    \end{align*}

    % \item \textit{Periodic:} In the special case of $A = \Lambda_n$ we consider free boundary conditions on $[-n, n+1]^d$ but identify vertices $u$ and $v$ as the same if $u_i = -(v_i + 1) = n+1$ for some $i \in \{1, \ldots, d\}$. In other words, we obtain the discrete $d$ dimensional torus of radius $N$ with weights originating from weights in $[-n, n+1]^d$. OR TAKE TORUS AND EDD WEIGHTED EDGES r!
\end{itemize}

We denote by $\T^F(A)$ and $\T^W(A)$ the set of spanning trees on $G^F_A$ and $G^W_A$, respectively. We use $\#$ to denote either $F$ (free) or $W$ (wired) boundary conditions. The edge weights on $G^\#_A$ are taken to be the corresponding weights in $\mathbb{L}^d$. As in \eqref{eq:defRSTRE}, this gives rise to the Gibbs Measure on $\T^\#(A)$ defined by
\begin{equation} \label{eq:def_bP_FWP}
    \bP_{A, \beta}^{\omega, \#}(\cT = T) := \frac{1}{Z_{A}^{\#}} \prod_{e \in T} \exp(-\beta \omega_e) = \frac{1}{Z_{A}^{\#}}  \exp\big(- \beta \sum_{e \in T}\omega_e \big) ,
\end{equation}
with normalization constant $Z^{\#}_A = Z_{A, \beta}^{\omega, \#}$. Occasionally we write $H(\cT, \omega) = \sum_{e \in \cT} \omega_e$ for the Hamiltonian of $\cT$. In Theorem \ref{T:RSTRE_limit}, we will show that there exists a unique limiting measure regardless of what boundary conditions we set on $\partial A$.

\subsection{Weak convergence}

It is useful to represent subgraphs $H$ of $\mathbb{L}^d$ as functions from the edge set $E^d$ to $\{0,1\}$, i.e.\ $H$ applied on an edge $e$ is $1$ if and only if $e$ is in the edge set of $H$. We view the set of all subgraphs as a topology on $\{0,1\}^{E^d}$ generated by cylinder sets of the form $C(A)$, where $A$ is a finite set of edges and $C(A)$ are all the subgraphs of $\mathbb{L}^d$ containing all edges of $A$. Occasionally we also identify $H$ with its edge set $E(H)$ and say for $A \subset E^d$ that $A \subset H$ if $A \subset E(H)$. Note that the Borel $\sigma$-algebra on $\{0,1\}^{E^d}$ is also generated by cylinder sets $C(A)$ for finite edge sets $A$. One may now define weak convergence of measures $\nu_n$ to $\nu$ in the usual way, however, the following (equivalent) characterization will be easier to work with.
\begin{definition}[Weak convergence]
    We say that a sequences of probability measures $(\nu_n)_{n \in \N}$ on $\{0,1\}^{E^d}$ converges weakly to $\nu$ if for all finite $A \subset E^d$ we have
\begin{equation} \label{eq:weak_conv}
    \nu_n( A \subset H) \xrightarrow{n \rightarrow \infty} \nu(A \subset H).
\end{equation}
\end{definition}

We remark that if $\nu_n$ is defined on cylinder sets of the form $\{ A \subset H\}$, then for any finite set of disjoint edges $A,B \subset E^d$, using the inclusion-exclusion principle, we can write
\begin{equation*}
    \nu_n(A \subset H, B \cap H = \emptyset) = \sum_{S \subseteq B} \nu_n( A \cup S \subset H) (-1)^{|S|}.
\end{equation*}
Hence, if the convergence in \eqref{eq:weak_conv} holds for all finite $A$, then the equality of measures also holds for any set depending on only finitely many edges. Furthermore, by Kolmogorov's extension theorem (see e.g.\ \cite[Theorem A.3.1]{Dur19}), the measure $\nu$ in \eqref{eq:weak_conv} is uniquely determined by its values on cylinder sets. The measures $\nu_n$ will be defined on (finite) exhaustions of $\mathbb{L}^d$, which we define as follows.
\begin{definition}[Exhaustion]
     An exhaustion $\langle G_n \rangle$ of an infinite graph $G = (V,E)$ is a sequence of finite subgraphs $G_n = (V_n, E_n)$ such that $V_n \subset V_{n+1}$, $E_n \subset E_{n+1}$ for all $n \in \N$, and $V = \cup_{n} V_n$ as well as $E = \cup_n E_n$. Given $\omega$, we define $\bP^\omega_{G_n, \beta}$ as in \eqref{eq:def_bP_FWP} with free boundary conditions.
\end{definition}
\noindent One particularly useful exhaustion of $\mathbb{L}^d$ is when $G_n$ is the induced subgraph $G^F_{\Lambda_n}$ with free boundary conditions.  Although the wired construction of $G^W_{\Lambda_n}$ is not an exhaustion of $\mathbb{L}^d$ (recall that we add an extra vertex $\dagger$), we may still view spanning trees on $G^W_{\Lambda_n}$ as elements of $\{0,1\}^{E^d}$ since every edge of the form $(u, \dagger)$, for $u \in \partial \Lambda_n$, corresponds to a unique edge in $\mathbb{L}^d$ so that we may identify the edges of $G^W_{\Lambda_n}$ with a subset of edges in $E^d$.

\subsection{Uniform and minimal spanning forests}

Suppose that $\langle G_n \rangle$ is an exhaustion of $\Z^d$, and let $\nu^{0}_n$ be the (unweighted) UST measure on $G_n$. Pemantle \cite{Pem91} showed that there exists a unique measure $\nu^{0}$ supported on spanning forests on $\Z^d$ to which the sequence $\nu^{0}_n$ weakly converges to. Furthermore, the resulting spanning forest is almost surely connected if and only if $d \leq 4$. As exhaustions of graphs have free boundary conditions, we call this measure the \textit{free uniform spanning forest} measure. In \cite{BLPS01}, a similar construction but with wired boundary conditions was considered, and it was shown that the limiting measure coincides with the case of free boundary conditions (for $\Z^d$). We shall therefore refer to this measure as the \textit{uniform spanning forest} (USF) measure. For arbitrary infinite graphs the wired and free uniform spanning forest may differ, see e.g.\ \cite[Example 9.9]{BLPS01} for a counter-example.

To define the limiting MST measure on $\Z^d$, assume that the edge weights are i.i.d.\ according to some continuous distribution as in Assumption \ref{as:atomless}. For simplicity, one often considers the case when the weights are uniformly distributed on $[0,1]$. The \textit{free minimal spanning forest} is obtained by deleting any edge whose value is largest in some finite cycle of $\mathbb{L}^d$. One may again consider a wired construction as in \cite{LPS06}, however, as these measures agree on $\Z^d$ we shall simply call the limiting measure the \textit{minimal spanning forest} (MSF) and denote it by $\nu^\infty$. Equivalently, one may also define the MSF using an exhaustion $\langle G_n \rangle$ of $\Z^d$. We state one observation in the following lemma and refer to \cite[Section 3]{LPS06} for more details.
\begin{lemma} \label{L:MSF_prop}
    Let $\langle G_n \rangle$ be an exhaustion of $\Z^d$ and let $\nu^\infty$ be the MSF measure on $\Z^d$. Then 
    \begin{align*}
        e \in \text{MSF} &\Longrightarrow \exists N \in \N \text{ with } e \in \text{MST of } G_n \text{ for } n \geq N, \\
        e \not\in \text{MSF} &\Longrightarrow \exists N \in \N \text{ with } e \not\in \text{MST of } G_n \text{ for } n \geq N.
    \end{align*}
    Furthermore, the edge $e$ is in the MSF if and only if there exists some finite vertex set $W$ such that $e$ has the minimal weight in the (finite) edge set $E(W,V \setminus W)$.   
\end{lemma}

\section{Limiting RSTRE measure} \label{S:RSTREZd}

We first recall a stochastic domination property between spanning trees on the graph with free boundary conditions $G^F_A$ and on the graph with wired boundary conditions $G^W_A$.

\begin{definition}[Stochastic domination]
Let $\nu_1, \nu_2$ be two probability measures supported on a countable set $E$. We say that $\nu_1$ stochastically dominates $\nu_2$ on $S \subseteq E$ if there exists a probability measure $\nu$ on $E \times E$, with marginals given by $\nu_1$ and $\nu_2$, such that $\nu$ is supported on the set
\begin{equation*}
    \big\{ (T_1, T_2) \subseteq E \times E : \ T_1 \cap S \supseteq T_2 \cap S \big\}.
\end{equation*}
\end{definition}

For any finite $A \subset \Z^d$, the graph $G^F_A$ contains fewer edges than $G^W_A$, so intuitively one might guess that the UST measure on $G^F_A$ stochastically dominates the UST measure on $G^W_A$. Using electric network arguments, we can indeed verify this intuition.

\begin{lemma}\label{L:UST_domination}
    For any finite $A \subset \Z^d$ and $\beta \geq 0$ the (weighted) free UST measure $\bP^{\omega, F}_{A, \beta}$ stochastically dominates the (weighted) wired UST measure $\bP^{\omega, W}_{A, \beta}$ on $ E(G^F_A)$.
\end{lemma}
\begin{proof}[Proof]
    For any $F = \{e_1, \ldots, e_k\} \subseteq E(G^F_A)$, using the spatial Markov property (Lemma \ref{L:USTmarkov}) and Kirchhoff's formula (Theorem \ref{T:Kirchhoff}), we may write the probability of $F$ being a subset of the tree as
    \begin{align}
        \bP^{\omega, \#}_{A, \beta}( F \subset \cT) &= \prod_{i=1}^k \bP^{\omega, \#}_{A, \beta}( e_i \in \cT \mid e_{i-1}, \ldots, e_1 \in \cT) \nonumber \\
        &= \prod_{i=1}^k \weight(e_i) \effR{G^{\#}_A / \{ e_1, \ldots e_{i-1}\}}{e_i^-}{e_i^+}. \label{eq:effR_decomp}
    \end{align}
    Rayleigh's monotonicity principle gives that
    \begin{equation*}
        \effR{G^{F}_A / \{ e_1, \ldots e_{i-1}\} }{e_i^-}{e_i^+} \geq \effR{G^{W}_A / \{ e_1, \ldots e_{i-1}\}}{e_i^-}{e_i^+},
    \end{equation*}
    so that with \eqref{eq:effR_decomp} we have
    \begin{equation*}
        \bP^{\omega, F}_{A, \beta}(F \subset \cT)  \geq \bP^{\omega, W}_{A, \beta}(F \subset \cT),
    \end{equation*}
    which (as it holds for any $F \subseteq E(G^F_A)$) is equivalent to stochastic domination.
\end{proof}

We are now ready to prove the existence of a limiting measure, where we do not require the atomless condition in Assumption \ref{as:atomless}.

\begin{theorem} \label{T:RSTRE_limit}
    Let $\beta \geq 0$. There exists a unique measure $\bP^\omega_{\infty, \beta}$, supported on cycle-free graphs on $\mathbb{Z}^d$, to which the measures $\bP^{\omega, F}_{\Lambda_n, \beta}$ and $\bP^{\omega, W}_{\Lambda_n, \beta}$ weakly converge.
\end{theorem}

\begin{proof}
    The proof for free boundary conditions follows analogously as the unweighted case in \cite{Pem91}, see also \cite{BLPS01}. Namely, let $A = \{e_1, \ldots, e_k\}$ and, for $n$ large enough, decompose the probability that $A$ is contained in $\cT$ into 
    \begin{equation*}
        \bP^{\omega, \#}_{\Lambda_n, \beta}( A \subset \cT) = \prod_{i=1}^k \bP^{\omega, \#}_{\Lambda_n, \beta}( e_i \in \cT \mid e_{i-1}, \ldots, e_1 \in \cT).
    \end{equation*}
    Further decomposing this into a product of weights and effective resistances (as in the proof of Lemma \ref{L:UST_domination}), we see that Rayleigh's monotonicity principle (Theorem \ref{T:Rayleigh}) gives the following inequalities
    \begin{align*}
        \bP_{\Lambda_{n}, \beta}^{\omega, W}(A \subset \cT) &\leq \bP_{\Lambda_{n+1}, \beta}^{\omega, W}(A \subset \cT), \\
        \bP_{\Lambda_{n}, \beta}^{\omega, F}(A \subset \cT) &\geq \bP_{\Lambda_{n+1}, \beta}^{\omega, F}(A \subset \cT).
    \end{align*}
    As these sequences are monotone and bounded, the limits
    \begin{align*}
        \bP_{\infty, \beta}^{\omega, F}(A \subset \cT) &:= \limn \bP_{\Lambda_{n}, \beta}^{\omega, F}(A \subset \cT), \\
        \bP_{\infty, \beta}^{\omega, W}(F \subset\cT) &:= \limn \bP_{\Lambda_{n}, \beta}^{\omega, W}(A \subset \cT)
    \end{align*}
    almost surely exist, so that the measures converge weakly. 
    % Furthermore, the stochastic domination in Lemma \ref{L:UST_domination} implies that
    % \begin{equation*}
    %     \bP_{\infty, \beta}^{\omega, F}(A \subset \cT) \geq \bP_{\infty, \beta}^{\omega, W}(A \subset \cT).
    % \end{equation*}
    For each $n$ the measures are supported on cycle-free graphs, and so the same must hold for the limiting measures. We proceed to show that the limiting measures are the same.

    We claim that any connected component, both under $\bP_{\infty, \beta}^{\omega, F}$ and $\bP_{\infty, \beta}^{\omega, W}$, is almost surely of infinite size. Assume otherwise, then there exists two finite disjoint edge sets $A,B \subset E^d$ such that $B$ is a cutset between the set of endpoints $V(A)$ of the edges in $A$ and $\Z^d \setminus V(A)$ (i.e.\ any path between $V(A)$ and $\Z^d \setminus V(A)$ uses an edge in $B$) and
    \begin{equation*}
        \bP_{\infty, \beta}^{\omega, \#}(A \subset \cT, B \cap \cT = \emptyset) > 0.
    \end{equation*}
    However, this is a contradiction to the fact that, for any $n$ large enough so that $\Lambda_n$ contains both $A$ and $B$, we have
    \begin{equation*}
        \bP_{\Lambda_n, \beta}^{\omega, \#}(A \subset \cT, B \cap \cT = \emptyset) = 0,
    \end{equation*}
    as the spanning tree must be connected for finite graphs.

    Furthermore, as $\mathbb{L}^d$ is an amenable graph (the ratio $|\partial \Lambda_n|/|\Lambda_n|$ tends to zero), \cite[Remark 6.1]{BLPS01} shows that for any spanning forest $\mathcal{F}$ whose connected components are of infinite size we have
    \begin{equation*}
        \limn \frac{\sum_{v \in \Lambda_n} \deg_{\mathcal{F}}(v)}{|\Lambda_n|} = 2,
    \end{equation*}
    where $\deg_{\mathcal{F}}(v)$ denotes the degree of $v$ in $\mathcal{F}$. In particular, the same must almost surely hold for the trees (or possibly disconnected forests - see Theorem \ref{T:RSTRE_connect_4}) both under $\bP_{\infty, \beta}^{\omega, F}$ and $\bP_{\infty, \beta}^{\omega, W}$. By \cite[Proposition 5.10]{BLPS01} this implies that the laws must coincide.
\end{proof}

\begin{remark}
    For any boundary condition $B$ on the edges in $E(\partial \Lambda_n, \Lambda_n^c)$, we have by Rayleigh's monotonicity principle, as in the proof of Lemma \ref{L:UST_domination}, that
    \begin{equation*}
        \bP^{\omega, W}_{\Lambda_n, \beta}( A \subseteq \cT) \leq \bP^{\omega, B}_{\Lambda_n, \beta}( A \subseteq \cT) \leq \bP^{\omega, F}_{\Lambda_n, \beta}( A \subseteq \cT).
    \end{equation*}
    In view of Theorem \ref{T:RSTRE_limit}, the limiting measure is thus independent of the choice of boundary. Hence, from now on we will often drop the superscript denoting the boundary conditions.
    \end{remark}

\begin{lemma} \label{L:translation_rotation}
    The averaged law 
    \begin{equation}
        \widehat{\P}_{\infty, \beta}(\cdot) := \E\big[ \bP_{\infty, \beta}^{\omega}(\cdot) \big]
    \end{equation}
    is invariant under translations as well as rotations of the axes. 
\end{lemma}
\begin{proof}
    First, notice that by the dominated convergence theorem we have
    \begin{equation*}
        \widehat{\P}_{\infty, \beta}(A \subset \cT) = \E \big[ \limn \bP_{\Lambda_{n}, \beta}^{\omega, F}(A \subset \cT) \big] = \limn \E \big[ \bP_{\Lambda_{n}, \beta}^{\omega, F}(A \subset \cT) \big]. %\label{eq:invariant_DMT}
    \end{equation*}
    Suppose first that $\tau$ is a rotation of an axis of $\mathbb{L}^d$. Then 
    \begin{equation*}
        \bP_{\Lambda_{n}, \beta}^{\omega, F}(A \subset \cT) = \bP_{\Lambda_{n}, \beta}^{\tau(\omega), F}(\tau(A) \subset \cT).
    \end{equation*}
    As the environment is i.i.d., the laws of the random variables
    \begin{equation*}
        \bP_{\Lambda_{n}, \beta}^{\tau(\omega), F}(\tau(A) \subset \cT) \qquad \text{and} \qquad \bP_{\Lambda_{n}, \beta}^{\omega, F}(\tau(A) \subset \cT)
    \end{equation*}
    must coincide. Hence, 
    \begin{equation*}
        \E \big[  \bP_{\Lambda_{n}, \beta}^{\tau(\omega), F}(\tau(A) \subset \cT) \big] = \E \big[ \bP_{\Lambda_{n}, \beta}^{\omega, F}(\tau(A) \subset \cT) \big],
    \end{equation*}
    from which the rotation invariance follows from the above.

    If $\tau$ is a translation with inverse $\tau^{-1}$, then
    \begin{equation*}
        \bP_{\Lambda_{n}, \beta}^{\omega, F}(\tau(A) \subset \cT) = \bP_{\tau^{-1}(\Lambda_{n}), \beta}^{\tau^{-1}(\omega), F}(A \subset \cT)
    \end{equation*}
    with
    \begin{equation*}
        \E \big[ \bP_{\tau^{-1}(\Lambda_{n}), \beta}^{\tau^{-1}(\omega), F}(A \subset \cT) \big] = \E \big[ \bP_{\tau^{-1}(\Lambda_{n}), \beta}^{\omega, F}(A \subset \cT) \big].
    \end{equation*}
    Whenever $n$ is large enough (depending on $\tau$), we have
    \begin{equation*}
        \Lambda_n \subseteq \tau^{-1}(\Lambda_{2n}) \subseteq \Lambda_{4n},
    \end{equation*}
    and by Rayleigh's monotonicity principle
    \begin{equation*}
        \bP_{\Lambda_{n}, \beta}^{\omega, F}(\tau(A) \subset \cT) \geq \bP_{\tau^{-1}(\Lambda_{2n}), \beta}^{\omega, F}(\tau(A) \subset \cT) \geq \bP_{\Lambda_{4n}, \beta}^{\omega, F}(\tau(A) \subset \cT).
    \end{equation*}
    From this one may deduce that
    \begin{equation*}
        \limn \bP_{\tau^{-1}(\Lambda_{n}), \beta}^{\omega, F}(\tau(A) \subset \cT) = \limn \bP_{\Lambda_{n}, \beta}^{\omega, F}(\tau(A) \subset \cT)
    \end{equation*}
    so that, after applying the dominated convergence theorem once again, we obtain the translation invariance.
\end{proof}

\subsection{Number of trees} 

In \cite{Pem91}, it was proven that $\nu^{0}$ (the unweighted USF measure) concentrates on a single connected tree if and only if $d \leq 4$. We extend this to weighted graphs where the ratio of any two weights is uniformly bounded from below and above. 

\begin{theorem} \label{T:RSTRE_connect_4}
    Fix $\beta \geq 0$ and $(\omega_e)_{e \in E^d}$. Suppose that there exists some $K \geq 0$ such that $\omega_e \in [-K, K]$ for all edges $e$ in $E^d$. Then the probability measure $\bP^\omega_{\infty, \beta}$ is supported on connected trees if and only if $d \leq 4$.
\end{theorem}

\noindent We obtain the following immediate corollary, where we do not require $\mu$ to be atomless.
\begin{corollary} \label{C:RSTRE_connect_4}
    Suppose that $\mu$ is supported on a bounded set, i.e.\ there exists some $K > 0$ such that $\mu( [-K, K]) = 1$. Then for any $\beta \geq 0$ the RSTRE is almost surely connected if and only $d \leq 4$.
\end{corollary}

The main tool to prove Theorem \ref{T:RSTRE_connect_4} is to count the number of intersections of a random walk and a loop erased random walk given the environment. For two trajectories $X = (X_0, X_1, \ldots)$ and $Y = (Y_0, Y_1, \ldots)$ we define the number of intersections between $X$ and $Y$ up to time $n$ as
\begin{equation} \label{eq:num_intersection}
    I_n := \sum_{i=0}^n \sum_{j=0}^n 1_{X_i = Y_j},
\end{equation}
and we are interested in whether $\limn I_n$ is finite or not. We refer to Proposition 9.1 and Theorem 9.4 of \cite{BLPS01} and \cite{Pem91} for the following theorem, and remark that this can be seen as a consequence of Wilson's algorithm rooted at infinity as presented in Section \ref{S:Wilson}.
\begin{theorem} \label{T:UST_connected}
    Let $(G, \weight)$ be a weighted graph. Then
    \begin{enumerate}[(i)]
        \item the wired spanning forest is almost surely a single tree if an independent random walk and a loop erased random walk started at every (or some) vertex intersect infinitely often almost surely; \label{enu:finite}
        
        \item and if the number of intersections of two independent random walks started at every (or some) vertex is almost surely finite, then the wired spanning forest has infinitely many components. \label{enu:infinite}
    \end{enumerate}
\end{theorem}

\begin{remark}\label{R:RW_LERW_intersction}
    Suppose we consider two independent (non-loop erased) random walks $X$ and $Y$. Theorem 1 of \cite{LPS03} (see also \cite[Theorem 10.21]{LP16}) shows that if $X$ and $Y$ intersect infinitely often almost surely, then the same holds for the loop erasure $\LE(X)$ and $Y$. We may therefore replace part \eqref{enu:finite} in Theorem \ref{T:UST_connected} with the condition that $X$ and $Y$ intersect infinitely often. Computing the intersection of $X$ and $Y$ is often considerably easier than calculating the intersections of $\LE(X)$ and $Y$.
\end{remark}

To count the number of random walk intersections, we will use the following heat kernel bounds due to Delmotte \cite{Del99}. We also refer to Chapter 8 of \cite{Kum14} for similar results not assuming ellipticity (that is boundedness of $\mu$) of the environment. 

\begin{lemma} \label{L:RWRE_heat}
    Suppose the edge weights are such that $c^{-1} \leq \weight(e) \leq c$ for some $c \geq 1$. Then there exist constants $c_1, c_2, c_3, c_4 > 0$ such that the lazy random walk transition kernel $q_n(\cdot, \cdot)$ on the weighted graph $(\mathbb{L}^d, \weight)$ satisfies
    \begin{equation} \label{eq:RWRE_heat}
        c_1 \frac{1}{n^{d/2}} e^{- c_2 \frac{\Vert x - y \Vert^2}{n}} \leq q_n(x,y) \leq c_3 \frac{1}{n^{d/2}} e^{- c_4 \frac{\Vert x - y \Vert^2}{n}},
    \end{equation}
    for any $x, y \in \Z^d$ and $n \geq \Vert x - y\Vert$.
\end{lemma}
\noindent We easily obtain the following corollary from the heat kernel bounds for the lazy random walk on $(\mathbb{Z}^d, \weight)$ with $\weight(e) = \exp(-\beta \omega_e)$ for fixed $\beta \geq 0$ and bounded disorder variables. % $\omega_e \in [-K, K]$ for some $K \geq 0$.

\begin{corollary} \label{Co:green_n}
    Denote by $g_n(x,y)$ the Green's function up to time $n$, i.e.\ 
    \begin{equation*}
        g_n(x,y) := \sum_{k=0}^n q^\omega_k(x,y),
    \end{equation*}
    where $q^\omega_k(x,y)$ is the $k$-step transition kernel of the lazy random walk on $\mathbb{L}^d$ with weights $\weight(e) = \exp(-\beta \omega_e)$. If $\beta \geq 0$ is fixed and $\omega_e \in [-K, K]$ for some $K \geq 0$ and all edges $e \in E^d$, then for all $x \in \Z^d$
    \begin{equation} \label{eq:sym_gn_diag}
        \sum_{z \in \Z^d} g_n(x,z)^2 \asymp 
         \begin{dcases}
            \sqrt{n} & \text{ if } d=3, \\
            \log n & \text{ if } d=4, \\
            1 & \text{ if } d \geq 5,
        \end{dcases}
    \end{equation}
    and for all $x,y \in \Z^d$ with $\Vert x - y \Vert^2 \leq n$
    \begin{equation} \label{eq:sym_gn_offdiag} 
        \sum_{z \in \Z^d} g_n(x,z) g_n(y,z) \gtrsim
         \begin{dcases}
            \sqrt{n} - \Vert x-y \Vert & \text{ if } d=3, \\
            \log n - 2 \log \Vert x-y \Vert & \text{ if } d=4.
        \end{dcases}
    \end{equation}
\end{corollary}
\begin{proof}
    By our assumption on $\omega_e \in [-K, K]$, the heat kernel estimates of Lemma \ref{L:RWRE_heat} hold (with some constants depending on $\beta$ and $K$). By reversibility of the random walk and the boundedness of the environment, we have that $q^\omega_i(x,z) \asymp q^\omega_i(z,x)$. Using the heat kernel bounds we then obtain
    \begin{align*}
        \sum_{z \in \Z^d} g_n(x,z)^2 &= \sum_{z \in \Z^d} \sum_{i=0}^n \sum_{j=0}^n q^\omega_i(x,z) q^\omega_j(x,z) \asymp  \sum_{i=0}^n \sum_{j=0}^n \sum_{z \in \Z^d} q^\omega_i(x,z) q^\omega_j(z,x) \\
        &= \sum_{i=0}^n \sum_{j=0}^n q^\omega_{i+j}(x,x) \asymp \sum_{i=1}^n \sum_{j=1}^n \frac{1}{(i+j)^{d/2}}.
    \end{align*}
    From here one can verify (e.g.\ by computing an integral approximation) that the order is as in \eqref{eq:sym_gn_diag}.
    Similarly, we have that
    \begin{equation*}
        \sum_{z \in \Z^d} g_n(x,z)g_n(y,z) \asymp \sum_{i=0}^n \sum_{j=0}^n q^\omega_{i+j}(x,y).
    \end{equation*}
    As there are $k$ many pairs $(i,j)$, $0 \leq i,j \leq n$, such that $i+j=k$, the heat kernel bounds give that
    \begin{equation*}
        \sum_{i=0}^n \sum_{j=0}^n q^\omega_{i+j}(x,y) \geq \sum_{k=\Vert x - y \Vert^2}^n k \cdot q^\omega_{k}(x,y) \gtrsim \sum_{k=\Vert x - y \Vert^2}^n k^{1-d/2},
    \end{equation*}
    which satisfies \eqref{eq:sym_gn_offdiag}. 
\end{proof}

We are now ready to prove Theorem \ref{T:RSTRE_connect_4}. 

\begin{proof}[Proof of Theorem \ref{T:RSTRE_connect_4}]
    For $d=1$ the spanning tree is trivially connected, while for $d=2$ the random weighted graph is almost surely recurrent. Indeed, as $\omega_e$ is bounded from below and above, the effective resistance between any two sets in the weighted graph differs by at most a constant multiple of the quantity in the unweighted graph. Hence, the effective resistance from $0$ to $\infty$ is infinite and the graph is recurrent. We refer to Section 2.5 of \cite{LP16} for more details. Theorem \ref{T:UST_connected} part (\ref{enu:finite}) gives that the tree is almost surely connected.

    As $\omega_e \in [-K, K]$, for some $K \geq 0$, the transition probabilities $q^\omega_n(\cdot, \cdot)$ satisfy the heat kernel bounds in \eqref{eq:RWRE_heat} for some constants (depending on $\beta$ and $K$). Suppose now that $d \geq 5$, and start two independent random walks $X,Y$ at $0$. Recall from \eqref{eq:num_intersection} that $I_n$ denotes the number of intersections of $X$ and $Y$ up to time $n$. Then by Corollary \ref{Co:green_n} (we use here $\bE$ for expectation w.r.t.\ the random walk measure $Q$)
    \begin{equation*}
        \bE[I_n ] = \bE \big[ \sum_{i=0}^n \sum_{j=0}^n 1_{X_i = Y_j} \big] = \sum_{z \in \Z^d} \sum_{i=0}^n \sum_{j=0}^n  q^\omega_{i}(0,z)q^\omega_{j}(0,z) = \sum_{y \in \Z^d} g_n(0,z)^2 \asymp 1
    \end{equation*}
    for all $n$. Markov's inequality together with item (\ref{enu:infinite}) of Theorem \ref{T:UST_connected} shows that there are infinitely many tree components.
    
    Assume now that $d = 3$ or $d=4$, and consider two independent random walks $X$ and $Y$. Using Remark \ref{R:RW_LERW_intersction}, to prove that the tree is connected, it suffices to show that $X$ and $Y$ intersect infinitely often almost surely. We follow the approach of \cite[Section 5]{LPS03} (see also \cite[Section 10.5]{LP16}), by showing that the expectation of $I_n^2$ is of the same order as the expectation squared. Instead of using transitivity as in \cite{LPS03} (which does not necessarily hold in our case for given $\omega$), we will use Corollary \ref{Co:green_n} for bounding the number of intersections.
    
    Suppose that $X$ and $Y$ start at $u$ and $v$ and denote by $I_n(u,v)$ the number of intersections up to time $n$. For $n$ large enough, we have by Corollary \ref{Co:green_n} that
    \begin{equation} \label{eq:intersect_uv}
        \bE[I_n(u,v)] = \sum_{z \in \Z^d} g_n(u,z) g_n(v,z) \gtrsim
         \begin{dcases}
            \sqrt{n} - \Vert u-v \Vert & \text{ if } d=3, \\
            \log n - 2 \log \Vert u-v \Vert & \text{ if } d=4.
        \end{dcases}
    \end{equation}
    Furthermore, by independence of $X$ and $Y$ we have
    \begin{align*}
        \bE[I_n(u,v)^2] &= \sum_{i,j,k,\ell =0}^n \sum_{x,y \in \Z^d} Q^\omega_{(u,v)}(X_i = Y_k = x, X_j = Y_\ell = y) \\
        &=  \sum_{x,y \in \Z^d} \sum_{i,j,k,\ell =0}^n Q^\omega_u(X_i = x, X_j = y) Q^\omega_v(Y_k = x, Y_\ell = y).
    \end{align*}
    By splitting the sum into terms where $i\leq j$ and $k \leq \ell$, and using the Markov property, the expectation of $I_n(u,v)^2$ is bounded by
    \begin{align}
         & \sum_{x,y \in \Z^d} \Big( \sum_{0 \leq i \leq j\leq n} Q^\omega_u(X_i = x)Q^\omega_x (X_{j-i} = y) + \sum_{0 \leq j \leq i\leq n} Q^\omega_y(X_{i-j} = x) Q^\omega_u(X_{j} = y) \Big) \nonumber \\
         &\qquad \quad \times \Big( \sum_{0 \leq k \leq \ell \leq n} Q^\omega_v(Y_k = x)Q^\omega_x (X_{\ell-k} = y) + \sum_{0 \leq \ell \leq k \leq n} Q^\omega_y(Y_{k-\ell} = x) Q^\omega_v(Y_{\ell} = y) \Big) \nonumber \\
        &\quad \leq \sum_{x,y \in \Z^d} \Big( g_n(u,x) g_n(x,y) + g_n(y,x) g_n(u,y) \Big) \Big( g_n(v,x) g_n(x,y) + g_n(y,x)g_n(v,y) \Big) \nonumber \\
        &\quad \leq \sum_{x,y \in \Z^d} \Big( g_n(u,x)^2 g_n(x,y)^2 + g_n(y,x)^2 g_n(u,y)^2 \nonumber \\
        & \qquad \qquad  \qquad \qquad \qquad + g_n(v,x)^2 g_n(x,y)^2  + g_n(y,x) ^2 g_n(v,y)^2 \Big). \label{eq:I^2_split}
        % &\quad \asymp \sum_{x,y \in \Z^d} g_n(x,y)^2 \Big( g_n(u,x) + g_n(u,y) \Big) \Big( g_n(v,x) + g_n(v,y) \Big) \nonumber \\
        %  &\quad \asymp \sum_{x,y \in \Z^d} g_n(x,y)^2 \Big( g_n(u,x)^2 + g_n(u,y)^2 + g_n(v,x)^2 + g_n(v,y)^2 \Big) \nonumber \\
        %&\quad \asymp \sum_{x,y \in \Z^d} g_n(x,y)^2 \Big( g_n(u,x)g_n(v,x) +  g_n(u,x) g_n(v,y) + g_n(u,y) g_n(v,x) +  g_n(u,y)g_n(v,y) \nonumber
        % &\qquad \qquad \times \Big( \sum_{k=0}^n Q^\omega_0(Y_k = x)Q^\omega_x (X_{\ell-k} = y) + \sum_{0 \leq \ell \leq k \leq n} Q^\omega_y(Y_{k-\ell} = x) Q^\omega_0(Y_{\ell} = y) \Big) \nonumber \\
         %
        %  &\quad 4 \sum_{x,y \in \Z^d} \Bigg( \Big( \sum_{i=0}^n q^\omega_i(0,x) \sum_{j=0}^n q^\omega_j(x,y) \Big)^2 + \Big( \sum_{i=0}^n q^\omega_i(y,x) \sum_{j=0}^n q^\omega_j(0,y) \Big)^2 \Bigg) \nonumber \\
         %
        %  &\qquad \quad \times \Big( \sum_{0 \leq k \leq \ell \leq n} Q^\omega_0(Y_k = x)Q^\omega_x (X_{\ell-k} = y) + \sum_{0 \leq \ell \leq k \leq n} Q^\omega_y(Y_{k-\ell} = x) Q^\omega_0(Y_{\ell} = y) \Big) \nonumber \\
        %  %
        %  &\quad \leq 2 \sum_{x,y \in \Z^d} \Bigg( \Big( \sum_{i=0}^n q^\omega_i(0,x) \sum_{j=0}^n q^\omega_j(x,y) \Big)^2 + \Big( \sum_{i=0}^n q^\omega_i(y,x) \sum_{j=0}^n q^\omega_j(0,y) \Big)^2 \Bigg) \nonumber \\
        %  &\quad = 4 \sum_{x \in \Z^d} \Big( g_n(0,x)^2  \sum_{y \in \Z^d} g_n(x,y)^2  \Big) \label{eq:g_nx_g_ny}.
    \end{align}
    % where in the last row we used the inequality 
    % \begin{equation*}
    %    (a+b)(b+c) \leq  a^2 + b^2 + c^2 +d^2.
    % \end{equation*}
    As the environment is bounded, it holds that $g_n(x,y) \asymp g_n(y,x)$, and so first summing over $y$ (or $x$) shows that, up to constants, \eqref{eq:I^2_split} is bounded by
    \begin{align*}
         2 \sum_{x \in \Z^d} \bE[I_n(x,x)] \big( g_n(u,x)^2 + g_n(v,x)^2 \big) &\asymp \bE[I_n(0,0)] \big( \bE[I_n(u,u)] + \bE[I_n(v,v)] \big) \\
         &\asymp \bE[I_n(0,0)]^2. 
    \end{align*}
     
    By Corollary \ref{Co:green_n} and \eqref{eq:intersect_uv}, we thus have for some universal constant $c^* > 0$ (depending on $\beta$ and $K$ but importantly not on $u$ or $v$) that
    \begin{equation*}
        \liminf\limits_{n \rightarrow \infty} \frac{\bE [I_n(u,v)]^2}{\bE[I_n(u,v)^2]} \geq c^*.
    \end{equation*}
    Therefore, by the Paley–Zygmund inequality (see Theorem \ref{T:Paley-Z}) we have
    \begin{equation}
        Q \Big( I_n(u,v) \geq \frac{1}{2} \bE[I_n(u,v)] \Big) \geq \frac{1}{4} \frac{\bE[I_n(u,v)]^2}{\bE[I_n(u,v)^2]} \geq \frac{c^*}{4} - o(1). \label{eq:intersection_paley}
    \end{equation}
    
    Denote by $\Lambda$ the event that $X$ and $Y$ have infinitely many intersections, and by $Q_{(u,v)}(\cdot)$ the law of two independent random walks started at $u$ and $v$.  Lévy's zero–one law (see \cite[Theorem 4.6.9]{Dur19}) shows that
    \begin{equation}
       \lim\limits_{m \rightarrow \infty} Q_{(0,0)}( \Lambda \mid X_1, \ldots, X_m, Y_1, \ldots, Y_m) = 1_{\Lambda}, \label{eq:Q_0-1_law}
    \end{equation}
    and furthermore, by the Markov property, we have
    \begin{equation*}
       Q_{(0,0)}( \Lambda \mid X_1, \ldots, X_m, Y_1, \ldots, Y_m) = Q_{(X_m, Y_m)}( \Lambda).
    \end{equation*}
    As $\bE[I_n(u,v)] \rightarrow \infty$ for any $u,v \in \Z^d$, inequality \eqref{eq:intersection_paley} shows that
    \begin{equation*}
        Q_{(X_m, Y_m)}( \Lambda) = Q_{(X_m, Y_m)} \Big( \limn I_n(X_m, Y_m) = \infty \Big) \geq \frac{c^*}{8}.
    \end{equation*}
    Combining this with \eqref{eq:Q_0-1_law} gives
    \begin{equation*}
        1_{\Lambda} \geq \frac{c^*}{8} > 0,
    \end{equation*}
    so that $\Lambda$ must almost surely hold.    
\end{proof}

\begin{question}
    For which class of environments $\mu$ and which $\beta$ does the result of Corollary \ref{C:RSTRE_connect_4} hold? To what extent can the boundedness condition on disorder distribution $\mu$ be relaxed? Is there a $\mu$ which does not satisfy the heat kernel bounds in \eqref{eq:RWRE_heat}, but the result of Corollary \ref{C:RSTRE_connect_4} still holds?
\end{question}

\subsection{Convergence to the MSF for diverging \texorpdfstring{$\beta$}{beta}}

Now consider an arbitrary exhaustion $\langle G_n \rangle$ of $\mathbb{L}^d$, and define the Gibbs measure as in \eqref{eq:def_bP_FWP}, but allow $\beta$ to depend on $n$. Recall from Assumption \ref{as:atomless} that we assume that $\mu$ is atomless. %and further recall from Section \ref{SS:CoupleRG} that we can couple the MST on $G_n$ to the random environment $\omega$ by letting the MST be the almost surely unique (as we assume that $\mu$ is atomless) spanning tree of $G_n$ that minimizes $H(T, \omega) = \sum_{e \in T} \omega_e$. 
We obtain the following theorem via a straightforward application of the electric network tools developed in Chapter \ref{ch:prelims}.

\begin{theorem}
    If $\langle G_n \rangle$ is an exhaustion of $\mathbb{L}^d$ and $\beta = \beta(n) \rightarrow \infty$, then $\E[ \bP^\omega_{G_n, \beta}(\cdot)] $ converges weakly to the MSF measure $\nu^\infty$.
\end{theorem}
\begin{proof}
    By \eqref{eq:weak_conv} it suffices to prove that
    \begin{equation} \label{eq:RSTRE_MST_conv}
        \limn \E[ \bP^\omega_{G_n, \beta}(F \subset \cT)] = \nu^\infty(F \subset M)
    \end{equation}
    for any finite set $F \subset E^d$. Denote by $M = M(\omega)$ the (almost surely unique) MSF on $\Z^d$ obtained by a coupling with $\omega$ similar to that of Section \ref{SS:CoupleRG}, i.e.\ to determine $M$ we assign to each edge $e$ the weight $\omega_e$. We claim that $\P$-almost surely we have
    \begin{equation*}
        \bP^\omega_{G_n, \beta}(F \subset \cT) \xrightarrow{n \rightarrow \infty} 1_{F \subset M(\omega)},
    \end{equation*}
    so that \eqref{eq:RSTRE_MST_conv} follows by the dominated convergence theorem. 
    
    Let $F = \{e_1, \ldots, e_k\}$ and assume first that $\omega$ is such that $F \subset M$. By a union bound, we then obtain
    \begin{equation} \label{eq:RSTRE_MST_union}
        \bP^\omega_{G_n, \beta}(F \subset \cT) \geq 1 - \sum_{i = 1}^k\bP^\omega_{G_n, \beta}(e_i \not\in \cT).
    \end{equation}
    By Lemma \ref{L:MSF_prop}, for each edge $e_i$ there exists some finite $W_i$ such that $e_i$ is the smallest outgoing edge from $W_i$ to $V \setminus W_i$. Let $N$ be large enough such that $W_i \subset \Lambda_N$ and $E(W_i, V \setminus W_i) \subset E_N$ for each $i =1, \ldots, k$. Then the Nash-Williams inequality (Lemma \ref{L:nash_williams}) with cutset $E(W_i, V_N \setminus W_i)$ gives that for $n \geq N$
    \begin{equation*}
        \effR{G_n}{e^-_i}{e^+_i} \geq \frac{1}{\sum_{e \in E(W_i, V_N \setminus W_i)} \exp(-\beta \omega_e)}.
    \end{equation*}
    Kirchhoff's formula then shows that
    \begin{align*}
        \bP^\omega_{G_n, \beta}(e_i \in \cT) &\geq \frac{\exp(-\beta \omega_{e_i})}{\sum_{e \in E(W_i, V_N \setminus W_i)} \exp(-\beta \omega_e)} \\
        &= \frac{1}{1 + \sum_{e \in E(W_i, V_N \setminus W_i) \setminus \{ e_i\} } \exp\big(-\beta (\omega_e - \omega_{e_i}) \big)}.
    \end{align*}
    As $\beta$ diverges with $n$ and $\omega_{e_i} < \omega_e$ for $e \in  E(W_i, V_N \setminus W_i) \setminus \{e_i \}$, we therefore obtain that the probability that each $e_i$ is in the MSF goes to one, and thus    \begin{equation*}
        \bP_{G_n, \beta}(F \subset \cT) \rightarrow 1.
    \end{equation*}

    Now assume that $\omega$ is such that $F \not\subset M$, i.e.\ there is an $e_i$ such that $e_i \not\in M$. Then by Lemma \ref{L:MSF_prop} there exists some $N \in \N$ such that $e_i$ is not contained in the MST of $G_n$ for all $n \geq N$. That is, there exists a path $P$ in $G_N$ between the endpoints of $e_i$ on which each edge has a weight less than the weight of $e_i$. Let $L = |V(G_N)|$, then the length of $P$ is at most $L$. Denote by $f$ the edge with the largest weight on $P$. By Kirchhoff's formula, Rayleigh's monotonicity principle, and the series law, we have
    \begin{equation*}
        \bP^\omega_{G_n, \beta}(e_i \in \cT) \leq L \exp \big(- \beta (\omega_{e_i} - \omega_f) \big) \xrightarrow{n \rightarrow \infty} 0,
    \end{equation*}
    completing the proof as $\omega_{e_i} > \omega_f$. %    where we used that $\beta$ diverges with $n$.
\end{proof}

\subsection{Infinite electric networks} \label{SS:InfElectric} 
We end this section by providing some properties of electric networks in relation to (weighted) uniform spanning forests on infinite graphs. Suppose that $(G, \weight)$ is an infinite (with locally finite outgoing conductance) weighted graph. From now we will assume that $\bP^\weight_G$ is the weak limiting measure with wired boundary conditions constructed similarly as in Theorem \ref{T:RSTRE_limit} (note that in general, the limiting measure of the free and wired constructions may differ). The electric network results from Chapter \ref{ch:prelims} also apply to the infinite setting, although the connection to random walks may break down for general graphs. We refer to e.g.\ Proposition 9.2 and Section 10.2 of \cite{LP16} for the following.

\begin{lemma}[Kirchhoff's formula for infinite graphs] \label{L:InfKirch}
    Suppose $(G, \weight)$ is an infinite connected weighted graph and let $\bP^\weight_G$ be the law of the (wired) USF, then
	\begin{equation*}
	    \bP^\weight_G(e \in \cT) = \weight(e) \effR{G, \weight}{e^-}{e^+},
	\end{equation*}
    where
    \begin{equation} \label{eq:InfEffR}
        \effR{G, \weight}{e^-}{e^+} = \left\{ \Vert\theta\Vert^2_r \, : \, \theta \text{ is a unit flow from } e^- \text{ to } e^+ \right\}.
    \end{equation}
\end{lemma}
\noindent See also \cite[Theorem 10.15]{LP16} for the transfer-impedance theorem in the infinite setting (compare to Theorem \ref{T:transfer-impedance} for the finite graph case).

\begin{lemma}[Spatial Markov property for infinite graphs] \label{L:InfSpatialMarkov}
    \sloppy % 
    Suppose $(G, \weight)$ is an infinite connected weighted graph and let $A, B\subset E$ be two finite disjoint sets of edges such that $\bP^\weight_G(A\subset \cT, B\cap \cT=\emptyset)>0$. Then for any set of finite edges $F\subset E$,
	\begin{equation*}
	    \bP^\weight_G(F \subset \cT \mid A\subset \cT, B\cap \cT =\emptyset) = \bP^{\weight}_{(G-B)/A}(F \subset \cT).
	\end{equation*}
    \fussy %
\end{lemma}

\begin{proof}
    This follows readily by taking exhaustion of $G$, applying the finite version of the spatial Markov property, and then taking limits, where each term converges by Rayleigh's monotonicity principle.
\end{proof}

\medskip

Using these infinite electric network tools, we may easily obtain that the measures $\bP^\omega_{\infty, \beta}$, as in Theorem \ref{T:RSTRE_limit}, are distinct for different choices of $\beta$.
\begin{lemma}
Let $\beta, \beta' \geq 0$ with $\beta < \beta'$. Then almost surely there exists a finite set of edges $F$ such that
\begin{equation*}
    \bP_{\infty, \beta}^{\omega}(F \subset \cT) \neq \bP_{\infty, \beta'}^{\omega}(F \subset \cT).
\end{equation*}
\end{lemma}

%[[DRAFT COMMENT: This is extremely weak and I don't like that. What I think should be the case is that it is different for every $F$ (?). Furthermore, the averaged laws should also be different (for sets containing at least 2 (?) vertices). I will try and prove it in the meantime, but else I leave it like it is....]]}

\begin{proof}
Let $F_1$ be the edge set of some spanning tree on $\partial \Lambda_1$ using only edges with both endpoints in $\partial \Lambda_1$. If $\bP^\omega_{\infty, \beta}(F_1 \subset \cT) \neq \bP^\omega_{\infty, \beta'}(F_1 \subset \cT)$, then we are done. So assume otherwise and let $e_1, \ldots, e_{2d}$ be the edges incident to $0$ ordered such that $\omega_{e_1} < \ldots < \omega_{e_{2d}}$. Consider now $F = F_1 \cup \{e_1\}$, then by the spatial Markov property for infinite graphs (Lemma \ref{L:InfSpatialMarkov}) 
\begin{equation*}
    \bP^\omega_{\infty, \beta}(F \subset \cT) = \bP^\omega_{H, \beta}(e_1 \in \cT) \bP^\omega_{\infty, \beta}(F_1 \subset \cT) = \bP^\omega_{H, \beta}(e_1 \in \cT) \bP^\omega_{\infty, \beta'}(F_1 \subset \cT),
\end{equation*}
where $H$ is the graph obtained by contracting all edges in $F_1$. Notice that the probability that $e_1$ is in the RSTRE of $H$ is the same as the probability that $e_1$ is in the RSTRE on the graph $H'$ consisting of two vertices $0$ and $v$ with $2d$ edges $e_1, \ldots e_{2d}$ between $0$ and $v$. For such a graph, we have using Wilson's algorithm that
\begin{equation*}
    \bP^\omega_{H', \beta}(e_1 \in \cT) = \frac{e^{-\beta \omega_{e_1}}}{\sum_{i=1}^{2d} e^{-\beta \omega_{e_i}}} = \frac{1}{ 1 + \sum_{i=2}^{2d} e^{-\beta (\omega_{e_i} - \omega_{e_1})}},
\end{equation*}
which is strictly increasing in $\beta$ as (almost surely) $\omega_{e_i} > \omega_{e_1}$ for all $i =2, \ldots, 2d$. Consequently, for such a choice of $F$, we have that
\begin{equation*}
    \bP_{\infty, \beta}^{\omega}(F \subset T) < \bP_{\infty, \beta'}^{\omega, \#}(F \subset T). \qedhere
\end{equation*}
\end{proof}

We pose here two more natural questions regarding the difference in $\bP_{\infty, \beta}^{\omega}(\cdot)$ and $\bP_{\infty, \beta'}^{\omega}(\cdot)$ for $\beta \neq \beta'$.

\begin{question}
    Suppose that $\beta \neq \beta'$. Is it true that for any finite set of edges $F \subset E^d$ one has
    \begin{equation*}
        \bP_{\infty, \beta}^{\omega}(F \subset \cT) \neq \bP_{\infty, \beta'}^{\omega}(F \subset \cT) \, ?
    \end{equation*}
\end{question}
\begin{question}
    Suppose that $\beta \neq \beta'$. Is it true that for any finite set of edges $F \subset E^d$, with $|F| \geq 2$, one has
    \begin{equation*}
        \E[ \bP_{\infty, \beta}^{\omega}(F \subset \cT)] \neq \E[ \bP_{\infty, \beta'}^{\omega}(F \subset \cT)] \, ?
    \end{equation*}
    Note that if $F$ consists of just a single edge, then by transitivity we clearly have an equality for the averaged law.
\end{question}

\section{Free energy} \label{S:free_energy}

In statistical mechanics, much of the information about the Gibbs measure is encoded in the \textit{free energy}. For example, in the directed polymer model, the free energy determines the value of $\beta$ at which a phase transition occurs (see \cite[Chapter 2]{Com17} for more details). We define the (finite volume) free energy as
\begin{equation} \label{eq:Free_energy}
    F^\omega_n(\beta) := \frac{\log Z_{\Lambda_n}}{|\Lambda_n|},
\end{equation}
where the dependence of the partition function $Z_{\Lambda_n}$ on $\beta$ and $\omega$ is implicit. When $\beta = 0$, the partition function $Z_{\Lambda_n}$ simply corresponds to $|\mathbb{T}(\Lambda_n)|$, the number of spanning trees on $\Lambda_n$. In \cite{TW00}, it was shown that $|\mathbb{T}(\Lambda_n)|$ has well-defined asymptotically exponential growth, namely for any dimension $d \geq 2$ there exists a finite constant $z_d$ such that
\begin{equation} \label{eq:USTgrowth}
    \lim\limits_{n \rightarrow \infty} \frac{1}{|\Lambda_n|} \log |\mathbb{T}(\Lambda_n)| = z_d.
\end{equation}

In this section, we will assume that $\mu$ is centered, i.e. $\E[\omega_e] = 0$. We remark that for the weighted UST measure this restriction is essentially equivalent to requiring the mean to be finite, as adding any fixed value to all of the disorder variables does not affect the measure. For certain classes of disorder distributions, we will prove the existence of the infinite volume free energy, following along the lines of Theorem 2.1 of \cite{Com17}. %STOP HERE

\begin{lemma} \label{L:E_lim_free_energy}
    Suppose that $\E[\omega_e] = 0$ and $\beta \geq 0$. If almost surely
    \begin{equation} \label{eq:free_energy_assumption}
        \Big| F^\omega_n(\beta) - \E\big[F^\omega_n(\beta)\big] \Big| \xrightarrow{n \to \infty} 0,
    \end{equation}
    then almost surely the infinite volume free energy
    \begin{equation*}
        F(\beta) := \limn F^\omega_n(\beta) = \limn \E \big[  F^\omega_n(\beta)\big].
    \end{equation*}
    exists and is finite.
\end{lemma}

\begin{proof}
By our assumption in \eqref{eq:free_energy_assumption}, it suffices to show that the limit
\begin{equation*}
    \limn \E \big[  F^\omega_n(\beta)\big]
\end{equation*}
is some finite constant. This will be done using a multiparameter super-additivity argument. To ease notation, using the translation invariance of the environment, we consider instead the box $\Lambda^+_n = [0,n]^d$. For boxes of the form $I = [a_1, b_1] \times [a_2, b_2] \times \ldots \times [a_d, b_d]$ with $a_i \leq b_i$ integers, we let $\mathbb{T}(I)$ be the set of spanning trees on the graph $G_{I}$ with free boundary conditions, i.e.\ $G_I$ is the induced subgraph of $\mathbb{L}^d$ with vertex set $I \cap \Z^d$. Denote by $Z_{I}$ the corresponding partition function. 
 
Suppose the box $I$ can be partitioned into $M$ disjoint smaller boxes $I_k$ such that $\cup_k I_k \cap \Z^d = I \cap \Z^d$. First, we construct an injective mapping
\begin{equation*}
    \phi: \mathbb{T}(I_1) \times \mathbb{T}(I_2) \times \ldots \times \mathbb{T}(I_M) \longrightarrow \mathbb{T}(I).
\end{equation*}
Let $H = H(I_1, \ldots I_M)$ be a graph with vertex set $\{1, \ldots, M\}$ and edges $(i,j)$ if there is an edge between a vertex of $I_i$ to a vertex of $I_j$ in $\mathbb{L}^d$. As $H$ is connected, there is a spanning tree of $H$ corresponding to $M-1$ edges $\{e_1, \ldots, e_{M-1}\}$ in $\mathbb{L}^d$ that connect the boxes $I_k$ to each other. For any $(\cT_1, \ldots, \cT_M) \in \mathbb{T}(I_1) \times \ldots \times \mathbb{T}(I_M)$ let $\phi(\cT_1, \ldots, \cT_M)$ be the subgraph with vertices $I \cap \Z^d$ and edges
\begin{equation*}
    \bigcup_{k=1}^M E(\cT_k)  \ \bigcup \ \{e_1, \ldots, e_{M-1}\}.
\end{equation*}
By construction, $\phi(\cT_1, \ldots, \cT_M)$ forms a spanning tree of $I \cap \Z^d$ and the edges $\{e_1, \ldots, e_{M-1}\}$ are not contained in any $\cT_i$. 

It follows that
\begin{equation*}
    Z_{I} \geq \sum_{T \in \phi^{-1}(\mathbb{T}(I) )} e^{ - \beta H(T)} = \prod_{i=1}^{M-1} e^{-\beta \omega_{e_i}} \prod_{j=1}^M Z_{I_j}.
\end{equation*}
As we assume that $\mu$ is centered, taking expectations gives 
\begin{equation*}
     \E [\log Z_{I}] \geq -\beta \sum_{i=1}^{M-1} \E[\omega_{e_i}] + \sum_{i=1}^k \E[ \log Z_{I_i}] = \sum_{i=1}^k \E[ \log Z_{I_i}].
\end{equation*}
Translation invariance of the environment together with a multiparameter version of Fekete's Lemma (see \cite[Theorem 1]{Cap08}) shows that the limit
\begin{equation*}
    \limn \frac{\E [\log Z_{[0,n]^d}]}{(n+1)^d}   = \sup_{n \geq 1} \frac{1}{|\Lambda^+_n|} \E [\log Z_{[0,n]^d}].
\end{equation*}
exists, but it may be infinite. However, using Jensen's inequality and the well-defined exponential growth \eqref{eq:USTgrowth} on the number of spanning, gives that
\begin{equation*}
    \frac{\E [\log Z_{[0,n]^d}]}{(n+1)^d} \leq \frac{\log \E[ Z_{[0,n]^d}]}{(n+1)^d} = \frac{\log | \mathbb{T}_{[0,n^d]} | }{(n+1)^d} \xrightarrow{n \rightarrow \infty} z_d,   
\end{equation*}
and thus the limit is finite.
\end{proof}

The following lemma gives two special cases where the limiting free energy concentrates around its mean as required in \eqref{eq:free_energy_assumption}. As an immediate consequence of Lemma \ref{L:E_lim_free_energy}, we obtain that the free energy converges to a constant in these special cases.

\begin{lemma} \label{L:free_concentrate}
    Suppose that $\mu$ is either bounded or that $\mu$ is a normal distribution with finite variance. Then, for any $\beta \geq 0$, we have that
     \begin{equation}
       \Big| F^\omega_n(\beta) - \E\big[F^\omega_n(\beta)\big] \Big| \xrightarrow{n \to \infty} 0,
    \end{equation}
    almost surely.
\end{lemma}

\begin{proof}
    We will follow the same notation as the proof of Lemma \ref{L:E_lim_free_energy} and consider the boxes $\Lambda^+_n = [0,n]^d$. Assume first that $\mu([a,b]) = 1$ for some $b \geq a$, and consider some edge $e$ contained in $[0,n]^d$. We may view $Z_{[0,n]^d}(\omega_e)$ as a function in $\omega_e$ by fixing all other variables. One may easily verify that
\begin{equation*}
    \big| \log Z^\omega_{[0,n]^d}(b) - \log Z^\omega_{[0,n]^d}(a) \big| \leq \beta (b - a).
\end{equation*}
Let $m = |\Lambda^+_n \cap E^d|$ (which is approximately $2d (n+1)^d$), then an application of McDiarmid's inequality (see e.g.\ \cite[Appendix A.2]{Com17}) gives that for any $\epsilon > 0$
\begin{align}
    \P \Big( \big| \log Z_{[0,n]^d} - \E [\log Z_{[0,n]^d}] \big| \geq \epsilon (n+1)^d \Big) &\leq 2\exp\Big( - \frac{2 \epsilon^2 (n+1)^{2d}}{m \beta^2 (b-a)^2} \Big) \nonumber \\
    &\leq 2\exp \Big( - \frac{c \epsilon^2 (n+1)^d}{\beta^2 (b-a)^2} \Big), \label{eq:logZ_concentrate}
\end{align}
for some universal constant $c = c(d) > 0$. As \eqref{eq:logZ_concentrate} is summable in $n$, the Borel-Cantelli lemma shows that
\begin{equation*}
    \P \Big( \limsup\limits_{n \rightarrow \infty} \frac{1}{(n+1)^d} \big| \log Z_{[0,n]^d} - \E[ \log Z_{[0,n]^d} ] \big| \geq \epsilon \Big) = 0.
\end{equation*}
Letting $\epsilon$ go to zero we obtain almost sure convergence of $F^\omega_n(\beta)$ to its expectation.

For the Gaussian case, we may w.l.o.g.\ restrict ourselves to the standard normal distribution by suitably rescaling $\beta$. We will show that $\log Z_{[0,n]^d}$ is a Lipschitz function with respect to the disorder variables. Calculating derivatives for $e \in \Lambda^+_n \cap E^d$ gives that
\begin{equation*}
    \frac{\partial \log Z_{[0,n]^d}}{\partial \omega_e} = \frac{1}{Z_{[0,n]^d}}   \sum_{\substack{T \in \mathbb{T}([0,n]^d), \\ e \in T}} -\beta e^{-\beta H(T)} = - \beta \bP_{[0,n]^d}(e \in \cT).
\end{equation*}
Hence,
\begin{equation*}
    \Vert \nabla \log Z_{[0,n]^d} \Vert^2 = \beta^2 \sum_{e \in \Lambda^+_n \cap E} \bP_{[0,n]^d}(e \in \cT)^2 \leq 2d \beta^2 (n+1)^d,
\end{equation*}
It follows using Theorem A.2 of \cite{Com17} that (by possibly enlarging $c$) the concentration of \eqref{eq:logZ_concentrate} holds with $b=1$ and $a=0$ for the Gaussian case. Applying Borel-Cantelli again gives the convergence.
\end{proof}

\section{Overlap} \label{S:Zd_overlap}
In this section, we study the tree and edge overlap as functions of $\beta$. First, for finite graphs $G$, we briefly collect basic facts about some of the derivatives of the Gibbs measure $\bP^\omega_{G,\beta}$, both in terms of $\beta$ and $\omega_e$.

\begin{lemma} \label{L:derivative}
    Let $(G, \weight)$, $|V(G)| = n$, be a finite weighted graph with $\weight(e) = \exp(-\beta \omega_e)$ and (weighted) UST law $\bP^\omega_{G,\beta}$ with partition function $Z^\omega_{G, \beta}$. Recall that the Hamiltonian of $T$ was defined as $H(T,\omega) = \sum_{e \in T} \omega_e$. We have the following (partial) derivatives:
    \begin{align}
        \frac{\partial }{\partial \beta} \bP^\omega_{G,\beta}(\cT = T) &= \bP^\omega_{G, \beta}(\cT = T) \big(\bE^\omega_{G, \beta}[H(\cT, \omega)] - H(T, \omega) \big), \label{eq:T_partial_beta}\\
        %
        % \frac{\partial }{\partial \omega_e} \bP^\omega_{G,\beta}(\cT = T) &= \beta \bP^\omega_{G, \beta}(\cT = T) \big(\bP^\omega_{G, \beta}(e \in \cT) - 1_{e \in T} \big), \label{eq:T_partial_omega}\\
        %
        \frac{\partial }{\partial \beta } \bP^\omega_{G,\beta}(e \in \cT ) &= \bE^\omega_{G, \beta}[ 1_{e \in \cT}] \bE^\omega_{G, \beta}[H(\cT, \omega)] - \bE^\omega_{G, \beta}[ 1_{e \in \cT} H(\cT, \omega)], \label{eq:e_partial_beta}\\
        \frac{\partial }{\partial \omega_f }  \bP^\omega_{G,\beta}(e \in \cT) &= \beta \big(  \bP^\omega_{G,\beta}(e \in \cT)  \bP^\omega_{G,\beta}(f \in \cT) - \bP^\omega_{G,\beta}(e,f \in \cT) \big). \label{eq:e_partial_f}
    \end{align}
    Furthermore, if $(\omega_e)_{e \in E}$ are distributed as independent standard Gaussian random variables, then the partition function $Z^\omega_{G, \beta}$ satisfies
    \begin{equation}
         \frac{\partial }{\partial \beta } \E \Big[ \log Z^\omega_{G, \beta} \Big] =  \beta (n-1) -  \beta \sum_{e \in E} \E \big[ \bP^\omega_{G,\beta}(e \in \cT)^2 \big] .\label{eq:Z_partial_beta} 
    \end{equation}
\end{lemma}
\begin{proof}
    The first equalities \eqref{eq:T_partial_beta}, \eqref{eq:e_partial_beta} and \eqref{eq:e_partial_f} are straightforward calculations using the chain rule. For \eqref{eq:Z_partial_beta}, we first compute
    \begin{align*}
        \frac{\partial }{\partial \beta } \log Z^\omega_{G, \beta} &= -  \bE^\omega_{G, \beta}[H(\cT, \omega)] = - \sum_{T \in \T(G)} H(T, \omega) \bP^\omega_{G,\beta}(\cT = T)\\
        &= - \sum_{e \in E} \omega_e \sum_{T \in \T(G)} 1_{e \in T} \bP^\omega_{G,\beta}(\cT = T) = - \sum_{e \in E} \omega_e \bP^\omega_{G,\beta}(e \in \cT).
    \end{align*}
    Swapping integration and differentiation, and applying Gaussian integration by parts (see e.g.\ \cite[Section A.3.2]{Com17}) to $\E[\omega_e \bP^\omega_{G,\beta}(e \in \cT)]$, gives us that
    \begin{align*}
        \E \big[ \frac{\partial }{\partial \beta } \log Z^\omega_{G, \beta} \big] &= - \sum_{e \in E} \E \big[ \frac{\partial}{\partial \omega_e} \bP^\omega_{G,\beta}(e \in \cT) \big] \\
        &= - \beta \sum_{e \in E} \E \big[ \bP^\omega_{G,\beta}(e \in \cT)^2 - \bP^\omega_{G,\beta}(e \in \cT) \big] \\
        &= \beta (n-1) - \beta \sum_{e \in E} \E \big[ \bP^\omega_{G,\beta}(e \in \cT)^2 \big], 
    \end{align*}
    completing the proof.
\end{proof}

\begin{remark}
    The negative correlation of Remark \ref{R:correlation} and \eqref{eq:e_partial_beta} imply that $\bP^\omega_{G,\beta}(e \in \cT)$ is decreasing in $\omega_e$ and increasing in all other $\omega_f$, $f \neq e$. 
\end{remark}

% \begin{remark}
%     In view of Lemma \ref{L:e_MST_increase} and \eqref{eq:e_partial_beta}, we have
%     \begin{equation*}
%         \bE^\omega_{G, \beta}[ 1_{e \in \cT}] \bE^\omega_{G, \beta}[H(\cT, \omega)] > \bE^\omega_{G, \beta}[ 1_{e \in \cT} H(\cT, \omega)]
%     \end{equation*}
%     for edges $e$ that are contained in the MST. That is, the Hamiltonian and the indicator that $e$ is in the tree are negatively correlated. This is not too surprising: trees containing $e$ typically have a smaller Hamiltonian than trees without $e$. However, as Example \ref{ex:counter_increase} shows, if $e$ is not in the MST, then this does not imply positive correlation (for all $\beta$) of the above quantities.
% \end{remark}

\subsection{Tree overlap}

Suppose now that we sample two independent trees $\cT$ and $\cT'$ under the same law $\bP^\omega_{G, \beta}$. The probability that the trees are identical is given by
\begin{equation} \label{eq:same-tree}
    g(\beta) := \sum_{T \in \T(G)} \bP^{\omega, \otimes 2}_{G, \beta}( \cT = \cT' = T) =  \sum_{T \in \T(G)} \bP^{\omega}_{G, \beta}( \cT = T)^2.
\end{equation}
A straightforward consequence of Lemma \ref{L:derivative} is the following.
\begin{lemma} \label{L:exact-same-tree}
Let $G = (V,E)$ be a finite graph with at least one cycle and distinct edge weights $\{\exp(-\beta \omega_e)\}_{e\in E}$. Then the quantity 
\begin{equation}
    g(\beta) =  \sum_{T \in \T(G)} \bP^{\omega}_{G, \beta}( \cT = T)^2
\end{equation}
is strictly increasing in $\beta$ and tends to one as $\beta \rightarrow \infty$.
\end{lemma}
%%%%%%%%% Proof %%%%%%%%%%%%%
\begin{proof}
Denote by $T_1, \ldots, T_m$, $m \geq 2$, all the spanning trees of $G$ ordered such that $H(T_1, \omega) < H(T_2, \omega) < \ldots <  H(T_m, \omega)$, where $H(T, \omega) = \sum_{e \in T}\omega_e$. %This implies that $\bP_{\beta}(T_i) \geq \bP_{\beta}(T_j)$ for $i \leq j$. 
By \eqref{eq:T_partial_beta} we have 
\begin{align*}
    \frac{\partial }{\partial \beta} g(\beta) &= 2 \sum_{i=1}^m \bP^\omega_{G, \beta}(\cT = T_i) \frac{\partial }{\partial \beta} \bP^\omega_{G, \beta}(\cT = T_i)\\
    &= 2 \sum_{i=1}^m \bP^\omega_{G, \beta}(\cT = T_i)^2 \big( \bE^\omega_{G, \beta}[H(\cT, \omega)] - H(T_i, \omega) \big) .
\end{align*}
Let $k < m$ be the largest index such that $H(T_k, \omega) \leq \bE^\omega_{G, \beta}[H(\cT, \omega)]$, and note that as $H(T_i, \omega) < H(T_k, \omega)$, for $i < k$, the same must hold for all $i < k$ as well. Observe that
\begin{align*}
    \sum_{i=1}^k \bP^\omega_{G, \beta}(\cT = T_i) \frac{\partial }{\partial \beta} \bP^\omega_{G, \beta}(\cT = T_i) & \geq \bP^\omega_{G, \beta}(\cT = T_k)  \sum_{i=1}^k \frac{\partial }{\partial \beta} \bP^\omega_{G, \beta}(\cT = T_i), \\
    \sum_{j=k+1}^m \bP^\omega_{G, \beta}(\cT = T_i) \frac{\partial }{\partial \beta} \bP^\omega_{G, \beta}(\cT = T_i) & > \bP^\omega_{G, \beta}(\cT = T_k)  \sum_{i=k+1}^m \frac{\partial }{\partial \beta} \bP^\omega_{G, \beta}(\cT = T_i),
\end{align*}
since the derivatives in the respective sums share the same sign and as
\begin{equation*}
    \bP^\omega_{G, \beta} (\cT = T_i) > \bP^\omega_{G, \beta} (\cT = T_k) > \bP^\omega_{G, \beta} (\cT = T_j)
\end{equation*}
for all $i < k$ and $j > k$. Therefore, using that trivially $\sum_{i=1}^m \bP_{ \beta}^{\omega}(\cT = T_i) = 1$, we obtain
\begin{equation*}
   \frac{\partial }{\partial \beta}g(\beta) %&= 2 \sum_{i = 1}^{k}  \bP_{ \beta}^{\omega} (T_i) \rho_{ \beta}^{\omega} (T_i) +2 \sum_{j = k+1}^{m} \bP_{ \beta}^{\omega} (T_j) \rho_{ \beta}^{\omega} (T_j) \\
   > 2 \bP^\omega_{G, \beta} (\cT = T_k) \sum_{i=1}^m \frac{\partial }{\partial \beta} \bP^\omega_{G, \beta}(\cT = T_i)= 0,
\end{equation*}
and so $g(\beta)$ is strictly increasing.

Lastly, for any $\beta \geq 0$ we clearly have that $g(\beta) < 1$, and furthermore
\begin{equation*}
    g(\beta) \geq \bP^\omega_{G, \beta}(\cT = T_1)^2 \geq \frac{e^{-\beta H(T_1)}}{e^{-\beta H(T_1)} + m e^{-\beta H(T_2)}} =\frac{1}{1 + m e^{-\beta (H(T_2) - H(T_1) )}},
\end{equation*}
which goes to one as $\beta \rightarrow \infty$.
\end{proof}

\begin{remark}
    In statistical physics, it is often useful to study quantities that are monotone in the model parameter. Currently, the tree overlap is the only quantity that we are aware of that is monotone in $\beta$. In fact, as shown later in Section \ref{SS:monotone}, many observables are not almost surely (w.r.t.\ $\omega$) monotone. %Nonetheless, it is natural to conjecture that under the averaged law the monotonicity should hold. 
\end{remark}

\subsection{Edge overlap}

Sample again two independent spanning trees $\cT$ and $\cT'$ of a finite graph $G$ in the same environment as in \eqref{eq:same-tree}, but now consider the number of edges that $\cT$ and $\cT'$ share. The expected edge overlap is defined as
\begin{align}
    \mathcal{O}_G(\beta) := \bE^{\omega, \otimes 2}_{G, \beta} \big[ |\cT \cap \cT'| \big] &:= \sum_{T, T' \in \mathbb{T}(G)} \sum_{e \in E} 1_{e \in T \cap T'} \bP^\omega_{G, \beta}(\cT = T) \bP^\omega_{G, \beta}(\cT = T') \nonumber \\
    &= \sum_{e \in E} \bP^\omega_{G, \beta}(e \in \cT)^2. \label{eq:def_overlap}
\end{align}
One motivation for studying the edge overlap stems from the derivative of the log partition function \eqref{eq:Z_partial_beta} in the Gaussian environment case. That is, the derivative of the expectation of the free energy, as defined in \eqref{eq:Free_energy}, satisfies
\begin{align*}
    \frac{\partial }{\partial \beta} \E[ F^\omega_{n}(\beta)]
    % &= \beta \frac{1}{|\Lambda_n|} \sum_{e \in \Lambda_n \cap E^d} \E \big[ \bP^\omega_{\Lambda_n, \beta}(e \in \cT) - \bP^\omega_{\Lambda_n, \beta}(e \in \cT)^2 \big] \\
    &= \beta \frac{|\Lambda_n| -1}{|\Lambda_n|} - \beta \frac{1}{|\Lambda_n|} \sum_{e \in \Lambda_n \cap E^d} \E \big[ \bP^\omega_{\Lambda_n, \beta}(e \in \cT)^2 \big] \\
    &= \beta \frac{1}{|\Lambda_n|}\big( |\Lambda_n| -1  - \E [\mathcal{O}_{\Lambda_n}(\beta)] \big).
\end{align*}
% As $\bP^\omega_{\Lambda_n, \beta}(e \in \cT) > \bP^\omega_{\Lambda_n, \beta}(e \in \cT)^2$ for any finite $\beta \geq 0$, we have that
% \begin{equation*}
%     \frac{\partial }{\partial \beta} \E[ F^\omega_{n}(\beta)] > 0,
% \end{equation*}
% for all $n \in \N$.
Taking the limit as $n \rightarrow \infty$, this naturally leads us to study the density of the overlap in $\Lambda_n$. Namely, let
\begin{equation*}
    \widehat{\mathcal{O}}_{\Lambda_n}(\beta) := \bE_{\infty, \beta}^{\omega, \otimes 2} \big[ |\cT \cap \cT' \cap \Lambda_n| \big] = \sum_{e \in \Lambda_n \cap E^d} \bP^\omega_{\infty, \beta}(e \in \cT)^2,
\end{equation*}
be overlap density of the infinite volume measure restricted to a finite box. Does the limit
\begin{equation*}
   \limn \frac{\widehat{\mathcal{O}}_{\Lambda_n}(\beta)}{|\Lambda_n| -1 } = \limn \frac{1}{|\Lambda_n| -1 }\sum_{e \in \Lambda_n \cap E^d} \bP^\omega_{\infty, \beta}(e \in \cT)^2
\end{equation*}
exist?  We show in Lemma \ref{L:overlapZd} that this is indeed the case. 

\begin{remark}  
    For the (unweighted) USF the probability of an edge being in the spanning tree is by symmetry equal to $1/d$. Hence,
    \begin{equation*}
        \limn \frac{\widehat{\mathcal{O}}_{\Lambda_n}(0)}{|\Lambda_n| -1 } = \frac{1}{d}.
    \end{equation*}
    On the other hand, for the MSF the overlap is simply the number of edges in $\Lambda_n$ in the MSF. As $\Z^d$ is amenable (see also the proof of Theorem \ref{T:RSTRE_limit}), the average degree in the MSF of a vertex in $\Lambda_n$ tends to $2$, so that the density of the overlap approaches $1$. 
\end{remark}

% Using symmetry, it is easy to see that for the periodic instance the limit of the expectation is well defined, but for other boundary conditions this is not as trivial. We will next prove that the edge overlap converges almost surely to the same value regardless of the boundary conditions.

\begin{lemma} \label{L:overlapZd}
The limit
\begin{equation*}
    \rho(\beta) := \limn \frac{\widehat{\mathcal{O}}_{\Lambda_n}(\beta)}{|\Lambda_n| -1 }
\end{equation*}
exists almost surely, and for any edge $e \in E^d$, is given by
\begin{equation*}
    \rho(\beta) = d \, \E \big[ \bP^\omega_{\infty, \beta}(e \in \cT)^2 \big] .
\end{equation*}
\end{lemma}

\begin{proof}
For notational purposes consider again boxes $\Lambda^+_n = [0,n]^d$, the result for $\Lambda_n$ follows by translation invariance. We wish to apply a version of Birkhoff's Ergodic theorem that holds in this more general multiparameter setting. Consider the shift operators $\tau_1, \ldots, \tau_d$ acting on $\omega$ by letting 
\begin{equation*}
    \tau_i\big( \{\omega_{(u,v)}\}_{(u,v) \in E^d} \big) = \{\omega_{(u - z_i,v -z_i)}\}_{(u,v) \in E^d},
\end{equation*}
where $z_i$ is the $i$-th standard basis vector. As the environment is i.i.d., each $\tau_i$ and any composition of $\tau_i$'s are ergodic. For $x \in \Z^d$ with $x_i \geq 0$, we define
\begin{equation*}
    \tau^x = \tau^{x_1}_1 \circ \tau^{x_2}_2 \circ \ldots \circ \tau^{x_d}_d,
\end{equation*}
where $\tau^{x_i}_i$ denotes the $x_i$ fold application of $\tau_i$.

Consider now the functions
\begin{equation*}
    f_{i}(\omega) = \bP^\omega_{\infty, \beta} \big( (0, z_i) \in \cT \big)^2 \qquad i=1, \ldots, d,
\end{equation*}
and notice that
\begin{equation*}
    f_{i} \circ \tau^x (\omega) = \bP^{\tau^x(\omega)}_{\infty, \beta} \big( (0, z_i) \in \cT \big)^2 = \bP^{\omega}_{\infty, \beta} \big( (x, x + z_i) \in \cT \big)^2.
\end{equation*}
Further, define the operator $S_{\Lambda_n^+}$ as
\begin{equation*}
    S_{\Lambda_n^+} f_i = \sum_{x \in \Lambda^+_n, x_i \neq n} f_i \circ \tau^x.
\end{equation*}
Then $(S_{\Lambda_n^+})_{n \in \N}$ is an additive process (in terms of the notation in \cite[Chapter 6]{Kre85}), and each $f_i$ is an integrable function. By a multiparameter ergodic theorem (see e.g.\ \cite[Theorem 2.8 in Chapter 6]{Kre85}), the limit of
\begin{align*}
    \frac{1}{n(n+1)^{d-1}} S_{\Lambda_n^+} f_i &= \frac{1}{n(n+1)^{d-1}} \sum_{x \in \Lambda^+_n, x_i \neq n} f_i \circ \tau^x (\omega) \\
    &= \frac{1}{n(n+1)^{d-1}} \sum_{e = (x, x +z_i) \in \Lambda^+_n \cap E^d} \bP^{\omega}_{\infty, \beta}( e \in \cT)^2
\end{align*}
almost surely exists for each $i=1, \ldots, d$. Furthermore, as each $\tau_i$ is ergodic, the limit must be given by
\begin{equation*}
    \E \big[ \bP^{\omega}_{\infty, \beta}\big( (0,z_1) \in \cT \big)^2 \big].
\end{equation*}
The translation and axis rotation invariance of the averaged law in Lemma \ref{L:translation_rotation} further gives that for any edge $e \in E^d$ we have
\begin{equation} \label{eq:ergodic_translate}
    \E \big[ \bP^{\omega}_{\infty, \beta} \big( (0,z_1) \in \cT \big)^2 \big] = \E \big[ \bP^{\omega}_{\infty, \beta}( e \in \cT)^2 \big].
\end{equation}
The lemma now follows by summing over all $d$ directions to give
\begin{align*}
    \rho(\beta) &= \limn \frac{1}{(n+1)^{d} - 1} \sum_{e \in \Lambda_n^+ \cap E^d} \bP^{\omega}_{\infty, \beta}( e \in \cT)^2. \\
    &= \limn \Big( \frac{n(n+1)^{d-1}}{(n+1)^{d} - 1} \Big) \sum_{i=1}^d \frac{1}{n(n+1)^{d-1}}  \sum_{x \in \Lambda^+_n, x_i \neq n} f_i \circ \tau^x (\omega) \\
    &= d \, \E[ \bP^{\omega}_{\infty, \beta}( e \in \cT)^2],
\end{align*}
where the last line used the argument leading up to \eqref{eq:ergodic_translate}.
\end{proof}

\subsection{Behavior of the edge overlap between the extremes}

In the directed polymer model, as introduced in the middle of Section \ref{S:def_RSTRE}, there is a close relation between the expected overlap and the critical parameter for the phase transition. For instance, if the environment is Gaussian and $d \geq 3$, then there exists a finite value $\beta_c > 0$ and a subset $\mathcal{I}$ of $(\beta_c, \infty)$\footnote{It is conjectured that $\mathcal{I} = (\beta_c, \infty)$, see \cite[Conjecture 6.1]{Com17}.} such that, for $\beta \in \mathcal{I}$, the normalized overlap is strictly positive. On the other hand, if $\beta < \beta_c$, then the normalized overlap converges to $0$, which is the same as the normalized overlap of the simple random walk on $\mathbb{Z}^d$. We refer to \cite[Section 6]{Com17} for more details. In Corollary \ref{Co:non-trivial-overlap}, we show that for the RSTRE, the normalized overlap for any $\beta > 0$ is always larger than in the case when $\beta = 0$ and smaller than the case when (formally) $\beta = \infty$.

\smallskip

First, we show that $\bP^{\omega}_{\infty, \beta}(e \in \cT)$ is a random function, i.e.\ it is not independent of $\omega$. See also \eqref{eq:e_partial_f} for the case when the graph is finite.

\begin{lemma} \label{L:increase_e}
The function $ \bP^{\omega}_{\infty, \beta}(e \in \cT)$ almost surely is in $(0,1)$ and it is strictly decreasing in $\omega_e$ for $\beta > 0$.
\end{lemma}
\begin{proof} %Maybe use notation from earlier
Applying Kirchhoff's formula (Theorem \ref{T:Kirchhoff}) to an exhaustion $\langle G_n \rangle$, gives that
\begin{equation*}
    \bP_{G_n, \beta} \big( (u,v) \in \cT \big) = \weight(u,v) \effR{G_n}{u}{v} = \frac{\weight(u,v)}{\weight(u,v) + \big( \effR{G_n - (u,v)}{u}{v} \big)^{-1}} .
\end{equation*}
Applying Rayleigh's monotonicity principle then shows that the limit exists and that it is given by
\begin{equation*}
    \bP^{\omega}_{\infty, \beta} \big( (u,v) \in \cT \big) = \frac{\weight(u,v)}{\weight(u,v) + \big( \effR{\Z^d - (u,v)}{u}{v} \big)^{-1}},
\end{equation*}
where the effective resistance on $\Z^d - e$ is defined as in \eqref{eq:InfEffR} (see also Lemma \ref{L:InfKirch}). Clearly $\effR{\Z^d - (u,v)}{u}{v}$ is independent of $\omega_e$, and further, as the weight of any edge is almost surely in $(0, \infty)$, the same must hold for the effective resistance. The assertion of the lemma follows by comparing $\bP^{\omega}_{\infty, \beta}((u,v) \in \cT)$ to the function
\begin{equation*}
    f(x) = \frac{e^{-\beta x}}{e^{-\beta x} + a} \quad a \in (0, \infty),
\end{equation*}
which, for $\beta > 0$, is strictly decreasing in $x$ and always in $(0,1)$.
\end{proof}

% \begin{figure}[ht]
%     \centering
%     \resizebox{!}{5cm}{\includegraphics[]{cage.tex}}
%     \caption{On the left we condition on a \textit{cage}, depicted by thick edges, to be in a spanning tree. Contracting these edges to a single vertex $v$ produces the graph on right. For $T$ to be a spanning tree on this new graph, exactly one of the 4 edges between $0$ and $v$ must be contained in $T$.}
%      \label{fig:cage}
% \end{figure}

%Tree overlap is increasing
The above lemma leads to the following corollary that $\rho(\beta)$ behaves unlike the UST or MST for $\beta \neq 0, \infty$.

\begin{corollary} \label{Co:non-trivial-overlap}
For any $\beta > 0$, we have
\begin{equation*}
    \frac{1}{d} < \rho(\beta) < 1.
\end{equation*}
\end{corollary}

\begin{proof}
Lemma \ref{L:increase_e} implies that $\bP^{\omega}_{\infty, \beta}(e \in T)$ is a random function in terms of $\omega$. The strict lower bound immediately follows from the Cauchy-Schwarz inequality. To show the upper bound, assume otherwise that $\rho(\beta) = 1$. Then by Lemma \ref{L:overlapZd}, we have for any edge $e \in E^d$ that
\begin{equation*}
    \E \big[ \bP_{\infty,\beta}(e \in \cT)^2 \big] = \frac{1}{d} =  \E \big[ \bP_{\infty,\beta}(e \in \cT) \big].
\end{equation*}
Therefore, as $\bP^{\omega}_{\infty, \beta}(e \in \cT) \in [0,1]$, it must almost surely be
\begin{equation*}
    \bP_{\infty,\beta}(e \in \cT)^2 = \bP_{\infty,\beta}(e \in \cT),
\end{equation*}
so that $\bP^{\omega}_{\infty, \beta}(e \in \cT) \in \{0,1\}$, which is a contradiction to Lemma \ref{L:increase_e}.
\end{proof}

\begin{question}
    Is the overlap density $\rho(\beta)$ increasing in $\beta$? At what rate does it approach $\frac{1}{d}$ and $1$ when $\beta \rightarrow 0$ or $\beta \rightarrow \infty$ ?
\end{question}

\subsection{Lack of monotonicity} \label{SS:monotone}

To observe a phase transition in statistical physics, it is often useful to study properties that are monotone in $\beta$. One natural question to ask is if the following holds: 
\begin{equation*}
    `` \, \bP^\omega_{G,\beta}(e \in \cT) \text{ is increasing } \Longleftrightarrow e \in \textrm{MST}(\omega) \, ".
\end{equation*}
In the following two examples we show that this is not the case, see also Example \ref{ex:counter_overlap} for a remark about the edge overlap. In particular, we will see that this implies a lack of monotonicity for some of the observables studied in this thesis. Note that we will choose exact, and somewhat arbitrary, values of $(\omega_e)_{e \in E}$ corresponding to measure zero events, however, using continuity, this implies that these counter-examples occur with positive probability.

\begin{example} \label{ex:MST_e_decrease}
    Consider the complete graph on the four vertices $\{1, 2,3, 4\}$ with one diagonal removed. We assign disorder variables $\omega_{(1,2)} = 0.1, \omega_{(1,3)} = 0.11, \omega_{(1,4)} = 10, \omega_{(2,4)} = 10.1, \omega_{(3,4)} = 10.2$. Then, the MST consists of the edges $(1,2)$, $(1,3)$ and $(1,4)$. One might expect that $\bP^\omega_{G,\beta}(e \in \cT)$ is increasing in $\beta$ for all edges in the MST, but it turns out that
    \begin{equation*}
        \bP^\omega_{G,\beta}\big( (1,4) \in \cT \big)
    \end{equation*}
    is decreasing in an interval containing $[0,0.4]$.
\end{example}

\begin{example} \label{ex:counter_increase}
    Let $G$ be the graph consisting of two vertices with three parallel edges, labeled $a,b$ and $c$, between them. Let $\omega_a = 0.1$, $\omega_b = 1$ and $\omega_c = 10$. Then the MST consists of the edge $a$, however, the probability that $b$ is in the weighted UST is given by
    \begin{equation*}
        \frac{e^{-\beta}}{e^{-0.1 \beta} + e^{-\beta} + e^{-10 \beta}},
    \end{equation*}
    which is strictly increasing in an interval containing $[0, 0.2]$. Hence, even though $b$ is not contained in the MST, the probability of $b$ being in the RSTRE may still be increasing. This example also shows that for trees $T \in \T(G)$ that are not the MST, we may still have that
    \begin{equation*}
        \bP^{\omega}_{G, \beta}(\cT = T)
    \end{equation*}
    is an increasing function for some values of $\beta$.
\end{example}

\noindent The following example shows that the overlap is not almost surely increasing in $\beta$ for finite graphs. Nonetheless, we still believe that averaged (over $\omega$) overlap $\E[\mathcal{O}_G(\beta)]$ is increasing.

\begin{example} \label{ex:counter_overlap}
    Consider the complete graph on 3 vertices (i.e.\ a cycle) with edges denoted by $a,b,c$. To obtain $G$, add an edge $d$ that is parallel to $c$, and let $\omega_a =1$, $\omega_b =2$, $\omega_c =1.01$, and $\omega_d =1.1$. Then, the overlap $\mathcal{O}_G(\beta)$ is strictly decreasing in $\beta$ for $\beta$ in an interval containing $[0,0.1]$. A heuristic explanation for this fact is the following. %The function $\sum_{i=1}^m x^2_i$, under the condition that $\sum_{i=1}^m x_i = 1$, 
    The overlap is large when the summands are either close to $1$ or close to $0$ (as is the case for the MST). In this example, the MST is given by the edges $a$ and $c$, however, due to the parallel edge $d$, for small $\beta$ we have that $\bP^\omega_{G, \beta}(b \in \cT) > \bP^\omega_{G, \beta}(c \in \cT)$. In particular, the increase in $\bP^\omega_{G, \beta}(b \in \cT)^2$ does not compensate for the decrease in $\bP^\omega_{G, \beta}(c \in \cT)^2$.
\end{example}

\begin{remark}
    Consider again the weighted graph from Example \ref{ex:MST_e_decrease}. Further computations show that the diameter of $\cT$ is an increasing function of $\beta$ in an interval containing $[0,0.01]$ and it is a decreasing function for $\beta$ large enough. In particular, due to the lack of almost sure monotonicity, it is not clear how to define a critical threshold for the observable of the diameter, or in fact, for any other observable studied in this thesis. 
\end{remark}

%The above example does not imply that $\E[\mathcal{O}_G(\beta)]$ is

\chapter{Local statistics on the complete graph} \label{ch:local}

We now return to the RSTRE on finite graphs. In particular, this chapter will focus only on the complete graph with the disorder variables being uniformly distributed on $[0,1]$ (see Assumption \ref{As:K_n_U}). We will consider the following two observables as functions of $\beta = \beta(n)$: (1) the edge overlap of two RTSREs sampled under the same environment; (2) the local limit of the RSTRE. One of the main results is that there is a sharp transition in the local limit when $\beta = n^{\gamma}$ and $\gamma$ crosses the critical threshold $\gamma_c = 1$, which stands in contrast to the global observable of the diameter in Chapter \ref{ch:diam} (see also Remark \ref{R:local_global}). This chapter is essentially a summary of \cite{Mak24}. 

\section{Setup and main results}

Throughout this whole chapter, we will work with the following setup.
\begin{assumption} \label{As:K_n_U}
    Let $G_n = (V_n, E_n) = (V,E)$ be the complete graph on $n$ vertices. For $\beta = \beta(n) \geq 0$, let $\cT^\omega_{n, \beta}$ be the RSTRE on $G_n$ with i.i.d.\ disorder variables $(\omega_e)_{e\in E_n}$ uniformly distributed on $[0,1]$.
\end{assumption}

We denote the (quenched) law of $\cT^\omega_{n, \beta}$ by $\bP^\omega_{n, \beta}$. When there is no ambiguity, we drop the sub and superscripts and write $\bP$ and $\bE$ instead of $\bP^\omega_{n, \beta}$ and $ \bE^\omega_{n, \beta}$, respectively. In this case, the weights are given by
\begin{equation*}
    \weight(e) = \exp( - \beta \omega_e),
\end{equation*}
with mean
\begin{equation*}
    \xi = \xi(\beta) = \frac{1 - e^{-\beta}}{\beta}.
\end{equation*}
Let $S_m = \sum_{i=1}^m \weight(e_i)$ be a sum of the weights of $m$ (distinct) edges. The concentration inequality in Lemma \ref{L:low_concentration} shows that  for $\beta \geq 1$ and any $ 0 < \delta \leq 1$ we have
\begin{equation} \label{eq:BernU_ch3}
	\P \big( |S_m - m \xi | \geq \delta m \xi \big) \leq 2 \exp \Big( - \frac{\delta^2 m}{9 \beta } \Big).
\end{equation}

\smallskip

The first local observable we study is the expected overlap $\mathcal{O}(\beta)$, as already introduced in \eqref{eq:def_overlap}, of two independent trees $\cT$ and $\cT'$ sampled under the same environment. Recall that this quantity is defined as
\begin{equation} \label{eq:def_overlap_ch4}
    \mathcal{O}(\beta) :=  \bE^{\omega, \otimes 2}_{n, \beta} \big[ |\cT \cap \cT'| \big] = \sum_{T, T' \in \T(G)} | T \cap T'| \bP^\omega_{n, \beta}(\cT = T) \bP^\omega_{n, \beta}(\cT = T').
\end{equation}
Applying Kirchhoff's formula (Theorem \ref{T:Kirchhoff}), shows that we may rewrite this as
\begin{equation} \label{eq:OV_effR}
    \mathcal{O}(\beta)  = \bE^{\otimes 2} \Big[ \sum_{e \in E} 1_{e \in \cT, e \in \cT'}  \Big] = \sum_{e \in E}\bP(e \in \cT)^2 = \sum_{e \in E} w(e)^2 \effR{}{e^-}{e^+}^2.
\end{equation}

\noindent The first result of this chapter gives the overlap as a function of $\beta$ under the setup as in Assumption \ref{As:K_n_U}.

\begin{theorem}[{\cite[Theorem 1.1]{Mak24}}] \label{T:Overlap}
    Let $\cT^\omega_{n, \beta}$ be the RSTRE as in Assumption \ref{As:K_n_U}. If $0 < \beta = \beta(n) \ll n/\log n$, then with high probability
    \begin{equation} \label{eq:overlap_low}
        \mathcal{O}(\beta) = (1 + o(1)) \beta \frac{1 - e^{-2\beta}}{(1 -e^{-\beta})^2}.
    \end{equation}
    If on the other hand $\beta \gg n (\log n)^2$, then with high probability
    \begin{equation} \label{ovelap_high}
        \mathcal{O}(\beta) = (1 - o(1)) n.
    \end{equation}
\end{theorem}
\begin{remark}
    The probability of an edge being in the (unweighted) UST on the complete graph (which corresponds to $\beta = 0$ in \eqref{eq:PomegaT}) is $2/n$ and is independent of $\omega$, and hence the expected overlap of the UST is $2 (n-1)/n$. On the other hand, for the MST the overlap is trivially $(n-1)$. Notice that the overlap in \eqref{eq:overlap_low} is approximately equal to $\beta$, and hence the overlap smoothly interpolates between the overlap in the UST and MST as $\beta$ increases from $0$ to $n$.
    %In fact, when $\beta \gg n (\log n)^2$ the RSTRE will be equal to the MST up to $o(n)$ number of edges.
\end{remark}

We remark that (independently of our work) Kúsz in \cite{K24} studied yet another observable called the total length, which in our notation is defined as
 \begin{equation*}
    L(\cT) := H(\cT, \omega) = \sum_{e \in \cT} \omega_e,
 \end{equation*}
 and is closely related to the edge overlap.
\begin{theorem}[cf. {\cite[Theorem 1.4]{Mak24}} and {\cite[Theorem 1.8]{K24}}] \label{T:length}
      Let $\cT^\omega_{n, \beta}$ be the RSTRE as in Assumption \ref{As:K_n_U}. If $0 < \beta = \beta(n) \ll n/\log n$, then with high probability (and in $\P$-expectation)
    \begin{equation} \label{eq:length_low}
       \bE^\omega_{n,\beta} \big[ L(\cT) \big]= (1 + o(1)) \frac{n}{\beta} \cdot \frac{1 - \beta e^{-\beta} - e^{-\beta}}{1- e^{-\beta}}.
    \end{equation}
    If on the other hand $\beta \gg n (\log n)^5$, then with high probability (and in $\P$-expectation)
    \begin{equation} \label{eq:length_high}
        \bE^\omega_{n,\beta} \big[ L(\cT) \big] = (1 + o(1)) \zeta(3),
    \end{equation}
    where $\zeta(3) = \sum_{k=1}^\infty k^{-3} = 1.202\dots$ is Apéry's constant.
\end{theorem}

\noindent The constant $\zeta(3)$ is precisely the limiting total length of the MST on the complete graph, as shown by \cite{Fri85} (see also \cite{Jan95} for a central limit result). We remark that in \cite{K24}, Lemma \ref{L:effR_lowdisorder_both} (see also \cite[Lemma 3.1 and 3.5]{Mak24}) is used to prove the asymptotics in the low disorder regime. Furthermore, Kúsz in \cite{K24} also shows that the asymptotics in the high disorder regime hold under the milder condition $\beta \gg n \log n$. Theorem \ref{T:length} will be a consequence of the tools we develop for the edge overlap, see also the end of Section \ref{SS:low_edge} for an alternative simpler proof of \eqref{eq:length_high} that has not appeared elsewhere and works as long as $\beta \gg n \log n$. 

\begin{question}
    What is the (typical) order of
    \begin{equation*}
        X := \bE^\omega_{n,\beta} \big[ L(\cT) - L(\textrm{MST})  \big]   \, ?
    \end{equation*}
    Heuristic arguments might suggest that $X$ should be of order $n^2/\beta^2$ whenever $\beta \gg n$.
\end{question}

\smallskip

The proofs of Theorem \ref{T:Overlap} and Theorem \ref{T:length} rely on giving sharp bounds on the effective resistance between any pairs of vertices. Using some further arguments, these bounds will also us to prove the following theorem, where we refer to Section \ref{S:local} for a definition of local (weak) convergence.

\begin{theorem}[{\cite[Theorem 1.3]{Mak24}}] \label{T:local_limit}
    If $\beta = \beta(n) \ll n/\log$, then under the averaged law $\widehat{\P}$ we have the local convergence
    \begin{equation} \label{eq:localRSTRE_UST}
        \cT(1) \xrightarrow{d} \mathcal{P}(o),
    \end{equation}
    where $\cT(1)$ is $\cT$ rooted at $1$ and $\mathcal{P}$ is a Poisson(1) branching process conditioned to survive forever. If on the other hand $\beta = \beta(n) \geq n (\log n)^\lambda$, where $\lambda = \lambda(n) \rightarrow \infty$ arbitrarily slowly, then there exists a (random) rooted tree $\mathcal{M}(o)$ such that under the averaged law the local convergence
    \begin{equation} \label{eq:localRSTRE_MST}
        \cT(1) \xrightarrow{d} \mathcal{M}(o)
    \end{equation}
    holds.
\end{theorem}
\noindent The rooted trees $(\mathcal{P}, o)$ and $(\mathcal{M}, o)$ are precisely the local limits of the UST and MST on the complete graph, respectively. See also Section \ref{S:local}. Note that \cite[Theorem 1.9]{K24} shows that the local convergence to $\mathcal{M}$ holds under the weaker assumption of $\beta \gg n \log n$.

\begin{remark} \label{R:local_global}
    Theorem \ref{T:local_limit} shows that there is a sharp cut-off for the local limit when $\beta = n^\gamma$ and $\gamma$ crosses the critical value $\gamma_c = 1$. We therefore observe two regimes: (1) a low disorder regime $\beta \ll n/\log n$, where we observe UST behavior; (2) a high disorder regime $\beta \gg n (\log n)^{\lambda}$, where we observe MST behavior. Later in Section \ref{S:high_and_low}, we will again encounter two different regimes, however, there the high disorder corresponds to $\beta \geq n^{4/3}$. Furthermore, we conjecture there to be a smooth transition of the diameter when $\beta = n^{\gamma}$ and $1< \gamma < 4/3$, see Conjecture \ref{C:Intermediate}.
\end{remark}

\begin{question}
    What is the local limit of the RSTRE when $\beta = Cn$ for some constant $C > 0$? Is it different for every value of $C$?
\end{question}

\section{Introduction to local limits} \label{S:local} 

The pair $(G, v)$, where $v \in V(G)$, is referred to as a rooted graph with $v$ as its root. Two rooted graphs $(G, v)$ and $(G', v')$ are isomorphic if there exists a graph isomorphism $\phi : V(G) \rightarrow V(G')$ between $G$ and $G'$ that fixes the roots, i.e.\ there exists a bijection $\phi$ between the vertex sets $V(G)$ and $V(G')$ such that $\phi(v) = v'$ and $(x,y) \in E(G)$ if and only if $(\phi(x), \phi(y)) \in E(G')$, and in this case we write $(G,v) \simeq (G',v')$. For $r = 1, 2, \ldots$ denote by 
\begin{equation*}
    B_G(v,r) = \Big( \big\{ u \in V(G) : d_G(u,v) \leq r \big\}, \ \big\{ (x,y) \in E(G) : d_G(x,v), d_G(y,v) \leq r \big\} \Big)
\end{equation*}
the closed ball (viewed as a subgraph) of radius $r$ centered at $v$ in $G$. We define the pseudometric on the set of rooted graphs by
\begin{equation*}
    d_0\big( (G,v), (G',v') \big) :=  2^{-k},
\end{equation*}
where $k$ is defined as
\begin{equation*}
    \sup \big\{ r \in \N \, : \, \big(B_G(v,r), v \big) \simeq \big( B_{G'}(v',r), v' \big) \big\},
\end{equation*}
with the notation $2^{-\infty} := 0$. Finally, let $\mathcal{G}_*$ be the isomorphism classes of (locally finite) rooted graphs, then $(\mathcal{G}_*, d_0)$ forms a complete and separable metric space. See also \cite{BS01} or \cite[Chapter 2]{vdH24} for a more thorough introduction. We will use the following definition of local (weak) convergence.

\begin{definition}[Local weak convergence]
    Let $(G,v)$ be a (random) rooted graph, and consider a random graph sequence $(G_n)_{n \in \N}$ with corresponding laws $(\mu_n)_{n \in \N}$, and let $(v_n)_{n \in \N}$ be vertices drawn uniformly from the vertex sets of $(V(G_n))_{n \in \N}$. We say that the sequence $(G_n)_{n \in \N}$ converges locally to $(G,v)$ if, for all $r \geq 1$ and every finite rooted graph $(H,o)$, we have
    \begin{equation*}
        \limn \mu_n \Big( \big(B_{G_n}(v_n,r), v_n \big) \simeq (H,o) \Big) = \mu \Big( \big(B_G(v,r), v \big) \simeq (H,o) \Big).
    \end{equation*}
\end{definition}

Under Assumption \ref{As:K_n_U}, the averaged law is invariant under any permutation of the vertices (the complete graph is transitive and the environment is i.i.d.), so we may as well fix the random vertex $v_n = 1$. In particular, the local weak convergence of the RSTRE to some rooted graph $(\mathcal{S},v)$ is equivalent to
\begin{equation} \label{eq:def_local_conv}
    \widehat{\P}_{n,\beta} \Big( B_{\cT}(1,r) \simeq (H,o) \Big) \xrightarrow{n \rightarrow \infty} \mu \Big( B_{\mathcal{S}}(v,r) \simeq (H,o) \Big),
\end{equation}
for all $r \geq 1$ and rooted graphs $(H,o)$. If the convergence in \eqref{eq:def_local_conv} holds, then we write $\cT(1) \xrightarrow{d} \mathcal{S}(v)$.

\medskip

We briefly describe how to construct $(\mathcal{P}, o)$, the Poisson(1) branching process conditioned to survive forever. Consider an infinite backbone $o, x_1, x_2, \ldots$ started at the root $x_0 = o$. Then attach to each vertex $x_i$ on the backbone an independent Poisson(1) branching processes tree, each of which is almost surely finite. The following local convergence on the complete graph is due to Grimmett \cite{Gri80}.
\begin{theorem}
    The UST on the complete graph converges locally to $(\mathcal{P}, o)$.
\end{theorem}
Since the work of \cite{Gri80}, the local convergence has been extended to other families of graphs in e.g.\ \cite{BP93} and \cite{NP22}. For instance, in \cite{NP22}, the authors show that the above theorem holds for any (almost) d-regular graph sequence with degrees diverging with $n$. 

\smallskip

\begin{figure}[tbh]
  \centering
  %\hspace*{-0.3cm}%
  \subcaptionbox{UST\label{fig:localUST}}{\includegraphics[width=0.5\textwidth]{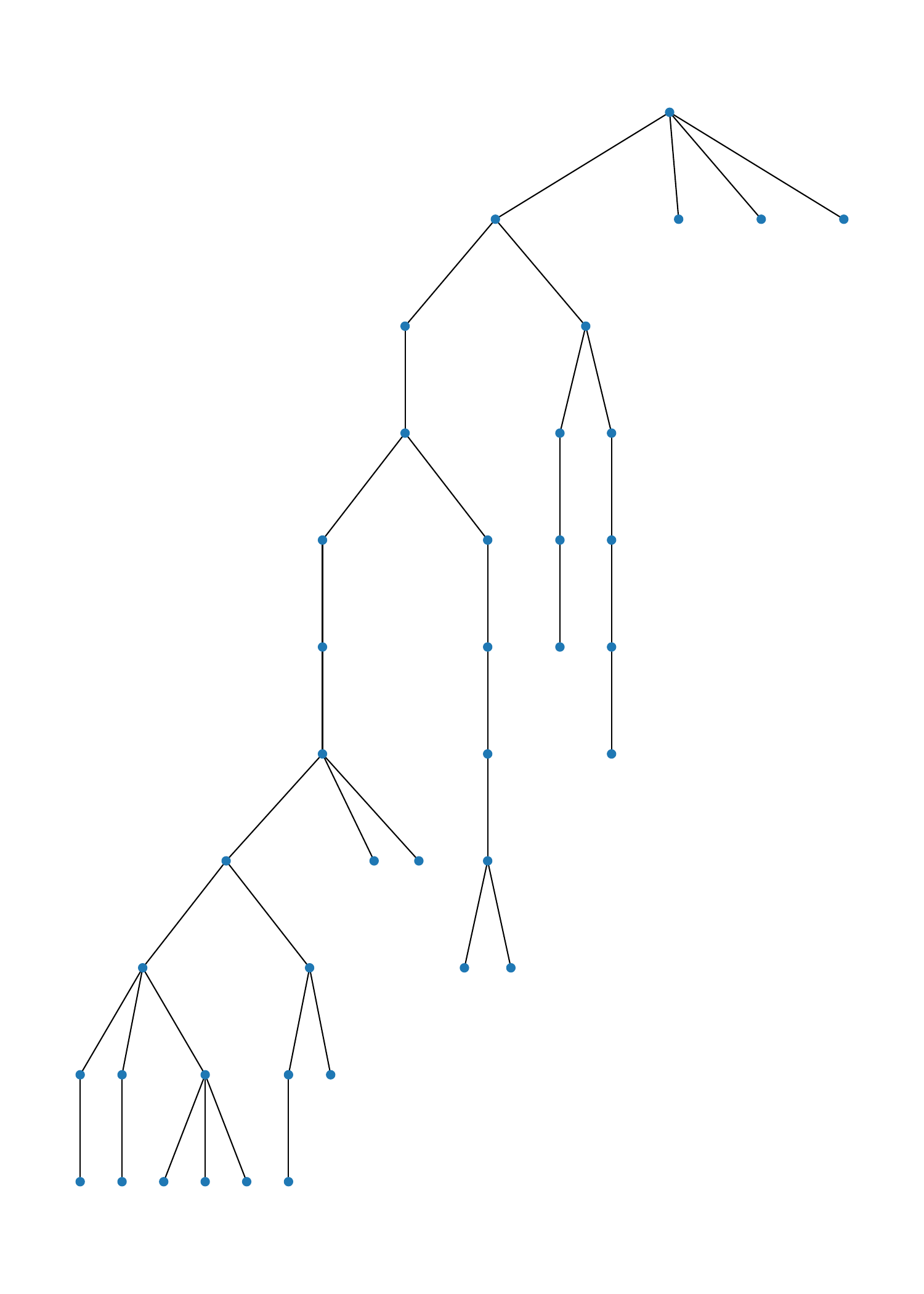}}\hspace*{-0.15cm}%
  \subcaptionbox{Random MST\label{fig:localMST}}{\includegraphics[width=0.5\textwidth]{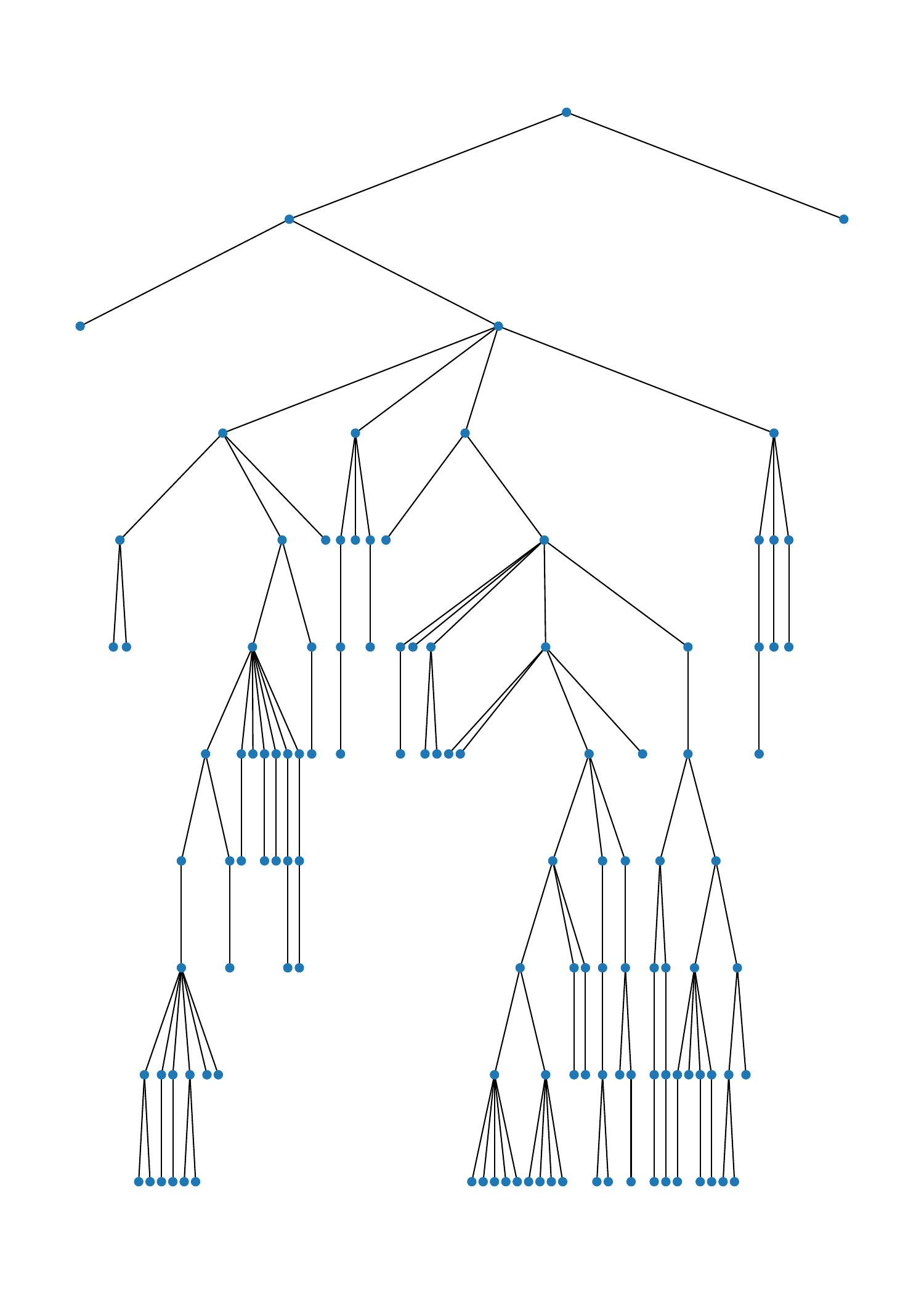}}
  \caption{A simulation of the ball of radius $10$ in the UST and (random) MST on 100,000 vertices. The MST typically has some vertices with large degrees, leading to a larger volume growth compared to the UST.}
  \label{fig:local}
\end{figure}

The local limit of the MST has been studied less intensively. Nonetheless, in \cite{Add13} it was proved that the MST on the complete graph converges to a local limit, and by showing that the limit has cubic growth compared to the quadratic growth of the  UST local limit, that the limiting object is different than $\mathcal{P}$. See also Figure \ref{fig:local} for a simulation of the local neighborhoods of the UST and the random MST. We state part of their result in the following theorem and note that this can also be deduced from \cite[Theorem 5.4]{AS04} and the more recent work of \cite{NT24}, which studies spectral and diffusive properties of the limiting object.

\begin{theorem} \label{T:local_MST}
    There exists a random tree rooted tree $(\mathcal{M} ,o)$ such that the MST on the complete graph rooted at $1$ locally converges to $(\mathcal{M}, o)$.
\end{theorem}

\noindent The construction of $\mathcal{M}$ is rather complicated and is linked to invasion percolation on the Poisson weighted infinite tree (PWIT). It is natural to ask the following question. 

\begin{question}
    Does there exist a ``back-bone'' construction of $\mathcal{M}$ similar to the Poisson(1) tree conditioned to survive forever? That is, to construct $\mathcal{M}$, take an infinite backbone $o=x_0, x_1, x_2, \ldots$ and independently attach to each $x_i$ a (finite) tree according to some measure.
\end{question}

%%%%%%%%%%%% NEW SECTION %%%%%%%%%%%%%%%%%%%%%%

\section{Low disorder} \label{S:localLow}

\subsection{Edge overlap}

To prove \eqref{eq:overlap_low}, the low disorder overlap of Theorem \ref{T:Overlap}, we show that the effective resistance concentrates around its mean. Denote by $|e \cap f|$ the number of endpoints that two edges $e,f \in E$ share. We have the following concentration bounds.

\begin{lemma}[cf.\ {\cite[Lemma 3.1 and Lemma 3.5]{Mak24}}] \label{L:effR_lowdisorder_both}
     Assume $\beta \ll n/ \log n$ and fix $K > 0$. Then the following holds with probability at least $1 - O(n^{-K})$:
     \begin{enumerate}[(i)]
        \item for all $u,v \in V$
            \begin{equation} \label{eq:effR_P_low}
                \effR{}{u}{v} = \frac{2(1 + o(1))}{n \xi};
            \end{equation}
        \item for all edges $e,f$ with $|e \cap f| = 0$
            \begin{equation} \label{eq:effR_ecapf=0_low}
                \effR{G / f}{u}{v} = \frac{2(1 + o(1))}{n \xi};
            \end{equation}
            
        \item and for all edges $e = (u,v)$ and $f$ with $|e \cap f| = 1$
            \begin{equation} \label{eq:effG_contract}
                \effR{G / f}{u}{v} = \frac{3(1 + o(1))}{2n \xi}.
            \end{equation}
     \end{enumerate}
\end{lemma}

\begin{proof}[Proof sketch of the lower bounds in \eqref{eq:effR_P_low}, \eqref{eq:effR_ecapf=0_low} and \eqref{eq:effG_contract}]
    \sloppy % Fix overflow
    The proof of the lower bounds is a simple application of the electric network theory developed in Section \ref{S:electric} together with a concentration inequality.
    Fix some edge $(u,v) \in E$. Rayleigh's monotonicity principle (Theorem \ref{T:Rayleigh}) gives that the effective resistance in $G$ is at least the effective resistance on the graph $G'$ where every edge between two vertices $x,y \not\in \{u,v\}$ is contracted. Notice that the graph $G'$ consists of the $3$ vertices $u,v$ and $z$, where $z$ represents all the contracted vertices. There is exactly one edge connecting the vertices $u$ and $v$, and there are $n-2$ parallel edges connecting $u$ to $z$ (similarly for $v$ to $z$). Using the series and parallel law on $G'$ shows that
    \begin{equation} \label{eq:up_con}
        \frac{1}{\weight(u,v) + \Big( \big(\sum_{y \neq u,v} \weight(u,y) \big)^{-1} + \big( \sum_{y \neq u,v} \weight(y,v) \big)^{-1} \Big)^{-1}}  = \effR{G'}{u}{v}.
    \end{equation}
    Applying the concentration inequality in Lemma \ref{L:low_concentration} yields \eqref{eq:effR_P_low}. The inequalities of \eqref{eq:effR_ecapf=0_low} and \eqref{eq:effG_contract} follow along the same lines.
    \fussy %
\end{proof}

For the upper bounds in \eqref{L:effR_lowdisorder_both}, we will require the following lemma about the effective resistance in the random graph when $p \gg \log n/ n$. We refer to \cite{Mak24} for a proof, and note that the inequality in \eqref{eq:Gnp_effR_upper} is essentially due to \cite{Jon98}.

\begin{lemma}[cf. {\cite[Lemma 3.2]{Mak24}}] \label{L:effGnp}
    Suppose $ \log n \ll np \ll \sqrt{n}$ and let $\epsilon = \epsilon(n) > 0$ be such that $\sqrt{\log n/ np} \ll \epsilon < 1/2$. There exists some universal constant $C > 0$ such that for any $K > 0$ the following holds with probability at least $1 - O(n^{-K})$: 
    \begin{enumerate}[(i)]
        \item for all $u,v \in V$
            \begin{equation} \label{eq:Gnp_effR_upper}
                \effR{G_{n,p}}{u}{v} \leq \frac{2(1 + C \epsilon)}{np};
            \end{equation}
        \item and for all edges $e = (u,v)$ and $f$ with $|e \cap f| = 1$
            \begin{equation} \label{eq:effGnp_contract}
                \effR{G_{n,p} / f}{u}{v} \leq \frac{3(1+C \epsilon)}{2np}.
            \end{equation}
    \end{enumerate}
\end{lemma}

Using Lemma \ref{L:effGnp} we may obtain a flow that estimates the effective resistance in the weighted graph up to constant factors. Namely, consider the environment coupled to the random graph as in Section \ref{SS:CoupleRG}. Let $p = 1/\beta \gg \log n/n$ and denote by $\theta$ the unit current flow (the flow minimizing the energy in \eqref{T:thompson}). Then lowering the weights of $p$-open edges to $\exp(-\beta p )$, gives with Rayleigh's monotonicity principle that
\begin{equation} \label{eq:flow_off_constant}
    \effR{}{u}{v} \leq \frac{(2+o(1))}{np} e^{\beta p} \approx \frac{2 \cdot e}{n \xi}.
\end{equation}

\noindent We will need one further straightforward lemma to remove the extra factor of $e$ in \eqref{eq:flow_off_constant}.
\begin{lemma} \label{L:disjoint_flow}
    Let $\theta_1, \ldots, \theta_m$ be unit flows from $u$ to $v$ that are supported on disjoint edges, that is for each edge $e$ we have $\theta_i(e) > 0$ for at most one $i \in \{1, \ldots, m\}$. Then, for any $\alpha_1, \ldots, \alpha_m \geq 0$
    \begin{equation*}
        \Big\Vert \sum_{i=1}^m \alpha_i \theta_i \Big\Vert^2_r = \sum_{i=1}^m \alpha^2_i \Vert\theta_i\Vert^2_r.
    \end{equation*}
\end{lemma}
\begin{proof}
    Using that the support of the flows are disjoint, we see that
    \begin{equation*}
        \Big\Vert \sum_{i=1}^m \alpha_i \theta_i \Big\Vert^2_r = \frac{1}{2}\sum_{e \in E} \Big( \sum_{i=1}^m \alpha_i \theta_i(e) \Big)^2 = \frac{1}{2}  \sum_{e \in E} \sum_{i=1}^m \alpha_i^2 \theta_i(e)^2 = \sum_{i=1}^m \alpha^2_i  \Vert\theta_i \Vert^2_r. \qedhere
    \end{equation*}
\end{proof} 

\sloppy % 
\begin{proof}[Proof sketch of the upper bounds in \eqref{eq:effR_P_low}, \eqref{eq:effR_ecapf=0_low} and \eqref{eq:effG_contract}]
    We couple snapshots of the random graph process $(G_{n,p})_{p \in [0,1]}$ with the environment $\omega$ as follows (see also Section \ref{SS:CoupleRG}). For $0 \leq p_0 \leq p_1 \leq 1$ let $G_n[p_0, p_1)$ be the subgraph consisting of vertices $V$ and edges satisfying
    \begin{equation*}
        e \in G_n[p_0, p_1) \iff \omega_e \in [p_0, p_1).
    \end{equation*}
    As $\beta \ll n/\log n$, we may take $p = \kappa(n) \log n /n$ for some function $\kappa = \kappa(n) \rightarrow \infty$ such that $\beta = o(p^{-1})$. Define $m := \lfloor p^{-1} \rfloor$, and consider the random graph at snapshots of the form
    \begin{equation*}
        p_j = j \cdot p \qquad j=1, \ldots, m.
    \end{equation*}
    
    Now for each $1 \leq j \leq m$, let $\theta_j$ be the unit current flow on the random graph $G_{n}[p_{j-1}, p_j]$, and define the scalars
    \begin{equation*}
        \tilde{\alpha}_j = e^{-\beta \cdot p_{j}} = e^{-j \beta p} \qquad j=1,\ldots, m.
    \end{equation*}
    Consider the unit flow
    \begin{equation*}
        \theta := \frac{1}{\sum_{j=1}^m \tilde{\alpha}_j} \sum_{j=1}^m \tilde{\alpha}_j \theta_j,
    \end{equation*}
    which is a convex combination of the flows $\theta_i$ giving much more weight to the ``early'' snapshots corresponding to small $j$. As the flows $\theta_i$ are supported on disjoint edges, using Lemma \ref{L:effGnp} and Lemma \ref{L:disjoint_flow}, one can show that with high $\P$-probability
    \begin{equation*}
        \effR{}{u}{v} \leq \frac{2(1 + o(1))}{n \xi}.
    \end{equation*}
    The upper bounds in \eqref{eq:effR_ecapf=0_low} and \eqref{eq:effG_contract} follow by the same argument.
\end{proof}

\fussy %

\noindent The results about the edge overlap in \eqref{eq:overlap_low} of Theorem \ref{T:Overlap} and the length \eqref{eq:length_low} in Theorem \ref{T:length}, for $\beta \ll n/\log n$, now follow readily.

\begin{proof}[Proof of \eqref{eq:overlap_low} in Theorem \ref{T:Overlap}]
    By Lemma \ref{L:effR_lowdisorder_both} we have with high probability that
    \begin{equation*}
        \sum_{e \in E} \weight(e)^2 \effR{}{e^-}{e^+}^2 = \frac{4(1 + o(1))}{n^2 \xi^2} \sum_{e \in E} \weight(e)^2.
    \end{equation*}
    Straightforward calculations give
    \begin{equation*}
        \E[ \weight(e)^2] = \frac{1-e^{-2\beta}}{2\beta},
    \end{equation*}
    so that the concentration inequality of Lemma \ref{L:low_concentration} for the above sum yield the equality in \eqref{eq:overlap_low}.
\end{proof}

\begin{proof}[Proof sketch of \eqref{eq:length_low} in Theorem \ref{T:length}]
    Notice that the expected (w.r.t.\ $\bE$) edge length satisfies
    \begin{equation*}
        \bE \big[ L(\cT) \big] = \sum_{(u,v) \in E} \omega_{(u,v)} \bP\big( (u,v) \in \cT \big) =  \sum_{(u,v) \in E} \omega_{(u,v)} e^{-\beta \omega_{(u,v)}} \effR{}{u}{v}.
    \end{equation*}
    Therefore, by Lemma \ref{L:effR_lowdisorder_both}, we have 
    \begin{equation}
        \bE \big[ L(\cT) \big] = (1+o(1)) \frac{2}{n \mu} \sum_{e \in E} \omega_{e} e^{-\beta \omega_e} \label{eq:low_length_concentrate}
    \end{equation}
    with probability, say, at least $1 - O(n^{-2})$. Let $X_{e} = \omega_e e^{-\beta \omega_e}$, then 
    \begin{equation*}
        \E[X_e] = \frac{1 - \beta e^{-\beta} - e^{-\beta}}{\beta^2},
    \end{equation*}
    and there exists a constant $C > 0$ such that for $\beta \geq 1$, we have 
    \begin{align*}
     \mathrm{Var} (X_e) \leq C \frac{1}{\beta^3}.
    \end{align*}
    One may complete the proof by applying Bernstein's inequality similarly to \eqref{eq:BernU}, with $m = n(n-1)/2$ and $m/\beta^2 \rightarrow \infty$, and using that the length is trivially bounded by $n$. See also Appendix \ref{AS:concentration}. 
\end{proof}

\subsection{Transfer-impedence}

Let $\bP_0 = \bP^\omega_{0,n}$ be the UST measure on the unweighted complete graph with corresponding transfer-impedence matrix $Y_0$. For two edges $e,f$ denote by $| e \cap f|$ the number of endpoints $e$ and $f$ share. Applying the series law and parallel law, one may easily show that 
\begin{equation*}
    | Y_0(e,f)| =
    \begin{dcases}
        \frac{2}{n} & \text{if } |e \cap f| = 2, \\
        \frac{1}{n} & \text{if } |e \cap f| = 1, \\
        0 & \text{if } |e \cap f| = 0.
    \end{dcases}
\end{equation*}
Using Lemma \ref{L:effR_lowdisorder_both}, we may obtain the following for the transfer-impedance matrix $Y$ with weights given by $\weight(e) = \exp(-\beta \omega_e)$.   

\begin{lemma}[cf.\ {\cite[Lemma 3.6 and Corollary 3.7]{Mak24}}] \label{L:Y(e,f)}
    Assume that $\beta = \beta(n) \ll n/\log n$ and let $K > 0$ be any constant. Then for all edges $e,f$ we have
    \begin{equation} \label{eq:Y(e,f)_abs}
        |Y(e,f)| = \frac{\weight(f)}{n \xi} \big( |e \cap f| + o(1) \big),
    \end{equation}
    with probability at least $1 - O(n^{-K})$.  In particular,
    \begin{equation} \label{eq:Y(e,f)-Y_0(e,f)}
        \Big| Y(e,f) - \frac{\weight(f)}{\xi} Y_0(e,f) \Big| = o(1) \cdot \frac{\weight(f)}{n\xi}
    \end{equation}
    with probability at least $1 - O(n^{-K})$. Furthermore, for fixed $k \geq 1$ and distinct edges $e_1, \ldots, e_k \in E$ we have
    \begin{equation} 
        \bP(e_1, \ldots, e_k \in \cT) = (1 + o(1)) \prod_{i=1}^k \frac{\weight(e_i)}{\xi} \bP_0(e_1, \ldots, e_k \in \cT). \label{eq:ek_in_T}
    \end{equation}
    with probability at least $1 - O(n^{-K})$.
\end{lemma}
\begin{proof}
    Recall from Section \ref{S:electric} that
    \begin{equation*}
        Y(e,e) = \bP(e \in \cT) = \weight(e) \effR{}{e^-}{e^+},
    \end{equation*}
    and for $e,f$ with $|e \cap f| \leq 1$
    \begin{align}
        Y(e,f)^2 &= \frac{\weight(f)}{\weight(e)} \big( \bP(e \in \cT) \bP(f \in \cT) - \bP(e,f \in \cT) \big) \nonumber \\
        &= \weight(f)^2 \effR{}{f^-}{f^+} \Big( \effR{}{e^-}{e^+}  - \effR{G / f}{e^-}{e^+} \Big). \label{eq:Y(e,f)^2}
    \end{align}
    %
    % Rayleigh's monotonicity principle shows that for edges with $|e \cap f| = 0$
    % \begin{equation*}
    %     \effR{G / f}{e^-}{e^+} \leq \effR{}{e^-}{e^+}.
    % \end{equation*}
    %  Therefore, 
    Lemma \ref{L:effR_lowdisorder_both} gives for $|e \cap f| \leq 1$ that
    \begin{align*}
        \effR{}{e^-}{e^+} &= (1 + o(1)) \frac{2}{ n \xi}, \\
        \effR{G / f}{e^-}{e^+} &= (1 + o(1)) \Big( \frac{4 - |e \cap f|}{2 n \xi} \Big),
    \end{align*}
    with probability at least $1 - O(n^{-K})$. Inserting this back into \eqref{eq:Y(e,f)^2} readily gives \eqref{eq:Y(e,f)_abs}.

    \smallskip 
    
    For equation \eqref{eq:Y(e,f)-Y_0(e,f)}, notice that the signs of $Y(e,f)$ and $Y_0(e,f)$ agree as long as $|e \cap f| \geq 1$. Namely, for any two edges $e = (e^-, e^+)$ and $f = (f^-, f^+)$, the flow along $f$ is non-negative if $f$ is pointing away from $e^-$ or towards $e^+$ (i.e.\ either $f^- = e^-$ or $f^+ = e^+$), while the flow along $f$ is  non-positive if $f$ is pointing towards $e^-$ or away from $e^+$ (i.e.\ either $f^+ = e^-$ or $f^- = e^+$). If on the other hand $|e \cap f| = 0$, then $Y_0(e,f) = 0$, so that the equality of \eqref{eq:Y(e,f)-Y_0(e,f)} readily follows.

    For the equality \eqref{eq:ek_in_T}, we may w.l.o.g.\ assume that $e_1, \ldots, e_k$ contains no cycles.  By the transfer-impedance theorem (Theorem \ref{T:transfer-impedance}) and the concentration of Lemma \ref{L:effR_lowdisorder_both}, we have
    \begin{align}
        \bP(e_1, \ldots, e_k \in \cT) &= \sum_{\sigma \in S_k} \Big( sgn(\sigma) \prod_{i=1}^k Y(e_i, e_{\sigma(i)})\Big) \nonumber \\
        &= \sum_{\sigma \in S_k} \bigg[ sgn(\sigma)  \prod_{i=1}^k \frac{\weight(e_i)}{\xi}\prod_{i=1}^k \Big(Y_0(e_i, e_{\sigma(i)}) + \frac{o(1)}{n} \Big) \bigg] \nonumber \\
        &= \prod_{i=1}^k \frac{\weight(e_i)}{\xi} \bigg[ \Big( \sum_{\sigma \in S_k}  sgn(\sigma) \prod_{i=1}^k Y_0(e_i, e_{\sigma(i)})  \Big) + k! \frac{k \cdot o(1)}{n^k}  \bigg] \nonumber \\
        &= \prod_{i=1}^k \frac{\weight(e_i)}{\xi} \Big( \bP_0(e_1, \ldots, e_k \in \cT) + k! \frac{k \cdot o(1)}{n^k} \Big), \label{eq:timp_intermediate}
    \end{align}
    with probability at least $1 - O(n^{-K})$. Using the spatial Markov property, we may obtain the bound
    \begin{align}
       \bP_0(e_1, \ldots, e_k \in \cT)  &= \prod_{i=1}^k \bP_0(e_i \in \cT \mid e_1, \ldots, e_{i-1} \in \cT) \nonumber \\
       &\geq \prod_{i=1}^k \frac{1}{i n} = \frac{1}{k! n^k} \gg k!\frac{k \cdot o(1)}{n^k}, \label{eq:k!_UST_bound}
    \end{align}   
    where the second line used that the probability of an edge being in the UST is lower bounded by the probability that a random walk started at one of the endpoints of $e_i$ (with at most $in$ neighbors) immediately uses the edge $e_i$. The equality in \eqref{eq:ek_in_T} now follows from \eqref{eq:timp_intermediate}.
\end{proof}

\subsection{Local limit via tree moments}

Consider a sequence of random variables $X_n$ taking values in $\R$. Suppose that $X_\infty$ is another random variable such that for all $r \geq 1$
\begin{equation*}
    \E[ X_n^r] \xrightarrow{n \rightarrow \infty } \E[X_\infty^r],
\end{equation*}
and that the moments of $X_\infty$ do not grow too fast, say, for instance, that
\begin{equation}
    \E[X_\infty^r] \leq e^{C r} \label{eq:method_of_moment}
\end{equation}
for some constant $C > 0$ and all $r \geq 1$. Then, the method of moments (see e.g.\ \cite[Section 3.3.5]{Dur19}) says that we must have weak convergence of the sequence $(X_n)_{n \in \N}$ to the random variable $X_\infty$. The authors in \cite{BP93} developed a similar approach using $t$\textsuperscript{th} tree moments to show the local convergence of trees, which we now describe.

%\footnote{ We remark that the convergence in \cite{BP93} considers fixed roots instead of uniformly chosen random roots as Section \ref{SS:local}, however, as the complete graph is transitive, these notions are equivalent.}. 

A \textit{tree-map} $f : V(t) \rightarrow V(G)$ is a function from the vertices of a rooted tree $(t,s_1)$ to the vertices of a rooted graph $(G,s_2)$ that satisfies the following properties:
\begin{enumerate}
    \item $f$ is injective;
    \item it maps the root $s_1$ to the root $s_2$;
    \item and if $u \sim v$ in $t$, then $f(u) \sim f(v)$ in $G$.
\end{enumerate}
Denote by $N(G,t)$ the number of distinct tree-maps between $(t,s_1)$ and $(G,s_2)$. Let $(\bP_n)_{n \geq 1}$ and $\bP_\infty$ be probability distributions on (rooted) trees with corresponding expectations $(\bE_n)_{n \geq 1}$ and $\bE_\infty$. The following is an analog of the method of moments for finite trees.

\begin{proposition}[Proposition 8.7 in \cite{BP93}] \label{P:tree_moment}
    Suppose that 
    \begin{equation*}
        \bE_n \big[ N(T, t) \big] \xrightarrow{n \rightarrow \infty} \bE_\infty \big[ N(T, t) \big] < \infty
    \end{equation*}
     for all finite trees $t$ and that $\bP_\infty$ is uniquely determined by the values of $\bE_\infty [N(T, t)]$. Then $\bP_n \xrightarrow{d} \bP_\infty$.
\end{proposition}

Similarly to \eqref{eq:method_of_moment}, if the expectation $\bE_\infty$ satisfies
\begin{equation*}
    \bE_\infty \big[ N(T, t) \big] \leq e^{C |t|},
\end{equation*}
for some $C > 0$ and all finite trees $t$, then $\bP_\infty$ is uniquely determined by the values of $\bE_\infty [N(T, t)]$ (see \cite[Section 8.3]{BP93}). In particular, as $\mathcal{P}$ the Poisson(1) tree conditioned to survive forever satisfies
\begin{equation*}
    \bE \big[ N(\mathcal{P}, t)\big]= |t|,
\end{equation*}
the law of $\mathcal{P}$ is uniquely determined by its $t$\textsuperscript{th} moments. Consequently, to prove the local convergence of $\cT$ in the case $\beta \ll n/\log n$ it suffices to prove convergence of the expectations of $N(\cT, t)$ to that of $\mathcal{P}$.

\begin{proof}[Proof of \eqref{eq:localRSTRE_UST} in Theorem \ref{T:local_limit}]
As in the proof of Theorem 5.3 in \cite{BP93}, for any finite rooted tree $t$ let $\phi(f)$ be the image (in $G$ rooted at $1$) of the tree map $f$. Let $K > |t|$, then by Lemma \ref{L:Y(e,f)}
\begin{equation*}
    \E\big[ \bP( \phi(f) \subseteq \cT) \big] = \big(1 + o(1) \big) \bP_0 \big( \phi(f) \subseteq \cT \big) + O(n^{-K}) = \big(1 + o(1) \big) \bP_0 \big( \phi(f) \subseteq \cT \big).
\end{equation*}
Hence, summing over all tree maps $f$
\begin{equation*}
    \E \Big[ \bE \big[ N(\cT, t) \big] \Big] = \E \Big[ \sum_f \bP \big( \phi(f) \subseteq \cT \big)  \Big] = \big(1+ o(1) \big) \sum_f \bP_0 \big( \phi(f) \subseteq \cT \big),
\end{equation*}
which converges to $\bE[ N(\mathcal{P}, t)]$ by Theorem 5.3 of \cite{BP93}. Applying Proposition \ref{P:tree_moment} finishes the proof.
\end{proof}

%%%%%%%%%%%%%%%%%%%%% NEW SECTION %%%%%%%%%%%%%%%%%%%%%%%

\section{High disorder} \label{S:localHigh}

In this regime, we assume that $\beta \gg n (\log n)^\lambda$, either for $\lambda = 2$ for the edge overlap, or $\lambda = \lambda(n) \rightarrow \infty$ arbitrarily slowly for the local limit.
\subsection{Edge overlap} \label{SS:low_edge}

First, we prove the overlap equality \eqref{eq:high_overlap} of Theorem \ref{T:Overlap}. The result is restated in the next lemma.
\begin{lemma}[{\cite[Lemma 4.1]{Mak24}}] \label{L:high_overlap}
    Let $\beta = \beta(n) \gg n (\log n)^2$, then with high probability
    \begin{equation} \label{eq:high_overlap}
        \mathcal{O}(\beta) = (1 - o(1)) n.
    \end{equation}
\end{lemma}

\begin{proof}
    We claim that
    \begin{equation*}
        \E[ \mathcal{O}(\beta)] = \frac{n(n-1)}{2} \E \big[ \bP(e \in \cT)^2 \big] = (1-o(1)) n,
    \end{equation*}
    where the first equality holds by transitivity for any edge in $E$. If the above is true, then, as the overlap squared is trivially upper bounded by $n^2$, the variance of $\mathcal{O}(\beta)$ must be $o(n^2)$. Chebyshev's inequality immediately gives \eqref{eq:high_overlap}. Hence, it suffices to prove that for an edge $e = (u,v)$
    \begin{equation}  \label{eq:expectation_overlap}
        \E \big[ \bP(e\in \cT)^2 \big] = \E \big[ \exp(-2\beta \omega_e) \effR{}{u}{v}^2 \big] = (1-o(1)) \frac{2}{n}. 
    \end{equation}

    To prove \eqref{eq:expectation_overlap}, we will define several events that imply that an edge in the MST (coupled to $\omega$ as in Section \ref{SS:CoupleRG}) is also an edge in the RSTRE with high probability. Consider an edge $e = (u,v) \in E$, and let $p = \omega_e$ and $p^- := \omega_e - e^{-n}$. It is useful to think of $p^-$ as the time in the random graph process just before the edge $e$ is added. We define
    \begin{equation} \label{eq:defG_e}
        \mathcal{G}_e := \big\{ G_{n, p} = G_{n, p^-} \cup \{ e \} \big\}
    \end{equation}
    as the event that between $p^-$ and $p$ no edge other than $e$ is added to the random graph. Notice that the length of the interval $[p^-, p]$ is $e^{-n}$ and that there are of order $n^2$ many edges, from which we may obtain that
    \begin{equation*}
        \P( \mathcal{G}_e^c) = O(n^{-4}).
    \end{equation*}

    Next, let $\epsilon \geq n^{-1/4}$ be a sequence tending to zero (e.g.\ $\epsilon = 1/\log \log n$), and consider the following two events
    \begin{align}
        % \mathcal{A}_e = \mathcal{A}_e(\epsilon) &:= \{ \omega_e \not\in [1 - \frac{1 - \epsilon}{n}, 1 - \frac{1 + \epsilon}{n}]\} \\
        \mathcal{M}_e &:=  \big\{ e \in \textrm{MST} \big\} \cap \big\{ \omega_e \leq \frac{5 \log n}{n} \big\}, \label{eq:defM_E}\\
        \mathcal{S}_e = \mathcal{S}_e(\epsilon) &:= \big\{ |\mathcal{C}_2(p^-)| \leq L \epsilon^{-2} \log n \big\}, \label{eq:defS_e}
    \end{align}
    for some $L$ large enough but fixed constant as in Theorem \ref{T:C2_size}. By Theorem \ref{T:Gnp_connect}, the graph $G_{n,5\log n/n}$ is connected with very high probability, so that
    \begin{equation*}
        \P( \mathcal{M}_e) = \frac{2}{n} - O(n^{-4}).
    \end{equation*}
    Furthermore, let 
     \begin{equation*}
        E^{*} := \Big\{ e \in E \, : \, \omega_e \not \in \Big[\frac{1 - \epsilon}{n}, \frac{1 + \epsilon}{n} \Big] \Big\},
    \end{equation*}
    then using bounds on the size of the 2nd largest cluster of Theorem \ref{T:C2_size} gives that
    \begin{equation}
        \P( \mathcal{S}_e) \geq \P\big( \mathcal{S}_e \mid e \in E^{*} \big) \big(1 - \frac{2 \epsilon}{n} \big)  = 1 - \frac{2 \epsilon}{n} - O(n^{-\eta}), \label{eq:S_e_bound}
    \end{equation}
    where $\eta$ is the constant from Theorem \ref{T:C2_size}.

    For the remainder of the proof, we restrict to the events $\mathcal{M}_e$, $\mathcal{S}_e$ and $\mathcal{G}_e$. Denote by $S(u)$ and $S(v)$ be the connected components of the random graph process at $p^-$ and w.l.o.g.\ assume that $|S(u)| \leq |S(v)|$. As $\mathcal{M}_e$ and $\mathcal{G}_e$ hold, the edge $e$ must be the edge with the lowest weight in the cutset $E(S(u), S(u)^c)$ separating $u$ and $v$. The Nash-Williams inequality (Lemma \ref{L:nash_williams}) together with Kirchhoff's formula then shows that
    \begin{equation} \label{eq:nash_component}
        \bP(e \in \cT) \geq \frac{\weight(e)}{\sum_{f \in E(S(u), S(u)^c)} \weight(f)} = \frac{1}{1 + \sum_{f \in E(S(u), S(u)^c) \setminus \{e\}} e^{-\beta (\omega_f - \omega_e)} }.
    \end{equation}
    We will use that on the event $\mathcal{S}_e$, the cutset is bounded in size by
    \begin{equation*}
        \big| E(S(u), S(u)^c) \big| \leq n |\mathcal{C}_2(p)| \leq L n \epsilon^{-2} \log n.
    \end{equation*}
    The edges in $E(S(u), S(u)^c) \setminus \{e\}$ are random variables uniformly distributed on $[0,1]$ conditioned to be larger than $\omega_e$, so that the random variables in the collection
    \begin{equation*}
        \big( \frac{\omega_f - \omega_e}{1- \omega_e} \big)_{f \in E(S(u), S(u)^c) \setminus \{e \}}
    \end{equation*}
    are i.i.d.\ uniform random variables on $[0,1]$. 
    
    Lastly, we define 
    \begin{align}
        \mathcal{B}_e  &:= \Big\{\forall f \in E \big( S(u), S(u)^c \big) \setminus \{e \} \, : \, e^{-\beta(\omega_f - \omega_e)} < \frac{1}{n^5} \Big\} \nonumber \\
        &= \Big\{ \forall f \in E \big( S(u), S(u)^c \big) \setminus \{e \} \, : \, \frac{\omega_f-\omega_e}{1-\omega_e} >  5 \frac{\log n}{ \beta( 1-\omega_e)} \Big\}, \label{eq:defB_e}
    \end{align}
    which is the event that there is a large enough gap between $\omega_e$ and any second smallest weight in $ E(S(u), S(u)^c)$. Let $\gamma > 0$ be such that
    \begin{equation} \label{eq:condition_beta_eps}
        \frac{n \epsilon^{-2} (\log n)^2 }{\beta} \leq \gamma,
    \end{equation}
    and let $U$ be a uniform random variable on $[0,1]$. It follows that on the event $\mathcal{S}_e$ and $\{ \omega_e \leq 5 \log n/n \}$, we have
     \begin{equation*}
        \P \big( \mathcal{B}_e \big) \geq \P \Big( U > 6 \frac{\log n}{ \beta} \Big)^{L n \epsilon^{-2} \log n} = \Big(1 - 6 \frac{\log n}{ \beta} \Big)^{L n \epsilon^{-2} \log n},
    \end{equation*}
    which is at least $1 - O(\gamma)$ if the condition \eqref{eq:condition_beta_eps} is true. We now use the inequality \eqref{eq:nash_component} under the restriction that the event $\mathcal{B}_e$ holds. Namely, on the event $\mathcal{B}_e$, the right hand side of \eqref{eq:nash_component} can be further lower bounded by $1/(1+ n^{-3})$, so that
    \begin{equation} \label{eq:e_in_cT_MST}
        \bP(e \in \cT) \geq \frac{1}{1 + n^{-3}} \geq 1 - n^{-3}.
    \end{equation}

    We can now finally complete the proof of \eqref{eq:expectation_overlap}. As we assume that $\beta \gg n (\log n)^2$, we may pick $\gamma = o(1)$ and $\epsilon = o(1)$ such that the condition \eqref{eq:condition_beta_eps} is satisfied. 
    Conditioning on the events $\mathcal{M}_e$, $\mathcal{S}_e$, $\mathcal{G}_e$ and $\mathcal{B}_e$, then gives for $\gamma \geq n^{-\eta}$ that
    \begin{align*}
        \E \big[ \bP(e \in \cT)^2 \big] &\geq (1 - n^{-3})^2 \P \big( \mathcal{B}_e \ \big| \ \mathcal{M}_e, \mathcal{S}_e, \mathcal{G}_e \big) \P\big( \mathcal{M}_e, \mathcal{S}_e, \mathcal{G}_e \big) \\
        &\geq \big(1 - O(\gamma) \big) \big(1 - \epsilon \big) \frac{2}{n},
    \end{align*}
     completing the proof.
\end{proof}

We briefly sketch how to obtain the equality \eqref{eq:length_high} in the high disorder regime of Theorem \ref{T:length}. We also refer to \cite[Theorem 1.8]{K24}, where they only require that $\beta \gg n \log n$ instead of $\beta \gg n (\log n)^5$, see also the proof below for yet another alternative approach.

\begin{proof}[Proof sketch of \eqref{eq:length_high} in Theorem \ref{T:length}]
    Denote by $M$ the MST coupled to the environment as in Section \ref{SS:CoupleRG}. Denote by $|\cT \setminus M|$ the number of edges that are in $\cT$ but not in $M$. Then, since the MST is the tree with the smallest Hamiltonian/length, we have
    \begin{equation*}
        L(M) \leq L(\cT) \leq L(M) + |\cT \setminus M| \max_{e \in \cT} \omega_e.
    \end{equation*}
    The proof of Theorem \ref{L:high_overlap} shows that, with high probability, the overlap between $M$ and $\cT$ satisfies
    \begin{equation*}
        \bE\big[ |M \cap \cT| \big] \geq \big(1 - O(\max \{ \epsilon, \gamma \}) \big) n,
    \end{equation*}
    where $\epsilon$ and $\gamma$ are as in \eqref{eq:condition_beta_eps}. Therefore,
    \begin{equation*}
        |\cT \setminus M| = O(\max \{ \epsilon, \gamma \}) n
    \end{equation*}
    with high probability. As $G_{n,p}$ is connected with high probability when $p = 2 \log n/n$, one can show (e.g.\ using Lemma \ref{L:gap}) that $\max_{e \in \cT} \omega_e \leq 3 \log n/n$ with high probability. Equation \eqref{eq:length_high} now follows by letting $\epsilon, \gamma = o( (\log n)^{-1})$, which is possible provided that $\beta \gg n (\log n)^5$.
\end{proof}

The following (easy) alternative proof only requires the milder condition $\beta \gg n \log n$ and has not appeared elsewhere.

\begin{proof}[Alternative proof of \eqref{eq:length_high} in Theorem \ref{T:length}.]
    Recall from Lemma \ref{L:derivative} that 
    \begin{equation*}
         \frac{\partial }{\partial \beta} \bP_{\beta}(\cT = T) = \bP_{\beta}(\cT = T) \big(\bE_{\beta}[H(\cT, \omega)] - H(T, \omega) \big),
    \end{equation*}
    and hence
    \begin{equation*}
         \bP_{2 \beta}(\cT = T) =  \bP_{\beta}(\cT = T) \exp \big( \int_\beta^{2\beta} \bE_{x}[H(\cT, \omega)] - H(T, \omega) dx \big).
    \end{equation*}
    Furthermore, notice that as 
    \begin{equation*}
        \frac{\partial}{\partial \beta} \bE_{\beta} [ H(\cT)] = - \textrm{Var}(H(\cT)) \leq 0,
    \end{equation*}
    the function $\bE_{\beta} [ H(\cT)]$ is non-increasing in $\beta$. Let $T_1$ be the minimum spanning tree coupled to $\omega$. Then, as Cayley's formula states that there are $n^{n-2}$ spanning trees on the complete graph, we obtain that
    \begin{align*}
       \exp \Big( \beta \big( \bE_{2 \beta}[H(\cT, \omega)] - H(T_1, \omega) \big) \Big) &\leq \exp \Big( \int_\beta^{2\beta} \bE_{x}[H(\cT, \omega)] - H(T_1, \omega) dx \Big)  \\
       &= \frac{\bP_{2\beta}(\cT = T_1)}{\bP_{\beta}(\cT = T_1)} \leq n^{n-2}.
    \end{align*}
    Taking the logarithm on both sides yields that
    \begin{equation*}
        \bE_{2 \beta}[H(\cT, \omega)] \leq \frac{(n-2) \log n}{\beta} + H(T_1, \omega),
    \end{equation*}
    from which \eqref{eq:length_high} follows whenever $\beta \gg n \log n$.
\end{proof}

In the above proof we used the rather trivial bound that the ratio of probabilities is upper bounded by $n^{n-2}$. We do not expect this to be sharp. If one were to obtain a better bound, it might be possible to improve the local limit result to the regime $\beta \gg n$ instead of $\beta \gg n \log n$. For instance, using \eqref{eq:length_low} one might be able to show that
\begin{equation*}
    \int_0^{o(n/\log n)} \bE_x[ H(\cT, \omega)] dx = (1 - o(1) ) n \log n,
\end{equation*}
from which it follows that
\begin{equation*}
    %\bP_{\beta}(\cT = \textrm{MST}) &\geq \frac{1}{(\log n)^{n-2}}, \\
    \bE_{\beta}[H(\cT, \omega)] = o \Big( \frac{n \log n}{\beta} \Big) + H(T_1, \omega),
\end{equation*}
whenever $\beta \geq n/\log n$. 
\begin{question}
    Of what order is $\bP_{\beta}(\cT = \textrm{MST})$ when $\beta \ll n$, $\beta = C n$, and $\beta \gg n$? In particular, can we give non-trivial bounds on the ratio
    \begin{equation*}
        \frac{\bP_{2\beta}(\cT = \textrm{MST})}{\bP_{\beta}(\cT = \textrm{MST})} \, ?
    \end{equation*}
\end{question}

\subsection{Local limit}

Recall that $B_\cT(u,r)$ and $B_M(u,r)$ are the closed balls of radius $r$ centered at $u$ in the RSTRE and MST, respectively. The proof of the local convergence in the high disorder regime follows by again constructing similar events as in the proof of Lemma \ref{L:high_overlap} that hold for all outgoing edges of a fixed vertex. We refer to \cite[Lemma 4.2]{Mak24} for the full proof details.

\begin{lemma}[{\cite[Lemma 4.2]{Mak24}}] \label{L:cT_M_deg}
    For any $u \in V$ we have that
    \begin{equation} \label{eq:cT_M_deg}
         \E \Big[ \bP \big( B_\cT(u,1 ) \neq B_M(u,1) \big) \Big]  = O\Bigg( \max\Bigg\{ \Big( \frac{n (\log n)^5}{\beta} \Big)^{1/3}, n^{-\eta/2} \Bigg\}\Bigg),
    \end{equation}
    where $\eta$ is the constant from Theorem \ref{T:C2_size}.
\end{lemma}

\begin{proof}[Proof sketch]
    We first show that the probability of $B_M(u,1) \subseteq B_\cT(u,1)$ is close to $1$ under the averaged law $\widehat{\P}$. For each outgoing edge $(u,v)$ of $u$ consider again the events $\mathcal{G}_{(u,v)}$, $\mathcal{M}_{(u,v)}$, $\mathcal{S}_{(u,v)}$ and $\mathcal{B}_{(u,v)}$ as in \eqref{eq:defG_e}, \eqref{eq:defM_E}, \eqref{eq:defS_e} and \eqref{eq:defB_e}. The arguments leading to up \eqref{eq:e_in_cT_MST} give that
    \begin{equation} \label{eq:uv_not_T_indicators}
        \sum_{v \neq u} \bP\big( (u,v) \not\in \cT \big) 1_{\mathcal{M}_{(u,v)}} 1_{\mathcal{G}_{(u,v)}} 1_{\mathcal{S}_{(u,v)}} 1_{\mathcal{B}_{(u,v)}}  \leq n^{-2}.
    \end{equation}
    
    Now let 
    \begin{equation*}
        \gamma = \epsilon = \max\Big\{ \Big( \frac{n (\log n)^2}{\beta} \Big)^{1/3}, n^{-2\eta/3}\Big\},
        %\epsilon &= \log n \gamma = \log n \Big( \frac{n}{\beta} \Big)^{1/3}
    \end{equation*}
    which are chosen in such a way that the condition in \eqref{eq:condition_beta_eps} holds (where we may implicitly assume that $\epsilon = o(1)$). Using suitable probability bounds on the events $\mathcal{G}_{(u,v)}$, $\mathcal{M}_{(u,v)}$, $\mathcal{S}_{(u,v)}$ and $\mathcal{B}_{(u,v)}$ similar to the proof of Lemma \ref{L:high_overlap}, one may show that
    \begin{equation}
        \E \Big[ \bP\big( B_M(u,1) \not\subseteq B_\cT(u,1) \big)  \Big]  = O(\gamma). \label{eq:M_sub_T}
    \end{equation}
    That is, with high probability all neighbors of $u$ in $M$ are also neighbors of $u$ in $\cT$. 

    Consider now the degrees $\deg_M(u) = |B_M(u,1)|$ and $\deg_\cT(u) = |B_\cT(u,1)|$ of $u$ in $M$ and $\cT$. The equality \eqref{eq:M_sub_T} shows that $\deg_\cT(u)$ almost stochastically dominates $\deg_M(u)$, however, by transitivity, we have
    \begin{equation*}
        \E \big[ \bE [\deg_\cT(u)] \big] = \sum_{v \neq u} \widehat{\P} \big( (u,v) \in \cT \big) = \frac{2(n-1)}{n} = \E [\deg_M(u)],
    \end{equation*}
    i.e.\ the degrees have the same expectation (under the joint law). This is only possible if we also have that $B_\cT(u,1) \subseteq B_M(u,1)$ with high probability. Indeed, using \eqref{eq:M_sub_T} we have
    \begin{align}
        \E \Big[ \bP \big( B_\cT(u,1) \not\subseteq B_M(u,1) \big) \Big] &\leq \E \Big[ \bP \big( B_\cT(u,1) \not\subseteq B_M(u,1), B_M(u,1) \subseteq B_\cT(u,1) \big) \Big] \nonumber\\
        &\qquad + \E \Big[ \bP \big( B_M(u,1) \not\subseteq B_\cT(u,1) \big) \Big] \nonumber \\
        &\leq \E \Big[ \bP \big( \deg_\cT(u) > \deg_M(u) \big) \Big] + O(\gamma). \label{eq:B_T_notsubset}  
    \end{align}
    To complete the proof, it thus suffices to show that $\E \big[ \bP \big( \deg_\cT(u) > \deg_M(u) \big) \big]$ is small, which we do by considering the cases whether $\deg_\cT(u)$ is larger than $\deg_M(u)$ or not. Namely,
    \begin{align}
        0 &= \E \Big[ \bE\big[ \big(\deg_\cT(u) - \deg_M(u)\big) 1_{ \deg_\cT(u) > \deg_M(u)} \big] \Big] \nonumber \\
        &\qquad \qquad + \E \Big[ \bE \big[ \big(\deg_\cT(u) - \deg_M(u) \big) 1_{ \deg_\cT(u) < \deg_M(u)}\big] \Big]  \nonumber \\
        &\geq \E \big[ \bP \big( \deg_\cT(u) > \deg_M(u) \big) \big] - 60 \log n \E \big[ \bP \big( \deg_\cT(u) < \deg_M(u) \big) \big] \nonumber \\
        &\qquad \qquad - n \P \big( \deg_M(u) > 60 \log n \big).  \label{eq:deg_split_case}
    \end{align}
    Rearranging \eqref{eq:deg_split_case} and using the bound on the largest degree of a vertex in $M$ from Theorem \ref{L:max_deg_MST}, finishes the proof.
\end{proof}

We now extend the $r=1$ case from Lemma \ref{L:cT_M_deg} to all $r \geq 1$ by using exchangeability and a union bound.

\begin{lemma} \label{L:cT_ball_M}
    Suppose that $\beta \gg n (\log n)^\lambda$ for $\lambda = \lambda(n) \rightarrow \infty$ arbitrarily slowly. Then for any fixed $r \geq 1$
    \begin{equation*}
         \E \Big[ \bP \big( B_\cT(1,r) \neq B_M(1,r) \big) \Big] = o(1).
    \end{equation*}
\end{lemma}

\begin{proof}
    Let $r \geq 1$ be given. If $B_M(1,r) \neq B_\cT(1,r)$, then there must exist at least one vertex $u$ in $B_M(1,r)$ such that $B_M(u,1) \neq B_\cT(u,1)$. Taking a union bound over all vertices shows that 
    \begin{align}
        \E \Big[ \bP \big( B_\cT(1,r) \neq B_M(1,r) \big) \Big] &\leq
        \E \Big[ \bP \big( B_\cT(1,1) \neq B_M(1,1) \big) \Big] \nonumber  \\
        &\quad + \sum_{u \neq 1} \E \Big[ \bP \big( B_\cT(u,1) \neq B_M(u,1) \big) 1_{u \in B_M(1,r)}\Big] .\label{eq:B_neq_M_union}
    \end{align}
    The following observations are key: (1) a vertex $u$ is in $B_M(1,r)$ if and only if $1$ is in $B_M(u,r)$; (2) by exchangeability, each term appearing in the sum of \eqref{eq:B_neq_M_union} is identical for any two vertices $u,v \neq 1$; and (3) again by exchangeability, we may replace the vertex $1$ in the sum of \eqref{eq:B_neq_M_union} by any other vertex $v \neq u$. Hence, for any vertex $u \neq 1$
    \begin{align*}
        & \sum_{u \neq 1 }\E \Big[ \bP \big( B_\cT(u,1) \neq B_M(u,1) \big) 1_{u \in B_M(1,r)}\Big] \\
        & \hspace{3cm} =  (n-1) \E \Big[ \bP \big( B_\cT(u,1) \neq B_M(u,1) \big) 1_{ 1 \in B_M(u,r)}\Big] \\
        &\hspace{3cm} = (n-1) \frac{1}{n-1} \sum_{v \neq u}\E \Big[ \bP \big( B_\cT(u,1) \neq B_M(u,1) \big) 1_{ v \in B_M(u,r)}\Big].
    \end{align*}
    Taking the sum inside the expectation we see that 
    \begin{multline*}
        \sum_{v \neq u}\E \Big[ \bP \big( B_\cT(u,1) \neq B_M(u,1) \big) 1_{ v \in B_M(u,r)}\Big] \\
        = \E \Big[ \bP \big( B_\cT(u,1) \neq B_M(u,1) \big) \big(\big|B_M(u,r)\big|-1 \big)\Big].
    \end{multline*}
    Using Lemma \ref{L:max_deg_MST} we can bound the size of $|B_M(u,r)|$ by $60^r (\log n)^r$ with high probability, so that with \eqref{eq:B_neq_M_union} and exchangeability we have
     \begin{multline}
         \E \Big[ \bP \big( B_\cT(1,r) \neq B_M(1,r) \big) \Big]  \\
         \leq 60^r (\log n)^r \E \Big[ \bP \big( B_\cT(1,1) \neq B_M(1,1) \big) \Big] + O(n^{-4}). \label{eq:B_M_size_bound}
    \end{multline}
    Applying Lemma \ref{L:cT_M_deg} together with the assumption that $\beta \gg n (\log n)^\lambda \gg n (\log n)^{5 + 3r}$ finishes the proof.
\end{proof}

The proof of the local limit in the high disorder regime of Theorem \ref{T:local_limit} now follows easily.

\begin{proof}[Proof of \eqref{eq:localRSTRE_MST} in Theorem \ref{T:local_limit}]
    Suppose $t$ is some finite rooted tree of height $r$. Lemma \ref{L:cT_ball_M} gives that
    \begin{align*}
        \widehat{\P}( B_\cT(1,r) \simeq t) &\leq \P( B_M(1,r) \simeq t) + \E \Big[ \bP \big( B_\cT(1,r) \simeq t) 1_{B_M(1,r) \not\simeq t} \Big] \\
        &\leq \P( B_M(1,r) \simeq t) + \E \Big[ \bP \big( B_\cT(1,r) \neq B_M(1,r) \big) \Big] \\
        &= \P( B_M(1,r) \simeq t) + o(1)
    \end{align*}
    and
    \begin{align*}
        \widehat{\P}( B_\cT(1,r) \simeq t) &\geq \E \Big[ \bP \big( B_\cT(1,r) \simeq t) 1_{B_M(1,r) \simeq t} \Big] \\
        &\geq \P( B_M(1,r) \simeq t) - \E \Big[ \bP \big( B_\cT(1,r) \neq B_M(1,r) \big) \Big] \\
        &= \P( B_M(1,r) \simeq t) - o(1).
    \end{align*}
    As the MST on the complete graph has a local limit (see Theorem \ref{T:local_MST}), the RSTRE must locally converge to the same limiting object.
\end{proof}

\begin{remark}
    To obtain \eqref{eq:B_M_size_bound} we used the neighborhoods $B_M(1,r)$ grow (with high probability) in size at most as fast $60^r (\log n)^r$. We expect this bound to be nowhere close to optimal, in fact, as shown in \cite{Add13}, the typical growth of $|B_M(1,r)|$ is of order $r^3$ which is independent of $n$. With this one might expect the same proof of the local limit to hold when $ \lambda = \lambda(n)$ is some fixed large fixed constant. However, currently we are not aware of any better tail bounds, say $|B_M(1,r)| = O(\log n)$, that hold with large enough probability to make \eqref{eq:B_M_size_bound} useful. 
\end{remark}

% MAYBE WOULD BE INTERESTING TO SAY THE FOLLOWING:
% include something about taking the highest weighted edge, but this does not generate local neighborhood

% We conclude this chapter with a remark on the edge weights selected by Prim's algorithm in constructing the MST.

% \begin{remark}
%     % To construct the local neighborhood of a vertex $v$ in the MST, suppose that we run  Prim's algorithm (see Section \ref{S:MSTintro}) for $T$ many steps. 
%     Adapting the proof of Theorem \ref{T:highDisorder}, we can show that, with high probability, the first $T \gg 1$ edges discovered by Prim's algorithm (see Section \ref{S:MSTintro}) are also contained in the RSTRE. However, this does not imply that the local neighborhoods of the MST and RSTRE around $v$ coincide with each other. Indeed, by Proposition 3.1 of \cite{Add13} (see also \cite{ABG12}), if $T = o(\sqrt{n})$, then the resulting local limit is that of invasion percolation on the Poisson weighted infinite tree, which has a quadratic volume growth, as opposed to the cubic volume growth of the local limit of the MST.
%     The following is a brief heuristic of why this is the case. Once Prim's algorithm discovers a vertex in the giant component, which it does quickly, everything in the giant component is explored before returning back to the local neighborhood of $v$.  WRITE BETTER!!!
% \end{remark}

\chapter{Diameter of RSTREs}
\label{ch:diam}

Having studied local observables of the RSTRE in the previous chapter, we now move to the global observable of the diameter. This chapter focuses on 3 different families of graphs: in Section \ref{S:toy_model}, we study bounded-degree expanders and finite boxes in $\Z^d$, whereas in Section \ref{S:high_and_low} onwards, we restrict ourselves to the complete graph. The main results of this chapter are Theorem \ref{T:toy_model} and Theorem \ref{T:main}, which correspond to our work in the papers \cite{MSS23} and \cite{Mak24}, respectively. Although we slightly generalize the results of \cite{MSS23} and \cite{Mak24}, we will often only sketch the proof ideas and refer to the relevant references for the full details. 

\section{Sufficient conditions for bounding the diameter} \label{SS:DiameterCondition}

Let $G$ be a finite graph and denote by $d_G(u, v)$ the graph distance (that is the length of the shortest path) between $u$ and $v$ in $G$. The diameter of $G$ is defined as
\begin{equation*}
    \diam(G) := \max\limits_{u,v \in G} d_G(u,v),
\end{equation*}
and we are interested in how $\diam(\cT^\omega_{G, \beta})$ behaves as we vary $\beta$ or choose a different environment distribution $\mu$. It is a classical fact \cite{Sze83} that the diameter of the UST on $K_n$ behaves as $\sqrt{n}$, and more recently in \cite{MNS21}, the authors showed that there exist general conditions on the lazy simple random walk on $G$ (recall the connection between random walks and USTs from Section \ref{S:Wilson}) that imply that the diameter is of order $\sqrt{n}$. In the following, we state an analog of the result in \cite{MNS21} for weighted graphs $(G, \weight)$.

A weighted graph $(G, \weight)$ is called
balanced, mixing and escaping with parameters $D, \alpha, \theta > 0$ if the following are satisfied:
\begin{enumerate}
    \item $(G, \weight)$ is balanced:
    \begin{equation} \label{eq:Nbalanced}
       \frac{\max_{u\in V} \pi(u)}{\min_{u\in V} \pi(u)}= \frac{\max_{u \in V} \sum_{v} \weight(u,v)}{\min_{u \in V} \sum_{v} \weight(u,v)} \leq D;
    \end{equation}
    \item $(G, \weight)$ is mixing:
    \begin{equation} \label{eq:Nmixing}
        \tmix(G, \weight) \leq n^{\frac{1}{2} - \alpha};
    \end{equation}
    \item $(G, \weight)$ is escaping:
    \begin{equation} \label{eq:Nescaping}
        \sum_{t=0}^{\tmix} (t+1) \sup_{v \in V}q_t(v,v) \leq \theta.
    \end{equation}
\end{enumerate}
The following theorem is an extension of the main theorem in \cite{MNS21}, where compared to \cite{MNS21}, we also allow weights on the edges and for the parameters $D$ and $\theta$ to depend on $n$. 

\begin{theorem}[{Theorem 2.3 of \cite{MSS23}}]\label{T:Nach}
For any $\alpha>0$, there exist $C, k, \gamma > 0$ such that if $(G, \weight)$ satisfies conditions \eqref{eq:Nbalanced}-\eqref{eq:Nescaping} with $\alpha$ and $D,\theta \leq n^{\gamma}$, then for any $\delta > n^{-\gamma}$,
\begin{equation}\label{eq:Nach1}
    \bP^\weight \big( (CD \theta \delta^{-1})^{-k} \sqrt{n} \leq {\diam}(\cT) \leq (CD \theta \delta^{-1})^{k} \sqrt{n} \big) \geq 1 - \delta.
\end{equation}
\end{theorem}

%%%%%%%%%%%%%%% NEW SECTION %%%%%%%%%%%%%%%%%%%%%%%

\section{Bounded-degree graphs with fixed weight distributions} \label{S:toy_model}

In this section, assume as usual that $\mu(-\infty, \infty) = 1$, however, we fix\footnote{In \cite{MSS23}, $\beta$ is implicit in the weight distribution which is supported on $(0, \infty)$. Furthermore, the definition of edge expansion is slightly different, however, as the graphs have bounded degrees, these agree up to some multiplicative constant.} $\beta = 1$. Denote by $\bP = \bP^\omega_{G, 1}$ the law of the RSTRE given the environment $(\omega_e)_{e \in E(G)}$. In \cite{MSS23}, we proved the following theorem, where we recall the notion of expander graphs from Definition \ref{D:expander}.

\begin{theorem}[Theorem 1.1 in \cite{MSS23}]\label{T:toy_model}
	Let $G = (V,E)$ be a graph with $|V|=n$, which is either of the following: %which belongs to one of the following two families:
	\begin{itemize}
		\item[(i)] a $b$-expander graph for some $b>0$ with maximum degree $\Delta<\infty$; 
		\item[(ii)] the box $[-L, L]^d\cap \Z^d$ of volume $n$, for some $d\geq 5$.
	\end{itemize}
        Let $\cT$ be the RSTRE on $G$ for some environment distribution $\mu$ with $\beta = 1$. Then there exist constants $c_1, c_2, \gamma>0$ (with $c_1 = c_1(\mu, b, \Delta)$ for case (i), and $c_1 = c_1(\mu, d), c_2 = c_2(d)$ for case (ii)), such that for all $n \geq 2$ and $\delta > n^{-\gamma}$, we have
        \begin{equation}\label{eq:probbdd}
            \widehat \P \left( \frac{\sqrt{n}}{c_1(\delta^{-1}\log n)^{c_2}} \leq \diam(\cT) \leq c_1 (\delta^{-1} \log n)^{c_2} \sqrt{n}  \right) \geq 1-\delta .
    \end{equation}
\end{theorem}

\subsection{Improvement via bumping}
For simplicity's sake, in this section, we will only focus on the case when the graphs are expander graphs and refer to Section \ref{SS:boxZd} for a brief explanation of the $\Z^d$ case. We will prove the following slight extension of Theorem \ref{T:toy_model}. %, where we provide more details on the constant $c_1$ depending on the disorder distribution $\mu$. 

\begin{theorem}\label{T:toy_model_bump}
        Let $(G_n)_{n\in \N} = ((V_n, E_n))_{n \in \N}$, with $|V_n| = n$, be a sequence of $b$-expanders with maximal degree $\Delta$, and let $(\mu_n)_{n \in \N}$ be a sequence of environment distributions. Then there exists $p^* = p^*(\Delta, b), c_1, c_2, \gamma > 0$ such that for $a_1 = a_1(n)$ and $a_2 = a_2(n) \geq a_1$ satisfying 
        \begin{align*}
            \mu_n(-\infty, a_1) &\leq 1/(\Delta - 1) - p^*/2,\\
            \mu_n(a_2, \infty) &\leq p^*, \\
            a_2 - a_1 &\ll \log n,
        \end{align*}
        for all $n \in \N$, we have for any $\delta > n^{-\gamma}$ that
        \begin{equation}\label{eq:toy_bump_thm}
            \widehat \P \left( \frac{\sqrt{n}}{c_1 \chi^{c_2} } \leq \diam(\cT) \leq c_1 \chi^{c_2} \sqrt{n}  \right) \geq 1-\delta,
        \end{equation}
        where
        \begin{equation*}
            \chi = \chi(a_1,a_2,\delta, n) = \Delta \exp(a_2-a_1) \delta^{-1} \log n .
        \end{equation*}
\end{theorem}
\noindent The proof of Theorem \ref{T:toy_model_bump} will follow along very similar lines of the argument in \cite{MSS23}, and hence for parts of the proof that do not differ from \cite{MSS23} too much, we will only provide a proof sketch. See also Remark \ref{R:toy_vs_bump} for a reason why this theorem is an improvement upon Theorem \ref{T:toy_model}.

\subsection*{Bumping procedure}

Let the environment $(\omega_e)_{e \in E}$ be given with corresponding weights $\weight(e) = \exp(- \omega_e)$. As there is no dependence on $\beta$, by abuse of notation, we write $(G, \omega)$ for the weighted graph corresponding to $(G, \weight)$ and we will occasionally refer to $\omega$ as the weights. Suppose that $|V(G)| = n$ and $|E(G)| = m$, and denote by $(G^0, \omega^0)$ the weighted graph $(G, \omega)$. For some $0 \leq k \leq m$, we shall construct a sequence of weighted graphs $(G^{i}, \omega^{i})$, $0 \leq i \leq k$, each depending on $(G^{i-1}, \omega^{i-1})$, by either: 
\begin{enumerate}[(1)]
    \item changing the weights $(\omega_e)_{e \in E}$, 
    \item or by removing and contracting some edges of $E(G)$.
\end{enumerate}
We extend the laws of $\bP^{(G^{i}, \omega^{i})}$ on $\{0,1\}^{E(G^{i})}$ to laws on $\{0,1\}^{E(G)}$ by letting contracted (resp.\ removed) edges be (resp.\ not be) in the spanning tree of $G$. We write 
\begin{equation} \label{eq:view_on_E}
    \tilde{\bP}^{(i)}(\mathcal{T} \in \cdot) := \bP^{(G^{i}, \omega^{i})}(\mathcal{T} \in \cdot) \text{ viewed as a measure on } \{0,1\}^E,
\end{equation}
for this new law. For some $D > 0$ to be chosen later (of order $\sqrt{|V(G)|}$ up to polylog factors), the goal is to obtain the following bound
\begin{equation*}
    \tilde{\bP}^{(0)}( \diam(\cT) \leq D) \geq \tilde{\bP}^{(k)}( \diam(\cT) \leq D) \geq 1- \delta,
\end{equation*}
where we will be able to apply Theorem \ref{T:Nach} on $(G^k, \omega^k)$. Furthermore, for a different sequence of $((G^i, \omega^i))_{0 \leq i \leq k}$ the same will hold for the event that the diameter is larger than some $D' > 0$.

\smallskip

Suppose we wish to ``bump'' some edges $e_1, \ldots e_k$ with $a_1 \leq \omega_i \leq a_2$ for $i=1, \ldots, k$. The iterative construction of the sequence $(G^i, \omega^i)$ is as follows. Assume we are given $(G^{i}, \omega^{i})$, then for $D > 0$ we may write
\begin{align}
    \tilde{\bP}^{(i)}( \diam(\cT) \leq D) &= \tilde{\bP}^{(i)}( \diam(\cT) \leq D \mid e_{i+1} \in \cT) \tilde{\bP}^{(i)}( e_{i+1} \in \cT) \nonumber \\
    &\quad + \tilde{\bP}^{(i)}( \diam(\cT) \leq D \mid e_{i+1} \not\in \cT) \tilde{\bP}^{(i)}( e_{i+1} \not\in \cT)  \nonumber \\
    &=: \lambda_1 \tilde{\bP}^{(i)}( e_{i+1} \in \cT) + \lambda_2 \tilde{\bP}^{(i)}( e_{i+1} \not\in \cT). \label{eq:bump_split} 
\end{align}
Notice that the value of the conditional laws $\lambda_1, \lambda_2$ in \eqref{eq:bump_split} are independent of $\omega_{e_{i+1}}$, and that trivially
\begin{equation*}
    \tilde{\bP}^{(i)}( e_{i+1} \in \cT) + \tilde{\bP}^{(i)}( e_{i+1} \not\in \cT) = 1.
\end{equation*}
This means that $\tilde{\bP}^{(i)}( \diam(\cT) \leq D)$ is a convex combination of $\lambda_1$ and $\lambda_2$ with $\tilde{\bP}^{(i)}( e_{i+1} \in \cT)$ decreasing in $\omega_{e_{i+1}}$ and $\tilde{\bP}^{(i)}( e_{i+1} \not\in \cT)$ increasing in $\omega_{e_{i+1}}$ (see also Lemma \ref{L:derivative}). If $\lambda_1 \leq \lambda_2$, then we let $\omega^{i+1}$ be $\omega^i$ except we replace the value of $\omega_{e_{i+1}}$ with $a_1$, and if $\lambda_1 > \lambda_2$, we replace the value with $a_2$. It follows that
\begin{align*}
    \tilde{\bP}^{(i)}( \diam(\cT) \leq D) &\geq \lambda_1 \tilde{\bP}^{(i+1)}( e_{i+1} \in \cT) + \lambda_2 \tilde{\bP}^{(i+1)}( e_{i+1} \not\in \cT) \\
    &= \tilde{\bP}^{(i+1)}( \diam(\cT) \leq D).
\end{align*}
In the special case when either $a_1 = - \infty$ or $a_2 = \infty$, we contract (resp. remove) the edge in $G^i$ to obtain $G^{i+1}$ instead. Here we implicitly assume that after contraction/deletion, we may still define a spanning tree measure on $G^{i+1}$, i.e.\ we do not create a loop in $\cT$ or disconnect the graph. Otherwise, it must have been that either $\tilde{\bP}^{(i)}( e_{i+1} \in \cT)$ or $\tilde{\bP}^{(i)}( e_{i+1} \not\in \cT)$ was zero, so that we may set $\omega_{e_{i+1}}$ to any value without changing the law of the spanning tree. The following lemma follows by construction of the bumping procedure.

\begin{lemma} \label{L:bump}
    Let $F$ be a set of edges in $(G, \omega)$ with $\omega_e \leq a_1$ for all $e \in F$, and fix $D > 0$. Then there exists a connected weighted graph $(G', \omega')$, such that an edge $e \in F$ is either contracted in $G'$ (without contracting a loop) or has $\omega'_e = a_1$, and 
    \begin{equation} \label{eq:bump_F}
        \bP( \diam(\cT) \leq D) \geq \tilde{\bP}^{(G', \omega')}( \diam(\cT) \leq D).
    \end{equation}
    If instead $a_2 \leq \omega_e$ for each $e \in F$, then an edge $e$ is either not contained in $G'$ or has $\omega'_e = a_2$ and \eqref{eq:bump_F} holds. The same conclusion (for a different $(G', \omega')$) also holds for the event that the diameter is larger than some $D'$.
\end{lemma}

\begin{figure}[ht]
	\centering
	\includegraphics[width=0.98\textwidth]{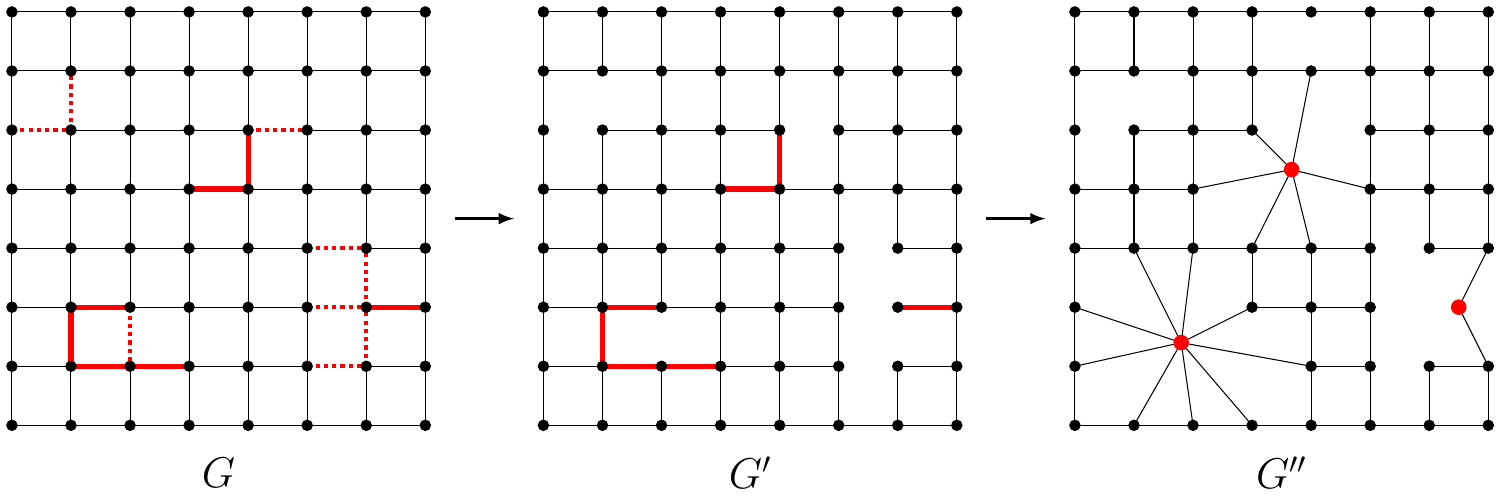}
	\caption{The thick (resp.\ dashed) red edges in the first graph $G$ correspond to edges with $\omega_e \leq a_1$ (resp.\ $\omega_e \geq a_2$). We obtain the graph $G''$ by first bumping the small weighted edges and then bumping the large weighted edges. Edges may be removed/contracted in the bumping procedure.}
	%Consider the graph $G$ on the left and perform percolation on $G$ by closing (indicated in red) edges with weight outside $[1/A, A]$ and call this set $K$. For any realization of $\cT$ on $K$ we can contract edges (drawn thick) in $\cT(K)$ and remove edges (drawn dotted) in $K \setminus \cT(K)$ such that we obtain a new graph on the right denoted by $G'$.}
\label{fig:ConDel}
\end{figure}

To prove Theorem \ref{T:toy_model_bump}, we apply Lemma \ref{L:bump} twice. First, we remove (or update the weight of) any edge with a small weight to obtain a graph $G'$, then we contract (or update the weight of) any edge which has a large weight in $G'$ to obtain $G''$. See also Figure \ref{fig:ConDel}. The final graph $G''$ will depend on $\omega$ and thus be random, however, the weight of any edge is bounded in some interval, and so we will be able to show that the conditions of Theorem \ref{T:Nach} hold with high $\P$-probability.

\begin{remark} \label{R:toy_vs_bump}
    The difference in Theorem \ref{T:toy_model_bump} compared to Theorem \ref{T:toy_model} is that in \cite{MSS23} we only contracted/removed edges and never updated their weight. Consequently, we can only condition on, say, a $p^*/2$ fraction of edges, while for Theorem \ref{T:toy_model_bump} we can deal with a $1/(\Delta - 1) + p^*/2$ fraction of edges. Furthermore, we allow the environment distribution $\mu$ to depend on the size of the graph\footnote{The exact same proof arguments of \cite{MSS23} can be very easily adapted to also include this last extension.}.
\end{remark}

\subsection{Bumping small edges}

Let $F^{-1}_\mu$ be the inverse CDF of $\mu$ and let
\begin{equation} \label{eq:toy_small_weight}
    E_{\geq 1 - p^*} := \big\{ e \in E : \omega_e \geq F^{-1}_\mu(1 - p^*) \big\}
\end{equation}
be the set of edges with ``small'' conductance for some $p^*$ to be determined later. We first bump all the edges in $E_{\geq 1 - p^*}$. Applying Lemma \ref{L:bump} to $F = E_{\geq 1 - p^*}$ with $a_2 = F^{-1}_\mu(1 - p^*)$, gives a connected weighted graph $(G', \omega')$ such that any edge $e$ in $G'$ has $\omega'_e \leq F^{-1}_\mu(1 - p^*)$ and
\begin{equation}
    \bP( \diam(\cT) \leq D) \geq \tilde{\bP}^{(G', \omega')}( \diam(\cT) \leq D). \label{eq:bump_small_G'}
\end{equation}
The following is essentially Lemma 3.6 of \cite{MSS23}, and says that, with high probability, any graph $G'$ produced by Lemma \ref{L:bump} (i.e.\ for any configuration of removing the edge or increasing its weight for edges in $E_{\geq 1 - p^*}$) the edge expansion $h_{G'}$ from Definition \ref{D:expander} is not too small.

\begin{lemma}[cf.\ {\cite[Lemma 3.6]{MSS23}}] \label{L:iso_G'}
    The (unweighted) graph $G'$ from \eqref{eq:bump_small_G'} satisfies
    \begin{equation*}
        h_{G'} \geq \frac{1}{8 \Delta^2 (\log n)^2}
    \end{equation*}
    with probability at least $1 - O(n^{-2})$, provided that $p^* = p^*(\Delta, b) > 0$ is small enough.
\end{lemma}

\begin{proof}[Proof sketch]
    The (unweighted) graph $G'$ is obtained by deleting a subset of edges with $\omega_e \geq F^{-1}_\mu(1- p^*)$ from $G$ under the condition that the graph remains connected. Denote by $\cC_1(G, 1- p^*)$ the largest connected component of the random graph coupled to $\omega$ as in Section \ref{SS:CoupleRG}. The key facts for the proof are that the following all hold with high probability:
    \begin{enumerate}
        \item The giant component $\cC_1(G, 1- p^*)$ has edge expansion $h_{\cC_1(G, 1- p^*)} \geq 1/\log n$ (see Lemma 3.4 of \cite{MSS23});
        \item The giant component has size $|\cC_1(G, 1- p^*)| \geq 3n/4$ (see Lemma 3.1 of \cite{MSS23});
        \item Any connected component in $V \setminus \cC_1(G, 1- p^*)$ has size at most $\log n$ (see Lemma 3.3 of \cite{MSS23}).
    \end{enumerate}
    Any set $S$ can then be decomposed into $S_1$ the part inside $\cC_1(G, 1- p^*)$ and $S_2$ the part outside of $\cC_1(G, 1- p^*)$. By the above facts, each connected component of $S_1$ or $S_2$ has an edge expansion not too small. It remains to show that the same holds for the union of the components. This is done by using that the giant component has a large size and the fact that the graph needs to stay connected to be able to support a spanning tree measure. We refer to \cite{MSS23} for the full details, and remark that there the expansion ratio
    \begin{equation*}
        \tilde{h}_{G'}(S) = \frac{|E(S,S^c)|}{|S|} \leq \frac{|E(S,S^c)|}{ \tfrac{1}{\Delta }|E(S,V)|} = \Delta \, h_{G'}(S)
    \end{equation*}
    was used, and hence we require an extra factor of $\Delta$ in Lemma \ref{L:iso_G'}.
\end{proof}

\subsection{Bumping large edges}

Similarly to \eqref{eq:toy_small_weight}, define
\begin{equation} \label{eq:toy_large_weight}
    E_{\leq 1/(\Delta - 1) - p^*/2} := \big\{ e \in E : \omega_e \leq F^{-1}_\mu(1/(\Delta - 1) - p^*/2) \big\} 
\end{equation}
as the set of edges with ``large'' conductance. Consider again $(G', \omega')$ produced by bumping edges in $E_{\geq 1 - p^*}$. Let $(G'', \omega'')$ be the weighted graph obtained from $(G', \omega')$ by applying Lemma \ref{L:bump} with $F = E_{\leq 1/(\Delta - 1) - p^*/2}$ and $a_1 = F^{-1}_\mu(1/(\Delta - 1) - p^*/2)$. In other words, for each edge in $G'$ with a large conductance, either update its edge weight to something smaller or contract it without creating loops. The weighted graph $(G'', \omega'')$ then satisfies
\begin{equation}
    \bP( \diam(\cT) \leq D) \geq \tilde{\bP}^{(G', \omega')}( \diam(\cT) \leq D) \geq \tilde{\bP}^{(G'', \omega'')}( \diam(\cT) \leq D). \label{eq:bump_large_G''}
\end{equation}
The following is the key observation about $(G'', \omega'')$ and is essentially given by Proposition 3.5 of \cite{MSS23}.

\begin{lemma}[cf.\ {\cite[Proposition 3.5]{MSS23}}] \label{L:bottle_G''}
    Let
    \begin{align*}
        a_1 &:= a_1(\Delta, b) = F^{-1}_\mu \big(\frac{1}{\Delta -1} - \frac{p^*}{2} \big), \\
        a_2 &:= a_2(\Delta, b) = F^{-1}_\mu(1 - p^*).
    \end{align*}
    Then there exists a constant $C_\Phi = C_\Phi(p^*) > 0$, such that the weighted graph $(G'', \omega'')$ from \eqref{eq:bump_large_G''} satisfies
    \begin{equation} \label{eq:bottle_G''}
        \Phi_{(G'', \omega'')} \geq C_\Phi \frac{\exp (a_1 -a_2) }{\Delta^3 (\log n)^3}
    \end{equation}
    with probability $1 - O(n^{-2})$, provided that $p^* = p^*(\Delta, b) > 0$ is small enough.
\end{lemma}
\begin{proof}[Proof sketch]
    We again only sketch the proof of Lemma \ref{L:bottle_G''} and refer to the proof of Proposition 3.5 in \cite{MSS23} for more details. The graph $G''$ is obtained by contracting a subset of edges in $E_{\leq 1/(\Delta - 1) - p^*/2}$ such that we never contract a cycle. By Lemma \ref{eq:bump_small_G'}, with high probability we may assume that
    \begin{equation*}
        h_{G'} \geq \frac{1}{8 \Delta^2 (\log n)^2}.
    \end{equation*}
    Furthermore, by Lemma \ref{L:perc_subcrit} as trivially $1/(\Delta -1) - p^*/2 < 1/(\Delta -1)$, we have for some $C_s = C_s(p^*) > 0$ that
    \begin{equation}
        \Big| \cC_1 \big(G, \frac{1}{\Delta - 1} - p^*/2 \big) \Big| \leq C_s \log n, \label{eq:contract_small}
    \end{equation}
    with high probability.
    
    Let $S \subset V(G'')$, and denote by $\tilde{S} \subset V(G')$ the set that was contracted in $G'$ to obtain $S$. By the expansion properties of $G'$, we have
    \begin{align*}
        |E_{G''}(S,S^c)| = |E_{G'}(\tilde{S}, \tilde{S}^c)| & \geq \frac{1}{8 \Delta^2 (\log n)^2} \min \big\{ \big|E(\tilde{S}, V(G')) \big|, \big|E(\tilde{S}^c, V(G') \big| \big\}\\
        &\geq \frac{1}{16 \Delta^2 (\log n)^2} \min \big\{ |\tilde{S}|, |\tilde{S}^c| \big\}.
    \end{align*}
    For the bottleneck ratio, we only need to consider $S$ with $0 < \pi(S) \leq 1/2$, however, this does not imply a uniform lower bound on $|\tilde{S}|$ and $|\tilde{S}^c|$. Nonetheless, by using \eqref{eq:contract_small}, which implies that $|\tilde{S}| \leq C_s \log n |S|$, we may obtain that
    \begin{equation*}
        \min \big\{ |\tilde{S}|, |\tilde{S}^c| \big\} \geq C \frac{1}{ \log n}
    \end{equation*}
    for some constant $C > 0$. As
    \begin{equation*}
         a_1 \leq \omega_e'' \leq a_2,
    \end{equation*}
    for all edges $e \in E(G'')$ and the maximum degree in $G''$ is bounded by $\Delta C_s \log n$, one may now deduce \eqref{eq:bottle_G''}.
\end{proof}

\subsection{Proof of Theorem \ref{T:toy_model_bump}}

To prove Theorem \ref{T:toy_model_bump} we will now apply Theorem \ref{T:Nach} to the graph $G''$ as constructed in Lemma \ref{L:bottle_G''}. We also refer to the proof of Theorem 1.1 in \cite{MSS23} for more details.

\begin{proof}[Proof of Theorem \ref{T:toy_model_bump}.]
    Fix $\delta > 0$ and let $D$ be some constant which we determine later in \eqref{eq:determine_D}. We first prove the upper bound on the diameter of $\cT$, the lower bound follows along the similar arguments.
    
    Assume that $p^* = p^*(\Delta, b)$ is small enough such that the inequalities in Lemma \ref{L:iso_G'} and Lemma \ref{L:bottle_G''} hold with probability at least $1 - O(n^{-2})$. As in Lemma \ref{L:bottle_G''}, we obtain a new graph $(G'', \omega'')$ (depending on $D$) by first bumping the small weighted edges in $E_{\geq p^*}$, and then by bumping the large weighted edges in $E_{\leq 1/(\Delta - 1) - p^*}$ such that $a_1 \leq \omega''_e \leq a_2$ for all edges $e$ in $G''$. By construction (see the argument before \eqref{eq:bump_large_G''}), the new graph $G''$ satisfies
    \begin{equation*}
        \bP( \diam(\cT) \leq D) \geq \tilde{\bP}^{(G'', \omega'')}( \diam(\cT) \leq D).
    \end{equation*}
    
    We verify the conditions in Theorem \eqref{T:Nach} for $(G'', \omega'')$:
    \begin{itemize}
        \item Condition \eqref{eq:Nbalanced}: Denote by $\pi''$ the stationary distribution of the random walk on $(G'', \omega'')$. For the remaining proof we may, with high probability by the inequality in \eqref{eq:contract_small}, restrict to the event that the largest contracted component has size at most $C_s \log n$ for some constant $C_s > 0$ (depending on $p^*$ and thus on $\Delta$ and $b$). This means that the largest degree in $G''$ is at most $\Delta  C_s \log n$, and as any vertex has degree at least one, we have
         \begin{equation*}
       \frac{\max_{u\in V} \pi''(u)}{\min_{u\in V} \pi''(u)} \leq \Delta C_s \log n \exp(a_2-a_1).
    \end{equation*}
    
        \item Condition \eqref{eq:Nmixing}:  Again, as each vertex has at least one incident edge with $\omega_e \leq a_2$, we have
        \begin{equation*}
            \min_{v \in V(G'')} \pi''(v) \geq \frac{\exp(-a_2)}{2 |E(G'')| \exp(-a_1)} \geq \frac{\exp(a_1- a_2)}{\Delta n}.
        \end{equation*}
        Furthermore, let $M = C_\Phi \frac{\exp (a_1 -a_2) }{\Delta^3 (\log n)^3}$, then using Theorem \ref{T:HeatCheeger}
    \begin{align*}
        \tmix(G'') &\leq 1 + \frac{4}{M^2} \int_{4/\min \pi''(v)}^8 \frac{1}{r} dr \\
        &\leq 1 + \frac{4}{M^2} \big( \log 8 - \log 4 + (a_2 - a_1) + \log(\Delta n) \big) \\
        &\leq C \Delta^6 \exp \big( 3(a_2-a_1) \big)  (\log n)^7  =: f(a_1, a_2, \Delta, b) (\log n)^7,
    \end{align*}
    for some constant $C > 0$ depending on $C_\phi$ and $p^*$, and thus on $\Delta$ and $b$. Therefore, in Theorem \ref{T:Nach} we may choose any $\alpha \in (0,1/2)$.

    \item Condition \eqref{eq:Nescaping}: We may use the trivial bound
    \begin{align*}\label{eq:tmixG'}
    	\sum_{t=0}^{\tmix(G'')} (t+1) \sup_{v \in V(G'')}q_t (v,v) &\leq (\tmix(G'')+1)^2 \\
     &\leq 2 f(a_1, a_2, \Delta, b)^2 (\log n)^{14}.
    \end{align*}
    \end{itemize}

    The conditions of Theorem \ref{T:Nach} are now verified for $(G'', \omega'')$ with parameters that are of order $n^{o(1)}$ (assuming that $a_2(n) - a_1(n) \ll \log n$). Let $m = |V(G'')|$, which we may assume satisfies $n/C_s\log n \leq m \leq n$ by \eqref{eq:contract_small}. Applying Theorem \ref{T:Nach} to the weighted UST $\cT_{G''}$ (viewed as an element of $\{0,1 \}^{E(G'')}$ with law $\bP_{G''}^{\omega''}$) gives that there exist constants $c_1, c_2 > 0$, where $c_2$ is independent of all parameters, such that for $\delta >0$,
    \begin{equation*}
         \chi = \frac{\delta}{2}f(a_1, a_2, \Delta, b) \log n,
    \end{equation*}
    we have
    \begin{equation} \label{eq:G''nachmias}
        \bP_{G''}^{\omega''} \Big( \frac{\sqrt{m}}{c_1 \chi^{c_2}}  \leq \diam(\cT_{G''}) \leq c_1 \chi^{c_2} \sqrt{m} \Big) \geq 1 - \frac{\delta}{2}.
    \end{equation}
    It remains to translate the bounds of $\diam(\cT_{G''})$ into bounds on $\diam(\cT)$.

    To construct $\cT$ from $\cT_{G''}$, we must undo all the contractions we did to obtain $G''$ from $G$. As the largest connected components at $p= 1/(\Delta - 1) - p^*/2$ are small, each contracted vertex in $G''$ can correspond to a maximum of $C_s \log n$ many vertices in $G$. Therefore, the length of any acyclic path in $\cT$ is at most $C_s \log n$ times longer than the corresponding path in $\cT_{G''}$. Consequently, with high probability (w.r.t.\ $\omega$) we have
        \begin{equation*}
        %\bP( \diam(\cT) \leq D) \geq 
        \tilde{\bP}^{(G'', \omega'')}( \diam(\cT) \leq D) \geq \bP_{G''}^{\omega''} \big( \diam(\cT_{G''}) \leq \frac{D}{C_s \log n} \big),
    \end{equation*}
    where we recall that $\tilde{\bP}$ is a measure on $\{0,1\}^{E(G)}$ instead of $\{0,1 \}^{E(G'')}$. The upper bound on the diameter of $\cT$ now follows with \eqref{eq:G''nachmias} by letting
    \begin{equation} \label{eq:determine_D}
        D = c_1 \chi^{c_2} \sqrt{n} C_s \log n,
    \end{equation}
    and absorbing the constant $C_s$ into $c_1$. The lower bound on the diameter follows along a similar argument, where we need to use the fact that with high probability $m = |V(G'')| \geq n/C_s \log n$.
\end{proof}

\begin{example} \label{ex:3-reg}
    Let $G_n$ be a random $3$-regular graph on $n$ vertices for some large $n \in \N$. It is known (cf.\ \cite{Bol88}) that, with high probability, $G_n$ is a $b$-expander for some $b > 0$ independent of $n$. Furthermore, the critical window of percolation is of the form $p=(1 + \lambda n^{-1/3})/2$ for $\lambda \in \R$, see \cite{NP10}. Theorem \ref{T:toy_model_bump} then says that, even if we ignore more than half of the mass of any edge weight distribution $\mu_n$, we may still see a tree with a diameter of order $\sqrt{n}$ (up to polylogarithmic factors). In this special case, using the well-understood critical window of the random $3$-regular graph, we believe that it might be possible to show that only a $o(1)$ portion of edges influence the UST measure (up to polylogarithmic factors). That is, if we ignore logarithmic factors, the diameter is determined by the values $F^{-1}_{\mu_n}$ takes on the interval $[1/2 - o(1), 1/2 + o(1)]$. 
\end{example}

% \begin{remark}
%     One can extend Theorem \ref{T:toy_model_bump} to weight distributions  $\mu_n$ that may depend on the size of the graph. Indeed, the only parameters in Theorem \ref{T:toy_model_bump} that depend on $\mu$ are $a_1$ and $a_2$. If $a_1 = a_1(n)$ and $a_2 = a_2(n)$ are such that $\exp(a_2 - a_1)$ is less than any power of $n$, then Theorem \ref{T:Nach} is still applicable. Therefore, we obtain upper and lower bounds that are, up to multiplicative polylogarithmic corrections terms, of order $\sqrt{n}$.
% \end{remark}

\subsection{High-dimensional bounded-degree graphs} \label{SS:boxZd}

We state here the analog of Theorem \ref{T:toy_model_bump} for boxes in $\Z^d$, $d \geq 5$, and give only a very brief sketch of the proof. We refer to Section 5 of \cite{MSS23} for more details for the case where we only contract and delete edges instead of updating their weight.

\begin{theorem}
    Let $(G_n)_{n \in \N}$,  $|G(V_n)| \asymp n^d$, be the sequence of boxes $[-n, n]^d$ in $\mathbb{L}^d$, $d \geq 5$. There exists constants $p^*, c_1, c_2, \gamma > 0$, all depending only on $d$, such that for a sequence of environment distributions $(\mu_n)_{n\geq 1}$ and $a_1 = a_1(n), a_2 = a_2(n)$ satisfying
    \begin{equation} \label{eq:boxes_a12}
    \begin{aligned}
        a_1 &:= F^{-1}_{\mu_n} \big( \frac{1}{2d - 1} - \frac{p^*}{2} \big), \\
        a_2 &:= F^{-1}_{\mu_n}(1 - p^* ), \\
        a_2 - a_1 &\ll \log n,
    \end{aligned}
    \end{equation}
    for all $n \in \N$, we have for any $\delta > n^{-\gamma}$ that
    \begin{equation}
            \widehat \P \left( \frac{n^{d/2}}{c_1 \chi^{c_2} } \leq \diam(\cT) \leq c_1 \chi^{c_2} n^{d/2}  \right) \geq 1-\delta,
    \end{equation}
        where
    \begin{equation*}
            \chi = \chi(a_1,a_2, \delta, n) = \exp(a_2-a_1) \delta^{-1} \log n.
    \end{equation*}
\end{theorem}

The main difference between the proof idea compared to Theorem \ref{T:toy_model_bump} is that instead of considering the bottleneck ratio as in Lemma \ref{L:bottle_G''}, one has to examine the whole bottleneck profile as in Definition \ref{D:bottleneck_profile}. The analogous statement of Lemma \ref{L:bottle_G''} is the following, where we again use the bumping procedure from Lemma \ref{L:bump} to obtain $(G_n'', \omega'')$ from $(G_n, \omega)$.  

\begin{lemma}[cf.\ {\cite[Lemma 5.2]{MSS23}}] \label{L:bottle_profile_G''}
    There exists constants $p^*, C_\Phi, c > 0$, all depending only on $d$, such that for $a_1 = a_1(n), a_2 = a_2(n)$ as in \eqref{eq:boxes_a12}, the weighted graph $(G_n'', \omega'')$ from \eqref{eq:bump_large_G''} satisfies
    \begin{equation} \label{eq:Zdbottle_G''}
        \Phi_{(G_n'', \omega'')}(r) \geq C_\Phi \frac{\exp ( a_1 -a_2) }{(\log n)^c} \Big( \frac{\pi''_{\mathrm{min}}}{r} \Big)^{\frac{1}{d}},
    \end{equation}
    with high probability.
\end{lemma}

\noindent We will not prove Lemma \ref{L:bottle_profile_G''}, instead, we refer to the previous section and the proof of Lemma 5.2 in \cite{MSS23} for the details. Using Lemma \ref{L:bottle_profile_G''} and Theorem \ref{T:HeatCheeger}, assuming $a_2 - a_1 \ll \log n$, one may again use Theorem \ref{T:Nach} to recover Theorem \ref{T:toy_model} for the case of boxes in $\Z^d$, $d \geq 5$. 

In fact, we conjecture that a similar bound as in Lemma \ref{L:bottle_profile_G''} should hold for a large class of ``high-dimensional'' graphs and not only boxes in $\Z^d$. The main hurdle in extending Theorem \ref{T:toy_model} to other graphs is showing the following implication:
\begin{equation} \label{eq:ConjIso}
`` \;\Phi_{G}(r) \geq c (\pi_{\mathrm{min}}^{-1} r)^{-\frac{1}{\beta}} \;\Longrightarrow\; \Phi_{G''} (r) \geq c (\log n)^{-k} (\pi_{\mathrm{min}}^{-1} r)^{-\frac{1}{\beta}} \text{ with high probability}",
\end{equation}
for some $k \geq 0$, see also \cite[Section 6.2]{MSS23}.

\subsection{Counter example for unbounded degrees}

The key ingredient of the proof of Theorem \ref{T:toy_model_bump} is that for a random graph with a parameter that is uniformly bounded away from $0$, the bottleneck ratio (or profile) is not drastically different from that of the original graph. We can not expect this to hold as the degrees are growing with the size of the graph. Indeed, the following gives a counterexample to Theorem \ref{T:toy_model_bump} when we no longer assume uniform boundedness of the degrees.

\begin{proposition}[cf.\ {\cite[Proposition 6.1]{MSS23}}] \label{P:counterexample}
Let $K_n$ be the complete graph with $n$ vertices and let $\mu$ be the law of the random variable $- \exp(U^{-1})$, where $U$ is uniformly distributed on $[0,1]$. Fix $\beta =1$, and let $M^{(n)}$ be the MST coupled to the environment as in Section \ref{SS:CoupleRG}, then
\begin{equation} \label{eq:convergeT1}
	\widehat \P( \cT = M^{(n)}) \xrightarrow{n \rightarrow \infty} 1.
\end{equation}
As a consequence, with high probability, the diameter of the uniform spanning tree $\cT$ is of order $n^{1/3}$ as $n \rightarrow \infty$.
\end{proposition}
\begin{proof}
Denote by $T_1$ and $T_2$ the spanning trees with the smallest and second smallest weight $\weight(T)=\sum_{e \in T} \omega_e$, respectively. We have the following facts:
\begin{enumerate}
	\item The number of spanning trees on $K_n$ is $n^{n-2}$ (Cayley's formula). \label{enu:Cayley}
 
	\item $T_1$ and $T_2$ differ by a single edge (otherwise we can swap an edge in $T_2$ with an edge in $T_1$ to obtain a spanning tree $T_3$ with $w_n(T_1) < w_n(T_3) < w_n(T_2)$). \label{enu:T1T2}  
 
	\item If $(U_i)_{1\leq i \leq m}$ is a collection of i.i.d.\ uniform random variables on $[0,1]$ and $X = \min_{i \neq j} |U_i - U_j|$, then
	\begin{equation} \label{eq:uniform_gap_law}
		\P( X > t) = (1 - (m-1) t)^m \qquad \text{for } 0 \leq t \leq \frac{1}{m-1}.
	\end{equation}
        This implies that $\P( X \leq m^{-2} (\log m)^{-1}) \rightarrow 0$. \label{enu:minGapU}
 
	\item By Theorem \ref{T:Gnp_connect}, every edge $e$ in $T_1$ has $U_e \leq 2 \log n/n$ with high probability. \label{enu:ConnectionP} %for some reason adding a label creates some vertical space below....?
\end{enumerate}
It thus suffices to show that
\begin{equation*}
	\frac{n^{n-2} \weight(T_2)}{\weight(T_1)} \xrightarrow{n \rightarrow \infty} 0\,.
\end{equation*}
To this end, let $e \in T_1 \setminus T_2$ and $f \in T_2 \setminus T_1$ be the two edges in which $T_1$ and $T_2$ differ, and define the gap $g := U_f - U_e$. Then by item \ref{enu:minGapU} with $m=\binom{n}{2}$, with high probability
\begin{equation*}
	\frac{g}{U_e U_f} \geq n^{-4} \log^{-1} n.
\end{equation*}
Hence,
\begin{align*}
	\frac{n^{n-2} \weight(T_2)}{\weight(T_1)} &= n^{n-2} \exp \Big( \exp(U_f^{-1}) - \exp(U_e^{-1})\Big) \\
	&= n^{n-2} \exp \Bigg( -\exp(U_e^{-1}) \Big(1 - \exp \big(- \frac{g}{U_e U_f} \big) \Big)\Bigg) \\
	&\leq n^{n-2} \exp \Bigg( -\exp(U_e^{-1}) \Big(1 - \exp \big(- n^{-4} \log^{-1} n \big) \Big)\Bigg) \\
	&\leq n^{n-2} \exp \Big( -\exp\big(\frac{n}{2 \log n} \big ) \frac{1}{2n^4 \log n} \Big) \longrightarrow 0  ,
\end{align*}
where we used the bound $(1-\e^{-x}) \geq x/2$ for $0 < x \leq 1$.

It is known from \cite{ABGM17} that $n^{-1/3} M^{(n)}$ converges in distribution to a random compact metric space. Therefore, with high probability as $n\to\infty$, the diameter of $T_1$, and hence $\cT$, is of order $n^{1/3}$. %It's a bit overkill... 
\end{proof}

%%%%%%%%%%%%%%% NEW SECTION %%%%%%%%%%%%%%%%%%%%%%%

\section{Diameter of the RSTRE on the complete graph} \label{S:high_and_low}

As demonstrated in Proposition \ref{P:counterexample}, once the degrees are growing with the number of vertices, we can no longer expect the conclusion of Theorem \ref{T:toy_model} to hold for such graphs. We now turn our focus to the complete graph, and for the rest of this section, we will make the following assumption on the environment distribution $\mu$, which essentially says that the tail of $\mu$ close to zero behaves as that of a Beta distribution $B(\alpha, \cdot)$.

\begin{assumption} \label{as:mu_tail}
    Let $\mu$ be supported on $[0, \infty)$, and assume that there exists constants $\alpha, \rho, c_\mu > 0$ such that
    \begin{equation*}
        F_{\mu}(t) = c_\mu t^{\alpha} \qquad \text{whenever } 0 \leq t \leq \rho,
    \end{equation*}
    where $F_\mu$ the CDF of $\mu$.
\end{assumption}

Under assumption \ref{as:mu_tail}, for $\beta \geq 0$, it is clear that $\weight(e) = \exp(-\beta \omega_e) \in [0,1]$, ensuring that both the mean and variance of $\weight(e)$ are finite. If the underlying graph $G$ has $n$ vertices and the environment distribution is given by $\mu$, we denote by $\cT^\omega_{n, \beta}$ the corresponding RSTRE with law $\bP^\omega_{n, \beta}$. To ease notation in the later sections, we will often drop the sub- and superscripts when there is no ambiguity. This section aims to present the main ideas behind the proof of the following theorem, which is a generalization of \cite[Theorem 1.1]{MSS24}

\begin{theorem} \label{T:main}
    Suppose $\mu$ satisfies Assumption \ref{as:mu_tail} for some $\alpha > 0$.  There exists a constant $C = C(\mu) >0$ such that for any $\delta > 0$, if $\beta \leq C (n/\log n)^{1/\alpha}$ and $n$ is sufficiently large, then 
    \begin{equation}\label{eq:low}
        \widehat \P \Big( C_1(\delta)^{-1} n^{1/2}\leq \diam(\cT^\omega_{n,                \beta})\leq C_1(\delta) n^{1/2} \Big) \geq 1-\delta,
    \end{equation}
    where $C_1(\delta)>0$ is a constant depending only on $\delta$. On the other hand, if $\beta \geq n^{1/\alpha + 1/3} \log n$ and $n$ is sufficiently large, then
    \begin{equation}\label{eq:high}
        \widehat \P \Big( C_2(\delta)^{-1} n^{1/3}\leq \diam(\cT^\omega_{n, \beta})\leq C_2(\delta) n^{1/3}   \Big) \geq 1-\delta,
    \end{equation}
    where $C_2(\delta)>0$ is a constant depending only on $\delta$.  
\end{theorem}

Theorem \ref{T:main} identifies two phases: (1) a \textit{low disorder} regime $\beta \leq C (n/\log n)^{1/\alpha}$, where the diameter of the RSTRE is of the same order as the diameter of the UST; (2) a \textit{high disorder} regime $\beta \geq n^{1/\alpha + 1/3} \log n$, where the diameter of the RSTRE is of the same order as the diameter of the MST. When $\mu$ is the uniform distribution on $[0,1]$, we have the following conjecture for the \textit{intermediate} regime, and refer to Chapter \ref{ch:conj} for a heuristic argument.

\begin{conjecture}[cf.\ {\cite[Conjecture 1.3]{MSS24}}] \label{C:Intermediate}
	Let $G=K_n$ be the complete graph with $n$ vertices and edge weights $\weight(e) = \exp(-\beta \omega_e)$, where $\omega_e$ are i.i.d.\ uniform on $[0,1]$. Then, with high probability with respect to the averaged law $\widehat{\P}$, we have that
	\begin{equation*}
	\diam \big( \cT^\omega_{n, \beta} \big) =
	\begin{dcases}
	n^{\frac{1}{2}+o(1)} & \beta \leq C n, \\
	n^{\frac{1}{2} - \frac{\gamma}{2}+o(1)} & \beta  = n^{1 + \gamma+o(1)}, \quad 0 < \gamma < \frac{1}{3}, \\
	n^{\frac{1}{3}+o(1)} & \beta \geq C n^{\frac{4}{3}}.
	\end{dcases}
	\end{equation*}
\end{conjecture}

Notice that compared to Theorem \ref{T:local_limit}, the high disorder regime (if $\alpha = 1$) corresponds to $\beta \gg n^{4/3}$ instead of $\beta \gg n$. That is, although there is a sharp cut-off in the local behavior around $\beta \approx n$, there is a conjectured smooth transition from the global behavior of the UST to the global behavior of the MST for $\beta = n^{1 + \gamma}$ with $0 < \gamma < 1/3$.

\begin{remark} \label{R:correspondence}
    As we multiply $\omega_e$ by $\beta = \beta(n)$, we are actually considering a sequence of RSTREs with environment distribution $\mu = \mu_n$ depending on $n$, the size of the graph. The same behavior of the diameter described in Theorem \ref{T:main} and Conjecture \ref{C:Intermediate} can also be obtained by choosing a fixed environment distribution. Namely, in the case of $\alpha = 1$, there is a correspondence between
    \begin{equation*}
         \weight(e) = \exp(-n^{1 + \gamma }U) \quad \text{ and } \quad \weight'(e) = \exp(-U^{-\gamma}),
    \end{equation*}
    where $U$ is uniform on $[0,1]$ and $0 < \gamma < 1/3$. See also Remark \ref{R:constant_dist_gap}, \cite[Remark 6.2]{MSS23}, and in particular \cite{K24}, in which the details of this correspondence are fully worked out. 
\end{remark}

\section{Low disorder} \label{SS:diamLow}

The proof of the low disorder result in Theorem \ref{T:main} will follow by using concentration arguments to show that the conditions of Theorem \ref{T:Nach} are satisfied for suitable parameters. The following theorem is a slight generalization of \cite[Theorem 3.1]{MSS24}

\sloppy
\begin{theorem} \label{T:GeneralLow}
	Let $G_n = (V_n, E_n)$, with $|V_n|=n$, be a sequence of $b$-expander graphs with minimum and maximum degree $\deg_{\mathrm{min}}(n) = \deg_{\mathrm{min}}$ and $\deg_{\max}(n) = \deg_{\max(n)}$, such that for some fixed $\Delta > 0$ we have $\deg_{\max}/\deg_{\mathrm{min}} \leq \Delta$ for all $n$. Suppose that $\mu$ satisfies Assumption \ref{as:mu_tail}. Then there exists constants $C_1 = C_1(b,\Delta, \mu), C_2(\mu)$ such that if $\deg_{\mathrm{min}} \geq C_2 \log n$, then we have for $\beta \leq C_1 (\deg_{\mathrm{min}}/ \log n)^{1/\alpha}$ and any $\delta > 0$ that
	\begin{equation*}
		\widehat \P \left( C_3^{-1} \sqrt{n} \leq \diam(\cT) \leq C_3  \sqrt{n}  \right) \geq 1-\delta,
	\end{equation*}
	for some $C_3 = C_2(\delta, b, \Delta) > 0$ and all $n = n(\delta)$ large enough.
\end{theorem} 

\fussy

% \begin{remark}
%     The constants $C_1$ and $C_2$ in Theorem \ref{T:GeneralLow} are independent of $n$, and dictate when the outgoing edge weights start to concentrate around their means. Furthermore, $C_3$ is independent of both $n$ and $\mu$. In particular, one may apply Theorem \ref{T:GeneralLow} (with the same constants $C_1, C_2$ and $C_3$) to a sequence of $b$-expanders whose ratio of maximal to minimal degree is uniformly bounded from above.
% \end{remark}

\noindent Several remarks regarding Theorem \ref{T:GeneralLow} are in order.

\begin{remark}
    In \cite{MSS24}, we have not assumed that $\deg_{\mathrm{min}} \geq C_2 \log n$, however, as we only considered the uniform distribution on $[0,1]$, the difference $\omega_e - \omega_{f}$ was uniformly bounded from above and below. Thus, for any fixed $\beta$, the ratio of weights was always bounded away from zero and infinity. This allows us to conclude the GHP (Gromov–Hausdorff–Prokhorov) convergence of the RSTRE in Proposition \ref{P:GHPbeta}. For $\mu$ which is only bounded from below by $0$, this is not necessarily the case, and consequently, the concentration in the proof of Theorem \ref{T:GeneralLow} may fail if $\deg_{\mathrm{min}}$ is small enough.
\end{remark}

\begin{remark} \label{R:mst_weight}
    Let $M_n$ be the MST on the complete graph on $n$ vertices with weights $\weight(e)$ that are i.i.d.\ distributed according to some $\lambda$ with CDF $F_\lambda$. Frieze in \cite{Fri85} showed that if $F_\lambda$ is differentiable at $0$, then
    \begin{equation*}
        L(M_n) = \sum_{e \in M_n} \weight(e) \xrightarrow{n \rightarrow\infty} \frac{\zeta(3)}{F_\lambda'(0)},
    \end{equation*}
    where $\zeta(3) = \sum_{k=1}^\infty k^{-3} = 1.202\dots$ and the convergence holds both in probability and in expectation, see also Theorem \ref{T:length}. Similarly, in view of Theorem \ref{T:GeneralLow}, the low disorder regime of the RSTRE is determined by the behavior of $\mu$ around $0$.
\end{remark}

\begin{remark}
    One might think that Theorem \ref{T:GeneralLow} implies that, given the environment, any two trees have a comparable probability of being picked as the RSTRE. This is far from the truth. For the moment consider the complete graph on $n$ vertices and let $\mu$ be the uniform distribution on $[0,1]$. Then by the results of \cite{Fri85} mentioned in Remark \ref{R:mst_weight}, the minimum and maximum weighted trees $T_{\mathrm{min}}, T_{\mathrm{max}}$ satisfy
    \begin{align*}
        H( T_{\mathrm{min}}, \omega) &\approx \zeta(3), \\
        H( T_{\mathrm{max}}, \omega) &\approx n-1 - \zeta(3),
    \end{align*}
    and, by the law of large numbers, for any fixed tree $T_{\textrm{typical}}$
    \begin{equation*}
        H( T_{\textrm{typical}}, \omega) \approx \frac{n-1}{2}.
    \end{equation*}
    That is, if $\beta \gg \frac{1}{n}$, then none of the trees $T_{\mathrm{min}}, T_{\mathrm{max}}$ and $T_{\textrm{typical}}$ have a comparable weight. A more realistic heuristic for Theorem \ref{T:GeneralLow} is that the addition of weights turns the graph into a $n/\beta^{1/\alpha}$ regular graph. As long as $\beta \ll n^{1/\alpha}$ we might expect a UST behavior.
\end{remark}

\noindent A natural follow-up question is the following.

\begin{question}
    Order the spanning trees $T_1, T_2, \ldots$ of the complete graph in increasing weight. Given $\beta, \delta > 0$, what value of $k = k(\beta, \delta)$ must we pick such that
    \begin{equation*}
        \sum_{i=1}^k \bP(\cT = T_i) \geq 1 - \delta \, ?
    \end{equation*}
\end{question}

\subsection{Concentration of random walk observables}
In \cite{MSS24}, % we proved Theorem \ref{T:GeneralLow} under the assumption that $\mu$ was the uniform distribution on $[0,1]$. 
we showed that the bottleneck ratio in Definition \ref{D:bottleneck} concentrates around its unweighted version, and then applied Theorem \ref{T:Nach} to the weighted graph. The same proof (with a slightly different concentration inequality) works for any $\mu$ satisfying Assumption \ref{as:mu_tail}. We will only sketch the details. %We shall require the following lemma, and from now on will assume that $\deg_{\mathrm{min}} \geq C_2 \log n$ for some large constant $C_2$ to be determined later.

We write
\begin{equation*}
        \xi = \xi(\mu, \beta) = \E[\weight(e)] = \E[e^{-\beta \omega_e}] 
\end{equation*}
for the mean of the edge weight. Under Assumption \ref{as:mu_tail}, we will show in Lemma \ref{L:mu_tail_concentrate} (see Appendix \ref{AS:concentration}) that there exists a constant $C_B = C_B(\mu)$, such that for $S_m = \sum_{i=1}^m \weight(e_i)$ and $0 < \delta < 1$, we have
\begin{equation} \label{eq:mu_tail_concentrate}
    \P \big( |S_m - m \xi | \geq \delta m \xi \big) \leq 2 \exp \left( - C_B \frac{\delta^2 m}{ \max\{\beta^\alpha, 1\} } \right).
\end{equation}
This concentration inequality will allow us to show the following lemma for the lazy random walk on the weighted graph $(G_n, \weight)$.

\begin{lemma}[cf. {\cite[Lemma 3.4 and Lemma 3.5]{MSS24}}] \label{L:deg_cheeger_concentrate}
    There exists a constant $C_1 = C_1(b, \Delta, \mu)$ such that if $\beta \leq C_1 (\deg_{\mathrm{min}}/\log n)^{1/\alpha}$, then the following holds with high probability:
    \begin{enumerate}
        \item For all $v \in V_n$ we have
        \begin{equation} \label{eq:stationary}
            \frac{1}{3 \Delta n} \leq \pi(v) \leq \frac{3\Delta}{n};
        \end{equation}
    \item and for all $S \subset V_n$ with $0 < \pi(S) \leq 1/2$, it holds that
        \begin{equation} \label{eq:Phi_concentrate}
            \Phi_{(G_n,\weight)}(S) \geq \frac{b}{72\Delta^2}.
        \end{equation}
    \end{enumerate}

\end{lemma}
\begin{proof}[Proof sketch.]
    Let $|E_n| = m$ with $\deg_{\mathrm{min}} \leq \frac{2m}{n} \leq \deg_{\mathrm{max}}$. For $v \in V_n$ denote by $\deg_v$ the degree of $v$, and recall that $\weight(v) = \sum_{u \sim v} \weight(u,v)$ and
    \begin{equation*}
        \pi(v) = \frac{\weight(v)}{2 \sum_{e \in E} \weight(e)}.
    \end{equation*}
    Applying the concentration inequality in \eqref{eq:mu_tail_concentrate} for $\beta \geq 1$ gives, with a union bound, that for some constant $C_B > 0$
    \begin{align}
        \P\Big( \big|\sum_{e \in E_n} \weight(e) - m \xi \big| \geq \frac{m}{2} \xi \Big) &\leq 2 \exp \left( - C_B \frac{ m}{4 \beta^\alpha}\right), \label{eq:concen_all_edges} \\
        \P\left(\exists v \in V_n : |\weight(v) - \deg_v \xi | \geq  \frac{\deg_v}{2} \xi \right) &\leq 2 n\exp\Big(- C_B\frac{ \deg_{\mathrm{min}}}{ 4 \beta^{\alpha}}\Big). \label{eq:concen_stationary}
    \end{align}
    Clearly if $\beta \leq C_1 (\deg_{\mathrm{min}}/ \log n)^{1/\alpha}$ for some small $C_1$, then \eqref{eq:stationary} holds with high probability. If $0 \leq \beta \leq 1$, according to the concentration in \eqref{eq:mu_tail_concentrate}, we must replace $\beta$ in \eqref{eq:concen_all_edges} and \eqref{eq:concen_stationary} with $1$, and in this case both \eqref{eq:concen_all_edges} and \eqref{eq:concen_stationary} still go to zero as long as $\deg_{\mathrm{min}} \geq C_2 \log n$ for $C_2 = 8/C_B > 0$.

    The proof of \eqref{eq:Phi_concentrate} follows similarly. First, one has to show (using that $G_n$ is a $b$-expander) that for $S$ with $0 < \pi(S) < 1/2$, the sets $E(S, S^c)$ and $E(S, V_n)$ both have a size of order $\deg_{\mathrm{min}} |S| \asymp \deg_{\mathrm{max}} |S|$. Then one applies a union bound over all choices of $S$, and finally, one uses the concentration inequality in \eqref{eq:mu_tail_concentrate}. We omit the details and refer to the proof of Lemma 3.5 in \cite{MSS23} for more details.
\end{proof}

Using Lemma \ref{L:deg_cheeger_concentrate} we may now prove Theorem \ref{T:GeneralLow}.

\begin{proof}[Proof of Theorem \ref{T:GeneralLow}.]
    We verify the three conditions of \eqref{eq:Nbalanced}-\eqref{eq:Nescaping} in Theorem \ref{T:Nach}. By Lemma \ref{L:deg_cheeger_concentrate}, with high probability we have that
    \begin{align*}
         \frac{1}{3 \Delta n} \leq \pi(v) &\leq \frac{3\Delta}{n} \quad \forall v \in V_n,\\
       \frac{b}{72\Delta^2} &\leq \Phi_{(G,\weight)} .
    \end{align*}
    Let $M = \frac{b}{72\Delta^2}$, then a consequence of Theorem \ref{T:HeatCheeger} is that
    \begin{align*}
        \tmix &\leq 1 + \frac{4}{M^2} \int_{4/(3 \Delta n)}^8 \frac{1}{r} dr \\
        &\leq 1 + \frac{4}{M^2} \big( \log 8 + \log( 3 \Delta n) ) \leq C_\mathrm{mix} \log n ,
    \end{align*}
    for some constant $C_\mathrm{mix} = C_\mathrm{mix}(b,\Delta) > 0$. Furthermore, for any $t \geq 2$ and $\chi = \exp\big( - (t-1)M^2/4 \big) \pi(v)^{-1}$, we again have by Theorem \ref{T:HeatCheeger} that
    \begin{equation*}
        q_t(v,v) \leq \pi(v) + \exp \big( - (t-1) \frac{M^2}{4} \big).
    \end{equation*}
    Hence, with high probability
	\begin{equation*}
		\sum_{t=0}^{\tmix} (t+1) \sup_{v \in V_n}q_t (v,v)\leq \sum_{t=0}^{ C_{\mathrm{mix}}\log n} (t+1) \big( \frac{3\Delta}{n} + e^{-(t-1) \frac{M^2}{4}} \big) \leq C_{\mathrm{esc}},
	\end{equation*}
    for some constant $C_{\mathrm{esc}} = C_{\mathrm{esc}}(b, \Delta)$. Fix $\delta > 0$ as in the statement of Theorem \ref{T:GeneralLow}. If $n$ is large enough (depending on $\delta$), then applying Theorem \ref{T:Nach} with $D = 9 \Delta$, $\alpha = 1/4$ (or any other $0 < \alpha < 1/2$) and $\theta = C_{\mathrm{esc}}$ completes the proof.
\end{proof}

\subsection{Lack of concentration above $n^{1/\alpha}$}  \label{SS:no_conc}

We do not expect the low disorder proof to work, at least not without major modifications, when $\beta$ exceeds $n^{1/\alpha}$. Namely, we will show that the conditions in \eqref{eq:Nbalanced} and \eqref{eq:Nmixing} can only satisfied by parameters that are of order at least $n^c$, for any $c >0$, so that Theorem \ref{T:Nach} will no longer be applicable. 

Theorem \ref{T:Gnp_connect} hints towards the fact that the connection threshold for the \ErdosRenyi random graph occurs at $p = \log n/n$. Indeed, if $p = \log n/ 2n$, then (see e.g.\ \cite[Theorem 5.8]{vdH17}) $G_{n,p}$ is disconnected with high probability. Consequently, using the coupling of $\omega$ to $G_{n,p}$ as described in Section \ref{SS:CoupleRG}, there exists a vertex $u \in V_n$ such that
\begin{equation*}
    \omega_{(u,v)} \geq F^{-1}\big( \frac{\log n}{2n} \big) = \big( \frac{\log n}{2 c_{\mu} n} \big)^{1/\alpha} \qquad \forall v \neq u,
\end{equation*}
see also \eqref{eq:inverse_CDF}. As $\sum_{e \in E} \weight(e)$ concentrates when $\beta \ll n^{2/\alpha}$ (see Lemma \ref{L:mu_tail_concentrate}), we have for $\beta \gg n^{1/\alpha} (\log n)^{1 - 1/\alpha}$ that
\begin{equation*}
    \pi(v) \leq (n-1)\frac{\exp\big(-\beta \big( \frac{\log n}{2 c_{\mu}n} \big)^{1/\alpha} \big)}{n(n-1) \xi} = O( n^{-c}),
\end{equation*}
for any $c > 0$. Condition \eqref{eq:Nbalanced} is therefore never satisfied for $D$ that is of polynomial order. 

For simplicity's sake, assume now that $\beta \ll \min\{ n^{4/3\alpha}, n^{1/\alpha + 1/3} \}$. If we let $S = \cC_1(1/n)$, then by Theorem \ref{T:C1critical} we have $|E(S,S)| \asymp n^{1/3}$ and $|E(S,S^c)| \asymp n^{4/3}$ with high probability. Notice that the edges in $E(S,S)$ and $E(S,S^c)$ are conditioned on $\{ \omega_e \leq F^{-1}_\mu(1/n) \}$ and $\{ \omega_e > F^{-1}_\mu(1/n) \}$, respectively. 
% The calculations in \eqref{eq:xi_conditional} show that, for $\beta \gg n^{1/\alpha}$, we have
% \begin{align*}
%     \E[e^{- \beta \omega_e} \mid \omega_e \leq F^{-1}_\mu(1/n)] &\asymp n \beta^{-\alpha},\\
%     \E[e^{- \beta \omega_e} \mid \omega_e > F^{-1}_\mu(1/n)] &\asymp n^{1/\alpha - 1}  \beta^{- 1} \exp\big( -\frac{\beta}{n^{1/\alpha}} \big).
% \end{align*}
In Appendix \ref{AS:concentration}, the concentration arguments leading to \eqref{eq:full_cond_conc} show that, if $\beta \ll \min\{ n^{4/3\alpha}, n^{1/\alpha + 1/3}\}$, then with high probability
\begin{equation} \label{eq:condition_concen}
    \begin{aligned}
    \sum_{e \in E(S,S)} \weight(e) &\asymp n^{4/3} \beta^{-\alpha}, \\
    \sum_{e \in E(S,S^c)} \weight(e) &\asymp n^{1/\alpha +1/3} \beta^{-1} \exp\big( -\frac{\beta}{(c_\mu n)^{1/\alpha}} \big).
\end{aligned}
\end{equation}
Therefore, with high probability, 
\begin{equation*}
    \Phi_{(G, \weight)} \leq \frac{\sum_{e \in E(S,S^c)} \weight(e)}{2 \sum_{e \in E(S,S)} \weight(e)} \lesssim (n\beta)^{\alpha - 1} \exp\big(-\frac{\beta}{(c_\mu n)^{1/\alpha}} \big) .
\end{equation*}
If $\beta \gg n^{1/\alpha} \log n$, then the order of the mixing time becomes larger than any polynomial so that condition \eqref{eq:Nmixing} can not be fulfilled.

\section{High disorder} \label{SS:diamHigh}

For the remainder of this Chapter, we will assume that $G_n = K_n$ is the complete graph on $n$ vertices. We restate the high disorder result of Theorem \ref{T:main} in the following theorem. 

\begin{theorem} \label{T:highDisorder}
    Assume that $\mu$ satisfies Assumption \ref{as:mu_tail}. If $\beta \geq n^{1/\alpha + 1/3} (\log n)^{1/\alpha}$, then
    \begin{equation*}
        \widehat{\E}[ \diam(\cT)] = \Theta(n^{1/3}).
    \end{equation*}
    Furthermore, given $\delta > 0$, there exists a constant $C = C(\delta) > 0$ such that
    \begin{equation*}
        \widehat{\P} \big( C^{-1} n^{\frac{1}{3}} \leq \diam(\cT) \leq C n^\frac{1}{3} \big)  \geq 1 - \delta.
    \end{equation*}
\end{theorem}

The proof of the high disorder regime is more involved than the low order regime and follows a similar approach of \cite{ABR09}\footnote{We refer to \cite[Remark 4.11]{MSS24} for the novel adaptations required to adapt the proof to our model.}, who showed that the MST on the complete graph has a diameter of order $n^{1/3}$. Recall the coupling from Section \ref{SS:CoupleRG} to an \ErdosRenyi random graph, where we call an edge $p$-open if
\begin{equation*}
    \omega_e \leq F_\mu^{-1}(p),
\end{equation*}
and $p$-closed otherwise. To construct $G_{n,p}$, we retain $p$-open edges while removing any $p$-closed edges, and we order the connected clusters $\cC_{\ell}(p)$, $\ell \geq 1$, in decreasing size. Note that, by Assumption \ref{as:mu_tail}, for $0 \leq p \leq c_\mu \rho^\alpha$ we have
\begin{equation} \label{eq:inverse_CDF}
    F^{-1}_\mu(p) = \frac{1}{c_\mu^{1/\alpha}} p^{1/\alpha}.
\end{equation}

\begin{definition} \label{D:ClCp}
    For a vertex set $A\subset V_n$ we define the minimal subtree of $\cT=\cT_{n,\beta}^\omega$ containing all vertices in $A$ as
	\begin{equation}  \label{eq:defCbar}
		\cT_A \ := \bigcup_{u,v \in A} \gamma_{\cT}(u, v),
	\end{equation}
	where $\gamma_{\cT}(u, v)$ is the unique path in $\cT$ (regarded as a graph with its own vertex and edge set) connecting $u$ and $v$.
\end{definition} 

Theorem \ref{T:highDisorder} will follow from the next proposition along with the fact that the giant component within the critical window is ``tree-like.''

\begin{proposition} \label{P:DiaC1}
There exists a constant $g_0 = g_0(\mu) \geq 1$ such that if $\epsilon \geq n^{-1/3}$, $\beta \epsilon \geq (n \log n)^{1/\alpha}$ and $p_0 = (1 + g_0 \epsilon )/n$, then
\begin{equation}
    \widehat{\E}[ \diam(\cT)] = \widehat{\E}\Big[ \diam\Big( \ClCp{1}{p_0} \Big) \Big] + \bigO{\frac{n^{1/3}}{\sqrt{\epsilon n^{1/3}}}}. \label{eq:E_DiaC1}
\end{equation}
Furthermore,
\begin{equation}
    \diam\Big(\ClCp{1}{p_0} \Big) \leq  \diam (\cT ) \leq \diam\Big( \ClCp{1}{p_0} \Big) + \bigO{\frac{n^{1/3}}{\sqrt{\epsilon n^{1/3}}} + (\log n)^6} \label{eq:P_DiaC1}
\end{equation}
with $\widehat{\P}$-probability at least
\begin{equation*}
    1 - \bigO{\exp \Big(- \frac{1}{4} \sqrt{g_0 \epsilon n^{1/3}} \Big)} - \bigO{\frac{1}{n}}. \label{eq:DiaC1_probbound}
\end{equation*}
\end{proposition}

\begin{figure}[ht]
\centering
    \includegraphics[width=0.7\textwidth]{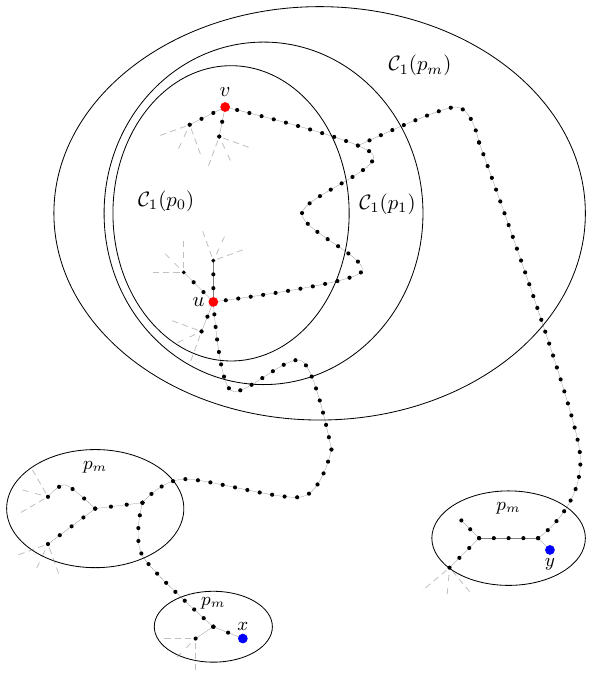}
    % \includestandalone[width=0.65\textwidth]{Write ups/glue_picture}%     without .tex extension
  % or use \input{mytikz}
  \caption{Typical vertices $x$ and $y$ are likely to be connected in $\cT$ via a short path from $x$ to $\ClCp{1}{p_m}$, a path inside $\ClCp{1}{p_m}$ and another short path from $\ClCp{1}{p_m}$ to $y$.}
  \label{fig:proofidea}
\end{figure}

The rough outline of the proof of Proposition \ref{P:DiaC1} is as follows:
\begin{enumerate} \label{page:proof_idea}
    \item Lemma \ref{L:gap}: For $p < q$, with $q - p$ ``large'', and any $u,v \in V_n$ in the same connected component of $G_{n,p}$, the path $\gamma_\cT(u,v)$ in the RSTRE between $u$ and $v$ typically uses only $q$-open edges. We will use this to show that with high probability $\ClCp{1}{p_0} \subseteq \cC_1(p_1)$, for $p=p_0$ and $q=p_1$ with $p_1-p_0\geq C \log n/\beta$.

    \item Lemma \ref{L:sizepm}: For $p_m$, chosen such that $\mathcal{C}_1(p_m)$ is of size of order $n/\log n$, the diameter of $\ClCp{1}{p_m}$ is comparable to the diameter of $\ClCp{1}{p_0}$ with $p_0 = (1+g_0 \epsilon)/n$ for some constant $g_0 \geq 1$.

    \item Lemma \ref{L:fasthitLERW}:  The LERW started at vertices in $V_n \setminus \mathcal{C}_1(p_m)$ ``quickly'' hits $\mathcal{C}_1(p_m)$ with high probability, which by Wilson's algorithm ensures that the diameter of $\cT$ remains of the same order as $\ClCp{1}{p_m}$.
\end{enumerate}
\noindent We refer to Figure \ref{fig:proofidea} for an illustration of the above steps.

\subsection{Gap in edge weights}

The first step is to show that for $p,q$, with $q - p$ large enough (depending on $\beta$), any two vertices that are in the same open $p$-cluster can not be connected in the UST via a path that uses a $q$-closed edge. For the proof of Theorem \ref{T:highDisorder}, it suffices to consider the case $p, q \approx 1/n$, where for $\mu$ satisfying Assumption \ref{as:mu_tail} the inverse CDF is given in \eqref{eq:inverse_CDF}.

\begin{lemma}[cf.\ {\cite[Lemma 4.3]{MSS24}}] \label{L:gap}
Suppose that $\epsilon  = o(1)$ and 
\begin{align*}
    p &= \frac{1 +  \epsilon}{n}, \\
    q &= \frac{1 +  t \epsilon}{n},
\end{align*}
for $t \ll \epsilon^{-1}$. For $u,v \in V_n$, define the event
\begin{equation*}
    F(u,v) := \big\{ u \xleftrightarrow{p} v, \,\gamma_{\cT}(u,v) \text{ contains an edge $e$ with }  \omega_e > q \big\}.
\end{equation*}
Then there exists a constant $c_g = c_g(\alpha, \rho, \mu) > 0$ such that for $n$ large enough
\begin{equation}
    \widehat{\P} \bigg( \bigcup_{u,v \in V_n, u \neq v} F(u,v) \bigg) \leq n^5 \exp \big( -c_g (t-1) \frac{\beta \epsilon }{n^\alpha}\big). \label{eq:gapUnion}
\end{equation}
\end{lemma}

\begin{proof}[Proof Sketch.]
    The event $F(u,v)$ is the union of the events
    \begin{equation}
    	F(u,v; e) :=  \big\{ u \xleftrightarrow{p} v, \, e \in \gamma_{\cT}(u,v), \, \omega_e > q \big\} .
    \end{equation}
    We show that $\widehat{\P}(F(u,v; e))$ is small so that a union bound will give \eqref{eq:gapUnion}. Assume first the special case that $e = (u,v)$. As any $p$-open path in $G$ has length at most $n$, Rayleigh's monotonicity principle and Kirchhoff's formula shows that
    \begin{equation} \label{eq:Kirch_gap}
        \bP^{\omega}_{\beta}(e \in \cT) = \weight(e) \effR{}{u}{v} \leq e^{-\beta F^{-1}_\mu(q) } \frac{n}{e^{-\beta F^{-1}_\mu(p)}} = n \exp \big(- \beta \frac{ q^{1/\alpha} - p^{1/\alpha} } {c_\mu^{1/\alpha}} \big),
    \end{equation}
    where we used \eqref{as:mu_tail} for $F^{-1}_\mu(p)$ and $p$ small enough, for which we require $n$ to be large enough.

    For the general case, if $e \in \gamma_{\cT}(u,v)$, then there must exist two (possibly empty) paths in $\cT$ from $u$ and $v$ to the endpoints of $e$. Conditional on this event, and using the spatial Markov property, we may contract these paths into single vertices. On this new graph, we may apply the same arguments as above to obtain an upper bound as in \eqref{eq:Kirch_gap}.
    
    To obtain \eqref{eq:gapUnion} from \eqref{eq:Kirch_gap}, notice that a Taylor expansion for $(1+x)^{1/\alpha}$ for $x$ close to zero, gives that
    \begin{equation*}
        (1 + t \epsilon)^{1/\alpha} - (1 + \epsilon)^{1/\alpha} = \frac{(t-1)}{\alpha} \epsilon + O(t^2\epsilon^2).
    \end{equation*}
    Hence, there exists some constant $c_g = c_g(\alpha, \rho, \mu)$, such that for $t \ll \epsilon^{-1}$ we have
    \begin{equation*}
        \beta \frac{ q^{1/\alpha} - p^{1/\alpha}}{c_\mu^{1/\alpha}} \geq c_g \beta \frac{(t-1)}{n^{1/\alpha}}.
    \end{equation*}
   This bound, together with a union bound over all edges and vertices for the events $F(u, v; e)$, establishes \eqref{eq:gapUnion}.
\end{proof}
    
\begin{figure}[t]
    \centering
    \includegraphics[width=0.9\textwidth]{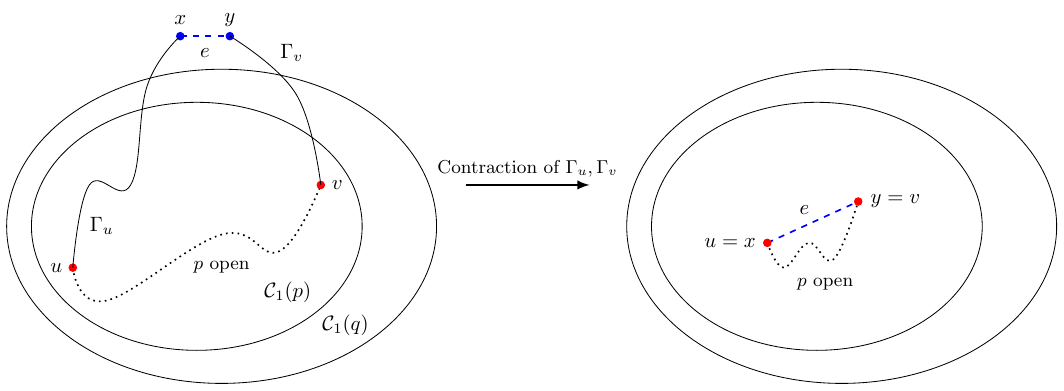}
    \caption{The path $\gamma_\cT(u,v)$ can be decomposed into a path $\Gamma_u$, an edge $e$ and another path $\Gamma_v$. After the contraction of $\Gamma_u$ and $\Gamma_v$, the edge $e$ becomes an edge between $u$ and $v$ in $(G', \weight)$.} 
    \label{fig:gap}
\end{figure}

We refer also to Figure \ref{fig:gap} and the proof of Lemma 4.3 in \cite{MSS24} for more details regarding the proof of Lemma \ref{L:gap}.

\begin{remark} \label{R:compare}
    When $t$ is fixed and $\beta \epsilon \geq C n^\alpha \log n$, for some large enough constant $C$, then we have a gap in the edge weight conductances:
    \begin{equation*}
        g(p,q) := \exp \big( -c_g (t-1) \frac{\beta \epsilon }{n^\alpha}\big) = O(n^{-6}).
    \end{equation*}
    For $\ClCp{1}{p}$ as in Definition \ref{D:ClCp}, Lemma \ref{L:gap} then implies that
    \begin{equation*}
        \ClCp{1}{p} \subseteq G_{n,q} ,
    \end{equation*}
    with high probability. If on the other hand $\beta \epsilon = O(n^{\alpha})$, then $g(p, q) = O(1)$. In this case, any two edges $e,f$ that are $p$-closed but $q$-open, satisfy that $\weight(e)/\weight(f)$ is uniformly bounded away from zero and infinity.
\end{remark}

\subsection{Building up to the giant}

We now consider the \ErdosRenyi random graph process $(G_{n,p})_{p \in [0,1]}$ at different snapshots $p_i$, where for some $g_0$ and $\epsilon = o(1)$ to be determined later, we set
\begin{equation} \label{eq:p_seq}
	\begin{aligned}
		g_i &= (5/4)^{i/2} g_0, \\
		p_i &= \frac{1 + g_i \epsilon}{n}. %\label{eq:p_seq}
	\end{aligned}
\end{equation}
Furthermore, denote by $m$ the integer satisfying  
\begin{equation*}
    m = m(n) = \min \big\{ i \geq 0 \, : \, g_i \epsilon \geq \frac{1}{\log n} \big\}.
\end{equation*}
We start with the tree $\ClCp{1}{p_0}$ on $\cC_1(p_0)$ (see Definition \ref{D:ClCp}) and over $m = O(\log n)$ steps extend this tree to the giant component at progressively larger snapshots until reaching $\ClCp{1}{p_m}$. The following lemma states that, for some suitable choice of $\epsilon$, the diameter of $\ClCp{1}{p_0}$ is of the same order as the diameter of $\ClCp{1}{p_m}$.

\begin{lemma}[cf. {\cite[Lemma 4.4]{MSS24}}] \label{L:sizepm}
Let $\epsilon \geq n^{-1/3}$ with $\beta \epsilon \geq n^\alpha \log n$. Then there exists some universal constant $g_0 > 0$ such that for $p_0$ as in \eqref{eq:p_seq}, we have
\begin{equation*}
    \widehat{\E} \Big[ \diam \big( \ClCp{1}{p_m} \big) \Big] \leq \widehat{\E} \Big[ \diam \big( \ClCp{1}{p_0} \big) \Big] +  \bigO{\frac{n^{1/3}}{\sqrt{\epsilon n^{1/3}}}}.
\end{equation*}
\end{lemma}

\begin{proof}[Proof sketch.]
    We briefly summarize the proof ideas from \cite[Lemma 4.4]{MSS24}. For connected graphs $H_1 \subseteq H_2$ denote by $\ell(H_2)$ the longest acyclic path in $H_2$, and by $\ell (H_2 \setminus H_1)$ the longest acyclic path in $H_2$ not using vertices of $H_1$. Any two vertices $u, v \in H_2$ can be connected by a path that consists of a path from $u$ to some $x_1 \in H_1$, a path from $v$ to some $x_2 \in H_1$, and a path from $x_1$ to $x_2$ within $H_1$. As observed in Lemma 2 of \cite{ABR09}, this gives us the following inequality
    \begin{equation}
        \diam (H_2) \leq \diam (H_1) + 2 \ell(H_2 \setminus H_1) + 2. \label{eq:long_path}
    \end{equation}
    
    Consider now the events 
    \begin{align}
        A(i) &:= \Big\{ |\mathcal{C}_1(p_i)| \in \Big[\frac{3}{2} g_i \epsilon n, \frac{5}{2} g_i \epsilon n\Big] \text{ and } \ell(\mathcal{C}_1(p_i)) \leq (g_i \epsilon n^{1/3})^4 n^{1/3} \Big\}, \label{eq:good1} \\
        B(i) &:= \Big\{  \ell \big( G_{n, p_{i+2}} \setminus \mathcal{C}_1(p_i) \big)\leq n^{1/3}/ \sqrt{g_i  \epsilon n^{1/3} } \Big\},\label{eq:good2} \\
        C(i) &:= \Big\{  \mathcal{C}_1(p_i) \subseteq \mathcal{C}_1(p_{i+1}) \Big\}.\label{eq:good3}
    \end{align}
    These events heuristically suggest that paths in $G_{n,p_i}$ are not too long, and any long path in $\cC_1(p_{i+1})$ must pass through $\cC_1(p_i)$. Assume now that for every $0 \leq i \leq m$ the events $A(i), B(i)$ and $C(i)$ hold, and that $\ClCp{1}{p_i} \subseteq G_{n,p_{i+1}}$. The latter holds with high probability by Lemma \ref{L:gap} if 
    \begin{equation*}
        \beta \epsilon \geq n^\alpha \log n,
    \end{equation*}
    and $g_0$ is large enough (see also Remark \ref{R:compare}). Denote by $D_i$ be the diameter of $\ClCp{1}{p_i}$, then by the definition of $B(i)$ and $C(i)$ we have
    \begin{align}
    D_{i+1} - D_i &\leq 2 \ell \left[ \ClCp{1}{p_{i+1}} \setminus \ClCp{1}{p_{i}} \right] + 2 \nonumber \\
    &\leq 2 \ell \left[ G_{n, p_{i+2}} \setminus  \mathcal{C}_1(p_i) \right] + 2 \nonumber \\
    &\leq \frac{2n^{1/3}}{\sqrt{g_i \epsilon n^{1/3}}} + 2. \label{eq:diff_in_D}
    \end{align}
    Repeating the same argument for every $0 \leq i \leq m$ gives that
    \begin{equation*}
        D_{m} - D_{0} \leq \sum_{i=0}^{m-1} \frac{2n^{1/3}}{\sqrt{g_i \epsilon n^{1/3}}} + 2 m = O \Big( \frac{2n^{1/3}}{\sqrt{\epsilon n^{1/3}}} \Big).
    \end{equation*}

    It remains to deal with the case that, for some $i^* > 0$, the events hold only for $i^* \leq i \leq m$ (here is where the event $A(i)$ comes into play) instead of for all $0 \leq i \leq m$. We leave the details to the proof of Lemma 4.4 in \cite{MSS24}.
\end{proof}

\subsection{Connections outside the giant}

The previous section established that $\diam(\ClCp{1}{p_0}) \asymp \diam(\ClCp{1}{p_m})$, the next step is to argue that 
\begin{equation*}
    \diam \big( \ClCp{1}{p_m} \big) \asymp \diam(\cT),
\end{equation*}
 i.e.\ the subtree $\ClCp{1}{p_m}$ captures the order of the diameter of the whole tree.

\begin{lemma}[cf.\ {\cite[Lemma 4.9]{MSS24}}] \label{L:fasthitLERW}
    Let $p_m$ be defined as in \eqref{eq:p_seq}. Then for all $v \in V_n$
    \begin{equation*}
        \widehat{\P} \Big( d_{\cT} \big( v, \ClCp{1}{p_m} \big) \geq (\log n)^6 \Big) = \bigO{ \frac{1}{n^2} }.
    \end{equation*}
\end{lemma}

The proof of Lemma \ref{L:fasthitLERW} is rather technical and we refer to \cite[Lemma 4.9]{MSS24} for the details. The heuristic argument is as follows. We can first construct $\ClCp{1}{p_m}$ by running Wilson's algorithm with starting points (and root) in $\cC_1(p_m)$. Given $\ClCp{1}{p_m}$, to find the path from $v \not\in \ClCp{1}{p_m}$ to $\ClCp{1}{p_m}$ we run a random walk started at $v$ and ended the first time it hits $\ClCp{1}{p_m}$, and loop erase the path. Each time the random walk explores a new $p_m$-open cluster other than $\cC_1(p_m)$, it can visit at most $|\cC_2(p_m)|$ many new vertices. Hence, we may bound the length of the path from $v$ to $\ClCp{1}{p_m}$ by the size of the 2nd largest component times the number of components the random walk visits. The key observation is that every time a new component is discovered, there is an approximately $|\mathcal{C}_1(p_m)|/n \asymp 1/\log n$ chance for the random walk to jump to $\mathcal{C}_1(p_m)$, so that at most a polylogarithmic number of clusters are visited. 

Combining the intermediate steps, we are now ready to sketch the proof of Proposition \ref{P:DiaC1}.

\begin{proof}[Proof sketch of Proposition \ref{P:DiaC1}]
    By Lemma \ref{L:sizepm}, we have
    \begin{align*}
        \widehat{\E} \Big[ \diam \big( \cT \big) \Big] &\leq \widehat{\E} \Big[ \diam \big( \ClCp{1}{p_m} \big) \Big] + 2\widehat{\E} \big[ \max_{v \in V_n} d_\cT \big( v,\ClCp{1}{p_m} \big) \big] \\
        &\leq \widehat{\E} \Big[ \diam \big( \ClCp{1}{p_0} \big) \Big] + 2\widehat{\E} \big[ \max_{v \in V_n} d_\cT \big( v,\ClCp{1}{p_m} \big) \big] + \bigO{\frac{n^{1/3}}{\sqrt{\epsilon n^{1/3}}}} .
    \end{align*}
    Furthermore, using Lemma \ref{L:fasthitLERW} we obtain that
    \begin{equation*}
         \widehat{\E} \big[ \max_{v \in V_n} d_\cT \big( v,\ClCp{1}{p_m} \big) \big] \leq  (\log n)^6 + n \widehat{\P} \Big( d_{\cT} \big( v, \ClCp{1}{p_m} \big) \geq (\log n)^6 \Big) = O\big((\log n)^6 \big).
    \end{equation*}
    We also refer to Section 4.5 of \cite{MSS24} for more details.
\end{proof}

\subsection{Proof sketch of Theorem \ref{T:highDisorder}}

As mentioned earlier, we now use Proposition \ref{P:DiaC1}, and the fact that for $\epsilon = n^{-1/3}$ and $p = (1 + g_0 \epsilon)/n$ the giant component is a tree with a bounded number of cycles, to prove Theorem \ref{T:highDisorder}. We will make use of the following straightforward lemma.
 
\begin{lemma}[Lemma 4.10 in \cite{MSS24}]\label{L:ratiodT1T2}
    Let $\cT_1, \cT_2$ be two spanning trees of a connected graph $G$. If $|E(\cT_2) \setminus E(\cT_1)| = k$, then
    \begin{equation}
        \frac{\diam(\cT_2)}{k+1} - 1 \leq \diam(\cT_1) \leq (k+1) \diam(\cT_2) + k   .\label{eq:ratio_dT1T2}
    \end{equation}
\end{lemma}

\begin{proof}[Proof sketch of Theorem \ref{T:highDisorder}]
     Let $\epsilon = n^{-1/3}$ and fix $\delta > 0$ as in the statement of Theorem \ref{T:highDisorder}. By Proposition \ref{P:DiaC1}, there exists a $g_0 = g_0(\mu, \delta) \geq 1$ such that, if $\beta \geq n^{1/3 + 1/\alpha} (\log n)^{1/\alpha}$ and $p_0 = (1+g_0 \epsilon)/n$, then 
     \begin{equation*}
         \diam(\cT) = \diam \big( \ClCp{1}{p_0} \big) + O(n^{1/3}),
     \end{equation*}
     both in $\widehat{\E}$-expectation and with $\widehat{\P}$-probability at least $1 - \delta/2$.
     
    Define the excess of a graph as $\mathrm{Exc}(G) = |E(G)| - |V(G)|$, and let $p_1 = (1 + (5/4)^{1/2}g_0 \epsilon)/n$. It can be shown that the following events
    \begin{enumerate}
        \item $\cC_1(p_0) \subseteq \ClCp{1}{p_0} \subseteq \cC_1(p_1)$ (see Lemma \ref{L:gap} and \eqref{eq:good3})
        \item $\diam(\cC_1(p_0)) \geq A^{-1} n^{1/3}$,  for some $A = A(\delta, g_0) > 0$ (see \cite[Theorem 1.1]{NP08})
        \item $\mathrm{Exc}(\cC_1(p_1)) \leq 200 g_0^3$ (see \cite[Lemma 2.9]{MSS24})
    \end{enumerate}
    all hold with large enough probability (possibly depending on $\delta$). For the remainder of the proof restrict to the above events. Suppose we are given $\cC_1(p_0)$ and $\ClCp{1}{p_0}$, then let $T^0$ be any spanning tree on $\cC_1(p_0)$. Using only edges that are $p_1$-open, we may extend $T^0$ to a spanning tree $T$ on the vertex set $V(\ClCp{1}{p_0})$. As the excess of $\cC_1(p_1)$ is small, the difference $|E(T) \setminus E(\ClCp{1}{p_0})|$ must also be small. Lastly, since the diameter of $T$ is lower bounded by the diameter of $T^0$, and therefore by the diameter of $\cC_1(p_0)$, Lemma \ref{L:ratiodT1T2} implies that the diameter of $\ClCp{1}{p_0}$ must be of order at least $n^{1/3}$.
\end{proof}

\begin{remark} \label{R:constant_dist_gap}
    The only place where the proof of Proposition \ref{P:DiaC1} depended on $\mu$ is in Lemma \ref{L:gap} and its applications. In particular, we require $\epsilon$ to be such that, for a fixed $t \geq 1$ (e.g.\ $t = (5/4)^{1/2}$ as in \eqref{eq:p_seq}), we have
    \begin{equation} \label{eq:depend_mu}
        \beta \Big( F_\mu^{-1}\big( \frac{1 + t \epsilon}{n} \big) - F_\mu^{-1}\big( \frac{1 +  \epsilon}{n} \big) \Big) \geq C \log n,
    \end{equation}
    where $C> 0$ is some large enough constant. If instead of taking, say, $\mu$ to be the law of $U$, a random variable uniformly distributed on $[0,1]$, if we let $\mu$ be the law of $-U^{-\gamma}$, then
    \begin{align*}
        \beta \Big( F_\mu^{-1}\big( \frac{1 + t \epsilon}{n} \big) - F_\mu^{-1}\big( \frac{1 + \epsilon}{n} \big) \Big) &= \beta \Big( \big( \frac{n}{1+ \epsilon} \big)^{-\gamma} - \big( \frac{n}{1+t \epsilon} \big)^{\gamma} \Big) \\
        &= \beta n^{\gamma} \big( (t-1) \epsilon + O(\epsilon^2) \big).
    \end{align*}
    If we fix $\beta =1$ and let $\epsilon = C n^{-\gamma} \log n $, then we obtain the same bound as required in \eqref{eq:depend_mu}. See also Remark \ref{R:correspondence}. 
\end{remark}

\subsection{Fixed and very large disorder}

When $\beta$ is either fixed in $n$ or $\beta = \beta(n) \gg n^5 \log n$, then we may use alternative methods to prove the behavior of the low and high disorder of Theorem \ref{T:main}. In this last section, we will only consider the case when $G$ is the complete graph and $\mu$ is the uniform distribution on $[0,1]$.

\subsection*{Fixed disorder strength}

In \cite{AS24}, the authors showed that if a sequence of weighted graphs converges to a connected graphon, then the rescaled weighted USTs converge to Aldous' Brownian continuum random tree. Graphons arise as continuum limits of sequences of dense graphs and can be encoded as symmetric measurable functions  $f:[0, 1]^2 \to [0,1]$. Roughly speaking, $f(x, y)$ can be interpreted as the probability that an edge exists between $x$ and $y$. We refer to \cite{LS06} or \cite[Section 2]{AS24} for a more formal introduction to graphons. The following proposition is an application of the results developed in \cite{AS24}, together with standard concentration inequalities. See Proposition 1.2 of \cite{MSS24} for more details. 

\begin{proposition}[{\cite[Proposition 1.2]{MSS24}}]\label{P:GHPbeta}
	Fix $\beta \geq 0$ and let $\cT_n$ be the RSTRE on the complete graph with $n$ vertices and uniform $[0,1]$ environment. Denote by $d_{\cT_n}$ the graph distance on $\cT_n$ and by $\nu_n$ the uniform distribution on the vertices of $K_n$. Then
	\begin{equation*}
	\big( \cT^\omega_{n, \beta}, \frac{1}{\sqrt{n}} d_{\cT^\omega_{n, \beta}}, \nu_n \big) \xrightarrow{(d)} \big( \cT, d_{\cT}, \nu \big),
	\end{equation*}
	where $(\cT, d_{\cT}, \nu)$ is Aldous' Brownian CRT equipped with its canonical mass measure $\nu$, and the convergence is in distribution with respect to the GHP (Gromov–Hausdorff–Prokhorov) distance between metric measure spaces.
\end{proposition}

\subsection*{Very large disorder strength}

We very briefly give an alternative argument for the fact that $\diam(\cT) = \Theta(n^{1/3})$ when $\beta$ is very large. The following essentially is an adaptation of Proposition \ref{P:counterexample}.

\begin{theorem}[cf.\  {\cite[Theorem 4.12]{MSS24}}] \label{T:SuperHighDisorder}
	Let $G_n = (V_n, E_n)$, $|V_n| = n$ and $|E_n| = m(n) = m$, be a sequence of graphs. Assign to each edge a weight $\weight(e) = \exp(-\beta \omega_e)$, where $(\omega_e)_{e \in E_n}$ are i.i.d.\ uniform on $[0,1]$. If $\beta \gg m^2 n \log n$, then w.h.p.\ the RSTRE concentrates on the MST.
\end{theorem} 
\begin{proof}[Proof sketch]
    We only need to verify that the trees $T_1, T_2$  with the lowest and second lowest Hamiltonian $H(T, \omega) = \sum_{e \in T} \omega_e$ satisfy
    \begin{equation*}
        \frac{n^{n-2} \weight(T_2)}{\weight(T_1)} \leq \exp\Big( (n-2) \log n - \beta( \min\limits_{e \neq f} | \omega_e - \omega_f|) \Big) \xrightarrow{n \rightarrow \infty} 0 .
    \end{equation*}
    The above may be proven in exactly the same way as in Proposition \ref{P:counterexample}.
\end{proof}
\noindent As the MST on the complete graph has a diameter of order $n^{1/3}$, the same must be true for the RSTRE when $\beta \gg n^5 \log n$.

\begin{remark}
    On the complete graph, we can improve the lower bound of $n^5 \log n$ in Theorem \ref{T:SuperHighDisorder} to $n^4 (\log n)^2$. Indeed, the random graph $G_{n,p}$ is, with high probability, connected at $p = 2 \log n / n$ (cf.\ Theorem \ref{T:Gnp_connect}), so Lemma \ref{L:gap} shows that any edges with, say, $\omega_e \geq 3 \log n / n$ are not contained in the RSTRE with high probability. Therefore, we only need to compare on the order of $n \log n$ random variables that are uniformly distributed on $[0, 3 \log n / n]$. The typical minimum gap between such random variables is of order $n^{-3} (\log n)^{-1}$.
\end{remark}
\noindent The requirement $\beta = n^{4 + o(1)}$ in the above remark for ensuring that the RSTRE concentrates on the MST is far from optimal. In fact, the following theorem due to K\'{u}sz shows that the correct threshold is $\beta = n^{2 + o(1)}$. 

\begin{theorem}[cf.\ {\cite[Theorem 1.5]{K24}}]
    Let $\cT^\omega_{n, \beta}$ be the RSTRE on the complete graph $K_n=(V_n, E_n)$ with $n$ vertices and i.i.d.\ disorder variables $(\omega_e)_{e\in E_n}$ uniformly distributed on $[0,1]$. Denote by $M_n$ the MST coupled to $\omega$ as in \ref{SS:CoupleRG}. If $\beta \ll n^2/\log n$, then
    \begin{equation*}
        \widehat{\P}( \cT^\omega_{n, \beta} = M_n) \rightarrow 0.
    \end{equation*}
    If on the other hand $\beta \gg n^2 \log n$, then
    \begin{equation*}
        \widehat{\P}( \cT^\omega_{n, \beta} = M_n) \rightarrow 1.
    \end{equation*}
\end{theorem}

\chapter{Graph reduction}
\label{ch:conj}

In this final chapter, we provide heuristics for the intermediate regime $n \leq \beta \leq n^{4/3}$ in Conjecture \ref{C:Intermediate}, where we assume that $\mu$ is the uniform distribution on $[0,1]$. The main argument is that the RSTRE on the giant component can be reduced to an RSTRE on a smaller, bounded-degree graph, where we expect to see UST behavior. As a byproduct of this reduction, we will show that the diameter of an unweighted UST on the slightly supercritical giant component is of order square root number of vertices (up to some logarithmic factors) of the giant component. The main results in this chapter are Theorem \ref{T:Kernel} and Theorem \ref{T:USTgiantdiam}. See also Proposition \ref{P:non-trivial} for a non-trivial upper bound on the diameter in the regime $\beta = n^{1 + \gamma}$ with $1/4 < \gamma < 1/3$.

\section{Outline}
We will make the following assumption on $\beta$, which ensures that the kernel graph is a bounded-degree graph.
\begin{assumption} \label{as:betaeps_s}
    There exists $s > 0$ such that $\beta = \beta(n)$ and $\epsilon = \epsilon(n)$ satisfy
    \begin{gather}
        n^{1+s} \ll \beta \ll n^{4/3} \log n,  \label{eq:beta_inter_s}\\
        \epsilon = g_0 \frac{n \log n}{\beta} \gg n^{-\frac{1}{3}}, \label{eq:eps_inter}
    \end{gather}
    where $g_0$ is the constant from Proposition \ref{P:DiaC1}.
\end{assumption}

\noindent The following steps provide a rough outline of our heuristic approach for Conjecture \ref{C:Intermediate} in the intermediate regime $\beta =n^{1+\gamma}$ with $0 < \gamma < 1/3$: 
\begin{enumerate}
    \item Remove all edges outside of $\cC_1(p)$ with $p = (1 + \epsilon )/n$ (see Section \ref{S:edge_remove}).
    \item Reduce the giant component to its kernel graph (see Section \ref{S:Kernel}).
    \item Argue that the RSTRE restricted to the kernel graph gives the correct order of the diameter (see Theorem \ref{T:Kernel}).
    \item Use the series law to calculate the distribution of edge weights (see Section \ref{S:macro_weight}).
\end{enumerate}

\noindent See in particular Section \ref{SS:heuristics} for a more detailed explanation of how these steps give rise to a heuristic argument. 

\section{Edge removal outside the giant} \label{S:edge_remove}

As we have seen in Proposition \ref{P:DiaC1}, the diameter of the random spanning tree $\cT$ is essentially determined by the diameter of the subtree connecting vertices in a (slightly) supercritical giant component. However, when constructing the connected component of the giant component in $\cT$, it is not immediately obvious that we may remove all the edges with small weights without changing the distribution too much. %unclear if the distribution chthat the distribution of removing an edge is equivalent to conditioning on the edge not appearing in the trajectory. %Indeed, the spatial Markov property only applies to the full tree, and for, e.g.\ for directed graphs, one may construct easy counterexamples showing that the distributions are drastically different. However, 
In Lemma \ref{L:TVgap}, we will strengthen the result of Lemma \ref{L:gap}, showing that this is indeed the case. %indeon the giant shows we may indeed remove edges outside of the giant component. We will work under the following assumption to ensure that the reduced graph has bounded degree.

Consider the random graph process $(G_{n,p})_{p \in [0,1]}$ at the following two parameters:
\begin{equation}
    p_0 = \frac{1 + \epsilon}{n} \qquad \text{and} \qquad p_1 = \frac{1 +  \sqrt{5/4} \epsilon}{n}, \label{eq:snapshots01}
\end{equation}
see also \eqref{eq:p_seq}. Recall the coupling of $\omega$ to $G_{n,p}$ from Section \ref{SS:CoupleRG}, and that $\ClCp{1}{p_0}$ is the minimal connected subgraph of the RSTRE containing all vertices of $\cC_1(p_0)$ (see Definition \ref{D:ClCp}).

\begin{definition}
    Define $\cT_{p_1}$ as the RSTRE on $\mathbb{T}(\cC_1(p_1))$ (the set of spanning trees of $\cC_1(p_1)$), where the disorder variables are inherited from $(\omega_e)_{e \in E(\cC_1(p_1))}$. If $\omega$ is such that $\cC_1(p_0) \subseteq \cC_1(p_1)$, further define $\cT_{p_1 \mid p_0}$ as the minimally connected subgraph of $\cT_{p_1}$ containing all vertices of $\cC_1(p_0)$.
\end{definition}

% CHAT GPT CORRECTION FROM HERE

\begin{remark}
    The trees $\cT_{p_1 \mid p_0}$ and $\ClCp{1}{p_0}$ differ in the sense that, for the construction of $\cT_{p_1 \mid p_0}$, we have already removed all vertices and edges outside of $\cC_1(p_1)$. In contrast, for $\ClCp{1}{p_0}$, the edges outside $\cC_1(p_1)$ only do not appear in $\ClCp{1}{p_0}$ w.h.p.\ (see Lemma \ref{L:gap}).
\end{remark}

The following lemma will allow us to work with $\cT_{p_1 \mid p_0}$ instead of $\ClCp{1}{p_0}$.

\begin{lemma} \label{L:TVgap}
    Assume $\omega$ is such that $\cC_1(p_0) \subseteq \cC_1(p_1)$. Denote by $\nu_1$ and $\nu_2$ the laws of $\ClCp{1}{p_0}$ and $\cT_{p_1 \mid p_0}$, respectively. Then
    \begin{equation} \label{eq:TVgap_eq}
        \Vert \nu_1 - \nu_2\Vert_{\mathrm{TV}} \leq n^4 e^{-\beta(p_1 - p_0)},
    \end{equation}
    where $\Vert \cdot \Vert_{\mathrm{TV}}$ is the total variation norm. 
\end{lemma}
\noindent The proof of Lemma \ref{L:TVgap} is not too difficult and follows in a similar spirit to the proof of Lemma \ref{L:gap}. Since the proof is somewhat technical, we defer it and the definition of the total variation distance $\Vert \cdot \Vert_{\mathrm{TV}}$ to Appendix \ref{A:TV}. Using Lemma \ref{L:TVgap} and the arguments in Section \ref{SS:diamHigh} we obtain the following result.

\begin{lemma} \label{L:Kernel_reduce}
   Let $p_1$ be as in \eqref{eq:snapshots01}. Then, with high $\widehat{\P}$-probability,
    \begin{equation*}
       \diam(\cT_{p_1}) - o(n^{1/3}) \leq \diam(\cT) \leq \diam(\cT_{p_1}) + o(n^{1/3}).
    \end{equation*}
\end{lemma}
\begin{proof}
    By Assumption \ref{as:betaeps_s}, the parameters satisfy
    \begin{gather*}
        \epsilon n^{1/3} \rightarrow \infty, \\
        n^5 e^{-\beta(p_1 - p_0)}  = \bigO{\frac{1}{n}},
    \end{gather*}
    so that the conclusions of Section \ref{SS:diamHigh} (where we absorbed $g_0$ into $\epsilon$) all hold with high probability. Proposition \ref{P:DiaC1} gives us, with high $\widehat{\P}$-probability, that
    \begin{equation*}
        {\diam}\big(\ClCp{1}{p_0} \big) \leq  {\diam} (\cT ) \leq {\diam}\big( \ClCp{1}{p_0} \big) + o(n^{1/3}).
    \end{equation*}
    Furthermore, as the event $\{\cC_1(p_0) \subseteq \cC_1(p_1) \}$ in \eqref{eq:good3} holds with high probability whenever $\epsilon n^{1/3} \rightarrow \infty$ (see \cite[Lemma 4.7]{MSS24} or \cite{Luc90}), Lemma \ref{L:TVgap} gives a coupling $\lambda$ with
    \begin{equation*}
        \lambda \Big( {\diam}(\ClCp{1}{p_0}) \neq {\diam}(\cT_{p_1 \mid p_0}) \Big) \leq n^4 e^{-\beta (p_1 - p_0)} = O\big( \frac{1}{n} \big),
    \end{equation*}
    so that we may equivalently bound the diameter of $\cT_{p_1 \mid p_0}$ instead of the diameter of $\ClCp{1}{p_0}$. Using again the probability bound of Lemma 4.7 in \cite{MSS24} for the event $B(i)$ in \eqref{eq:good2}, shows that with high probability the longest acyclic path in $\cC_1(p_1) \setminus \cC_1(p_0)$ is short. Thus,
    \begin{equation*}
        {\diam}(\cT_{p_1}) - \bigO{\frac{n^{1/3}}{\sqrt{\epsilon n^{1/3}}}} \leq {\diam}(\cT_{p_1 \mid p_0}) \leq {\diam}(\cT_{p_1}),
    \end{equation*}
    from which the assertion readily follows.
\end{proof}

\section{Kernel graph} \label{S:Kernel}

Given a finite graph $G$, the \textit{2-core} of $G$ is the maximal subgraph $H \subseteq G$ in which each vertex has a degree of at least $2$. We can construct the 2-core iteratively by removing all vertices of degree $1$ until all remaining vertices have a degree of at least $2$. The \textit{kernel} $\Ker(G)$ of $G$ is obtained by removing all disjoint cycles from the 2-core and contracting any maximal 2-path (a path in which every interior vertex has degree $2$) into a single edge. We denote the vertices of $\Ker(G)$ by capital letters and identify them with their corresponding vertices in the 2-core of $G$, i.e. those vertices that are not interior vertices in maximal 2-paths. Furthermore, we implicitly define a function $\phi$ that maps an edge $(U,V) \in E(\Ker(G))$ to edges in $E(G)$ via the relation
\begin{equation} \label{eq:edge_mapping}
    \phi(U,V) = \{ e_1, e_2, \ldots, e_\ell \} \subseteq E(G),
\end{equation}
if and only if the edges $e_1, e_2, \ldots, e_\ell$ in the 2-core of $G$ are contracted to obtain $(U,V)$ in the kernel.

Next, consider the giant component $\cC_1(p_1)$ and let $\mathcal{K} = \Ker(\cC_1(p_1))$ be the kernel graph of $\cC_1(p_1)$. We denote by $\cT^\Ker_{p_1}$ the RSTRE on $\mathbb{T}(\mathcal{K})$ with edge weights $\widehat{\weight}(\cdot)$ given by
\begin{equation} \label{eq:Kernel_weight_dist}
    \widehat{\weight}(U,V) = \Big( \sum_{e \in \phi(U,V)} \frac{1}{\weight(e)} \ \Big)^{-1}. 
\end{equation}
In Section \ref{S:macro_weight}, we will give bounds on the weights of typical edges in the kernel $\textrm{Ker}(\cC_1(p_1))$. Our main goal in this chapter is to prove the following theorem.
\begin{theorem} \label{T:Kernel}
    Suppose that $\beta$ and $\epsilon$ are, for some $s > 0$, as in Assumption \ref{as:betaeps_s}. Then there exists a constant $C_K = C_K(s) > 0$ such that, for $n$ large enough, we have with high $\widehat{\P}$-probability that
    % \begin{equation}
    %     \frac{\epsilon^{-1}}{C_K \log n} {\diam}(\cT^\Ker_{p_1}) - o(n^{1/3}) \leq {\diam}(\cT) \leq 3 \epsilon^{-1} {\diam}(\cT^\Ker_{p_1}) \log n + o(n^{1/3}) +O( \epsilon^{-1} \log n).
    % \end{equation}
    \begin{equation}
        \frac{\epsilon^{-1}}{C_K \log n} {\diam}(\cT^\Ker_{p_1}) - o(n^{1/3}) \leq {\diam}(\cT)
    \end{equation}
    and
    \begin{equation}
        {\diam}(\cT) \leq 3 \epsilon^{-1} {\diam}(\cT^\Ker_{p_1}) \log n + o(n^{1/3}) +O( \epsilon^{-1} \log n).
    \end{equation}
\end{theorem}

\noindent The main argument in the proof of Theorem \ref{T:Kernel} is that we can explicitly analyze the 2-core and kernel of the giant component in the slightly supercritical regime, and couple spanning trees on the kernel with spanning trees on the giant component. We remark that the above theorem leads to the following non-trivial bounds on the diameter of the RSTRE when $\beta = n^{1 + \gamma}$ and $1/4 < \gamma < 1/3$.

\begin{proposition} \label{P:non-trivial}
    Let $\beta =  n^{1 + \gamma}$ with $ 0 < \gamma < 1/3$. Then with high $\widehat{\P}$-probability we have
    \begin{equation*}
        \diam(\cT) = O \big( n^{1 - 2 \gamma} \log n \big) .
    \end{equation*}
\end{proposition}

\begin{proof}
    The kernel graph of $\cC_1(p_1)$ has, with high probability, of order $\epsilon^3 n = O( n^{1- 3 \gamma} (\log n)^2)$ many vertices (see e.g.\ Theorem \ref{T:C1structure}). Thus, by Theorem \ref{T:Kernel}, we have with high probability that 
    \begin{equation*}
        \diam(\cT) = O \Big( \epsilon^{-1}  n^{1- 3 \gamma} (\log n)^2  \Big) = O \big( n^{1 - 2 \gamma} \log n \big). \qedhere
    \end{equation*}
\end{proof}

\begin{remark}
    Observe that if $\gamma > 1/4$ in Proposition \ref{P:non-trivial}, then $1 - 2 \gamma < 1/2$, meaning the diameter of the RSTRE is of order much smaller than the diameter of the unweighted UST. In particular, Proposition \ref{P:non-trivial} shows that there cannot be a sharp transition in the behavior of the diameter around $\beta = n^{4/3}$.
\end{remark}

\subsection{Structure of the Giant Component}

Denote by $\mathcal{N}(\mu, \sigma^2)$ the normal distribution with mean\footnote{In the remainder of this chapter, $\mu$ will denote some real number instead of the disorder distribution, which is always assumed to be uniform on $[0,1]$.} $\mu$ and variance $\sigma^2$, by ${\mathrm{Geom}}(p)$ the geometric distribution with mean $1/p$, and by ${\mathrm{Poisson}}(\lambda)$ the Poisson distribution with mean $\lambda$.  The following theorem gives an explicit way to construct the 2-core and kernel of $\cC_1(p)$ for $p$ in the slightly supercritical regime.

\begin{theorem}[Theorem 2 in \cite{DKLP11}] \label{T:C1structure}
    Let $\cC_1(p)$ be the largest component of $G_{n,p}$ for $p= \frac{1+\epsilon}{n}$, where $\epsilon^3 n \rightarrow \infty$ and $\epsilon \rightarrow 0$. Let $\mu < 1$ denote the conjugate of $1+\epsilon$, that is, $\mu e^{-\mu} = (1+\epsilon)e^{-(1+\epsilon)}$. Then $\cC_1(p)$ is contiguous to the following model $\tilde{\cC}_1(p)$:
    \begin{enumerate}
        \item Let $\Lambda \sim \mathcal{N}(1+\epsilon - \mu, \frac{1}{\epsilon n})$ and assign i.i.d.\ variables $D_u \sim {\mathrm{Poisson}}(\Lambda)$ ($ u \in \{1,2, \ldots, n\})$ to the vertices, conditioned that $\sum D_u 1_{D_u \geq 3}$ is even. Let
        \begin{equation}
            N_k = \# \{ u : D_u = k\} \qquad \text{and} \qquad N = \sum_{k \geq 3} N_k.
        \end{equation}
        Select a random multigraph $\tilde{\mathcal{K}}$ on $n$ vertices, uniformly among all multigraphs with $N_k$ vertices of degree $k$ for $k\geq 3$.

        \item Replace the edges of $\tilde{\mathcal{K}}$ by paths of lengths i.i.d.\ ${\mathrm{Geom}}(1 - \mu)$.

        \item Attach an independent ${\mathrm{Poisson}}(\mu)$-Galton-Watson tree to each vertex.
    \end{enumerate}
    That is, $\P(\tilde{\cC}_1(p) \in \mathcal{A}) \rightarrow 0$ implies $\P(\cC_1(p) \in \mathcal{A}) \rightarrow 0$ for any set of graphs $\mathcal{A}$.
\end{theorem}
The above theorem states that we can construct $\cC_1(p)$, or more precisely an asymptotic approximation of $\cC_1(p)$, in three steps. First, we consider the kernel graph as a configuration model, then we extend the paths by multiplying the length with a geometric random variable to obtain the $2$-core, and finally, we attach trees to each vertex. Importantly, there is a lot of independence between the different steps in this model. We remark that the reason for requiring $\beta \gg n^{1 + s}$, $s> 0$, in Assumption \ref{as:betaeps_s} is that the kernel graph in step 1 will have a bounded degree.

\subsection{Coupling to kernel graph}

To prove Theorem \ref{T:Kernel}, we will need a coupling between the RSTRE on $\cC_1(p_1)$ and the RSTRE on the kernel graph $\mathcal{K}$ of $\cC_1(p_1)$. For a weighted graph $(G, \weight)$, we write $\cT(G, \weight)$ for a spanning tree distributed according to the weighted UST measure on $G$, denoted by  $\bP^\weight_G$. Furthermore, we define $r(e) := 1/\weight(e)$ as the resistance of an edge. Recall from \eqref{eq:edge_mapping} that $\phi : E(\Ker(G)) \rightarrow E(G)$ maps an edge $(U,V)$ in the kernel to the edges in $G$ that are contracted to obtain $(U,V)$.

\begin{lemma} \label{L:couple_mac_mic}
    Let $(G, \weight)$ be a finite connected weighted graph and let $\mathcal{K}$ be the kernel of $G$. Assign each edge $(U,V) \in E(\mathcal{K})$ the weight
    \begin{equation} %\label{eq:Ker_series}
        \widehat{\weight}(U,V) = \frac{1}{\sum_{e \in \phi(U,V)} \frac{1}{\weight(e)}}. 
    \end{equation}
    Then there exists a coupling between $\cT(\mathcal{K},\widehat{\weight})$ and $\cT(G, \weight)$ such that
    \begin{equation} \label{eq:coupleK_G}
        (U, V) \in \cT(\mathcal{K},\widehat{\weight}) \iff \phi(U,V) \subset \cT(G, \weight). 
    \end{equation}
\end{lemma}
\begin{proof}
    % \lu?ca{Can we make this proof quicker/nicer?}
    We first consider the simpler case in which, to obtain the kernel from the 2-core, we only contract a path of length $2$ into a single edge. Namely, let $(G_1,\weight)$ be a weighted (multi) graph, and let $v \in V(G_1)$ have degree $2$ with incident edges $e_1 = (u, v)$ and $e_2 = (v,x)$ having resistances $r_1$ and $r_2$, respectively. Denote by $Q_u^{G_1}$ the random walk measure on $(G_1, \weight)$ started at $u$. Consider the weighted graph $(G_2, \widehat{\weight})$ where we remove the vertex $v$ and add a new edge $e =(u,x)$ with resistance $r = 1/\widehat{\weight}(e) = r_1 + r_2$. All other weights remain the same. We construct $\cT(G_1, \weight)$ by running Wilson's algorithm with root $x$ and first vertex $u$. Clearly, $e_1$ and $e_2$ are both in the spanning tree if and only if the first LERW step is from $u$ to $v$. By the Laplacian random walk representation of the LERW (see e.g.\ Exercise 4.1 in \cite{LP16}), we may calculate this probability as
    \begin{align*}
        \bP^\weight_{G_1} \big( e_1, e_2 \in \cT(G_1,\weight) \big) &= Q_u^{G_1}(X_1 = v \mid \tau_x < \tau_u^+) \\
        &= \frac{\weight(u,v)}{\weight(u)} \frac{\weight(v,x)}{\weight(u,v) + \weight(v,x)} \frac{1}{Q_u^{G_1}(\tau_x < \tau_u^+)} \\
        &= \frac{\weight(u,v) \weight(v,x)}{\weight(u,v) + \weight(v,x)} \effR{G_1}{u}{x}.
    \end{align*}
    Notice that
    \begin{align*}
        \frac{\weight(u,v) \weight(v,x)}{\weight(u,v) + \weight(v,x)} &= \frac{1}{\frac{1}{\weight(u,v)} + \frac{1}{\weight(v,x)}}, \\
        \effR{G_1}{u}{x} &= \effR{G_2}{u}{x},
    \end{align*}
    where the second equation is an application of the series law, so that by Kirchhoff's Formula (Theorem \ref{T:Kirchhoff}) we have
    \begin{align*}
        \bP^\weight_{G_1} \big( e_1, e_2 \in \cT(G_1,\weight) \big) &= \frac{1}{\frac{1}{\weight(u,v)} + \frac{1}{\weight(v,x)}} \effR{G_1}{u}{x} = \frac{1}{r_1 + r_2} \effR{G_2}{u}{x} \\
        &= \frac{1}{r} \effR{G_2}{u}{x} = \bP^{\widehat{\weight}}_{G_2}(e \in \cT(G_2,\widehat{\weight})).
    \end{align*}
    
    We proceed to the case where the contracted paths have lengths longer than $2$. Let $U,V \in V(\mathcal{K})$ be given. Repeating the above argument for every interior vertex of degree 2 along the induced path of $\phi(U,V)$, and noting that $r(U,V) = \sum_{e \in \phi(U,V)} r(e)$, gives
    \begin{equation} \label{eq:Marco_micro}
        \bP^{\widehat{\weight}}_{\mathcal{K}} \big( (U,V) \in \cT(\mathcal{K}, \widehat{\weight}) \big) = \bP^\weight_{G}\big( \phi(U,V) \subset \cT(G, \weight) \big).
    \end{equation}
    Arbitrarily order the edges of $\mathcal{K}$ as $(U_1, V_1), \ldots (U_m, V_m)$. We sample $\cT(\mathcal{K}, \widehat{\weight})$ in the following way: in the first step, include the edge $(U_1, V_1)$ with probability $\bP^{\widehat{\weight}}_{\mathcal{K}}\big( (U_1, V_1) \in \cT(\mathcal{K}, \widehat{\weight}) \big)$ and then either contract or remove $(U_1, V_1)$, depending on whether $(U_1, V_1)$ was included or not, to obtain a new graph $\mathcal{K}_1$. In the $i$-th step, $2 \leq i \leq m$, include the edge $(U_i, V_i)$ with probability $\bP^{\widehat{\weight}}_{\mathcal{K}_{i-1}}\big( (U_i,V_i) \in \cT(\mathcal{K}_{i-1}, \widehat{\weight}) \big)$ and generate a new graph $\mathcal{K}_i$ by either contracting the edge or removing it. By the spatial Markov property (Lemma \ref{L:USTmarkov}), this construction has the same law as the weighted UST measure on $\mathcal{K}$.

    We can similarly sample paths of $\cT(G, \weight)$ by including all edges in $\phi(U_i,V_i)$ with probability $\bP^\weight_{G_{i-1}}\big( \phi(U_i,V_i) \subset \cT(G_i, \weight) \big)$, and then contracting all edges in $\phi(U_{i-1},V_{i-1})$, or deleting an arbitrary edge in $\phi(U_i, V_i)$ if the whole path is not in $\cT(G_i, \weight)$, to obtain $G_i$, with $G_0 = G$. We note that deleting an edge in $\phi(U_i, V_i)$ has the same effect on the effective resistance on $G_i$ as deleting $(U_i, V_i)$ has on the effective resistance on $\mathcal{K}_i$. By the argument around \eqref{eq:Marco_micro}, in each step
    \begin{equation*}
        \bP^{\widehat{\weight}}_{\mathcal{K}_{i-1}} \big( (U_i,V_i) \in \cT(\mathcal{K}_{i-1}, \widehat{\weight}) \big) = \bP^\weight_{G_{i-1}}\big( \phi(U_i,V_i) \subset \cT(G_{i-1}, \weight) \big),
    \end{equation*}
    so that we may couple $\cT(\mathcal{K},\widehat{\weight})$ and $\cT(G, \weight)$ as required in \eqref{eq:coupleK_G}. 
\end{proof}
%\lu%ca{[[DRAFT COMMENT: The last part of the proof can do with some improvements... I'll try rewriting it a bit.]]}

%%%%%%%%%%%% NEW SECTION %%%%%%%%%%%%%

\section{Proof of Theorem \ref{T:Kernel}}

Using Lemma \ref{L:Kernel_reduce}, Theorem \ref{T:C1structure} and Lemma \ref{L:couple_mac_mic}, we may now prove Theorem \ref{T:Kernel}. The upper bound on the diameter will follow from the fact that the geometric random variables in step 2 of Theorem \ref{T:C1structure} can not get much larger than a logarithmic factor of their mean, and that the attached ${\mathrm{Poisson}}(\mu)$ trees in step 3 have a short diameter that only gives an additive (and not multiplicative) factor. The lower bound will use that not all geometric random variables in step 2 of Theorem \ref{T:C1structure} can be too short, so that a long path in the kernel will also give a long path in the 2-core. We will prove the following theorem, and remark that this together with Lemma \ref{L:Kernel_reduce} immediately implies Theorem \ref{T:Kernel}.

\begin{theorem} \label{T:general_kernel}
    Let $p = (1 + \epsilon)/n$ with $n^{1/3} \ll \epsilon \ll n^{-s}$, for some $s > 0$. %Let $\mu_n$ be any disorder distribution, and 
    Denote by $\cT_p$ and $\cT^\Ker_p$ the RSTREs on $\cC_1(p)$ and on $\Ker(\cC_1(p))$, respectively, where the weights on the edges in the kernel are given via the series law as in \eqref{eq:Kernel_weight_dist}. Then, there exists a constant $C_K = C_K(s) > 0$ such that there is a coupling between $\cT_p$ and $\cT^\Ker_p$ satisfying
    \begin{equation}
        \frac{\epsilon^{-1}}{C_K \log n} {\diam}(\cT^\Ker_{p})  \leq {\diam}(\cT_p) \leq 3 \epsilon^{-1} {\diam}(\cT^\Ker_{p}) \log n  +O( \epsilon^{-1} \log n),
    \end{equation}
    with high $\P$-probability, provided that $n$ is large enough.
\end{theorem}

\begin{proof}
    We first prove the upper bound on the diameter. For a graph $G$, we denote by $L_{\mathcal{K}}(G)$ the length of the longest contraction of a path in the 2-core of $G$ (denoted by $\textrm{2-core}(G)$) to obtain the kernel $\mathcal{K}$ of $G$. Consider the events
    \begin{align}
       \mathcal{B}_1 = \mathcal{B}_1(G) &:= \big\{ L_{\mathcal{K}}(G) > 3 \epsilon^{-1} \log n \big\}, \label{eq:2-core_len} \\
        \mathcal{B}_2 = \mathcal{B}_2(G) &:= \big\{ \max_{u \in V(G)} \min_{v \in \textrm{2-core}(G) }d_{G}(u,v) >  \epsilon^{-1} \log n \big\}.\label{eq:attachedTree}
    \end{align}
    Notice that for the model $\tilde{\cC}_1(p)$ in Theorem \ref{T:C1structure}, the event $\mathcal{B}_1$ is a condition on independent geometric random variables, and the event $\mathcal{B}_2$ is a condition on the height (or equivalently extinction time) of independent Galton-Watson trees. 
    
    A Taylor expansion gives that for $\mu$ as in Theorem \ref{T:C1structure}
    \begin{equation} \label{eq:taylor_nu}
        \mu = 1 - \epsilon + \frac{2}{3} \epsilon^2 + O(\epsilon^3) = 1 - (1+o(1))\epsilon,
    \end{equation}
    and so, by using a union bound with index size maximally $n^2$, we have for $\tilde{\P}(\cdot)$ the law of $\tilde{\cC}_1(p)$ that
    \begin{equation*}
            \tilde{\P}( \mathcal{B}_1) \leq n^2 (1 - \mu)^{3 \epsilon^{-1} \log n} = O\big( \frac{1}{n} \big).
    \end{equation*}
    Similarly, as the probability of a ${\mathrm{Poisson}}(\mu)$-Galton-Watson surviving for at least $k$ generations is bounded by $\mu^k$ (see e.g.\ Section 12.3 of \cite{LP16}), we have
    \begin{equation*}
        \tilde{\P}( \mathcal{B}_2) = O \big( \frac{1}{n} \big).
    \end{equation*}
    Applying Theorem \ref{T:C1structure} then gives that both $\mathcal{B}_1^c$ and $\mathcal{B}_2^c$ hold for $\cC_1(p)$ with high probability, and for the remainder of the proof we restrict to these events.

    Let $\mathcal{P}(x,y)$ be some (acyclic) path in $\cT_{p}$ between two vertices $x$ and $y$. Then, we can decompose $\mathcal{P}(x,y)$ into a path from $x$ to a vertex $u$ in the 2-core, a path from $u$ to a vertex $U$ that is kept in the kernel, a path between two vertices $U$ and $V$ that are in the kernel, a path from $V$ to a vertex $v$ in the 2-core, and a path from $v$ to $y$. That is, we can write $\mathcal{P}(x,y)$ as
    \begin{equation} \label{eq:path_decomp_kernel}
        \mathcal{P}(x,y) = \mathcal{P}(x,u) \cup \mathcal{P}(u,U) \cup \mathcal{P}(U,V) \cup \mathcal{P}(V,v) \cup \mathcal{P}(v,y),
    \end{equation}   
    where some of the paths may possibly be empty if $x$ (or $y$) is already in the 2-core. Write $|\mathcal{P}|$ for the length of the path $\mathcal{P}$, then, restricted to $\mathcal{B}_1^c$ and $\mathcal{B}_2^c$, we have
    \begin{equation} \label{eq:couple_upper}
        |\mathcal{P}(x,y)| \leq  2 \big( \epsilon^{-1} \log n + 3 \epsilon^{-1} \log n \big) + |\mathcal{P}(U,V)|.
    \end{equation}

    Next, consider the coupling from Lemma \ref{L:couple_mac_mic} between $\cT_{p}$ and $\cT^\Ker_{p}$ that guarantees that
    \begin{equation}
        (U, V) \in \cT^\Ker_{p} \iff \phi(U,V) \subset \cT_{p}.
    \end{equation} 
    This coupling, together with the event $\mathcal{B}_1^c$, gives that
    \begin{equation*}
        | \mathcal{P}(U,V) | \leq 3 \epsilon^{-1} {\diam}(\cT^\Ker_{p}) \log n,
    \end{equation*}
    so that the upper bound of the diameter in Theorem \ref{T:general_kernel} follows from \eqref{eq:couple_upper}.

    \medskip

    For the lower bound, consider again the coupling from Lemma \ref{L:couple_mac_mic}, and let $U,V$ be vertices in $\cT^\Ker_{p}$ which are at distance ${\diam}(\cT^\Ker_{p})$ apart, i.e.\ the path $\mathcal{P}(U, V)$ in $\Ker(\cC_1(p))$ has maximal length among paths in $\cT^\Ker_{p}$.  By abuse of notation, we write
    \begin{equation*}
        \phi \big( \mathcal{P}(U, V) \big) = \bigcup_{(X,Y) \in \mathcal{P}(U, V)} \phi(X,Y)
    \end{equation*}
    for the corresponding (uncontracted) path in $\cC_1(p)$.  We will show that there exists constants $C^*, C_s > 0$, such that, with high probability,
    \begin{equation} \label{eq:Kern_lower_bound}
        \big| \phi \big( \mathcal{P}(U, V) \big) \big| \geq C_s \epsilon^{-1}\bigg\lfloor \frac{|\mathcal{P}(U, V)|}{ C^* \log n} \bigg\rfloor.
    \end{equation}
    
     Let $U_0, U_1, \ldots, U_m$, with $U_0 = U$ and $U_m = V$, be the vertices (in order) on the path $\mathcal{P}(U, V)$. We claim that, with high probability, there are at least $r := \lfloor m/(C^* \log n)\rfloor$ indices $i$, $1 \leq i \leq m$, with
    \begin{equation} \label{eq:kernel_short_paths}
        |\phi(U_{i-1}, U_i)| > C_s \epsilon^{-1}.
    \end{equation}
    Assume for the moment that this is the case. Then this would imply the equality in \eqref{eq:Kern_lower_bound}, and consequently, the lower bound in Theorem \ref{T:general_kernel} follows with $C_K = C_s/2C^*$.
    
    Assume otherwise that \eqref{eq:kernel_short_paths} does not hold for $r$ many indices. Then there must exist a subsequence $U_{i}, U_{i+1}, \ldots U_{i+k}$ of length at least $C^* \log n$, such that for each $1 \leq j \leq k$ we have that $|\phi(U_{i+j-1}, U_{i+j})|$ is at most $C_s \epsilon^{-1}$. 
    % In other words, the event
    % \begin{equation*}
    %     \mathcal{B} = \bigg\{
    %     \begin{aligned}
    %         \ \exists &\text{ subpath $\mathcal{P}'$ of $\mathcal{P}(U,V)$ with $|\mathcal{P}'| \geq \log n$ } \\ \vspace*{-0.1cm}
    %         &\text{and $|\phi(X,Y)| \leq C_s \epsilon^{-1}$ } \forall (X,Y) \in \mathcal{P}'
    %     \end{aligned}
    %     \bigg\}
    % \end{equation*}
    % must hold
    In other words, there is a subpath of length at least $ C^* \log n$ such that each corresponding path of an edge in the 2-core has a length at most $C_s \epsilon^{-1}$. Denote this event by $\mathcal{B}$.  To show that $\mathcal{B}$ only occurs with small probability, consider again the construction of $\tilde{\cC}_1(p)$ with corresponding $\Lambda$ and $D_u$ as in Theorem \ref{T:C1structure}. By the approximation in \eqref{eq:taylor_nu} and the requirement that $\epsilon^3 n$ diverges, we see that $\epsilon \leq \Lambda \leq 3 \epsilon$ with high probability. Applying Chernoff's inequality (see e.g.\ \cite[Section 2.3]{Ver18}),
    %with $t = \log(2/s\Lambda)$ 
    furthermore shows that
    \begin{equation*}
        \P( D_u \geq 2/s) \leq \Big( \frac{e \Lambda }{2/s} \Big)^{2/s} e^{-\Lambda} \leq \Big( \frac{3 e s \epsilon }{2} \Big)^{2/s}.
    \end{equation*}
    As $\beta$ was chosen such that $\beta \gg n^{1+s}$ and thus $\epsilon \ll n^{-s +o(1)}$, a union bound gives that
    \begin{equation} \label{eq:bounded_degree}
        \P \big( \exists u \in \{1, \ldots, n\} \, : \, D_u \geq 2/s \big) = O \big( \frac{1}{n} \big).
    \end{equation}
    In particular, with high probability, the kernel $\tilde{\mathcal{K}}$ of $\tilde{\cC}_1(p)$ has maximum degree bounded by $2/s$. 
    
     For $k \geq 1$, we perform a percolation-like process on the kernel graph $\tilde{\mathcal{K}}$ by keeping each edge $(U,V)$ if and only if the replacement path $\phi(U,V)$ in step 2 of Theorem \ref{T:C1structure} has length at most $k$. If we denote by $F$ the CDF of a random variable distributed as ${\mathrm{Geom}}(1-\mu)$, then the probability of keeping an edge in this procedure is $F(k)$ and is independent of all other edges. As a consequence, we may equivalently percolate $\tilde{\mathcal{K}}$ with $p = F(k)$ as in Section \ref{S:ERintro}. Since we may assume that the kernel graph has degrees bounded by $2/s$, Lemma \ref{L:perc_subcrit} gives that there exists a $p_s > 0$ such that if $p \leq p_s$, then with probability at least $1- o(1)$ every $p$-open cluster in $\tilde{\mathcal{K}}$ has size less than $C^* \log n$, for some fixed constant $C^* > 0$.
     
     Let
     \begin{equation*}
         k^* = \frac{\log\big( \frac{1}{1-p_s} \big)}{2 \epsilon} \in \R,
     \end{equation*}
     then, using the approximation in \eqref{eq:taylor_nu} for $n$ large enough, we have for all $k \leq k^*$ that
     \begin{equation*}
         F(k) = 1 - \mu^{\lfloor k \rfloor} \leq p_s.
     \end{equation*}
    In particular, if $k = k(p_s) = \lfloor k^* \rfloor =: C_s \epsilon^{-1} > 0$, then, with high probability, every $F(k)$-open cluster contains at most $C^* \log n$ many vertices, and any $F(k)$-open path can have length at most $C^* \log n - 1$. Therefore, every path in $\tilde{\mathcal{K}}$ of length greater than $C^* \log n $ contains at least one edge for which the corresponding path in the 2-core of $\tilde{\cC}_1(p)$ has length at least $C_s \epsilon^{-1}$.  By Theorem \ref{T:C1structure}, the same must hold in $\cC_1(p)$ with high probability. It now follows that \eqref{eq:Kern_lower_bound} holds with high probability. %
    %Small issue with size being number of vertices and not number of edges... Should be ok - a path has length = #vertices - 1. So if length path $\geq \log n$ - so does the number of vertices
\end{proof}

\begin{remark}
    The proof of Theorem \ref{T:general_kernel} did not require the fact that the edge weights were i.i.d.\ uniforms on $[0,1]$. In fact, the theorem holds for arbitrary weights on the edges in $\cC_1(p)$ as long as the series law is used to determine the edge weights in the kernel.
\end{remark}

%%%%%%%%%%%%% NEW SECTION %%%%%%%%%%%%%%%

\section{Kernel graph with edge weights} \label{S:macro_weight}
In this section, we give some bounds on the edge weights in the kernel. Consider the construction of $\tilde{\cC}_1(p_1)$ as in Theorem \ref{T:C1structure}, and let $L = L(U,V) \sim {\mathrm{Geom}}(1 - \mu)$, where $\mu = (1+o(1)) \epsilon$. Recall that each edge $e$ in the 2-core has a disorder variable $\omega_e$ that is conditioned on $\{ \omega_e \leq p_1 \}$, i.e.\ the law of $\omega_e$ is the same as $p_1 U_e$ for $U_e \sim U[0,1]$. Therefore, by the series law (see also Section \ref{S:Kernel}), the weight of an edge $(U,V)$ in the kernel $\tilde{\mathcal{K}}$ of $\tilde{\cC}(p_1)$ satisfies
\begin{equation}
    \widehat{\weight}(U,V) \stackrel{(d)}{=} \sum_{k=1}^L e^{-\beta p_1 U_k}, \label{eq:law_macro}
\end{equation}
where $U_k$, $k = 1, \ldots, L$, are i.i.d.\ uniforms on $[0,1]$. Furthermore, the weight of any other edge in the kernel is independent of $\widehat{\weight}(U,V)$, but has the same distribution. %Denote the CDF of the sum \eqref{eq:law_macro} by $F$ (which depends both on $L$ and the $U_k$'s). 
Suppose that $L$ is given, then we let 
\begin{equation}
    M_L := \min_{k =1, \ldots, L} U_k \label{eq:min_geo_U}
\end{equation}
be the minimum of the uniform random variables. We will prove the following bounds on $\widehat{\weight}(U,V)$.

\begin{lemma} \label{L:macro_weight_dist}
    There exists some constant $C > 0$, such that, with probability at least $1 - O(n^{-4})$, it holds that
    \begin{equation*}
       e^{-\beta p_1 M_L} \leq \widehat{\weight}(U,V) \leq C  e^{-\beta p_1 M_L} \log n.
    \end{equation*}
\end{lemma}

To prove Lemma \ref{L:macro_weight_dist}, we will make use of the following classical order statistics result and refer to \cite[Theorem 6.6]{Das11} for a proof.
\begin{theorem} \label{T:order_stat}
Let $U_1, U_2, \ldots, U_m$ be independent $U[0, 1]$ variables, and denote by $U(1), U(2), \ldots , U(m)$ the order statistics, i.e.\ $U(1) \leq U(2) \leq \ldots \leq U(m)$. If $X_1, X_2, \ldots, X_{m+1}$  are $(m + 1)$ independent standard exponentials (also independent of the $U_i$, $1 \leq i \leq m$), then
\begin{equation*}
    \big( U(1), U(2), \ldots, U(m) \big) \stackrel{(d)}{=} \Big( \frac{X_1}{\sum_{i=1}^{m+1} X_i}, \frac{X_1 + X_2}{\sum_{i=1}^{m+1} X_i}, \ldots, \frac{\sum_{i=1}^m X_i}{\sum_{i=1}^{m+1} X_i}\Big).
\end{equation*}
\end{theorem}

\begin{proof}[Proof of Lemma \ref{L:macro_weight_dist}]
    The lower bound trivially holds. For the upper bound, consider the sum in \eqref{eq:law_macro} and denote by $U(1), \ldots, U(L)$ the order statistics of the uniforms (with $U(j) := \infty$ for $j > L$). We claim that for some large $C \geq 1$, we have with high probability that
    \begin{equation} \label{eq:U_order_bound}
        \sum_{k=1}^{L} e^{-\beta p_1 U_k} \leq (C \log n) e^{-\beta p_1 U(1) } + L e^{-\beta p_1 U(C \log n)} \leq (2C \log n) e^{-\beta p_1 U(1) }.
    \end{equation}
    The first inequality is trivial, while for the second inequality, it suffices to show that
    \begin{equation*}
        U(1) \leq  U(C \log n) - \log L \frac{ n }{\beta}. %U(C \log n) - \frac{ \log L}{\beta p_1} \leq 
    \end{equation*}
   Hence, to prove \eqref{eq:U_order_bound} using Theorem \ref{T:order_stat}, it is enough to show that for i.i.d.\ standard exponentials $X_i$ we have
    \begin{equation*}
        \log L \frac{n}{\beta} \sum_{i=1}^{L+1} X_i \leq \sum_{i=2}^{C \log n} X_i
    \end{equation*}
    with high probability.
    
    Recall that $L$ is a geometric random variable with mean $(1 +o(1)) \epsilon^{-1}$, so that we have $L \leq 4 \epsilon^{-1} \log n$ with probability at least $1 - O(n^{-4})$. If we restrict to this event, it is sufficient to show that
    \begin{equation} \label{eq:exp_i_bound}
         \frac{\epsilon}{g_0} \sum_{i=1}^{4 \epsilon^{-1} \log n + 1} X_i \leq \sum_{i=2}^{C \log n} X_i,
    \end{equation}
    where we used Assumption \ref{as:betaeps_s} on the relation between $\beta$ and $\epsilon$, and the fact that $\log(4\epsilon^{-1} \log n) \leq \log n$. An application of the sub-exponential Bernstein's inequality (see \cite[Section 2.8]{Ver18}) for i.i.d.\ exponentials $X_i$ gives that there exists a constant $c > 0$ such that for any $m \in \N$ and $ 0 < \delta < 1$ we have
    \begin{equation*}
        \P \Big( \Big| \sum_{i=1}^m X_i - m \Big| \geq \delta m \Big) \leq 2 \exp( -c \delta^2 m).
    \end{equation*}
    If $C > 0$ is sufficiently large, then with probability at least $1 - O(n^{-4})$, the bound in \eqref{eq:exp_i_bound} holds, completing the proof of Lemma \ref{L:macro_weight_dist}.
\end{proof}

\subsection{Approximation of edge weight distribution}

Suppose now that we replace the weight of the edges $(U,V)$ by
\begin{equation*}
    \weight^*(U,V) = e^{-\beta p_1 M_L},
\end{equation*}
where $M_L$ is as in \eqref{eq:min_geo_U}. Then by Lemma \ref{L:macro_weight_dist}, we have with high probability that 
\begin{equation*}
    \weight^*(U,V) \leq \widehat{\weight}(U,V) \leq C \weight^*(U,V)  \log n 
\end{equation*}
for all edges of $\tilde{\mathcal{K}}$ simultaneously. Denote by $F$ the CDF of $\weight^*(U,V)$. Then for $t \in [\exp(-\beta p_1),1]$,  it holds that
\begin{align*}
    F(t) &= \P \Big( M_L \geq - \frac{\log t}{\beta p_1} \Big) = (1 -\mu) \sum_{i=1}^\infty \P \Big(U_1 \geq - \frac{\log t}{\beta p_1} \Big)^i \mu^{i-1} \\
    &= (1 - \mu) \frac{1 + \frac{\log t}{\beta p_1}}{1 - \mu \big(1 + \frac{\log t}{\beta p_1} \big)} = (1+o(1)) \frac{\epsilon \big(1 + \frac{\log t}{\beta p_1} \big)}{\epsilon  - \frac{\log t}{\beta p_1} + \epsilon \frac{\log t}{\beta p_1} },
\end{align*}
and for the inverse CDF
\begin{align*}
    F^{-1}(q) = \exp \Big( - \beta p_1 \frac{(1-\mu)(1-q)}{1 - \mu + \mu q} \Big) &= \exp \Big(-(1+o(1))\beta p_1 \frac{\epsilon (1-q)}{q+ \epsilon - q\epsilon} \Big) \\
    &= \exp \Big(-(1+o(1))g_0 \log n \frac{(1-q)}{q+ \epsilon - q\epsilon} \Big) .
\end{align*}
In particular, for $q_1, q_2 \in [0,1]$ with $q_1, q_2 \gg \epsilon$, by further absorbing the $\epsilon$ and $q_i \epsilon$ terms into the $o(1)$ factor, we obtain
\begin{align}
    \frac{F^{-1}(q_1)}{F^{-1}(q_2)} &= \exp \Big(-(1+o(1))g_0 \log n \big( \frac{(1-q_1)}{q_1} -  \frac{(1-q_2)}{q_2}\big) \Big) \nonumber \\
     &= \exp \Big(-(1+o(1))g_0 \log n \frac{q_2-q_1}{q_1q_2} \Big). \label{eq:invCDF_kernel}
\end{align}

\medskip

For some parameter $q \in [0,1]$, call an edge $(U,V)$ $q$-open if $\weight^*(U,V) \leq F^{-1}(q)$; otherwise, call it $q$-closed. We end this section with three conclusions:
\begin{enumerate}
    \item In view of \eqref{eq:invCDF_kernel} and Lemma \ref{L:macro_weight_dist}, the ratio of weights of two typical edges in $\tilde{\mathcal{K}}$ is of polynomial size instead of exponential size as on the complete graph;
    
    \item If two edges $(U,V)$, $(X,Y)$ in $\tilde{\mathcal{K}}$ are $q_1$-closed but $q_2$-open with $|q_1 - q_2| = \bigO{\frac{1}{ \log n}}$ and $q_1, q_2 \gg \epsilon$, then
    \begin{equation*}
        \frac{\weight^*(U,V)}{\weight^*(X,Y)} \asymp 1;
    \end{equation*}

    \item If instead $(U,V)$ is $q_1$-open and $(X,Y)$ is $q_2$-closed with $q_2 - q_1 \geq 4/g_0$, then
    \begin{equation*}
        \frac{\weight^*(U,V)}{\weight^*(X,Y)} \geq n^4.
    \end{equation*}
\end{enumerate}
That is, the gap we encountered in Remark \ref{R:compare} (see also Remark \ref{R:constant_dist_gap}) is, for the RSTRE on $\tilde{\mathcal{K}}$, of logarithmic order instead of order $n^{-c}$ for some $c > 0$. This is one of the key ingredients for the heuristics in the next section.

\section{Heuristics and diameter of UST on the giant component} \label{SS:heuristics}

Theorem \ref{T:Kernel} gives rise to the following heuristic argument for the intermediate regime $\beta = n^{1 + \gamma}$, $0 < \gamma < 1/3$, of Conjecture \ref{C:Intermediate}. If $\beta \gg n^{1+s}$, then the kernel graph $\mathcal{K}$ is a bounded-degree expander (see \eqref{eq:bounded_degree}) with 
\begin{equation*}
     |V(\mathcal{K})| \asymp n^4 \log^3 n/ \beta^3 \approx n^{1 - 3 \gamma}, %\\
     %&  {\diam}(\cT) \approx n^{\gamma} {\diam}(\cT^\Ker_{p_1}), 
\end{equation*}
where $\approx$ ignores logarithmic factors. In Section \ref{S:macro_weight}, we have seen that on the kernel graph, the ratio of two typical edge weights is only of polynomial size, and not exponentially large in $n$. Furthermore, with a probability of order $1/\log n$, the weights of two randomly chosen edges are within logarithmic multiples of each other. This might lead one to believe that the diameter of $\cT^\Ker_{p_1}$ behaves as the diameter of the unweighted UST on the kernel of $\cC_1(p)$. If this were true, Theorem \ref{T:Kernel} shows that (up to polylogarithmic factors)
\begin{equation} \label{eq:heuristic_kern}
    {\diam}(\cT) \approx \frac{\beta}{n} \sqrt{|V(\mathcal{K})|} \approx \frac{\beta}{n} \sqrt{\frac{n^4 }{\beta^3}} = \frac{n}{\sqrt{\beta}} = n^{1/2 - \gamma/2}.
\end{equation}
This is precisely the statement of Conjecture \ref{C:Intermediate} for the intermediate regime.

\subsection{3-regular graph}

In the case of $1/4 < \gamma < 1/3$, we may further strengthen this heuristic. Namely, for $\gamma$ in this range, the kernel graph from Theorem \ref{T:C1structure} is a random $3$-regular graph on roughly $n^{1 - 3 \gamma}$ many vertices. As seen in Example \ref{ex:3-reg}, a lot is known about critical and near-critical percolation on the random $3$-regular graph. In particular, if we denote by $G^3_{m}$ a realization of such a graph on $m$ vertices, then for $\epsilon = o(1)$ we have with high probability that
\begin{align*}
    \big| \cC_1( G^3_m, \frac{1 - \epsilon}{2}) \big| &\asymp 2 \epsilon^{-1} \log(\epsilon^3 m) \\
    \big| \cC_1( G^3_m, \frac{1 + \epsilon}{2}) \big| &\asymp 2 \epsilon m,
\end{align*}
which mimics the behavior of the \ErdosRenyi random graph. We refer to \cite{NP10} for more general results of this form for $d$-regular graphs.

We make the following conjecture about the expansion properties of $\cC_1(G^3_m, p)$, where we remark that the analogous statement is true for the \ErdosRenyi random graph.
\begin{conjecture} \label{C:expand3}
    There exists a constant $k > 0$ such that with high probability
    \begin{equation}
        \cC_1 \big( G^3_m, \frac{1 + \frac{1}{\log m}}{2} \big) \text{ is a } (\log m)^{-k} \text{ expander}.
    \end{equation}
\end{conjecture}

\noindent Let $G^3_m$ be the 3-regular kernel graph of the giant component $\cC_1( (1+\epsilon)/n)$ in the regime $\epsilon = n^{-\gamma}$, $1/4 < \gamma < 1/3$, with $m \asymp \epsilon^3 n$. Denote by $p^{\pm} = (1 \pm 1/\log n)/2$, and let $\cT^3$ and $\cT^3_{p^+}$ be RSTREs on $G^3_m$ and $\cC_1 \big( G^3_m, p^+)$, respectively. As shown in Section \ref{S:macro_weight}, the ratio of edge weights of any two $p^+$-open but $p^-$-closed edges is bounded from below and above by logarithmic factors. If Conjecture \eqref{C:expand3} is true, then bumping the $p^-$-open edges as in Section \ref{S:toy_model} gives, with arguments similar to the proof of Theorem \ref{T:toy_model_bump}, that
\begin{equation} \label{eq:bump3}
    \diam(\cT^3_{p^+}) \approx \sqrt{ \Big|  \cC_1 \big( G^3_m, p^+ \big) \Big|} \approx \sqrt{\epsilon^3 n} \approx n^{1/2 - 3\gamma/2},
\end{equation}
where we again ignored any logarithmic factors. 

The last step in this heuristic, and the step in which we believe the heuristic might break down, is to compare the diameter of $\cT^3$ and $\cT^3_{p^+}$. Namely, if we can show that, similarly to Lemma \ref{L:Kernel_reduce}, we have
\begin{equation} \label{eq:compare3_giant}
    \diam(\cT^3_{p^+}) - o(n^{1/3}) \leq \diam(\cT^3) \leq \diam(\cT^3_{p^+}) + o(n^{1/3}),
\end{equation}
then the argument around \eqref{eq:heuristic_kern} and \eqref{eq:bump3} gives the conjectured order of the diameter in the intermediate regime $1/4 < \gamma < 1/3$. One reason to believe that \eqref{eq:compare3_giant} might hold is that $p^+$ is almost supercritical, in the sense that the giant component nearly contains a linear number of vertices. This ensures that a random walk started outside the giant hits the giant component within a small number of steps, which is very similar to the last argument in the proof of Proposition \ref{P:DiaC1} (see Lemma \ref{L:fasthitLERW}). %However, paths inside $\cT^3_{p^+}$ may leave the giant component $\cC_1(G^3_m, p^+)$.

\subsection{Unweighted UST on the giant component}

In this last section, we use Theorem \ref{T:general_kernel} and Theorem \ref{T:toy_model_bump} to prove the following theorem about the diameter of the (unweighted) UST on the giant component.

\begin{theorem} \label{T:USTgiantdiam}
    Let $ p  = p(n) = (1 + \epsilon)/n$ with $n^{-1/3} \ll \epsilon \ll n^{-s}$, for some $s > 0$. Denote by $\cT_p$ a realization of the unweighted ($\beta = 0$) UST on $\cC_1(p)$. Then there exists some constant $c =c(s)> 0$ such that
    \begin{equation} \label{eq:C1USTd}
        \widehat{\P} \Big( \frac{1}{(\log n)^c} \sqrt{\epsilon n} \leq \diam(\cT_p) \leq \sqrt{\epsilon n} (\log n)^c \Big) \xrightarrow{n \rightarrow \infty} 1.
    \end{equation}
\end{theorem}
\begin{proof}
    Let $c > 0$ be sufficiently large, to be determined later. In view of Theorem \ref{T:C1structure}, we may work with $\tilde{\cC}_1(p)$ instead of $\cC_1(p)$. Indeed, if we let $\mathcal{A}$ be the following set of graphs
    \begin{equation*}
        \Big\{ G = (V,E) : \bP_{G,0}\Big( \frac{1}{(\log |V|)^c} \sqrt{|V|} \leq \diam(\cT) \leq \sqrt{|V|} (\log |V|)^c  \Big) < 1 - \frac{1}{\log |V|}  \Big\},
    \end{equation*}
    and show that $\P( \tilde{\cC}_1(p) \in \mathcal{A}) \rightarrow 0$, then the fact that $|\cC_1(p)| \asymp |\tilde{\cC}_1(p)| \asymp \epsilon n \rightarrow \infty$ implies \eqref{eq:C1USTd} by Theorem \ref{T:C1structure}.

    \smallskip

    We wish to apply Theorem \ref{T:Nach} to $\tilde{\cC}_1(p)$, however, as the mixing time on $\cC_1(p)$ (and $\tilde{\cC}_1(p)$) is of order $\epsilon^{-3} (\log \epsilon^3 n)^2$ (see \cite{DLP12}), the mixing condition in \eqref{eq:Nmixing} is not necessarily satisfied for all choices of $\epsilon$. Nonetheless, Lemma 3.5 of \cite{DKLP11} states that the kernel $\tilde{\mathcal{K}}$, from step 1 in Theorem \ref{T:C1structure}, is an $\alpha$-expander for some constant $\alpha > 0$. Furthermore, as we assume $\epsilon \ll n^{-s}$ for some $s > 0$, the calculations around \eqref{eq:bounded_degree} show that the kernel $\tilde{\mathcal{K}}$ has degrees bounded by $2/s$ with high probability. This means that the graph $\tilde{\mathcal{K}}$ falls into the family of graphs covered by Theorem \ref{T:toy_model_bump}. In the last step of the proof, we will apply Theorem \ref{T:general_kernel} so that for the moment it suffices to consider the RSTRE on the (weighted) kernel graph.

    %Next, we identify the weight distribution on the edges of $\tilde{\mathcal{K}}$. 
    In the construction of $\tilde{\cC}_1(p)$ in Theorem \ref{T:C1structure}, each edge in the kernel corresponds to a $\textrm{Geom}(1-\mu)$ number of edges in the 2-core. In particular, if we denote by $\widehat{\weight}(U,V)$ the weight on the kernel edge $(U,V)$ induced by the series law, then
    \begin{equation*}
        \widehat{\weight}(U,V) = L_{(U,V)} = \exp\big( \log L_{(U,V)} \big),
    \end{equation*}
    where $(L_{(U,V)})_{(U,V) \in E(\tilde{\mathcal{K}})}$ are i.i.d.\ geometric random variables with mean $1 - \mu = (1 +o(1)) \epsilon$. Denote by $\mu_n$ the law of the random variable\footnote{This parametrization ensures that the weight is given by $\widehat{\weight}(e) = \exp( - \omega_e) = L_e$, which is in line with the parametrization of the RSTRE in \eqref{eq:defRSTRE}.} $-\log L_{(U,V)}$, and by $\cT^\Ker_p$ the RSTRE on $\tilde{\mathcal{K}}$ with disorder variables distributed according to $\mu_n$. We verify the conditions on $a_1,a_2, \mu_n$ in Theorem \ref{T:toy_model_bump} by noticing that for $a_1, a_2 \in \R$ with
    \begin{align*}
        a_1 &= a_1(n) = - \log( \epsilon^{-1}) - \log\big( \log \frac{1}{q_1} \big) - \log 2, \\
        a_2 &= a_2(n) = - \log( \epsilon^{-1}) - \log \big( \log \frac{1}{1-q_2} \big) + \log 2,
    \end{align*}
    we have, for $n$ large enough, 
    \begin{align*}
        \mu_n( - \infty, a_1) &\leq q_1, \\
        \mu_n( a_2, \infty) &\leq q_2.
    \end{align*}
    If $q_1, q_2$ are fixed in $n$, then such a choice of $a_1, a_2$ satisfies
    \begin{equation*}
        a_2 - a_1 \asymp 1 \ll \log n.
    \end{equation*}
    Applying Theorem \ref{T:toy_model_bump} to $\cT^\Ker_p$, with $\Delta = 2/s$, $b = \alpha$ and $\delta = (\log n)^{-1}$, gives that there exists a constant $c = c(s) > 0$ such that 
    \begin{equation*}
        \widehat{\P} \Big( \frac{1}{(\log n)^c} \sqrt{\epsilon^3 n} \leq \diam(\cT^\Ker_p) \leq \sqrt{\epsilon^3 n} (\log n)^c \Big) \xrightarrow{n \rightarrow \infty} 1,
    \end{equation*}
    where we used that the kernel graph contains on the order of $\epsilon^3 n$ many vertices. Theorem \ref{T:general_kernel} applied on $\cT_p$ and $\cT^\Ker_p$, after possibly enlarging $c$,  completes the proof.  
\end{proof}

We give a few more remarks about the diameter of the UST on $\cC_1(p)$ for values of $p$ that are not in the slightly supercritical regime.

\begin{remark}
    When $p= (1 + \lambda n^{-1/3})/n$ is in the critical window, then the number of excess edges is with high probability bounded, so that the UST has a diameter of order $n^{1/3} \asymp |\cC_1(p)|$. For this, see also the proof of Theorem \ref{T:highDisorder}. Furthermore, if $p = (1 - \epsilon)/n$ and $\epsilon \gg n^{-1/3}$, then the largest component is with high probability already a tree. In this case, it is uniformly distributed among all spanning trees on $|\cC_1(p)|$ many vertices, i.e.\ $\cC_1(p)$ is already a UST with a diameter of order $\sqrt{|\cC_1(p)|}$.
\end{remark}

\begin{remark}
    Assume that $p= (1 + \epsilon)/n$ for $\epsilon = \epsilon(n)$, with $c < \epsilon = O(\sqrt{\log n})$ for some fixed constant $c>0$. In \cite[Theorem 1.2]{FR08}, it was shown that the mixing time of the lazy random walk is of order $(\log n)^2$. Furthermore, it is easy to verify that the ratio of maximum to minimum degree is of order at most $\log n$. Therefore, since the size of $\cC_1(p)$ is of order $n$, naively applying Theorem \ref{T:Nach} already gives that the diameter of the UST on $\cC_1(p)$ is of order $\sqrt{n}$, up to some logarithmic corrections.
\end{remark}

\noindent The most important regime not covered by Theorem \ref{T:USTgiantdiam} and the above remarks is when $p = (1 + \epsilon)/n$ with $n^{-s} \ll \epsilon = o(1)$ for every $s > 0$. We believe that the analogous statement of Theorem \ref{T:USTgiantdiam} in this regime should also follow from Theorem \ref{T:toy_model_bump}, as long as one takes more care with the maximum degree $\Delta$ and the dependence of $a_1$ on $\Delta$. Additionally, since by the results of \cite{DLP12} the mixing time in this regime is of order $n^{o(1)}$, we can apply Theorem \ref{T:Nach} directly to the giant component to obtain similar bounds as \eqref{eq:C1USTd}. However, in this case, the correction terms are of order $n^{o(1)}$, which may not be of polylogarithmic size.

% \begin{question}
%     Let $p = (1 + \epsilon)/n$ with $n^{-s} \ll \epsilon = o(1)$ for every $s > 0$. What is the order of the diameter of the (unweighted) UST on $\cC_1(p)$?
% \end{question}

\appendix

% \begin{theorem}[McDiarmid..]
    
% \end{theorem}

% \begin{definition}[Sub-gaussian norm]
%     \begin{equation*}
%         \Vert X \Vert_{\psi_2} = \inf \{ t > 0 : \E[ \exp( X^2/t^2) \leq 2]
%     \end{equation*}
% \end{definition}

% There exists some universal constant $c > 0$, such that if $X$ is a sub-gaussian random variable, then  
% \begin{equation*}
%     \P ( |X| \geq t) \leq 2 \exp( -\frac{c t^2}{\Vert X \Vert_{\psi_2}^2}) \qquad \forall t \geq 0.
% \end{equation*}
% High dimensional probability section 2.5
% Bernstein...

\chapter{Concentration of measure} 

\section{Concentration inequalities}

We first recall the following standard concentration inequality for positive random variables.

\begin{theorem}[Paley–Zygmund inequality] \label{T:Paley-Z}
    If $Z \geq 0$ is a random variable with finite expectation, then for any $0 \leq t \leq 1$
    \begin{equation*}
        \P \big( Z > t \E[Z] \big) \geq (1- t)^2 \frac{\E[Z]^2}{\E[Z^2]}.
    \end{equation*}
\end{theorem}
\begin{proof}
    Observe that
    \begin{equation*}
        \E[Z] = \E[Z 1_{Z \leq t \E[Z]}] + \E[Z 1_{Z > t \E[Z]}] \leq  t \E[Z] + \E[Z 1_{Z > t \E[Z]}].
    \end{equation*}
    Applying the Cauchy–Schwarz inequality gives that 
    \begin{equation*}
        \E[Z 1_{Z > t \E[Z]}] \leq \E[Z^2]^{1/2} \P( Z > t \E[Z])^{1/2},
    \end{equation*}
    from which the theorem readily follows. 
\end{proof}

The main concentration inequality we use in this thesis is the following Bernstein's inequality, and we refer to Theorem 2.8.4 in \cite{Ver18} for a proof. The more general versions of Bernstein's inequalities for sub-exponential random variables can be found in Section 2.8 of \cite{Ver18}. 

\begin{theorem}[Bernstein's Inequality] \label{T:Bern}
	Let $X_1, \ldots, X_m$ be i.i.d.\ random variables with $|X_i| \leq K$ almost surely, mean $\xi$ and variance $\sigma^2$. Then for $S_m = \sum_{i=1}^n X_i$ and any $\delta > 0$, we have
	\begin{equation*}
		\P \big( |S_m - m \xi| \geq \delta m \big) \leq 2 \exp\left(-\frac{m \delta^2}{2 \sigma^2 + 2K\delta/3}\right).
	\end{equation*}
\end{theorem}

\section{Edge weight concentration} \label{AS:concentration}

As usual, we assume that $\weight(e) = \exp(-\beta \omega_e)$, where $(\omega_e)_{e \in E}$ are i.i.d.\ with common distribution $\mu$. In the special case when $\mu$ is the uniform distribution on $[0,1]$ and $\beta > 0$, it is easy to verify that the probability density function is given by
\begin{equation*}
    f_{\beta}(x) = \frac{1}{\beta x}, \qquad \exp(-\beta) \leq x \leq 1,
\end{equation*}
and that the mean and variance satisfy
\begin{align*}
    \xi &= \xi(\beta) = \E[\weight(e)] = \frac{1 - \exp(-\beta)}{\beta} \leq \frac{1}{\beta}, \\
    \sigma^2 &= \sigma^2(\beta) = \frac{1}{2\beta} \Big(1 - \frac{2}{\beta} + e^{-\beta} \Big( \frac{4}{\beta} - e^{-\beta} - \frac{2 e^{-\beta}}{\beta} \Big) \Big) \leq \frac{1}{\beta}. % =\frac{(1 - e^{-\beta}) \beta - 2 + \beta e^{-\beta} + 2 e^{-\beta}}{2 \beta^2
\end{align*}
Notice that the weights are trivially bounded in absolute value by $K = 1$, and that for $\beta \geq 1$ and $ 0 < \delta_* \leq 1$, we have
\begin{equation} \label{eq:ratio_xi_sigma}
%\frac{ \delta^2}{2 \sigma^2 + 2K\delta/3} \geq  
	\frac{ \xi^2}{2\sigma^2 + 2\delta_*\xi/3 }  
	\geq \frac{ \xi^2}{\frac{2}{\beta} + \frac 1{\beta} } 
	= \frac{ (1-\e^{-\beta})^2 }{3\beta}
	\geq  \frac{1}{9 \beta}. 
\end{equation}

\begin{lemma} \label{L:low_concentration}
Let $S_m = \sum_{i=1}^m \weight(e_i)$ for distinct edges $e_1, \ldots, e_m$. If $\mu$ is the uniform distribution on $[0,1]$, then for any $ 0 < \delta_* \leq 1$, we have
\begin{equation} \label{eq:BernU}
	\P \big( |S_m - m \xi | \geq \delta_* m \xi \big) \leq 2 \exp \left( - \frac{\delta_*^2 m}{9 \max\{ \beta, 1 \}}\right).
\end{equation}
\end{lemma}
\begin{proof}
    Let $\delta = \delta_* \xi$, then \eqref{eq:BernU} follows directly from Bernstein's inequality and \eqref{eq:ratio_xi_sigma} whenever $\beta \geq 1$. For $\beta \in [0,1]$, we remark that $\xi(\beta)$ is a decreasing function in $\beta$, while $\sigma^2(\beta)$ is an increasing function in $\beta$. We may thus verify that for $\beta \in [0,1]$
    \begin{equation*}
        \frac{\xi(\beta)^2}{2 \sigma^2(\beta) + 2 \xi(\beta)/3} \geq \frac{\xi(1)^2}{2 \sigma^2(1) + 2 \xi(0)/3} \geq \frac{1}{9},
    \end{equation*}
    so that \eqref{eq:BernU} holds for $\beta \leq 1$ as well.
\end{proof}
\noindent We remark that Bernstein's inequality is useful for $\beta = o(m)$, whereas Hoeffding's inequality only provides good bounds for $\beta = o(\sqrt{m})$.

\subsection{Assumption \ref{as:mu_tail}}

When $\mu$ is no longer the uniform distribution, it may not be easy to compute all the moments; nonetheless, we may obtain a similar result as in Lemma \ref{L:low_concentration} at the cost of a smaller constant in the exponent. Recall that Assumption \ref{as:mu_tail} says that $\mu$ is supported on $[0,\infty)$ and that there exist constants $\alpha, \rho, c_\mu > 0$ such that, for $F_\mu$ the CDF of $\mu$, we have
\begin{equation*}
    F_{\mu}(t) = c_\mu t^{\alpha} \qquad 0 \leq t \leq \rho.
\end{equation*}
The probability density function $g_{\mu, \beta}(x) = g_\beta(x)$ of $\omega_e$ then satisfies
\begin{equation*}
    g_\beta(x) = \alpha c_\mu x^{\alpha - 1} \qquad \text{for } 0 \leq x \leq \rho.
\end{equation*}
We write
\begin{align*}
    \xi = \xi(\mu, \beta) &= \E[\weight(e)] = \E[e^{-\beta \omega_e}], \\
    \sigma^2 = \sigma^2(\mu, \beta) &= \E[\weight(e)^2] - \E[\weight(e)]^2,
\end{align*}
for the mean and variance of $\weight(e)$. For fixed $r \geq 1$, we have the following (asymptotic) bounds
 \begin{align*}
        \E[ e^{- r \beta \omega_e} ] &\leq \E[ e^{- r \beta \omega_e} 1_{\omega_e \leq \rho}] + e^{- r \beta \rho} = c_\mu \alpha \int_{0}^{\rho} e^{- r \beta t} t^{\alpha - 1} dt + e^{- r \beta \rho} \\
        &=  c_\mu \alpha  \frac{1}{(r \beta)^\alpha} \Big[ - \Gamma(\alpha,  r\beta x) \Big]^{x=\rho}_{x=0} +  e^{- r \beta \rho} \\
        &\leq c_\mu \alpha \frac{\Gamma(\alpha)}{(r \beta)^\alpha} + e^{- r \beta \rho} 
        \asymp (r \beta)^{-\alpha},
\end{align*}
and if $\beta$ is large enough (depending on $\mu$)
\begin{equation*}
        \E[ e^{-r \beta \omega_e} ] \geq  c_\mu \alpha  \frac{1}{(r\beta)^\alpha} \Big[ - \Gamma(\alpha, r \beta x) \Big]^{x=\rho}_{x=0} \geq \frac{1}{2}  c_\mu \alpha \frac{\Gamma(\alpha)}{(r \beta)^\alpha} \asymp (r \beta)^{-\alpha},
\end{equation*}
where $\Gamma(\cdot, \cdot)$ is the upper incomplete Gamma function and $\Gamma(\cdot) = \Gamma(\cdot, 0)$ is the Gamma function. We may conclude that
\begin{equation} \label{eq:tail_xi_bound}
    \begin{aligned}
         \xi &\asymp \beta^{-\alpha}, \\
    \sigma^2 &\asymp \beta^{-\alpha}.
    \end{aligned}
\end{equation}
Using these estimates we may prove the following concentration inequality.

\begin{lemma} \label{L:mu_tail_concentrate}
    If $\mu$ satisfies Assumption \ref{as:mu_tail}, then there exists a constant $C_B = C_B(\mu)$ such that for $\beta \geq 0$ and $S_m$ the sum of $m$ edge weights, we have
    \begin{equation} \label{eq:AP_mu_tail_concentrate}
        \P \big( |S_m - m \xi | \geq \delta m \xi \big) \leq 2 \exp \left( - C_B \frac{\delta^2 m}{ \max\{\beta^\alpha, 1\} } \right).
    \end{equation}
\end{lemma}
\begin{proof}
    By \eqref{eq:tail_xi_bound} and the bounds leading up to it, if $N = N(\mu)$ is some fixed large constant and $\beta \geq N$, then there exists constants $c_1, c_2, c_3 > 0$ such that
    \begin{align*}
        c_1 \beta^{-\alpha} &\leq \xi \leq c_2 \beta^{-\alpha}, \\
        c_3 \beta^{-\alpha} &\leq \sigma^2.
    \end{align*}
    Furthermore, as $\mu$ is supported on $[0, \infty)$ we have that
    \begin{equation*}
        \inf_{\beta \in [0, N]} \frac{\xi^2}{2\sigma^2 + \frac{2}{3} \xi} > 0,
    \end{equation*}
    so that if $C_B = C_B(\mu)$ is chosen small enough, then for all $\beta \geq 0$
    \begin{equation*}
        \frac{\xi^2}{2\sigma^2 + \frac{2}{3} \xi} \geq  C_B \frac{1}{\max \{ \beta^\alpha, 1\}}.
    \end{equation*}
    Applying Bernstein's inequality with $K = 1$ gives the inequality in \eqref{eq:AP_mu_tail_concentrate}.
\end{proof}

Lastly, we compute some conditional expectations required for Section \ref{SS:no_conc}. Since for $x \rightarrow \infty$, we have
\begin{equation*}
    \Gamma(\alpha, x) = x^{\alpha - 1} e^{-x} ( 1 + O(x^{-1}) ),
\end{equation*}
see e.g.\ \cite[Section 8.11]{OLBC10}, it follows from the calculations leading up to \eqref{eq:tail_xi_bound} that, for $\beta,y \geq 0$ with $y \rightarrow 0$ but $\beta y \rightarrow \infty$, we have
\begin{equation} \label{eq:xi_conditional} 
    \begin{aligned}
        \E[ e^{- r \beta \omega_e} \mid \omega_e \leq y ] &\asymp \frac{( r \beta)^{- \alpha}}{F_\mu(y)}, \\
        \E[ e^{- r \beta \omega_e} \mid \omega_e > y ] &\asymp e^{-r \beta y} (r \beta)^{ -1 }y^{\alpha - 1 }.
    \end{aligned}
\end{equation}
Denote by $\xi_1$ and $\sigma^2_1$ the mean and variance of $\omega_e$ conditioned on $\{\omega_e \leq y\}$. If $\beta y \rightarrow \infty$, then as a consequence of \eqref{eq:xi_conditional}, we obtain that 
 \begin{equation} \label{eq:xi_condition_large}
        \frac{\xi^2_1}{2 \sigma^2_1 + \frac{2}{3} \xi_1} \asymp \frac{\beta^{- \alpha}}{F_\mu(y)}.
\end{equation}
Similarly, for $\xi_2$ and $\sigma^2_2$, the mean and variance of $\omega_e$ conditioned on $\{\omega_e > y\}$, we have for $\beta y \rightarrow \infty$ that
 \begin{equation} \label{eq:xi_condition_small}
        \frac{\xi^2_2}{2 \sigma^2_2 + \frac{2}{3} e^{-\beta y} \xi_2} \asymp \beta^{-1} y^{\alpha -1},
\end{equation}
where we remark that we added the factor $ e^{-\beta y}$ corresponding to $K$ in Bernstein's inequality.

Using \eqref{eq:xi_condition_large} and \eqref{eq:xi_condition_small}, we briefly outline how to obtain the asymptotic inequalities in \eqref{eq:condition_concen}. In this case, as $p = 1/n$, we have
\begin{gather*}
    y =  F^{-1}_\mu(1/n) = \frac{1}{(c_\mu n)^{1/\alpha}}, \\
    F_\mu(y) = 1/n,
\end{gather*}
see also \eqref{eq:inverse_CDF}. Therefore, by \eqref{eq:xi_condition_large} and \eqref{eq:xi_condition_small} 
\begin{equation} \label{eq:xi1_xi2_ratio}
    \begin{aligned}
    \frac{\xi^2_1}{2 \sigma^2_1 + \frac{2}{3} \xi_1} &\asymp n \beta^{- \alpha}, \\
    \frac{\xi^2_2}{2 \sigma^2_2 + \frac{2}{3} e^{-\beta y}\xi_2} &\asymp n^{1/\alpha -1} \beta^{- 1}.
\end{aligned}
\end{equation}

\sloppy %
Assume now that $\beta \ll \min\{ n^{4/3\alpha}, n^{1/\alpha + 1/3 } \}$ and let $S = \cC_1(1/n)$, where $|E(S,S)| \asymp n^{1/3}$ and $|E(S,S^c)| \asymp n^{4/3}$ with high probability. Since we have
\begin{align*}
    &|E(S,S)| n \beta^{-\alpha}  \longrightarrow \infty, \\
    &|E(S,S^c)| n^{1/\alpha -1} \beta^{-1} \longrightarrow \infty,
\end{align*}
 it follows from \eqref{eq:xi1_xi2_ratio}, together with Bernstein's inequality (Theorem \ref{T:Bern}), that with high probability
\begin{equation} \label{eq:full_cond_conc}
    \begin{aligned}
    \sum_{e \in E(S,S)} \weight(e) &\asymp n^{4/3} \beta^{-\alpha}, \\
    \sum_{e \in E(S,S^c)} \weight(e) &\asymp n^{1/\alpha +1/3} \beta^{-1} \exp\big( -\frac{\beta}{(c_\mu n)^{1/\alpha}} \big),
\end{aligned}
\end{equation}
as required in \eqref{eq:condition_concen}.

\fussy %

\chapter{Total variation and the proof of Lemma \ref{L:TVgap}} \label{A:TV}

\section{Total variation distance}

We define the \textit{total variation distance} $\Vert\cdot \Vert_{\mathrm{TV}}$ between two probability measures $\nu_1$ and $\nu_2$, on some space $\mathcal{X}$, as
\begin{equation*}
    \Vert\nu_1 - \nu_2 \Vert_{\mathrm{TV}} = \sup\limits_{A \subseteq \mathcal{X}} \big| \nu_1(A) - \nu_2(A) \big|.
\end{equation*}
% As our choice of $\nu_1 = \nu_1^\omega$ and $\nu_2 = \nu_2^\omega$ will be probability measures on finite sets, that is,  the finite set of all spanning trees on $G$, the total variation norm may also be written as
% \begin{equation*}
%      \Vert\nu_1 - \nu_2 \Vert_{\mathrm{TV}} = \frac{1}{2} \sum_{x \in \mathcal{X}} \big| \nu_1(x) - \nu_2(x) \big|,
% \end{equation*}
% see e.g.\ Proposition 4.2 of \cite{LP17}. 
Another equivalent characterization is given by the following proposition.

\begin{proposition}[{\cite[Proposition 4.7]{LP17}}] \label{P:TV_fact}
    Let $\nu_1$ and $\nu_2$ be two probability measures on $\mathcal{X}$. Then
   \begin{equation*}
     \Vert\nu_1 - \nu_2 \Vert_{\mathrm{TV}} = \min \big\{ \nu( X \neq Y) \, : \, (X,Y) \text{ is a coupling of $\nu_1$ and $\nu_2$} \big\}.
\end{equation*}
\end{proposition}
\noindent We will therefore often consider couplings between the measures, which by Proposition \ref{P:TV_fact} gives a connection to the total variation norm.%By abuse of notation, we extend the probability measures $\nu_1, \nu_2$ to also support the couplings we will define later.

\section{Proof of Lemma \ref{L:TVgap}}

To prove Lemma \ref{L:TVgap}, we will work with the Aldous-Broder algorithm instead of Wilson's algorithm\footnote{Almost the same proof also works for Wilson's algorithm, however, the notation becomes heavier as we must consider many different random walks.}. Discovered independently by Aldous \cite{Ald90} and Broder \cite{Bro89}, the Aldous-Broder algorithm constructs the weighted UST by running a random walk (started at any vertex) on the weighted graph until every vertex is visited. Whenever a previously unvisited vertex is seen, the most recent edge that the random walk used is added to the tree. The resulting (random) spanning tree is distributed as a weighted UST, see also \cite[Corollary 4.9]{LP16}.

\smallskip

We will make use of the following bound for the random walk measure $Q_v$ on $(G, \weight)$ started at the vertex $v$.

\begin{lemma}
    Let $A,B \subset V$ and suppose that $v \not\in A \cup B$. Then
    \begin{equation} \label{eq:escape_electric}
    Q_v(\tau_A < \tau_B) \leq \frac{\effR{}{v}{A \cup B}}{\effR{}{v}{A}} \leq \frac{\effR{}{v}{B}}{\effR{}{v}{A}}.
\end{equation}
\end{lemma}
\begin{proof}
    For the first inequality, see \cite[Exercise 2.36]{LP16} for more details. Since 
    \begin{equation*}
        Q_v( \tau_{B} < \tau_v^+) \leq Q_v( \tau_{A \cup B} < \tau_v^+),
    \end{equation*}
     the second inequality follows essentially from the definition of the effective resistance in \eqref{eq:def_effR}.
\end{proof}

Recall the notation of $\cT_{p_1 \mid p_0}$ and $\ClCp{1}{p_0}$ from Section \ref{S:edge_remove}. We are now ready to prove the following lemma, where we remark that many of the proof ideas are similar in spirit to those in the proof of Lemma 2.10 in \cite{K24}.
\begin{lemma*}[Lemma \ref{L:TVgap}]
     Assume $\omega$ is such that $\cC_1(p_0) \subseteq \cC_1(p_1)$. Denote by $\nu_1$ and $\nu_2$ the laws of $\ClCp{1}{p_0}$ and $\cT_{p_1 \mid p_0}$, respectively. Then
    \begin{equation} %\label{eq:TVgap_eq}
        \Vert \nu_1 - \nu_2\Vert_{\mathrm{TV}} \leq n^4 e^{-\beta(p_1 - p_0)}.
    \end{equation}
\end{lemma*}

\begin{proof}
    Suppose that $(G, \weight)$ is given with $\cC_1(p_0) \subseteq \cC_1(p_1)$. Denote by $(H, \weight)$ the weighted graph where $H = \cC_1(p_1)$ and the weights are induced by $(\omega_e)_{e \in E(\cC_1(p_1))}$. Let $v$ be any vertex in $\cC_1(p_0)$, and let $X^{(1)}$ and $X^{(2)}$ denote random walks on $(G, \weight)$ and $(H, \weight)$, respectively, both starting at $v$. For $x \in V(\cC_1(p_0))$, define $\tau_x$ as the first hitting time of $x$ for the random walk $X^{(1)}$. Furthermore, define
    \begin{equation*}
        E_{\geq p_1} = \{ e \in E(G) : \omega_e \geq p_1 \},
    \end{equation*}
    and, for $X^{(1)}$, consider the stopping times
    \begin{align*}
        \tau_{p_0\textrm{-cover}} &= \max_{x \in \cC_1(p_0)} \tau_x \, , \\
        \tau_e &= \min \{ t \geq 1 : (X_{t-1}, X_t) = e \}, \\
        \tau_{\geq p_1} &= \min_{e \in E_{\geq p_1} } \tau_e \, .
    \end{align*}
    
    Notice that if we stop the Aldous-Broder algorithm at time $\tau_{p_0\textrm{-cover}}$, we will have constructed a tree that contains $\ClCp{1}{p_0}$ as a subgraph (there may also be dangling branches containing no vertices of $\cC_1(p_0)$). Moreover, as long as $t < \tau_{\geq p_1}$, we can let $X^{(2)}$ follow $X^{(1)}$ without changing its law, that is, we may couple $X^{(1)}$ and $X^{(2)}$ in such a way that
    \begin{equation*}
        X^{(1)}_t = X^{(2)}_t \qquad  0 \leq t < \tau_{\geq p_1}.
    \end{equation*}
    In particular, we may couple the random walks such that, on the event
    \begin{equation}
        \big\{ \tau_{p_0\textrm{-cover}} < \tau_{\geq p_1} \big\},
    \end{equation}
    the trees $\ClCp{1}{p_0}$ and $\cT_{p_1 \mid p_0}$ coincide. Consequently, by Proposition \ref{P:TV_fact}, we have that
    \begin{equation} \label{eq:TV_bound}
        \Vert \nu_1 - \nu_2\Vert_{\mathrm{TV}} \leq Q_v( \tau_{p_0\textrm{-cover}} > \tau_{\geq p_1}).
    \end{equation}

    To show that the right-hand side of \eqref{eq:TV_bound} is small, we follow a similar approach to the proof of Lemma \ref{L:gap}. For $x \in \cC_1(p_0)$ and $e \in E_{\geq p_1}$ consider the event
    \begin{equation*}
        F(x; e) = \big\{ \tau_e < \tau_x \big\}.
    \end{equation*}
    Suppose that $e = (y,z)$, then we construct a new graph $G^e$ with vertex set $V(G^e) = V(G) \cup \{ yz \}$, edge set 
    \begin{equation*}
        E(G^e) =  \big( E(G) \setminus \{e\} \big) \cup \{ (y, yz), (yz, z) \},
    \end{equation*}
    and weights
    \begin{equation*}
        \weight(f) = 
        \begin{dcases}
            \weight(f) & \text{if } f \in E(G), \\
            \weight(e) & \text{if } f \in \big\{ (y,yz), (yz, z) \big\}.
        \end{dcases}
    \end{equation*}
    That is, we replaced the edge $e = (y,z)$ with a vertex $yz$ that has incident edges going to the original endpoints of $e$. Importantly, if $X^{(3)}$ is a random walk on $G^e$, then $\tau_e$ is precisely equal to the first hitting time of $yz$ for $X^{(3)}$.

    As $v$ and $x$ are in the same $p_0$-open component, the effective resistance (in $G^e$) between $v$ and $x$ is trivially upper bounded by $n \exp(\beta p_0)$. On the other hand, using the parallel law and Rayleigh's monotonicity principle, it is easy to see that the effective resistance in $G^e$ between $v$ and $yz$ is lower bounded by $2 \exp(\beta p_1)$. Therefore, if we denote by $Q_v^{G^e}$ the random walk transition kernel on $(G^e, \weight)$, then we have by \eqref{eq:escape_electric} that
    \begin{equation*}
        Q_v\big( F(x; e) \big) = Q_v^{G^e} \big( \tau_{yz} < \tau_x \big) \leq \frac{\effR{G^e}{v}{x}}{ \effR{G^e}{v}{yz} } \leq 2n e^{-\beta(p_1 - p_0)}.
    \end{equation*}
    A union bound together with \eqref{eq:TV_bound} gives that
    \begin{equation*}
         \Vert \nu_1 - \nu_2\Vert_{\mathrm{TV}} \leq \sum_{x \in \cC_1(p_0)} \sum_{e \in E_{\geq p_1}} Q_v\big( F(x; e) \big) \leq n^4 e^{-\beta(p_1 - p_0)},
    \end{equation*}
    completing the proof.
\end{proof}

\clearpage

\bookmarksetup{startatroot}
\printbibliography[heading=bibintoc]

\end{document}